\documentstyle[11pt,twoside,amsfonts,amssymb]{article}
  
\date{}
\begin{titlepage}

\title{Cohomology of convex cocompact groups and invariant distributions on limit sets}

\author{Martin
Olbrich\thanks{Mathematisches Institut, Universit\"at G\"ottingen, Bunsenstr. 3-5, 37073 G\"ottingen, GERMANY, E-mail: olbrich@uni-math.gwdg.de}
}

\end{titlepage} 

\newcommand{\proof}{{\it Proof.$\:\:\:\:$}}
\newcommand{\D}{\displaystyle}
\newcommand{\dist}{{\mathrm{dist}}}
\newcommand{\kaaa}{{\frak k}}
\newcommand{\paaa}{{\frak p}}
\newcommand{\taaa}{{\frak t}}
\newcommand{\haaa}{{\frak h}}
\newcommand{\Q}{{\Bbb Q}}
\newcommand{\R}{{\Bbb R}}
\newcommand{\Z}{{\Bbb Z}}
\newcommand{\C}{{\Bbb C}}
\newcommand{\F}{{\Bbb F}}
\newcommand{\HH}{{\Bbb H}}
\newcommand{\OO}{{\Bbb O}}

\newcommand{\gaaa}{{\frak g}}
\newcommand{\maaa}{{\frak m}}
\newcommand{\aaaa}{{\frak a}}
\newcommand{\naaa}{{\frak n}}
\newcommand{\qaaa}{{\frak q}}
\newcommand{\laaa}{{\frak l}}

\newcommand{\res}{{\mathrm{ res}}}

\newcommand{\cZ}{{\cal Z}}
\newcommand{\cH}{{\cal H}}

\newcommand{\cC}{{\cal C}}
\newcommand{\mod}{{\mathrm{ mod}}}
\newcommand{\cM}{{\cal M}}
\newcommand{\cO}{{\cal O}}

\newcommand{\cA}{{\cal A}}

\newcommand{\cU}{{\cal U}}
\newcommand{\Hom}{{ \mathrm{ Hom}}}
\newcommand{\vol}{{\mathrm{ vol}}}

\newcommand{\vp}{\varphi}
\newcommand{\ve}{\varepsilon}
\newcommand{\p}{\varrho}

\newcommand{\End}{{ \mathrm{ End}}}

\newcommand{\im}{{ \mathrm{ im}}}
\newcommand{\sign}{{\mathrm{ sign}}}

\newcommand{\Ree}{{\mathrm{ Re }}}

\newcommand{\ee}{{\mathrm{ e}}}

\newcommand{\tr}{{ \mathrm{ tr}}}

\newcommand{\Ad}{{\mathrm{ Ad}}}
\newcommand{\ad}{{ \mathrm{ ad}}}

\newcommand{\coker}{{\mathrm{ coker}}}
\newcommand{\id}{{ \mathrm{ id}}}
\newcommand{\ord}{{ \mathrm{ ord}}}
\newcommand{\nat}{{\Bbb  N}}
\newcommand{\supp}{{ \mathrm{ supp}}}
\newcommand{\spec}{{ \mathrm{ spec}}}

\newcommand{\aca}{{\aaaa_\C^\ast}}

\def\hB{\hspace*{\fill}$\Box$\newline\noindent}

\newcommand{\Fin}{{ \mathrm{ Fin}}}

\newtheorem{prop}{Proposition}[section]
\newtheorem{lem}[prop]{Lemma}
\newtheorem{ddd}[prop]{Definition}
\newtheorem{theorem}[prop]{Theorem}
\newtheorem{kor}[prop]{Corollary}

\newtheorem{con}[prop]{Conjecture}

\def\imath{i}
\newcommand{\ii}{i}

\begin{document}
\setcounter{page}{1}
\maketitle

\vspace{1.5cm}
\begin{abstract}
This paper contains a thorough investigation of invariant distributions
supported on limit sets of discrete groups acting convex cocompactly
on symmetric spaces of negative curvature. It can be considered as
a continuation of \cite{bunkeolbrich000}. Based on this investigation
we
\begin{itemize}
\item provide proofs of the Hodge theoretic results for the cohomology
of real hyperbolic manifolds
announced in \cite{olbrich000},
\item improve the bounds for the critical exponents obtained by Corlette 
for the quaternionic and the Cayley case,
\item compute the $L^2$-cohomology for the corresponding locally symmetric
spaces,
\item prove a version of the Harder-Borel conjecture for real hyperbolic
manifolds, and
\item compute higher cohomology groups with coefficients in hyperfunctions
supported on the limit set.
\end{itemize}

\end{abstract}
\newpage

\tableofcontents
\newpage
 
\parskip3ex

\section{Introduction}
\subsection{The central result}\label{least}

The present paper is an outgrowth of our long lasting attempt to prove a conjecture
of Patterson \cite{patterson93} concerning a Hodge theory for noncompact hyperbolic manifolds. It is a slightly revised version of my habilitation thesis
deliverd at the University of G\"ottingen.

Patterson's conjecture was stated for quotients $Y=\Gamma\backslash X$ of the real hyperbolic space $X=\R H^n$ by a convex cocompact discrete group
of isometries with the option to generalize to more general geometrically
finite hyperbolic manifolds. $\Gamma$ acts by conformal transformations
on the sphere $S^{n-1}=\partial X$ which appears as the geodesic boundary of $X$. Let $\Lambda$ be the limit
set of this action. We consider the space $\Omega^p_{-\infty}(\partial X)$ of complex valued differential $p$-forms on
$X$ with distributional coefficients, i.e., $p$-currents. The conformal group
$G$ acts on it. It fits into the de Rham complex 
\begin{equation}\label{00}
0\rightarrow
\Omega_{-\infty}^0(\partial X)
\stackrel{d}{\longrightarrow}
\Omega_{-\infty}^1(\partial X)
\stackrel{d}{\longrightarrow}
\dots \stackrel{d}{\longrightarrow} 
\Omega_{-\infty}^{n-1}(\partial X)
\stackrel{\int_{\partial X}}{\longrightarrow} \C
\rightarrow 0
\end{equation}
which we have completed by $\C$ in order to make it acyclic.
(\ref{00}) is a complex of $G$-representations. In particular, the center
$\cZ(\gaaa)$ of the universal enveloping algebra $\cU(\gaaa)$ of the Lie algebra
$\gaaa$ of $G$ acts on it. Let $F$ be an irreducible finite-dimensional representation of $G$. Let $(\Omega^p_{-\infty}(\partial X)\otimes F)^{\chi_F}
\subset \Omega^p_{-\infty}(\partial X)\otimes F$ be the largest subrepresentation
on which $\cZ(\gaaa)$ acts by the same character $\chi_F$ as on $F$.
It turns out (see Proposition \ref{trans}) that $(\Omega^p_{-\infty}(\partial X)\otimes F)^{\chi_F}$ is isomorphic to the space of distribution sections
$C^{-\infty}(\partial X,V(F,p))$ of a homogeneous vector bundle $V(F,p)\rightarrow \partial X$.
Thus (\ref{00}) induces the acyclic complex
\begin{equation}\label{01}
0\rightarrow
C^{-\infty}(\partial X,V(F,0))
\rightarrow
C^{-\infty}(\partial X,V(F,1))
\rightarrow
\dots \rightarrow 
C^{-\infty}(\partial X,V(F,n-1))
\rightarrow F
\rightarrow 0
\end{equation}
which appears
in the literature under various names like BGG-resolution or \v Zelobenko complex (see e.g. \cite{zelo76}, \cite{bastoneastwood89}, \cite{juhl01}).
The $G$-representations $C^{-\infty}(\partial X,V(F,p))$, $p=0,\dots,n-1$,
constitute all principal series representations of $G$ with infinitesimal character $\chi_F$.
Note that in case of the trivial representation $F=\C$ the complex (\ref{01})
coincides with (\ref{00}).
By $Z^p_{F,\Lambda}$ we denote the space of such $\Gamma$-invariant $p$-cocycles
of (\ref{01}) which are supported on the limit set $\Lambda$. 
Let $H^*(\Gamma,F)$ be the group cohomology of $\Gamma$ with coefficients in
the $\Gamma$-representation $F$. Patterson conjectured that
for all $p$
\begin{equation}\label{02}
H^p(\Gamma,F)\cong Z^{n-p}_{F,\Lambda}\ .
\end{equation}
$H^p(\Gamma,F)$ is equal to the de Rham cohomology group
$H^p(\Omega^*(Y,E_F))$ of the complex of differential forms on $Y$
with values in the flat vector bundle $E_F\rightarrow Y$ induced by
the $\Gamma$-representation $F$. On the other hand, one can interpret
elements of
$Z^{n-p}_{F,\Lambda}$ as a kind of boundary values of the $\Gamma$-invariant lifts from $Y$ to $X$ of very special closed and coclosed
$E_F$-valued $p$-forms (for noncompact $Y$ the space of all
closed and coclosed $p$-forms is infinite-dimensional).    
Thus (\ref{02}) could be considered as a version of Hodge theory for the present situation.

We are mainly interested in the case of noncocompact $\Gamma$, since
for cocompact $\Gamma$ the validity of (\ref{02}) is a rather direct consequence of classical
Hodge theory. Then the general receipt for producing closed and coclosed forms on $Y$ or, more
generally, eigenforms of the Laplacian $\Delta$ is given by the theory of Eisenstein series
(see \cite{bunkeolbrich000} and the references cited therein). Set $\Omega:=\partial X\setminus\Lambda$ and $B:=\Gamma\backslash \Omega$. $B$ should be considered as the boundary at infinity of $Y$. Let us
for a moment assume that $F=\C$. Then there is a distinguished Eisenstein series
$E^p_s(\phi)\in \Omega^p(Y)$, $s\in\C$, $\phi\in \Omega^{n-p}_{-\infty}(B)$, such that
$\Delta E^p_s(\phi)=((p-\frac{n+1}{2})^2-s^2) E^p_s(\phi)$, $d  E^p_s(\phi)=0$.
We are especially interested in the point $s=p-\frac{n+1}{2}$, which corresponds
to closed harmonic forms. The possible poles of $E^p_s$ at
$s=p-\frac{n+1}{2}$ will prevent us to prove (\ref{02}) in full generality.
Indeed, there are situations, where (\ref{02}) is not true.
Assume that $E^p_s$ has a pole at $s=p-\frac{n+1}{2}$ of order $k>0$.
The boundary value of the leading singular term $E^p(-k,\phi)$ of $E^p_s(\phi)$ at $s=p-\frac{n+1}{2}$ is an $(n-p)$-current with support on the limit set. On the other hand, 
\begin{eqnarray*} 
\lefteqn{d\delta
\left(\frac{d}{ds}_{|s=p-\frac{n+1}{2}}(s-(p-\frac{n+1}{2}))^k E^p_s(\phi)\right)}\\
&=&\frac{d}{ds}_{|s=p-\frac{n+1}{2}}(s-(p-\frac{n+1}{2}))^k 
\Delta E^p_s(\phi)\\  
&=&-\frac{d}{ds}_{|s=p-\frac{n+1}{2}}(s-(p-\frac{n+1}{2}))^{k+1} (s+p-\frac{n+1}{2})E^p_s(\phi)\\
&=&-(k+1)(2p-(n+1))E^p(-k,\phi)\ .
\end{eqnarray*}
Thus $E^p(-k,\phi)$ is exact (if $p\ne \frac{n+1}{2}$). If in addition 
$E^p(-k,\phi)$ is coclosed, then the boundary value of $E^p(-k,\phi)$ 
belongs to $Z^{n-p}_{\C,\Lambda}$. Therefore, (\ref{02}) could
not be true. An example, where such a situation occurs, is discussed at the
end of Section \ref{twist}.

For general $F$ there are Eisenstein series $E^p_{F,s}$ and a special parameter
$s=s_{F,p}$ analogous to $E^p_s$ and $s=p-\frac{n+1}{2}$. 
In order to deal with singularities of the Eisenstein series at $s=s_{F,p}$ we
consider certain successive non-split extensions of the homogeneous bundles $V(F,p)=:V^1(F,p)$
$$ 0\rightarrow V^{k-1}(F,p)\rightarrow V^{k}(F,p)\rightarrow V(F,p)\rightarrow
0\ .$$
By ${}^\Gamma C^{-\infty}(\Lambda,V^{k}(F,p))$ we denote the space of $\Gamma$-invariant
distribution sections of $V^k(F,p)$ which are supported on the limit set.
Set
$$ {}^\Gamma C^{-\infty}(\Lambda,V^{+}(F,p)):=\bigcup_{k\in\nat} {}^\Gamma C^{-\infty}(\Lambda,V^{k}(F,p))\ .$$
We have $Z^{p}_{F,\Lambda}\subset {}^\Gamma C^{-\infty}(\Lambda,V(F,p))\subset
{}^\Gamma C^{-\infty}(\Lambda,V^{+}(F,p))$.
The whole principal part of the Laurent expansion of $E^p_{F,s}(\phi)$ at
$s=s_{F,p}$ defines via the boundary value map an element of 
${}^\Gamma C^{-\infty}(\Lambda,V^{+}(F,p))$. By 
$E^+_\Lambda(F,p)\subset {}^\Gamma C^{-\infty}(\Lambda,V^{+}(F,p))$ 
we denote the space of all these boundary
values. In particular, $E^+_\Lambda(F,p)=\{0\}$ if and only if $E^p_{F,s}$ is regular at
$s=s_{F,p}$. In the case of cocompact $\Gamma$ we set $E^+_\Lambda(F,p):=0$.

Our result can now be formulated as follows 
(see Proposition \ref{bartII}, Theorem \ref{main2} and Corollary \ref{reg}).
It has been announced in \cite{olbrich000}.

\begin{theorem}\label{main00}
\begin{enumerate}
\item The space ${}^\Gamma C^{-\infty}(\Lambda,V^{+}(F,p))$ is finite-dimensional.
\item For $p=1,\dots,n-1$, $p\ne\frac{n+1}{2}$ there is an exact sequence
\begin{equation}\label{es}
0\rightarrow E^+_\Lambda(F,n-p)\rightarrow {}^\Gamma C^{-\infty}(\Lambda,V^{+}(F,n-p))\rightarrow H^p(\Gamma,F)\rightarrow 0\ .
\end{equation}
For $p=\frac{n+1}{2}$ there is a slightly modified exact sequence.
\item If $p\ge\frac{n+1}{2}$, then $E^+_\Lambda(F,n-p)=0$. 
\item If $E^+_\Lambda(F,n-p)=0$, then
$Z^{n-p}_{F,\Lambda}={}^\Gamma C^{-\infty}(\Lambda,V^{+}(F,p))$.\\ 
In particular,
$H^p(\Gamma,F)\cong Z^{n-p}_{F,\Lambda}$ by (\ref{es}).
\item For $p\ge\frac{n}{2}$ we have
$H^p(\Gamma,F)\cong Z^{n-p}_{F,\Lambda}$.
\end{enumerate}
\end{theorem}

Assertion 5 says in particular that (\ref{02}) is true in dimension $n=2$. 
Since (\ref{es}) is a sequence of finite-dimensional spaces
which are related to (generalized) harmonic forms the theorem can be considered as a variant of Hodge
theory. Up to now there is no example for which we can prove that the
map $Z^{n-p}_{F,\Lambda}\rightarrow H^p(\Gamma,F)$ is not surjective. 
We always have $\dim Z^{n-p}_{F,\Lambda}\ge \dim H^p(\Gamma,F)$ (Proposition \ref{irr2}).

As already remarked by Patterson \cite{patterson93} and independently observed
by Lott \cite{lott00}, Assertion 5 for $F=\C$
and $p\ne\frac{n}{2}$ follows from the results of Mazzeo and Phillips \cite{mazzeophillips90} on
comparison of $L^2$-cohomology with cohomology with compact support.
Indeed, also our proof of Assertions 3-5 relies heavily on
$L^2$-methods which do not apply to $p<\frac{n}{2}$. 

The theorem has an obvious generalization for $\Gamma$-representations
of the form $F\otimes V$, where $V$ is finite-dimensional and unitary (see
Theorem \ref{main2}).
We will also consider the
case of non-unitary $V$. In contrast to Assertions 3-5, Assertions 1 and 2 
also hold in this case (for noncocompact $\Gamma$).

Strictly speaking, Theorem \ref{main01} and all the results reviewed in this introduction are proved under the additional assumption
that $\Gamma$ is torsion-free since we like to work in the category
of smooth manifolds. Using the fact that $\Gamma$ always has a torsion-free
normal subgroup of finite index the results can be easily extended to the
general case (compare the discussion at the beginning of Subsection \ref{nancy}).

The reader should be warned that this introduction is the only
place in the present paper, where Eisenstein series are mentioned. Instead
we will use the extension map $ext_s$ which can be considered as a direct
construction of the boundary values of Eisenstein series. 
The Eisenstein series are then
obtained by composing $ext_s$ with a suitable Poisson transform (see \cite{bunkeolbrich000}).  
   
\subsection{Higher cohomology groups}

The exact sequence (\ref{es}) can also be viewed as the determination
of the space of invariants of the $\Gamma$-module of distribution sections
of $V^+(F,p)$ with support on the limit set in terms
of spectral data (boundary values of residues of Eisenstein series) and
topological data (the cohomology groups $H^{n-p}(\Gamma,F)$). 
This space of invariants constitutes
by definition the cohomology group $H^0(\Gamma, C^{-\infty}(\Lambda,V^{+}(F,p)))$. As it is often the case a full understanding
of the situation requires to know the higher cohomology groups
$H^j(\Gamma, C^{-\infty}(\Lambda,V^{+}(F,p)))$, $j>0$, too. Indeed, in \cite{patterson93} Patterson spelled out a second conjecture stating
that the dimensions of the spaces $H^j(\Gamma, C^{-\infty}(\Lambda,V(F,p)))$
are finite and
should be related to the order of the singularity of a twisted Selberg
zeta function at a certain integer point. We will discuss this conjecture in little
more detail at the beginning of Section \ref{sjp}. It turns out that it is more appropriate
to consider the spaces $C^{-\infty}(\Lambda,V^{+}(F,p))$, instead.
Moreover, for technical reasons (see Subsection \ref{connes}) we have to work with the
$\Gamma$-modules $C^{-\omega}(\Lambda,V^{+}(F,p))$ of hyperfunction sections
supported on the limit set. It is relatively easy to see (Theorem \ref{bart})
that
$$ H^0(\Gamma, C^{-\infty}(\Lambda,V^{+}(F,p)))\cong H^0(\Gamma, C^{-\omega}(\Lambda,V^{+}(F,p)))\ .$$
One expects that this isomorphism remains to be true in higher degrees
$j>0$. Up to now this has been proved for cocompact $\Gamma$, only (\cite{bunkeolbrich970}, \cite{bunkeolbrich980}). Results of this kind
strongly depend on the fact that $\Gamma$ does not contain parabolic elements.
While we will say nothing substantial on the Selberg zeta function we
will prove the following result (see Theorem \ref{louvre}).

\begin{theorem}\label{main01}
For all $j\ge 1$, $p\ne\frac{n-1}{2}$, there is a natural isomorphism
$$ H^j(\Gamma, C^{-\omega}(\Lambda,V^{+}(F,p)))\cong H^{n-p+j}(\Gamma, F)\ .$$
For $p=\frac{n-1}{2}$ there is a splitting 
$ V(F,\frac{n-1}{2})=V(F,+)\oplus V(F,-)$. Then
$$ H^j(\Gamma, C^{-\omega}(\Lambda,V^{+}(F,+)))\cong H^j(\Gamma, C^{-\omega}(\Lambda,V^{+}(F,-)))\cong H^{\frac{n+1}{2}+j}(\Gamma, F)\ .$$
In particular, the spaces $H^j(\Gamma, C^{-\omega}(\Lambda,V^{+}(F,p)))$
are finite-dimensional.
\end{theorem}

Again, it is allowed to incorporate twists by finite-dimensional $\Gamma$-representations into the theorem. We also gain information on the
spaces $H^j(\Gamma, C^{-\omega}(\Lambda,V^{k}(F,p)))$. In particular, they are
finite-dimensional. A result
which implies Theorem \ref{main01} in the special case of spherical $F$ and $p=n-1$
was already obtained in \cite{bunkeolbrich990}. 

The bundles $V(F,p)$ sit inside families of irreducible
$G$-homogeneous bundles $V(\sigma_\lambda)$ para\-metrized by irreducible
representations of the isotropy group $P$ of the chosen base point in 
$\partial X$. We call such a representation very special if $V(\sigma_\lambda)$
is one of the bundles $V(F,p)$, $p\ne 0$, $V(F,\pm)$. Moreover, it is natural
to deal with the more general framework of a linear real rank one Lie group
$G$ and a convex cocompact discrete subgroup $\Gamma\subset G$. Again we
have the symmetric space $X$, the limit set $\Lambda$ sitting in its
geodesic boundary $\partial X$, and homogeneous vector bundles $V(\sigma_\lambda)\rightarrow
\partial X$. The spaces $C^{-\omega}(\partial X,V(\sigma_\lambda))$ of their hyperfunction sections constitute the
principal series representations of $G$. There is also a notion of very special parameters $\sigma_\lambda$ (Definition \ref{ere}). Unless stated otherwise,  we will
exclude the exceptional symmetric space $X=\OO H^2$ from the considerations. We  then have $\Gamma$-modules $C^{-\omega}(\Lambda,V^{+}(\sigma_\lambda))$ and a finite-dimensional
subspace $E^+_\Lambda(\sigma_\lambda)\subset {}^\Gamma C^{-\infty}(\Lambda,V^{+}(\sigma_\lambda))=H^0(\Gamma,C^{-\omega}(\Lambda,V^{+}(\sigma_\lambda)))$ generated by the singular part of $ext_s$ at $s=\lambda$. Again it is conjectured that the order of the singularity at $\lambda$
of the Selberg zeta function associated to $\sigma$ is related to the dimensions
of the cohomology groups of $\Gamma$ with coefficients in $C^{-\omega}(\Lambda,V^{+}(\sigma_\lambda))$.

\begin{theorem}\label{main02}
Let $G$ be as above, $\Gamma\subset G$ convex cocompact, and $\sigma_\lambda$
not very special. Then
$$ H^j(\Gamma,C^{-\omega}(\Lambda,V^{+}(\sigma_\lambda)))=\{0\}\quad \mbox{ for all }
j\ge 1\ .$$
If, in addition, $\Gamma$ is not cocompact, then
$$ H^0(\Gamma,C^{-\omega}(\Lambda,V^{+}(\sigma_\lambda)))= E^+_\Lambda(\sigma_\lambda)\ .$$
\end{theorem}

This theorem combines the assertions of Proposition \ref{zwirbel} and Theorem \ref{impo}.
It would be desirable to extend also Theorems \ref{main00} and \ref{main01}
to the general rank one situation. However, for complex and
quaternionic hyperbolic spaces the structure of analogous results
is expected to be much more involved.    

The proof of each of the above theorems is based on three essential
ingredients: the surjectivity of Laplace-like operators on smooth
sections of vector bundles over connected noncompact Riemannian
manifolds, a geometric version of scattering theory (developed in
joint work with Ulrich Bunke mainly in \cite{bunkeolbrich000}), where the above
mentioned  operator $ext_s$ plays an essential role, and the theory
of Poisson transforms for homogeneous vector bundles over symmetric
spaces (which is essentially the theory of asymptotic expansions
of matrix coefficients of admissible $G$-representations).
The latter two ingredients allow us to follow a strategy which became
more and more popular and promising during the last years, 
even in the classical case of finite volume quotients $Y=\Gamma\backslash X$, namely,
study automorphic forms via the behaviour of their boundary distributions
(compare \cite{schmid00}, \cite{bernsteinreznikov98}, \cite{bruggeman97}, \cite{lewiszagier01}, 
\cite{lott00}, \cite{miatellowallach89}, \cite{Gai1}, \cite{juhl88}, and the joint papers of the author with U. Bunke
\cite{bunkeolbrich980}, \cite{bunkeolbrich991}, \cite{bunkeolbrich000}, \cite{bunkeolbrich01}).   

\subsection{The structure of the paper and side results}\label{connes}

Section \ref{sur} contains a couple of remarks on the surjectivity of the
$p$-form Laplacian on noncompact connected Riemannian manifolds.
We obtain a weak form of Hodge theory which in particular says that any
de Rham cohomology class has a coclosed harmonic representative.
These facts, though based on classical results, do not seem to have been
systematically studied so far. See, however, the note \cite{gaillard00}.

Sections \ref{hyper} and \ref{round} appear as a longish
commentary on \cite{bunkeolbrich000} and parts of \cite{bunkeolbrich990}.
Since we only know the surjectivity of the Laplacian on the space of
all smooth forms (or smooth sections of more general vector bundles)
and not on the space of forms of moderate growth (see the discussion in Subsection \ref{nancy}) we
are forced to redevelop the theory of the extension operator $ext_s$
for convex cocompact groups acting on symmetric spaces of negative curvature in a hyperfunction setting (note that the boundary value of an arbitrary eigensection is a hyperfunction while eigensections of moderate growth have distribution boundary values). This is entirely parallel to the corresponding
theory for distributions treated in \cite{bunkeolbrich000}.
In Section \ref{round} we introduce and study the spaces   
${}^\Gamma C^{-\infty}(\Lambda,V^{+}(\sigma_\lambda))$ and $E^+_\Lambda(\sigma_\lambda)$. In particular, we show that they
are finite-dimensional. For $\Ree(\lambda)\ge 0$ we draw some consequences of
unitarity which go beyond the results of \cite{bunkeolbrich000}, Section 7.
In particular, we obtain the following new results which are of independent 
interest:
\begin{enumerate}
\item[(i)] The operator $ext_\lambda$ is regular at such non-negative $\lambda$
which correspond to regular integral infinitesimal character (Proposition
\ref{newer}). This in particular implies Assertion 3 of Theorem \ref{main00}.
\item[(ii)] If $X=\R H^n$ and $\sigma$ is a faithful representation of $Spin(n-1)$,
then ${}^\Gamma C^{-\infty}(\Lambda,V^{+}(\sigma_0))=\{0\}$ (Proposition \ref{new}).
\item[(iii)] If $X=\HH H^n$ or $X=\OO H^2$ and $\Gamma$ is not cocompact, then
the Hausdorff dimension of $\Lambda$ is strictly less then $4n$ or $16$,
respectively (Corollary \ref{corlette}). This slightly improves a result of Corlette \cite{corlette90}.
Note that in these cases the Hausdorff dimension of $\partial X$ is
equal to $4n+2$ and $22$, respectively.
\end{enumerate}

The reader who is mainly interested in the consequences of these
results for spectral theory and $L^2$-cohomology can directly jump to
Section \ref{mazzeo}. There we give a slightly refined version of the
Plancherel Theorem for $L^2(\Gamma\backslash G)$ obtained in \cite{bunkeolbrich000}. For instance, thanks to (ii) we now know that in the real hyperbolic case limits of discrete series representations do not contribute to the discrete spectrum. Based on the Plancherel Theorem we compute the $L^2$-cohomology of $Y=\Gamma\backslash X$ with coefficients in the flat bundle $E_F$ induced by a finite-dimensional
$G$-representation by standard relative Lie algebra cohomology
methods (Theorem \ref{l2co}). Again, the result will be mainly in terms of invariant distributions supported on the limit set. It follows that these cohomology
groups are finite-dimensional except for the middle degree.

In Sections \ref{heart} and \ref{twist} we have to assume that $X=\R H^n$.
They culminate in the proof of Theorem \ref{main01}. The reason for the
restriction to the real hyperbolic case is that for the remaining cases
the theory of Poisson transforms for harmonic differential forms is neither 
developed sufficiently far nor would it be possible to use it in order to
translate the Hodge theoretic results of Section \ref{sur} directly to the boundary
(e.g., an appropriate Hodge theory for quotients of the complex hyperbolic
space should involve a filtration by bidegrees). For the real hyperbolic
case the theory of Poisson transforms of differential forms is well understood \cite{Gai1}. In Section \ref{heart} we adjust
this theory according to our needs. Combining it with the results of
Section \ref{round} we obtain a proof of Theorem \ref{main01} for the
special case $F=\C$ (except for Assertion 5 in case $p=\frac{n}{2}$), but with arbitrary twists allowed. The general case
follows by an application of the translation functor (Section \ref{twist})
which is well-known in representation theory. The proof
of Assertion 5 in case $p=\frac{n}{2}$ is then a consequence of the theory
of discrete series representations combined with various scalar product
formulas as e.g. Proposition \ref{highgreen}.

In Subsection \ref{plum} we derive some vanishing results for cohomology
based on Theorem \ref{main00}. In Subsection \ref{supp} we study of the 
relation of the spaces 
$Z^{p}_{F,\Lambda}$ to cohomology with compact support. 

For any Riemannian locally symmetric
space there is the notion of automorphic forms which involves the condition of moderate growth. Gaillard \cite{gaillard01} formulated a conjecture which
he calls the Harder-Borel conjecture and which states that
the subcomplexes of the de Rham complex consisting of coclosed harmonic
automorphic forms and of automorphic forms, respectively, are quasi-isomorphic
to the full de Rham complex. The problem is widely investigated in the special and most interesting case of finite
volume spaces. In our situation ($G$
of rank one, $\Gamma\subset G$ convex cocompact)
the conjecture would immediately follow if we would know that the Laplacian acting on forms
of moderate growth is surjective. Such an assertion, however, seems to
be very difficult to prove. Theorem \ref{main00} provides
representatives of cohomology classes by (boundary values of) very special
automorphic forms. This does not imply the conjecture yet. In Subsection \ref{nancy} we provide some additional arguments which together with
Assertion 2 of Theorem \ref{main00} show that the Harder-Borel conjecture is
true for quotients of the real hyperbolic space by a convex cocompact subgroup.

Section \ref{sjp} is devoted to the proofs of Theorem \ref{main02} and
Theorem \ref{main01}. Here we make use of several ideas already employed
in \cite{bunkeolbrich950} and \cite{bunkeolbrich990}. The main step is to compute the higher cohomology
groups of $C^{-\omega}(\partial X,V^{+}(\sigma_\lambda))$. The first
observation we make is that generalized eigenspaces of Laplacians
on homogeneous vector bundles over $X$ are acyclic $\Gamma$-modules (Lemma \ref{ac}). Here again the surjectivity of elliptic operators on analytic noncompact manifolds
(see Section \ref{sur}) becomes crucial. Then we show that for parameters
$\sigma_\lambda$ which are not very special suitably chosen
Poisson transforms map $C^{-\omega}(\partial X,V^{+}(\sigma_\lambda))$
isomorphically to a direct summand of such a generalized eigenspace.
Especially for $\lambda=0$ this requires a thorough discussion of various
cases. It follows that all higher cohomology groups of $C^{-\omega}(\partial X,V^{+}(\sigma_\lambda))$ vanish. If $X=\R H^n$ and
$\sigma_\lambda$ is very special, then the theory developed in Sections \ref{sur}, \ref{heart}, and \ref{twist} implies that certain de Rham 
complexes of generalized harmonic forms provide acyclic resolutions of  
$C^{-\omega}(\partial X,V^{+}(\sigma_\lambda))$. This establishes the connection
to the cohomology groups $H^p(\Gamma,F)$.  
In the final step we use the exact sequence
$$0\rightarrow C^{-\omega}(\Lambda,V^{+}(\sigma_\lambda))\rightarrow
C^{-\omega}(\partial X,V^{+}(\sigma_\lambda))\rightarrow 
C^{-\omega}(\Omega,V^{+}(\sigma_\lambda))\rightarrow 0 $$
in order to conclude that for all $p\ge 1$
$$ H^p\left(C^{-\omega}(\Lambda,V^{+}(\sigma_\lambda))\right)\cong 
H^p\left(C^{-\omega}(\partial X,V^{+}(\sigma_\lambda))\right)\ .$$
For $p=1$ this uses the meromorphy of $ext_s$.

We expect that Theorem \ref{main00} has a natural generalization to arbitrary
geometrically finite groups $\Gamma$ acting on $\R H^n$. Note, however,
that in this generality the spaces ${}^\Gamma C^{-\infty}(\Lambda,V(F,p))$ and ${}^\Gamma C^{-\infty}(\Lambda,V^+(F,p))$ can be infinite-dimensional. 
In order to obtain finite-dimensional spaces one should replace the
condition ``supported on the limit set'' by ``strongly supported on the limit
set'' (see \cite{bunkeolbrich01}). Indeed, the meromorphy
of $ext_s$ in the distribution setting has been established in \cite{bunkeolbrich991}, and the singular part of $ext$ produces invariant distributions which are ``strongly supported on the limit
set''. However, it is far from being obvious how to transfer
this theory to a hyperfunction setting.

In the special case of finite volume quotients there is an alternative
approach. In this case $\Lambda=\partial X$ and $\dim {}^\Gamma C^{-\infty}(\partial X,V(F,p))<\infty$. Thus there is no need
to ``move'' the (strong) support of the boundary value of a harmonic form to the
limit set and therefore things can be done without considering the
extended bundles $V^k(F,p)$ and $V^+(F,p)$. Moreover, in this case
Hodge theory on forms of moderate growth
is established (see the discussion in Subsection \ref{nancy}, in particular Proposition \ref{smom}). Therefore one can work from the very beginning
in the distribution setting. Arguing as in the proof of Corollary \ref{preth}
we obtain

\begin{prop}\label{last}
Let $\Gamma$ be such that $Y=\Gamma\backslash \R H^n$ has finite volume.
Let $F$ be a finite-dimensional $G$-representation as above. Let $E(F,p)
\subset {}^\Gamma C^{-\infty}(\partial X,V(F,p))$
be the subspace spanned by the boundary values of the (not necessarily
singular) leading terms of  
the Eisenstein series $E_{F,s}^p(\phi)$ at $s=s_{F,p}$ (see \ref{least}).
By $Z^p_F$ we denote the space of $\Gamma$-invariant cocycles of (\ref{01}).
Then there are exact sequences
\begin{eqnarray*} 
0&\rightarrow& E(F,p)\rightarrow {}^\Gamma C^{-\infty}(\partial X,V(F,p))
\rightarrow H^p(\Gamma,F)\rightarrow 0\ ,\quad p\ne\frac{n+1}{2}\ ,\\
0&\rightarrow& E(F,p)\cap Z^p_F\rightarrow Z^p_F
\rightarrow H^p(\Gamma,F)\rightarrow 0\ .
\end{eqnarray*}
\end{prop}
\hspace{3cm}

\noindent
{\it Acknowledgements}: First of all I want to thank Ulrich Bunke.
Most of the results of this paper strongly depend on ideas developed
in previous joint work with him. Stimulating discussions with him
accompanied also the work on the present paper. I am indebted to S. J.
Patterson whose conjectures and insightful remarks were a kind of orientation
guide not only for this work. Discussions with A. Juhl, who also carefully
read previous versions of the manuscript, and P.-Y. Gaillard
substantially influenced my view on the subject, too.

\newpage
\section[Hodge theory]{A residue of Hodge theory for noncompact Riemannian manifolds}\label{sur}

Given a smooth manifold $Y$ and a flat finite-dimensional complex vector 
bundle $(E,\nabla)$ on it
we can consider the associated de Rham complex 
$$ (\Omega^*(Y,E), d)$$
of smooth, $E$-valued differential forms, i.e., smooth sections of the bundles
$\Lambda^*(T^*Y)\otimes E$. Here the differential 
$$ d:\Omega^p(Y,E)\rightarrow \Omega^{p+1}(Y,E)\ .$$
is induced by the flat connection 
$$\nabla: \Omega^0(Y,E)\rightarrow \Omega^1(Y,E)\ .$$ 
We are interested in its cohomology groups $H^p(Y,E)$. The most important case is the one-dimensional trivial bundle equipped with the trivial connection. 
The cohomology
of the corresponding de Rham complex then becomes canonically isomorphic
to the usual (say singular) cohomology $H^*(Y,\C)$ of $Y$ with coefficients in $\C$. For us a manifold always has a countable base of the topology. For convenience we assume $Y$ to be connected.

Let $g$ be a Riemannian metric on $Y$, and let $h$ be a Hermitian metric 
on $E$. We do not require $h$ to be parallel with respect to $\nabla$.
These structures induce the Riemannian measure $\mu_g$ on $Y$ and Hermitian forms
$(.,.)$ on $\Lambda^p(T^*Y)\otimes E$. Then we can form the formal adjoint of $d$
$$\delta: \Omega^{p+1}(Y,E)\rightarrow \Omega^p(Y,E)$$
which is characterized by 
$$ \int_Y (\delta \omega_y,\eta_y)\:d\mu(y) 
=\int_Y ( \omega_y,d\eta_y)\:d\mu(y)\quad \mbox{ for all } \omega\in\Omega^{p+1}(Y,E),\eta\in \Omega_c^p(Y,E)\ .$$ 
Here the subscript $c$ means compact support. We have $\delta^2=0$. The corresponding $p$-form Laplacian is given by
$$ \Delta=\delta d+d\delta: \Omega^{p}(Y,E)\rightarrow \Omega^p(Y,E)\ .$$
$d$ and $\delta$ commute with the elliptic and formally
selfadjoint operator $\Delta$.

These operators give rise to three distinguished sub{\em complexes}
$$ \Omega^*(Y,E)_{\Delta,\delta}\subset\Omega^*(Y,E)_{\Delta}\subset\Omega^*(Y,E)_{(\Delta)} $$
of the de Rham complex $(\Omega^*(Y,E),d)$, namely
\begin{itemize}
\item $\Omega^p(Y,E)_{\Delta,\delta}:=\{\omega\in \Omega^{p}(Y,E)\:|\: \Delta\omega=0, \delta\omega=0\}\ $ (coclosed harmonic forms),
\item $\Omega^p(Y,E)_\Delta:=\{\omega\in \Omega^{p}(Y,E)\:|\: \Delta\omega=0\}\ $
(harmonic forms),
\item $\Omega^p(Y,E)_{(\Delta)}:=\{\omega\in \Omega^{p}(Y,E)\:|\: \mbox{ there exists } k\in\nat\mbox{ s.th. }\Delta^k\omega=0\}\ $
(generalized harmonic forms).
\end{itemize}
If $Y$ is compact, then by classical Hodge theory all these complexes
coincide with $\ker d\cap \ker\delta$, are finite-dimensional and provide
canonical representatives of cohomology classes by harmonic forms.
For noncompact $Y$ all these statements are false. However, one of the main results of the present section is that in general at least the following is true.

\begin{theorem}\label{ho}
The inclusions of complexes
$$\Omega^*(Y,E)_{\Delta,\delta}\hookrightarrow \Omega^*(Y,E)_{(\Delta)}\hookrightarrow\Omega^*(Y,E)$$
are quasi-isomorphisms, i.e., they induce isomorphisms in cohomology.
\end{theorem}

For noncompact manifolds the inclusion $\Omega^*(Y,E)_\Delta\stackrel{i}{\longrightarrow}\Omega^*(Y,E)$ is far from being
a quasi-isomorphism in general. In fact, one can show that
$i\oplus (-1)^*\delta:  \Omega^*(Y,E)_{\Delta}\longrightarrow \Omega^*(Y,E)\oplus\Omega^{*-1}(Y,E)$ is a quasi-isomorphism, i.e., for any $p\in\nat_0$ there is a canonical isomorphism
$$ H^p(\Omega^*(Y,E)_{\Delta})\cong H^p(Y,E)\oplus H^{p-1}(Y,E)\ .$$

The theorem is a simple consequence of the surjectivity of the Laplacian on $\Omega^p(Y,E)$
for noncompact manifolds,
a classical result which goes back to Malgrange \cite{malgrange56}
but which does not seem to have received much attention among global analysts. 
We will discuss this result in a moment. It has been used in several joint
papers of U. Bunke and the author (\cite{bunkeolbrich950}, \cite{bunkeolbrich990}). The validity of Theorem \ref{ho} was independently observed by P.-Y. Gaillard \cite{gaillard00} to whome
I am also indebted for providing adequate references.

First we need the following

\begin{ddd}
Let $E_1$, $E_2$ be vector bundles over $Y$. We say that a differential operator
$$ D: C^{\infty}(Y,E_1)\rightarrow C^{\infty}(Y,E_2) $$
has the weak unique continuation property if for any connected open subset 
$U\subset Y$, any section $f\in C^{-\infty}(U,E_1)$ with $Df=0$ the condition
$f_{|U_0}=0$ for a nonempty open subset $U_0\subset U$ implies $f\equiv 0$ on $U$.
\end{ddd}

If $Y$, $E_1$, $E_2$ are analytic, and $D$ is an elliptic operator with analytic
coefficients, then $D$ has the weak unique continuation property by analytic elliptic regularity (see e.g. \cite{hoermander831}, Thm. 8.6.1 or \cite{kashiwarakawaikimura86}, Thm. 3.4.4). In fact, in this case $Df=0$ implies
that $f$ is analytic. But for operators like $\Delta$ this analyticity condition
is not necessary for unique continuation. Indeed, in the setting introduced above the form Laplacians
$\Delta: \Omega^{p}(Y,E)\rightarrow \Omega^p(Y,E)$
have the weak unique continuation property (\cite{aronszajnkrzywickiszarski62}).
For an overview concerning unique continuation theorems we refer to \cite{kenig86}.
We now have

\begin{theorem}[\cite{malgrange56}, p.341]\label{masur}
Let $Y$ be a connected noncompact smooth manifold, and let 
$$D: C^{\infty}(Y,E_1)\rightarrow C^{\infty}(Y,E_2)$$ 
be an elliptic
operator which has the weak unique continuation property. Then  
$D$ is surjective.
\end{theorem}
In view of the above discussion this implies

\begin{kor}\label{susi}
The $p$-form Laplacian
$\Delta: \Omega^{p}(Y,E)\rightarrow \Omega^p(Y,E)$ on a noncompact connected
smooth manifold  
is surjective.
\end{kor}

\noindent
{\it Proof of Theorem \ref{ho}.$\:\:\:\:$} 
For compact manifolds the theorem is covered by classical Hodge theory (see e.g. \cite{wells73}, Sect. IV.5).
Thus we may assume that $Y$ is noncompact. We first show that the inclusion $\Omega^*(Y,E)_{\Delta,\delta}\hookrightarrow \Omega^*(Y,E)$ induces a surjection in cohomology. Let $\omega\in \Omega^p(Y,E)$
be closed. By Corollary \ref{susi} we find $\eta\in \Omega^p(Y,E)$ solving
the equation $\Delta \eta=\omega$.
Set 
$$\omega_0:=\omega-d\delta\eta=\delta d\eta\ .$$
Then $\omega_0$ is cohomologous to $\omega$ and coclosed. In particular, $\omega_0\in\Omega^p(Y,E)_{\Delta,\delta}$. This proves the desired surjectivity.

For injectivity we consider $\phi\in\Omega^{p-1}(Y,E)$, $\omega=d\phi$.
Solving $\Delta\eta=\phi$ in $\Omega^{p-1}(Y,E)$ we can form $\phi_0:=\phi-d\delta\eta=\delta d\eta$. Then $d\phi_0=\omega$, $\delta\phi_0=0$, $\Delta\phi_0=\delta d\phi_0=\delta\omega$. 

If $\omega\in \Omega^p(Y,E)_{\Delta,\delta}$, then 
$\Delta\phi_0=\delta\omega=0$, i.e., $\phi_0\in\Omega^{p-1}(Y,E)_{\Delta,\delta}$.
Thus in this case $\omega$ is a coboundary in $\Omega^*(Y,E)_{\Delta,\delta}$.
If $\Delta^k\omega=0$ for some $k\ge 1$, then
$$\Delta^{k+1}\phi_0=\Delta^k\delta\omega=\delta \Delta^k\omega=0\ .$$
Hence $\omega$ is a coboundary in $\Omega^*(Y,E)_{(\Delta)}$. Thus we have
shown that both injections of complexes induce injective maps in cohomology.
The proof is now complete.
\hB

There is an alternative way of expressing
cohomology in terms of generalized harmonic forms which will 
turn out to be crucial in the proof of the main theorem of the present paper.
By $Z^p(Y,E)_{(\Delta)}$ we denote the space of $p$-cocycles in $\Omega^*(Y,E)_{(\Delta)}$.

\begin{prop}\label{infi}
A $p$-form $\omega\in \Omega^p(Y,E)$ is exact if and only
if there exists a closed $p$-form $\eta$ such that $\Delta\eta=\omega$.
In particular we have
\begin{equation}\label{ik} H^p(Y,E)\cong Z^p(Y,E)_{(\Delta)}/\Delta\left(Z^p(Y,E)_{(\Delta)}\right)\cong \ker d\cap \ker\delta/ \Delta\left(Z^p(Y,E)_{(\Delta)}\right)\cap \ker\delta\ .
\end{equation}
Moreover, 
\begin{equation}\label{ak}
\Delta\left(Z^p(Y,E)_{(\Delta)}\right)=\Delta^k\left(Z^p(Y,E)_{(\Delta)}\right)
\qquad \mbox{for any }k\in\nat \ .
\end{equation}
\end{prop} 
\proof
For compact manifolds the proposition is a consequence of classical Hodge theory. Thus we can assume that $Y$ is noncompact. If $\eta$ is closed, then $\Delta\eta=d\delta\eta$ is exact.
On the other hand, for exact $\omega=d\phi$ we use Corollary \ref{susi} in order to solve the equation 
$\Delta \psi=\phi$ 
and set $\eta=d\psi$. Then
$$ \Delta\eta=d\Delta\psi=\omega\ .$$ 
This proves the first assertion. 
Since according to Theorem \ref{ho} any cohomology class can be represented by a generalized harmonic form, or even by a closed and coclosed form,
Equation (\ref{ik}) follows, too. Let now $\omega=d\delta \eta\in \Delta\left(Z^p(Y,E)_{(\Delta)}\right)$. By Corollary \ref{susi} we can solve $ \Delta^k \psi=\delta \eta$. Then $\omega=\Delta^k d\psi\in \Delta^k\left(Z^p(Y,E)_{(\Delta)}\right)$. This proves (\ref{ak}) and finishes the proof of the proposition.
\hB

There is another canonical codifferential $\hat\delta$ acting on $\Omega^*(Y,E)$
which only depends on the Riemannian metric $g$ but not on any chosen
Hermitian metric on $E$. Indeed, let
$$ \delta_1: \Omega^p(Y)\rightarrow\Omega^{p-1}(Y)$$
be the usual codifferential for the trivial bundle. Then we define
$$ \hat\delta : \Omega^p(Y,E)\cong C^\infty(Y,\Lambda^pT^*Y\otimes E)\rightarrow \Omega^{p-1}(Y,E)$$
by
\begin{equation}\label{Hut}
\hat\delta (\omega\otimes e):=\delta_1 \omega\otimes e- i_g(\nabla e)\omega\ , \end{equation}
where 
$$ i_g: \Omega^1(Y,E)\otimes \Omega^p(Y)\rightarrow \Omega^{p-1}(Y,E) $$
is the insertion operator defined by the Riemannian metric $g$.
Again we have $\hat\delta^2=0$. $\hat\delta$ coincides with $\delta$ if and only if the Hermitian metric
$h$ is parallel. We form
$$ \hat\Delta=\hat\delta d+d\hat\delta: \Omega^{p}(Y,E)\rightarrow \Omega^p(Y,E)\ .$$
$\hat\Delta$ is elliptic and commutes with $d$ and $\hat\delta$.
But in general there is no $L^2$-scalar product on $\Omega^{p}(Y,E)$ such
that $\hat\Delta$ becomes formally selfadjoint. Replacing in the above
definitions $\Delta$ by $\hat\Delta$ and $\delta$ by $\hat\delta$ we obtain new subcomplexes of $(\Omega^*(Y,E),d)$   
$$ \Omega^*(Y,E)_{\hat\Delta,\hat\delta}\subset\Omega^*(Y,E)_{\hat\Delta}\subset\Omega^*(Y,E)_{(\hat\Delta)}\ . $$

While classical Hodge theory does not work if $\hat\Delta$ is not selfadjoint the same arguing as above yields

\begin{kor}\label{ellen}
The Laplacian
$\hat\Delta: \Omega^{p}(Y,E)\rightarrow \Omega^p(Y,E)$ on a noncompact connected
smooth manifold  
is surjective.
\end{kor}

Thus the same proofs as above show that for
noncompact manifolds we are allowed to put a hat on the operators appearing in Theorem \ref{ho} and Proposition \ref{infi}.

\begin{theorem}\label{hop}
Let $Y$ be noncompact. Then the inclusions of complexes
$$\Omega^*(Y,E)_{\hat\Delta,\hat\delta}\hookrightarrow \Omega^*(Y,E)_{(\hat\Delta)}\hookrightarrow\Omega^*(Y,E)$$
are quasi-isomorphisms, i.e., they induce isomorphisms in cohomology.
\end{theorem}

\begin{prop}\label{infip}
Let $Y$ be noncompact. 
A $p$-form $\omega\in \Omega^p(Y,E)$ is exact if and only
if there exists a closed $p$-form $\eta$ such that $\hat\Delta\eta=\omega$.
In particular we have
\begin{equation}\label{ikp} H^p(Y,E)\cong Z^p(Y,E)_{(\hat\Delta)}/\hat\Delta\left(Z^p(Y,E)_{(\hat\Delta)}\right)\cong \ker d\cap \ker\hat\delta/ \hat\Delta\left(Z^p(Y,E)_{(\hat\Delta)}\right)\cap \ker\hat\delta\ .
\end{equation}
Moreover, 
\begin{equation}\label{akp}
\hat\Delta\left(Z^p(Y,E)_{(\hat\Delta)}\right)=
\hat\Delta^k\left(Z^p(Y,E)_{(\hat\Delta)}\right)
\qquad \mbox{for any }k\in\nat \ .
\end{equation}
\end{prop}

For compact manifolds Theorem \ref{hop} and Proposition \ref{infip} are not
true in general. However, we have
\begin{prop}\label{comp}
Let $Y$ be compact. There exists a $k\in\nat$ such that $\Omega^p(Y,E)_{(\hat\Delta)}=
\ker \hat\Delta^k$ for all $p$. In particular, $\dim \Omega^*(Y,E)_{(\hat\Delta)}<\infty$. There is a decomposition of complexes
\begin{equation}\label{spl}
\Omega^*(Y,E)=\Omega^*(Y,E)_{(\hat\Delta)}\oplus \hat\Delta^k\left(\Omega^*(Y,E)\right)\ .
\end{equation}
The embedding $\Omega^*(Y,E)_{(\hat\Delta)}\hookrightarrow\Omega^*(Y,E)$ is a
quasi-isomorphism.
\end{prop}
\proof We choose a Hermitian metric $h$ on $E$. Then we can form the codifferential $\delta$, the Laplacians $\Delta$ as above, and the Hilbert space completion $L^2(Y,\Lambda^pT^*Y\otimes E)$ of $\Omega^p(Y,E)$. The difference
$\hat\delta-\delta$ is an operator of order zero, hence $T:=\hat\Delta-\Delta$
has at most first order. $\Delta$ becomes an unbounded selfadjoint operator
with the Sobolev space $H^2(Y,\Lambda^pT^*Y\otimes E)$ as its domain of definition. 
$T$ as well as $\hat\Delta=\Delta+T$ are defined on $H^2(Y,\Lambda^pT^*Y\otimes E)$,
and $T(1+\Delta)^{-1}$ defines a compact operator on  $L^2(Y,\Lambda^pT^*Y)$. As a relative compact perturbation of $\Delta$ the
operator $\hat\Delta$ inherits from $\Delta$ the property that 
\begin{equation}\label{glik} 
\dim \bigcup_l\ker (\hat\Delta^l)<\infty\ ,\qquad 
L^2(Y,\Lambda^pT^*Y)=\bigcup_l\ker (\hat\Delta^l)\oplus \bigcap_l\im (\hat\Delta^l)\ 
\end{equation}
(see \cite{gohbergkrein69}, Ch. V, Lemma 10.1). In fact, this becomes a simple consequence
of the compactness of the resolvent $(\hat\Delta-\lambda)^{-1}$ as soon one knows that the resolvent set of $\hat\Delta$ is not empty which can be ensured
by
$$ \lim_{\lambda\to i\infty} \|T(\Delta-\lambda)^{-1})\|=0\qquad 
\mbox{(\cite{gohbergkrein69}, Ch. V, Lemma 7.2) .}$$
By elliptic regularity $\bigcup_l\ker (\hat\Delta^l)\subset \Omega^p(Y,E)$.
Equation (\ref{spl}) now follows from (\ref{glik}).

As in the proof of Theorem \ref{ho} one shows that any cohomology class in
$H^p(\hat\Delta^k\left(\Omega^*(Y,E)\right)$ has a $\hat\Delta$-harmonic representative. Now (\ref{spl}) implies that it is zero. This finishes the proof of the proposition.
\hB

\newpage
\section[Hyperfunction scattering]{Scattering theory for convex cocompact groups: the hyperfunction setting}\label{hyper}

The main purpose of this section is to transfer the results of \cite{bunkeolbrich000}, Sections 4-5, from distributions to hyperfunctions.
First we have to recall the setting of that paper.

Let $G$ be a connected, linear, real simple Lie group of rank one. We fix a
maximal compact subgroup $K\subset G$ and an Iwasawa decomposition  $G=KAN$ 
of $G$. Note that the rank one assumption is equivalent to $\dim A=1$. Let
$M:=Z_K(A)$ be the centralizer of $A$ in $K$. We denote the Lie algebras
of these groups by $\gaaa$, $\kaaa$, $\aaaa$, $\naaa$, and $\maaa$, respectively. We form 
the corresponding minimal parabolic subgroup $P:=MAN$ of $G$. 
The group $G$ acts isometrically and orientation-preserving on the rank-one symmetric space $X:=G/K$.
Let $\partial X:=G/P=K/M$ be its geodesic boundary. We consider $\bar{X}:=X\cup\partial X$ as a compact manifold with boundary.

Let $\Gamma\subset G$ be a discrete torsion-free subgroup. Its
limit set $\Lambda\subset \partial X$ is defined to be the set of accumulation
points in $X\cup \partial X$ of the $\Gamma$-orbit of an arbitrary point $x\in X$. The complement
$\Omega:=\partial X\setminus\Lambda$ is called the domain of discontinuity
of $\Gamma$. Indeed, $\Gamma$ acts properly discontinuously on $\Omega$ as
well as on $X\cup \Omega$. Hence, $\overline Y:=\Gamma\backslash (X\cup \Omega)$
is a manifold with boundary $B:=\Gamma\backslash \Omega$. Its interior
$Y:=\Gamma\backslash X$ carries the structure of a locally symmetric space of
negative sectional curvature.
We assume $\Gamma$ to be convex cocompact which means by definition that $\overline Y$ is compact. In particular,
a cocompact subgroup is convex cocompact. In this case $\Lambda=\partial X$,
$\Omega=\emptyset$, and $\overline Y=Y$ is a compact locally symmetric space.
We are mostly interested in the opposite case $\Omega\ne\emptyset$.
Then $Y$ is a locally symmetric space of infinite volume without cusps which
we will call a Kleinian manifold generalizing
the corresponding notion for three-dimensional hyperbolic manifolds.

By $\aca$ we denote the complexified dual of $\aaaa$.
For $a\in A$ and $\lambda\in \aca$ we set $a^\lambda:=\ee^{\langle \lambda,\log(a)\rangle}\in\C$. Let $\alpha$ be the short root of
$\aaaa$ in $\naaa$. We set 
$A_+:=\{a\in A\:|\:a^\alpha\ge 1\}$.
Define $\rho\in \aaaa^*$ as usual by $\rho(H):=\frac{1}{2}\tr(\ad(H)_{|\naaa})$, $\forall H\in\aaaa$. We have \vspace{0.2cm}\\  
\centerline{\begin{tabular}{|c||c|c|c|c|}
\hline
$X$&$\R H^n$&$\C H^n$&$\HH H^n$&$\OO H^2$\\
\hline
$\rho$&$\frac{n-1}{2}\alpha$&$n\alpha$&$(2n+1)\alpha$&$11\alpha$\\
\hline
\end{tabular} . }\vspace{0.3cm}
Any element $g\in G$ has a Cartan decomposition $g=k_ga_gh$, $k_g,h\in K$, $a_g\in A_+$, where $a_g$ and $k_gM\in K/M$ are uniquely determined by $g$.

\begin{ddd}\label{dodo}
The critical exponent $\delta_\Gamma\in \aaaa^*$ of $\Gamma$ is the smallest element such that the series $\sum_{g\in\Gamma} a_g^{-(\lambda+\rho)}$ converges for all $\lambda\in\aaaa^*$ with $\lambda>\delta_\Gamma$. If $\Gamma$ is the trivial group,
then we set $\delta_\Gamma:=-\infty$.
\end{ddd}

If $(\tau,V_\tau)$ is a finite-dimensional representation of $P$, then we denote by $V(\tau):=G\times_P V_{\tau}$ the associated
homogeneous vector bundle over $\partial X=G/P$. 
It induces a bundle on $B=\Gamma\backslash\Omega$ defined by $V_B(\tau):=\Gamma\backslash V(\tau)_{|\Omega}$.
These bundles are defined in the analytic category. Therefore we
can consider not only
smooth and distribution sections but also analytic and hyperfunction sections
of these bundles. Of main interest are the $P$-representations $\sigma_\lambda$ which arise as follows.
Let $(\sigma,V_\sigma)$ be a finite-dimensional unitary representation of $M$. For $\lambda\in \aca$  we  form the representation $\sigma_\lambda$
of $P$ on $V_{\sigma_\lambda}:=V_\sigma$, which is given by
$\sigma_\lambda(man):=\sigma(m)a^{\rho-\lambda}$.
There is a chain of inclusions of continuous $G$-representations
(called principal series representations) on complete locally convex Hausdorff topological vector spaces
$$C^{\omega}(\partial X,V(\sigma_\lambda))\subset 
C^{\infty}(\partial X,V(\sigma_\lambda))\subset 
C^{-\infty}(\partial X,V(\sigma_\lambda))\subset
C^{-\omega}(\partial X,V(\sigma_\lambda))$$
which we all denote by $\pi^{\sigma,\lambda}$. There is a corresponding
chain of inclusions of spaces of sections over $B$. Let $\tilde{\sigma}$ be the dual representation to $\sigma$. Since $\Lambda^{max} T^*\partial X\cong V(1_{-\rho})$ and $\partial X$ and $B$
are compact we have  
\begin{eqnarray*}
C^{-\sharp}(\partial X,V(\sigma_\lambda))&=&C^{\sharp}(\partial X,V(\tilde{\sigma}_{-\lambda}))^\prime\\
C^{-\sharp}(B,V_B(\sigma_\lambda))&=&C^{\sharp}(B,V_B(\tilde{\sigma}_{-\lambda}))^\prime\ ,\qquad \sharp\in\{\infty,\omega\}\ .
\end{eqnarray*}
As explained in \cite{bunkeolbrich000}, p. 86, we can equip
$$\bigcup_{\lambda\in\aca} V(\sigma_\lambda)\rightarrow \aca\times \partial X
\quad 
\mbox{ and }\quad \bigcup_{\lambda\in\aca} V_B(\sigma_\lambda)\rightarrow \aca\times B $$
with the structure of holomorphic families of bundles. 

The structure of a holomorphic family of bundles allows us to consider holomorphic and meromorphic families
of sections $U\ni\mu\mapsto f_\mu\in C^{\pm\sharp}(\partial X,V(\sigma_\mu))$,
$U\ni\mu\mapsto f_\mu\in C^{\pm\sharp}(B,V_B(\sigma_\mu))$, $U\subset\aca$ 
open, as well
as meromorphic families of continuous operators between such section spaces in
the sense of meromorphic functions with values in topological vector spaces.
Here we equip the space of continuous linear operators between two such
spaces with the topology of uniform convergence on bounded sets.
For a discussion of these notions we refer to
\cite{bunkeolbrich000}, p. 87, and \cite{bunkeolbrich990}, Section 2.2.
Let us recall that a meromorphic family of operators is said to
have a finite-dimensional singularity at $\lambda\in\aca$ if the coefficients
of the principal part of the Laurent expansion at $\lambda$ are operators
of finite rank.

Let $(\vp,V_\vp)$ be a finite-dimensional representation of $\Gamma$.
We form the bundle $V(\tau,\varphi):=V(\tau)\otimes V_\vp$ on $\partial X$ carrying the tensor product action of $\Gamma$ and define $V_B(\tau,\varphi):=\Gamma\backslash (V(\tau)\otimes V_\vp)_{|\Omega}$. In particular, we have the spaces of sections $C^{\pm\sharp}(\partial X,V(\sigma_\lambda,\varphi))$ and $C^{\pm\sharp}(B,V_B(\sigma_\lambda,\varphi))$ as
well as the various notions of $\aca$-parametrized families of sections and operators.

We chose some norm on $V_\vp$. Since $\Gamma$ is finitely generated we can find an element $\mu\in\aaaa^*$, $\mu \ge 0$, and a constant $C$ such that 
\begin{equation}\label{wau}
\|\varphi(g)\|\le Ca_g^{\mu} \quad \mbox{ for all } g\in \Gamma\ .
\end{equation}
 
\begin{ddd}
Let $\delta_\varphi\in \aaaa^*$ be the infimum of all $\mu\in\aaaa^*$ satisfying Equation (\ref{wau}) for some $C$. It is independent of the chosen norm. We call $\delta_\varphi$ the exponent of $(\varphi,V_\vp)$.
\end{ddd}

The goal of the present section is to develop the technical means needed  for the understanding of the spaces of
$\Gamma$-invariant hyperfunction sections ${}^\Gamma C^{-\omega}(\partial X,V(\sigma_\lambda,\vp))$ and its subspaces ${}^\Gamma C^{-\omega}(\Lambda,V(\sigma_\lambda,\vp))$ of those $\Gamma$-invariant hyperfunction sections on $\partial X$ which are supported on the limit set.
The starting point is the restriction map
$$ res: {}^\Gamma C^{-\omega}(\partial X,V(\sigma_\lambda,\vp))\rightarrow
C^{-\omega}(B,V_B(\sigma_\lambda,\vp))$$
which is given by the restriction $res_\Omega$ of a hyperfunction section of $V(\sigma_\lambda,\vp)$ to the open subset $\Omega\subset\partial X$ followed
by the identification ${}^\Gamma C^{-\omega}(\Omega,V(\sigma_\lambda,\vp))
\cong C^{-\omega}(B,V_B(\sigma_\lambda,\vp))$.
Then ${}^\Gamma C^{-\omega}(\Lambda,V(\sigma_\lambda,\vp))=\ker res$.
If $\Omega\ne\emptyset$, then we are going to construct a meromorphic family 
$\aca\ni\lambda\mapsto ext_\lambda$, 
$$ ext_\lambda: C^{-\omega}(B,V_B(\sigma_\lambda,\vp))\rightarrow {}^\Gamma C^{-\omega}(\partial X,V(\sigma_\lambda,\vp))\ ,$$  
of right inverses of $res$. This has already been done for the case of
the real hyperbolic space and trivial representations $\sigma=1$ and $\vp=1$ in \cite{bunkeolbrich990}. For
general rank one spaces (with the exception of the Cayley hyperbolic plane
$\OO H^2$) 
and bundles the analogous construction for distribution sections has
been carried out in \cite{bunkeolbrich000}. We shall follow these references
quite closely.

Let us
begin with a comparison theorem which is independent of the theory of $ext_\lambda$, and which also holds in the cocompact case $\Omega=\emptyset$. 

\begin{theorem}\label{bart}
Let $\tau$ and $\vp$ be finite-dimensional representations of $P$ and $\Gamma$,
respectively.
Then a $\Gamma$-invariant hyperfunction section $f\in {}^\Gamma  C^{-\omega}(\partial X, V(\tau,\vp))$ is a distribution section if and only
if
$res(f)\in C^{-\infty}(B, V_B(\tau,\vp))$. In particular,
$$ {}^\Gamma  C^{-\omega}(\Lambda, V(\tau,\vp))={}^\Gamma  C^{-\infty}(\Lambda, V(\tau,\vp))\ .$$
\end{theorem}
\proof We fix an $M$-invariant Hermitian scalar product $(.,.)$ on $V_{\tilde\tau}$.
Completing $C^{\infty}(\partial X,V(\tilde\tau\otimes 1_{-\rho}))$ with respect to
the inner product
$$ (\phi,\psi):=\int_K (\phi(k),\psi(k))\: dk $$
we obtain the Hilbert space $H=L^2(\partial X,V(\tilde\tau\otimes 1_{-\rho}))$. Here we view $\phi,\psi \in C^{\infty}(\partial X,V(\tilde\tau\otimes 1_{-\rho}))$ as functions on $G$ with values in $V_{\tilde\tau}$. $H$ carries a continuous 
representation $\pi$ of $G$, and its subspace $H_\infty$ of smooth vectors coincides
with $C^{\infty}(\partial X,V(\tilde\tau\otimes 1_{-\rho}))$ as a smooth
Fr\'echet representation. Thus the
dual space $(H_\infty)^\prime$ coincides with  $C^{-\infty}(\partial X, V(\tau))$. We therefore want to show that $f\in {}^\Gamma  C^{-\omega}(\partial X, V(\tau,\vp))$ defines an element in 
$(H_\infty)^\prime\otimes V_\vp$, whenever $res(f)\in C^{-\infty}(B, V_B(\tau,\vp))$.

Let $H_K$ be the subspace of $K$-finite vectors of $H$. Then according
to \cite{wallach92}, Lemma 11.6.1 and Proposition 11.6.2 a functional
$f$ on $H_K$ extends continuously to $H_\infty$ if and only if for all
$\phi\in H_K$ the matrix coefficients $c_{f,\phi}$, defined by
$$ c_{f, \phi}(g):=\langle f,\pi(g)\phi \rangle\ , $$
are analytic functions on $G$ having moderate growth, i.e., there exist
constants $C_{f,\phi}\in\R$, $d_f\in \aaaa^*$ such that for all $g\in G$
$$
|c_{f, \phi}(g)|\le C_{f,\phi} a_g^{d_f}\ .
$$

Since $\bar Y$ is compact we find a compact set $F\subset X\cup\Omega$
such that $\bigcup_{\gamma\in \Gamma} \gamma F= X\cup\Omega$. Set $F_\infty:=F\cap\Omega$, $F_G:=FK\subset G$.

If $f\in C^{-\omega}(\partial X, V(\tau,\vp))$ then we can form matrix
coefficients $c_{f, \phi}$ with elements in $\phi\in H_K$ as above which
are now $V_\vp$-valued functions. If $f\in {}^\Gamma C^{-\omega}(\partial X, V(\tau,\vp))$ such that $res(f)\in C^{-\infty}(B, V_B(\tau,\vp))$, then we 
can find $f_1\in C^{-\omega}(\partial X, V(\tau,\vp))$ supported on some
compact subset $Q\subset \partial X\setminus F$ such that $f_2:=f-f_1 \in C^{-\infty}(\partial X, V(\tau,\vp))$. By the above we find constants
$C_{2,\phi}$, $d_2\in \aaaa^*$ such that for all $g\in G$
\begin{equation}\label{mg}
|c_{f_2, \phi}(g)|\le C_{2,\phi} a_g^{d_2}\ .
\end{equation}

Note that 
$H_K\subset C^{\omega}(\partial X,V(\tilde\tau\otimes 1_{-\rho}))$.
We now need the following lemma (compare \cite{bunkeolbrich000}, Equation (37)).

\begin{lem}\label{birth}
Let $Q\subset \partial X$, $L\subset G$ be compact such that $gP\not\in Q$ for all $g\in L$. Let $f\in C^{-\omega}(\partial X, V(\tau,\vp))$ with $\supp f\subset Q$ and $\phi\in C^{\omega}(\partial X,V(\tilde\tau\otimes 1_{-\rho}))$. Then there
exists a constant $C$ such that for all $g\in L$, $a\in A_+$, $k\in K$
$$  |c_{f,\phi}(gak)|\le C a^{-2\rho}\|\tau (a)\| \ .$$
\end{lem}
\proof
Let $w\in N_K(\aaaa)\setminus M$ be a representative of the nontrivial element of the Weyl group $W(\gaaa,\aaaa)\cong N_K(\aaaa)/M\cong \Z_2$.
Then by the Bruhat decomposition of $G$ we have
$\partial X\setminus eP=Nw P$. Thus we can view $\phi$ as a $V_{\tilde\tau}$
valued function on $Nw$. We have
$$ (\pi(ak)\phi)(nw)=(\pi(k)\phi)(a^{-1} n w)=a^{-2\rho}\tilde\tau(a)(\pi(k)\phi)(a^{-1} nw a^{-1})
=a^{-2\rho}\tilde\tau(a)(\pi(k)\phi)(a^{-1} naw)\ .$$
Since $gP\not\in Q$ for all $g\in L$ there is a compact subset 
$N_0\subset N$ such that
for all $g\in L$ the support of the hyperfunctions 
$nw\mapsto (\tilde\pi(g^{-1}) f)(nw)=f(gnw)$ is contained in $N_0w$.
We obtain
\begin{eqnarray}
|c_{f,\phi}(gak)|&=&|\langle f, \pi(gak)\phi \rangle| \nonumber\\
&=&a^{-2\rho}|\langle \tilde\pi(g^{-1})f,\tilde\tau(a)(\pi(k)\phi)(a^{-1}.aw)\rangle |\nonumber\\
&\le& a^{-2\rho}\|\tau(a)\|\: |\langle \tilde\pi(g^{-1})f,(\pi(k)\phi)(a^{-1}.aw)\rangle |\label{still}\ .
\end{eqnarray}
Since
$$ \lim_{a\to\infty}(\pi(k)\phi)(a^{-1}.aw)\equiv (\pi(k)\phi)(w) \qquad 
\mbox{in } C^\omega(N_0w,V_{\tilde\tau})$$
we see that the pairing in (\ref{still}) defines a continuous function
on the compact set $L\times(A_+\cup \{\infty\})\times K$ and is therefore
uniformly bounded by some constant $C$. This finishes the proof of the lemma.
\hB

We continue the proof of Theorem \ref{bart}. Let us assume for a moment that
$\Gamma$ is not cocompact, i.e., $F_\infty\ne \emptyset$.
Then $F$ can be covered by the sets 
$$U_n:=\{gaK\in X, gP\in \partial X \:|\: a_g^\alpha< n,\: gP\in F_\infty,\: a\in A_+\}
\ ,\qquad n\in \nat\ .$$
Thus by compactness of $F$ we can find a compact $L\subset G$ such that $LP=F_\infty$ and 
$F_G\subset LA_+K$. Since there are constants $C^\prime$, $d_1>0$ such that for 
$g\in L$, $a\in A_+$ 
$$ a\le C^\prime a_{ga}\ ,\qquad a^{-2\rho}\|\tau(a)\|\le a^{d_1}$$
Lemma \ref{birth} implies that
$$ |c_{f_1, \phi}(g)|\le C_{1} a_g^{d_1}\qquad \mbox{for all } g\in F_G\ .$$
Combining this with (\ref{mg}) we obtain 
$$ |c_{f, \phi}(g)|\le C_{0} a_g^{d_0}\qquad \mbox{for all } g\in F_G\ .$$
Obviously, this inequality remains true in the cocompact case (with $d_0=0$).
Since $f$ is $\Gamma$-invariant, we have
$$ c_{f, \phi}(\gamma g)=\vp(\gamma)c_{f, \phi}(g)\ .$$
We arrive at 
$$ |c_{f, \phi}(\gamma g)|\le C a_\gamma^{\delta_\vp+\ve} a_g^{d_0} 
\qquad \mbox{for all } g\in F_G, \gamma\in\Gamma\ .$$
Now by \cite{bunkeolbrich000}, Corollary 2.4, there is some $C^{\prime\prime}$
such that
$$a_\gamma a_g \le C^{\prime\prime} a_{\gamma g} \qquad \mbox{for all } g\in F_G, \gamma\in\Gamma\ .$$
Setting $d_f=\max(d_0,\delta_\vp+\ve)$ we obtain
$$ |c_{f, \phi}(g)|\le C_{f,\phi} a_g^{d_f}  
\qquad \mbox{for all } g\in G\ .$$
Thus $c_{f, \phi}$ has moderate growth which implies $f\in C^{-\infty}(\partial X, V(\tau,\vp))$. This finishes the proof of Theorem \ref{bart}.
\hB

\begin{kor}\label{bartI}
If $\Ree(\lambda)>\delta_\Gamma+\delta_\varphi$, then   
${}^\Gamma C^{-\omega}( \Lambda ,V(\sigma_\lambda,\varphi))=\{0\}$.
\end{kor}
\proof If $\Omega\ne\emptyset$, then Theorem 4.7 in \cite{bunkeolbrich000}
states that for $\Ree(\lambda)>\delta_\Gamma+\delta_\varphi$
$${}^\Gamma C^{-\infty}( \Lambda ,V(\sigma_\lambda,\varphi))=\{0\}\ .$$
This result can also be proved for cocompact $\Gamma$ by a similar
but much easier argument. The corollary now follows from Theorem \ref{bart}.
\hB
 
From now on we assume $\Omega\ne\emptyset$. 
As a space of hyperfunction sections over a noncompact manifold $C^{-\omega}(\Omega,V(\sigma_\lambda,\vp))$ does not carry a natural
topology. For this reason the continuity of $res$ is not obvious.
But the proof of the continuity in case of trivial $\sigma$ and $\vp$
given in \cite{bunkeolbrich990}, Lemma 2.11, carries over to the general situation.
Thus we have

\begin{lem}\label{bor}
The map
$$ res: {}^\Gamma C^{-\omega}(\partial X,V(\sigma_\lambda,\vp))\rightarrow
C^{-\omega}(B,V_B(\sigma_\lambda,\vp))$$ 
is continuous.
\end{lem}

We want to define $ext_\lambda$ as the adjoint of a push-down map
$$ \pi_{*,-\lambda}: C^\omega(\partial X,V(\tilde\sigma_\lambda,\tilde\vp))
\rightarrow C^\omega(B,V_B(\tilde\sigma_{-\lambda},\tilde\vp))\ ,$$
which should be given by
\begin{equation}\label{summ}
\pi_{*,-\lambda}(f)(kM)= \sum_{g\in\Gamma} (\pi(g)f)(kM),\quad kM\in\Omega \ ,
\end{equation}
if the sum converges. Here we have used the identification $$C^\omega(B,V_B(\tilde\sigma_{-\lambda},\tilde\varphi))\cong{}^\Gamma C^\omega(\Omega,V(\tilde\sigma_{-\lambda},\tilde\varphi))\ .$$ 
$\pi(g)$ is the action
given by $\pi^{\tilde\sigma,-\lambda}(g)\otimes\tilde\varphi(g)$. In \cite{bunkeolbrich000}
it is shown that if we replace on both sides analytic by smooth sections,
then $\pi_{*,-\lambda}$ converges for $\Ree(\lambda)>\delta_\Gamma+\delta_\vp$
and defines a holomorphic family of continuous maps in this half-plane. Adapting
the proofs of \cite{bunkeolbrich990}, Lemma 2.6 and Lemma 2.7, in a straightforward way to our more general situation we obtain 

\begin{lem}\label{anal1}
If $\Ree(\lambda)>\delta_\Gamma+\delta_\varphi$, then the sum (\ref{summ}) converges and defines
a  continuous map 
$$\pi_{\ast, -\lambda}: C^\omega(\partial X,V(\tilde\sigma_{-\lambda},\tilde\varphi))\rightarrow  C^\omega(B,V_B(\tilde\sigma_{-\lambda},\tilde\vp))\ .$$
For any $g\in\Gamma$ it satisfies
\begin{equation}\label{iwik}
\pi_{\ast, -\lambda}\circ \pi(g)=\pi_{\ast, -\lambda}\ .
\end{equation}
Moreover, $\pi_{*,-\lambda}$ depends holomorphically on $\lambda$.
\end{lem}

\begin{ddd}\label{defofext}
For $\Ree(\lambda)>\delta_\Gamma+\delta_\varphi$ we define 
the extension map 
$$ext_\lambda: C^{-\omega}(B,V_B(\sigma_\lambda,\varphi))\rightarrow {}^\Gamma C^{-\omega}(\partial X,V(\sigma_\lambda,\varphi))$$
to be the adjoint of 
$$\pi_{*,-\lambda}:C^{\omega}(\partial X,V(\tilde{\sigma}_{-\lambda},\tilde\varphi)) \rightarrow C^{\omega}(B,V_B(\tilde{\sigma}_{-\lambda},\tilde\varphi))\ .$$
\end{ddd}
In fact, 
by Lemma \ref{anal1} the extension exists, is continuous, and  by \cite{bunkeolbrich990}, Lemma 2.3, it  depends holomorphically on $\lambda$ as an
operator with values in $C^{-\omega}(\partial X,V(\sigma_\lambda,\varphi))$. It follows from (\ref{iwik}) that the range of $ext_\lambda$ consists of $\Gamma$-invariant vectors. Moreover, the restriction of $ext_\lambda$ to
distribution sections coincides with the extension map
$$ext_\lambda: C^{-\infty}(B,V_B(\sigma_\lambda,\varphi))\rightarrow {}^\Gamma C^{-\infty}(\partial X,V(\sigma_\lambda,\varphi))$$
considered and shown to be meromorphic on $\aca$ in \cite{bunkeolbrich000}.

\begin{lem}\label{basic}
For $\Ree(\lambda)>\delta_\Gamma+\delta_\varphi$ the extension map
is the inverse of $res$:
\begin{eqnarray}
res\circ ext_\lambda&=&\id\ , \label{paul}\\
ext_\lambda\circ res&=&\id\ . \label{paula}
\end{eqnarray}
\end{lem}
\proof By Lemma 4.5 of \cite{bunkeolbrich000} the
relation (\ref{paul}) holds on the dense
subspace $C^{-\infty}(B,V_B(\sigma_\lambda,\varphi))\subset C^{-\omega}(B,V_B(\sigma_\lambda,\varphi))$. It extends by continuity
to the whole space of hyperfunction sections. It follows that 
$$res\circ ext_\lambda \circ res= res\ .$$
Since by Corollary \ref{bartI} the restriction
map is injective for $\Ree(\lambda)>\delta_\Gamma+\delta_\varphi$ this implies
(\ref{paula}). 
\hB

We want to construct a meromorphic continuation of $ext_\lambda$ to
all of $\aca$. For this reason we will study the scattering matrix
acting on analytic sections. First we have to recall the Knapp-Stein
intertwining operators.

Let $w\in N_K(\aaaa)\setminus M$ be a representative of the nontrivial element of the Weyl group $W(\gaaa,\aaaa)\cong N_K(\aaaa)/M\cong \Z_2$. If $\sigma$ is a representation of $M$, then its Weyl-conjugate $\sigma^w$,
acting on the same vector space $V_\sigma$, is defined by 
$\sigma^w(m):=\sigma( w^{-1} mw)$. If $\sigma$ is equivalent to $\sigma^w$, then we say that $\sigma$ is Weyl-invariant. 
Unless indicated otherwise $\sigma$ shall from now on denote a Weyl-invariant representation
of $M$ which is either irreducible or of the
form $\sigma^\prime\oplus \sigma^{\prime w}$ with $\sigma^\prime$ irreducible 
and not Weyl-invariant. In both cases the representation of $M$ on $V_\sigma$ can be extended to a representation of $N_K(\aaaa)$ which we also denote by $\sigma$. This extension is unique up to a character of the Weyl group, i.e.,
the two possible choices of $\sigma(w)$ can differ by a sign, only.
Let us fix such an extension.

For $\Ree(\lambda)<0$ the (unnormalized) $G$-intertwining operator
$$\hat{J}_{\sigma,\lambda}:C^{\infty}(\partial X,V(\sigma_\lambda)) \rightarrow C^{\infty}(\partial X ,V(\sigma_{-\lambda}))$$ 
is defined by the
convergent integral
\begin{equation}\label{furunkel}
(\hat{J}_{\sigma,\lambda}f)(g):=\sigma(w)\int_{N} f(gnw)\: dn\ .
\end{equation}
Here we consider $f\in C^{\infty}(\partial X,V(\sigma_\lambda))$
as a function on $G$ with values
in $V_{\sigma_\lambda}$ satisfying $f(gp)=\sigma_\lambda(p)^{-1}f(g)$ for all $p\in P$. 
The operator $\hat{J}_{\sigma,\lambda}$ does not depend on the choice of $w$.
For $\sigma=\sigma^\prime\oplus \sigma^{\prime w}$ as above we denote
by
$\hat J_{\sigma^\prime,\lambda}$ the restriction of $\hat{J}_{\sigma,\lambda}$ to $C^{\infty}(\partial X,V(\sigma^\prime_\lambda))$.

The holomorphic family $\lambda\mapsto \hat{J}_{\sigma,\lambda}$ has a meromorphic continuation to all of $\aca$ with poles of at most first order. There is a meromorphic function
$p_\sigma: \aca\rightarrow \C$, called the Plancherel density, such that
\begin{equation}\label{rob}
\hat{J}_{\sigma,-\lambda}\circ\hat{J}_{\sigma,\lambda}=
\frac{1}{p_\sigma(\lambda)}\id\ .
\end{equation}
For all this see \cite{knappstein71} or \cite{wallach92}, Ch. 10. 
Since $C^{\omega}(\partial X,V(\sigma_\lambda))$ is the space of analytic vectors of the continuous Fr\'echet representation $C^{\infty}(\partial X,V(\sigma_\lambda))$ the intertwining operator restricts to a continuous operator
$$\hat{J}_{\sigma,\lambda}:C^{\omega}(\partial X,V(\sigma_\lambda)) \rightarrow C^{\omega}(\partial X ,V(\sigma_{-\lambda}))\ .$$
The argument of the proof of \cite{bunkeolbrich990}, Lemma 2.16, shows that 
the latter operator indeed comes as a meromorphic family in the strong sense 
of this paper. The restriction of the adjoint of $\hat J_{\tilde\sigma,\lambda}$
to $C^{\omega}(\partial X,V(\sigma_\lambda))$ coincides with $\hat{J}_{\sigma,\lambda}$ (\cite{knappstein71}, Lemma 24). Therefore this adjoint can be used to define the continuous extension to hyperfunction sections
$$\hat{J}_{\sigma,\lambda}:C^{-\omega}(\partial X,V(\sigma_\lambda)) \rightarrow C^{-\omega}(\partial X ,V(\sigma_{-\lambda}))\ ,$$
which is again a meromorphic family of operators 
(\cite{bunkeolbrich990}, Lemma 2.3).

For any open $U\subset \partial X$ we introduce the space
$$ C^{-\omega}_U(\partial X,V(\sigma_\lambda)):=\{f\in C^{-\omega}(\partial X,V(\sigma_\lambda))\:|\: res_Uf\in C^{\omega}(U,V(\sigma_\lambda))\}\ .$$ 
We equip $C^{-\omega}_U(\partial X,V(\sigma_\lambda))$ with the weakest
topology such that the embedding $C^{-\omega}_U(\partial X,V(\sigma_\lambda))$
$\hookrightarrow C^{-\omega}(\partial X,V(\sigma_\lambda))$ and the restrictions
$ C^{-\omega}_U(\partial X,V(\sigma_\lambda))\rightarrow  C^{\omega}(W,V(\sigma_\lambda))$ to the spaces of germs of real analytic sections along any compact $W\subset U$ are continuous. 
We also consider the analogously defined twisted versions $ C^{-\omega}_U(\partial X,V(\sigma_\lambda,\vp))$.
The argument of the proof of \cite{bunkeolbrich990}, Lemma 2.19, shows
that the intertwining operators are off-diagonally smoothing in the following
strong sense.

\begin{lem}\label{offD}
For any open $U\subset\partial X$ the intertwining operators induce a meromorphic family of continuous
operators
$$\hat{J}_{\sigma,\lambda}:C^{-\omega}_U(\partial X,V(\sigma_\lambda)) \rightarrow C^{-\omega}_U(\partial X ,V(\sigma_{-\lambda}))\ .$$
\end{lem}

Tensoring with a finite-dimensional representation $(\varphi,V_\vp)$ of $\Gamma$ we obtain a
meromorphic family of $\Gamma$-intertwining operators which we denote
by the same symbol
$$\hat J_{\sigma,\lambda}:=\hat{J}_{\sigma,\lambda}\otimes\id: 
C^{-\omega}_U(\partial X,V(\sigma_\lambda,\varphi)) \rightarrow C^{-\omega}_U(\partial X ,V(\sigma_{-\lambda},\varphi))\ .$$

\begin{ddd}\label{scatdef}
For $\Ree(\lambda)>\delta_\Gamma+\delta_\varphi$ we define the unnormalized
scattering matrix
$$\hat{S}_{\sigma, \lambda}: C^{-\omega}(B,V_B(\sigma_\lambda,\varphi))
\rightarrow C^{-\omega}(B,V_B(\sigma_{-\lambda},\varphi))$$
as the continuous operator given by the composition
$$\hat{S}_{\sigma, \lambda}:=res\circ \hat{J}_{\sigma,\lambda}\circ ext_\lambda\ .$$
\end{ddd}
From now on we assume $X\ne \OO H^2$. Theorem 5.10 in \cite{bunkeolbrich000} tells us that the restriction of the scattering matrix
to smooth and distribution sections, respectively,
\begin{equation}\label{owl}
\hat{S}_{\sigma, \lambda}: C^{\pm\infty}(B,V_B(\sigma_\lambda,\varphi))
\rightarrow C^{\pm\infty}(B,V_B(\sigma_{-\lambda},\varphi))
\end{equation}
has a meromorphic continuation to all of $\aca$.
We are now able to prove the corresponding statement in the analytic category.

\begin{prop}\label{coocoo}
The restriction of the scattering matrix (\ref{owl}) to
$C^{\omega}(B,V_B(\sigma_\lambda,\varphi))$ provides a meromorphic family on
$\aca$ of continuous operators
$$\hat{S}_{\sigma, \lambda}: C^{\omega}(B,V_B(\sigma_\lambda,\varphi))
\rightarrow C^{\omega}(B,V_B(\sigma_{-\lambda},\varphi))\ .$$
The meromorphic continuation of the scattering matrix acting on hyperfunctions,
initially defined for $\Ree(\lambda)>\delta_\Gamma+\delta_\varphi$, is
given by the adjoint of
\begin{equation}\label{muh}
\hat{S}_{\tilde\sigma, \lambda}: C^{\omega}(B,V_B(\tilde\sigma_\lambda,\tilde\varphi))
\rightarrow C^{\omega}(B,V_B(\tilde\sigma_{-\lambda},\tilde\varphi))\ .
\end{equation}
We have
\begin{equation}\label{moik}
\hat{S}_{\sigma,-\lambda}\circ\hat{S}_{\sigma,\lambda}=
\frac{1}{p_\sigma(\lambda)}\id\ .
\end{equation}
\end{prop}
\proof \cite{bunkeolbrich000}, Theorem 5.10, asserts that the restriction 
of $ext_\lambda$ to distribution sections has a meromorphic continuation to $\aca$. This together with (\ref{paul}) implies that
$$ ext_\lambda: C^{\omega}(B,V_B(\sigma_\lambda,\varphi))\rightarrow 
{}^\Gamma C^{-\omega}_\Omega(\partial X,V_B(\sigma_\lambda,\varphi))$$
is meromorphic. Now 
$$\hat{J}_{\sigma,\lambda}:C^{-\omega}_\Omega(\partial X,V(\sigma_\lambda,\vp)) \rightarrow C^{-\omega}_\Omega(\partial X ,V(\sigma_{-\lambda},\vp))$$
is meromorphic by Lemma \ref{offD}, while
$$ res_\Omega: C^{-\omega}_\Omega(\partial X ,V(\sigma_{-\lambda},\vp))
\rightarrow C^{\omega}(\Omega,V(\sigma_{-\lambda},\vp))$$
is holomorphic by definition. We conclude the meromorphy of
$$\hat{S}_{\sigma, \lambda}=res\circ \hat{J}_{\sigma,\lambda}\circ ext_\lambda: C^{\omega}(B,V_B(\sigma_\lambda,\varphi))
\rightarrow C^{\omega}(B,V_B(\sigma_{-\lambda},\varphi))\ .$$
The adjoint of (\ref{owl}) is given by
$$ \hat{S}_{\tilde\sigma, \lambda}: C^{\mp\infty}(B,V_B(\tilde\sigma_\lambda,\tilde\varphi))
\rightarrow C^{\mp\infty}(B,V_B(\tilde\sigma_{-\lambda},\tilde\varphi))$$
(\cite{bunkeolbrich000}, Lemma 5.8). This implies that the adjoint of (\ref{muh}) indeed provides a meromorphic continuation of the scattering matrix
as defined in Definition \ref{scatdef}. Concerning the functional equation (\ref{moik}) we also could refer to \cite{bunkeolbrich000}, Theorem 5.10.
However, it might be instructive to verify (\ref{moik}) directly. Note
that (\ref{paula}) remains true on all of $\aca$, if we restrict both sides
to $C^{\omega}(B,V_B(\sigma_\lambda,\varphi))$. Let us compute $\hat{S}_{\sigma,-\lambda}\circ\hat{S}_{\sigma,\lambda}$ on $C^{\omega}(B,V_B(\sigma_\lambda,\varphi))$ for $\Ree(\lambda)<\delta_\Gamma+\delta_\varphi$  
$$\hat{S}_{\sigma,-\lambda}\circ\hat{S}_{\sigma,\lambda}=
res\circ \hat J_{\sigma,-\lambda}\circ ext_{-\lambda}\circ res \circ
\hat J_{\sigma,\lambda}\circ ext_{\lambda}
= res\circ \hat J_{\sigma,-\lambda}\circ 
\hat J_{\sigma,\lambda}\circ ext_{\lambda}
=\frac{1}{p_\sigma(\lambda)}\id\ .$$
Here we have also used (\ref{paul}) and (\ref{rob}). By
meromorphy (\ref{moik}) holds on all of $\aca$.
\hB 

We now come to the meromorphic continuation of $ext_\lambda$. First we treat the
case $\delta_\Gamma+\delta_\varphi<0$.

\begin{lem}\label{col} 
If $\delta_\Gamma+\delta_\varphi<0$, then
$$ext_\lambda: C^{-\omega}(B,V_B(\sigma_\lambda,\varphi))\rightarrow {}^\Gamma C^{-\omega}(\partial X,V(\sigma_\lambda,\varphi))\ ,$$
initially defined for $\Ree(\lambda)>\delta_\Gamma+\delta_\varphi$, admits
a meromorphic continuation to all of $\aca$ with at most finite-dimensional
singularities.
\end{lem}  
\proof
For $\Ree(\lambda)<-(\delta_\Gamma+\delta_\varphi)$ we set
$$ \widetilde{ext}_\lambda:=p_\sigma(\lambda)\hat{J}_{\sigma,-\lambda}\circ
ext_{-\lambda}\circ\hat{S}_{\sigma,\lambda}\ .$$
Using (\ref{paula}) and the functional equation (\ref{moik})
we obtain for $\Ree(\pm\lambda)<-(\delta_\Gamma+\delta_\varphi)$
\begin{eqnarray*}
\widetilde{ext}_\lambda=ext_\lambda\circ res\circ \widetilde{ext}_\lambda
&=&p_\sigma(\lambda) ext_\lambda\circ res\circ \hat{J}_{\sigma,-\lambda}\circ
ext_{-\lambda}\circ\hat{S}_{\sigma,\lambda}\\
&=&p_\sigma(\lambda) ext_\lambda\circ
\hat{S}_{\sigma,-\lambda}\circ\hat{S}_{\sigma,\lambda} = ext_\lambda\ .
\end{eqnarray*}
Thus $\widetilde{ext}_\lambda$ provides a meromorphic continuation of $ext_\lambda$ to all of $\aca$. Since $ext_\lambda$ restricted to distribution
sections has finite-dimensional singularities (\cite{bunkeolbrich000}, Theorem 5.10), and since distribution sections are
dense in $C^{-\omega}(B,V_B(\sigma_\lambda,\varphi))$
its continuous extension to hyperfunction sections has the same finite-dimensional singularities.
\hB

As in \cite{bunkeolbrich990} and \cite{bunkeolbrich000} we now apply the embedding trick in order to remove the assumption $\delta_\Gamma+\delta_\varphi<0$. 

\begin{lem}\label{muehsam}
If $\sigma=1$ is the trivial representation of $M$, then Lemma \ref{col} holds true without the assumption
$\delta_\Gamma+\delta_\varphi<0$.
\end{lem}
\proof The proof is almost the same as the one of Lemma 5.13 in \cite{bunkeolbrich000}. For the convenience of the reader we reproduce it here.

$X$ belongs to a series of rank-one symmetric
spaces.
First we assume  that $G^n$ belongs to the list
$\{Spin(1,n),SO(1,n)_0,SU(1,n),Sp(1,n)\}$.
Then we have a sequence of real, semisimple, linear Lie groups $\dots \subset G^n\subset G^{n+1}\subset \dots$ inducing embeddings of the corresponding
Iwasawa constituents $K^n\subset K^{n+1}$,
$N^n \subset  N^{n+1}$, $ M^n\subset M^{n+1}$
such that $A=A^n = A^{n+1}$.
Furthermore, there are totally geodesic embeddings of the symmetric spaces
$X^n\subset X^{n+1}$
inducing embeddings of their boundaries
$\partial X^n\subset \partial X^{n+1}$.
If $\Gamma\subset G^n$ is convex cocompact, then it is still convex cocompact
viewed  as a subgroup of $G^{n+1}$.
We obtain  embeddings
$\Omega^n\subset  \Omega^{n+1}$ inducing
$B^n\subset  B^{n+1}$
while the limit set $ \Lambda^n$ is identified with $ \Lambda^{n+1}$.
Let $\rho^n(H)=\frac{1}{2}\tr(\ad(H)_{|\naaa^n})$, $H\in\aaaa$.

The exponent of $\Gamma$ now depends on $n$ and is denoted by
$\delta_\Gamma^n$.
We have the relation $\delta_\Gamma^{n+1}=\delta_\Gamma^n-\zeta$, where $\zeta:=\rho^{n+1}-\rho^n>0$.
Thus $\delta_\Gamma^{n+m}\to -\infty$ as $m\to \infty$. Hence,
taking $m$ large enough we obtain $\delta_\Gamma^{n+m}+\delta_\varphi<0$.
The aim of the following discussion is to show how the meromorphic continuation of $ext^{n+1}_\lambda$ leads to the continuation of  $ext_\lambda^n$.

Let $P^{n}:=M^nA^nN^n$, $V(1_\lambda,\varphi)^n:=G^n\times_{P^n} V_{1_\lambda}\otimes V_\vp$, and $V_{B^n}(1_\lambda,\varphi)=\Gamma\backslash V(1_\lambda,\varphi)^n_{|\Omega^n}$. Here as always $(\varphi,V_\vp)$ is a finite-dimensional representation of $\Gamma$.
The representation $V_{1_\lambda}$ of $P^{n+1}$ restricts to the representation
$V_{1_{\lambda-\zeta}}$ of $P^n$.
This induces isomorphisms of bundles
$$V(1_\lambda,\varphi)^{n+1}_{|\partial X^n}\cong
V(1_{\lambda-\zeta},\varphi)^n,\quad  V_{B^{n+1}}(1_\lambda,\varphi)_{|B^n}\cong V_{B^n}(1_{\lambda-\zeta},\varphi)\ .$$

Let 
\begin{eqnarray*}
i^*:C^\omega(B^{n+1},V_{B^{n+1}}(1_\lambda,\tilde\varphi))&\rightarrow& C^\omega(B^{n },V_{B^{n}}(1_{\lambda-\zeta},\tilde\varphi))\ ,\\
j^*:C^\omega(\partial X^{n+1},V(1_\lambda,\tilde\varphi)^{n+1})&\rightarrow& C^\omega(\partial X^n,V(1_{\lambda-\zeta},\tilde\varphi)^n)
\end{eqnarray*}
denote the maps given by restriction of sections.
Note that $j^*$ is $G^n$-equivariant.
The adjoint maps define
the push-forward of hyperfunction sections
\begin{eqnarray*}
i_*:C^{-\omega}(B^n,V_{B^n}(1_{\lambda},\varphi))&\rightarrow& C^{-\omega}(B^{n+1},V_{B^{n+1}}(1_{\lambda-\zeta},\varphi))\ ,\\
j_*:C^{-\omega}(\partial X^n,V(1_{\lambda},\varphi)^n)&\rightarrow& C^{-\omega}(\partial X^{n+1},V(1_{\lambda-\zeta},\varphi)^{n+1})\ .
\end{eqnarray*}

If $\phi\in C^{-\omega}(B^n,V_{B^n}(1_\lambda,\varphi))$, then the push forward
$i_\ast \phi $ has support in $B^n\subset B^{n+1}$.
Since $res^{n+1}\circ ext^{n+1}_\lambda=\id$
we have 
\begin{equation}\label{ikmo}
\supp(ext^{n+1}_\lambda\circ i_*)(\phi)\subset 
\Lambda^{n+1}\cup \Omega^{n}=\partial X^n\ .
\end{equation}
Assume that $ext^{n+1}_\lambda$ is meromorphic on $\aca$.
We are now going to continue $ext_\lambda^n$ using $i_*$, $ext^{n+1}_{\lambda-\xi}$ and a left inverse of $j_*$.  
We identify $C^{\omega}(\partial X^{n+1},V(1_\lambda,\tilde\varphi)^{n+1})$
with $C^{\omega}(\partial X^{n+1})\otimes V_{\tilde\vp}$ for all $\lambda\in\C$.
 Denote by $\bar B(0,1)$ the closed unit ball in $\F$, where $\F=\R$,
$\C$ or $\HH$, respectively. We choose an analytic diffeomorphism $T:\bar B(0,1) \times \partial X^n\stackrel{\cong}{\rightarrow} U$ to a tubular neighborhood  $U$ 
of $\partial X^n$ in $\partial X^{n+1}$. Then we define
a continuous extension $t:C^{\omega}(\partial X^n)\otimes V_{\tilde\vp}\rightarrow C^\omega(U)\otimes V_{\tilde\vp}$
by 
$$T^*tf(r,x):=f(x)\ .$$ 
Let
$t^\prime:\:C^{-\omega}(U,V(1_{\lambda-\zeta},\varphi)^{n+1})\rightarrow C^{-\omega}(\partial X^n,V(1_{\lambda},\varphi)^n)$ be the adjoint of $t$. Then $t^\prime\circ j_*=\id$. Because of (\ref{ikmo}) we can define
$$\widetilde{ext}^n_\lambda \phi:= (t^\prime \circ ext^{n+1}_{\lambda-\xi}\circ i_*)(\phi)\ .$$
Then
$$\widetilde{ext}^n_\lambda : C^{-\omega}(B^n,V_{B^n}(1_\lambda,\varphi)) \rightarrow C^{-\omega}(\partial X^n,V(1_\lambda,\varphi)^n)$$
is a meromorphic family on $\aca$ of continuous maps with at most finite-dimensional singularities. 

In order to prove that $\widetilde{ext}^n_\lambda$ provides the desired
meromorphic continuation it remains to show that it coincides with ${ext}^n_\lambda$ 
in the region $\Ree(\lambda)>\delta_\Gamma^{n}+\delta_\varphi$. If $\Ree(\lambda)>\delta_\Gamma^{n}+\delta_\varphi$, then $\Ree(\lambda)-\zeta>\delta_\Gamma^{n+1}+\delta_\varphi$,
and the push-down maps $\pi_{*,-\lambda}^n$, $\pi_{*,-\lambda+\zeta}^{n+1}$ are defined.
It is easy to see from the definition of the push-down that in the domain of convergence
$$i^*\circ \pi_{*,-\lambda+\xi}^{n+1}=\pi_{*,-\lambda}^n\circ j^*\ .
$$
Taking adjoints we obtain $ext^{n+1}_{\lambda-\xi}\circ i_*
=j_*\circ {ext}^n_\lambda$.
Therefore we have
$$\widetilde{ext}^n_\lambda
=t^\prime \circ ext^{n+1}_{\lambda-\xi}\circ i_*
=t^\prime \circ j_*\circ {ext}^n_\lambda
={ext}^n_\lambda\ .$$
It follows by meromorphy that $\im(\widetilde{ext}^n_\lambda)$ consists of $\Gamma$-invariant sections for all $\lambda\in\aca$.
 
If $G^n$ does not belong to the list
$\{Spin(1,n),SO(1,n)_0,SU(1,n),Sp(1,n)\}$, then there is a finite covering
$p:\tilde G^n\rightarrow G^n$ with $\tilde G^n\in 
\{Spin(1,n),SO(1,n)_0,SU(1,n),Sp(1,n)\}$.
In this case one can find a normal subgroup
$\Gamma^0\subset \Gamma$ of finite index and a discrete torsion-free subgroup $\tilde \Gamma^0\subset \tilde G^n$ such that $p$ induces an isomorphism from $\tilde \Gamma^0$ to $\Gamma^0$. Indeed, using Selberg's Lemma (see e.g.
\cite{ratcliffe94}) we can take
a torsion-free subgroup $\tilde\Gamma^0$  of $p^{-1}(\Gamma)$ of finite index
and set $\Gamma^0:=p(\tilde\Gamma^0)$.
We can apply the concept of embedding to
the subgroup $\tilde\Gamma^0$. In order to transfer results for
$\tilde\Gamma^0$ to $\Gamma$ we use averages over the finite group $\Gamma/\Gamma^0$.
This finishes the proof of Lemma \ref{muehsam}.
\hB 

Now we are able to prove the main result of this section in full generality.

\begin{theorem}\label{extmer}
If $X\ne \OO H^2$, then the extension map
$$ext_\lambda: C^{-\omega}(B,V_B(\sigma_\lambda,\varphi))\rightarrow {}^\Gamma C^{-\omega}(\partial X,V(\sigma_\lambda,\varphi))\ ,$$
initially defined for $\Ree(\lambda)>\delta_\Gamma+\delta_\varphi$, admits
a meromorphic continuation to all of $\aca$ with at most finite-dimensional
singularities and satisfies
\begin{equation}\label{paul'}
res\circ ext_\lambda=\id\ .
\end{equation}
\end{theorem}
\proof 
In view of Lemma \ref{muehsam} it is enough to reduce the statement to
the case $\sigma=1$. As in \cite{bunkeolbrich000} we use tensoring with finite-dimensional $G$-representations. We can assume that $\sigma$ is irreducible. Then there exists an irreducible finite-dimensional representation $\pi_{\sigma,\mu}$ of $G$ with highest $\aaaa$-weight $\mu\in\aaaa^*$ such that the representation of $M$ on the highest weight space is equivalent to $\sigma$. The embedding of $P$-representations
$\sigma_\lambda\hookrightarrow 1_{\lambda+\mu}\otimes\pi_{\sigma,\mu}$ induces
a $\Gamma$-equivariant embedding of bundles
\begin{equation}\label{str} 
i_{\sigma,\mu}: V(\sigma_\lambda,\vp)\hookrightarrow V(1_{\lambda+\mu},\pi_{\sigma,\mu}\otimes\vp)\ ,
\end{equation}
where we have used the $\Gamma$-equivariant isomorphism $$V({\pi_{\sigma,\mu}}_{|P})\cong V(1_\rho, {\pi_{\sigma,\mu}}_{|\Gamma})\ ,\qquad 
[g,v]\mapsto (gP, \pi_{\sigma,\mu}(g)v)\ .$$
Therefore we have the corresponding embedding of bundles over $B$
$$ i^B_{\sigma,\mu}: V_B(\sigma_\lambda,\vp)\hookrightarrow V_B(1_{\lambda+\mu},\pi_{\sigma,\mu}\otimes\vp)\ .$$
These embeddings induce corresponding embeddings of the spaces of sections
which we denote by the same symbols.
One now checks that for $\Ree(\lambda)>\delta_\Gamma+\delta_\varphi$
\begin{equation}\label{hope}
ext_{\lambda+\mu}\circ i_{\sigma,\mu}^B=i_{\sigma,\mu}\circ ext_\lambda\ .
\end{equation}
By Lemma \ref{muehsam} the left hand side of (\ref{hope}) is meromorphic on $\aca$ with 
finite-dimensional singularities. This provides the meromorphic
continuation of the right hand side, thus of $ext_\lambda$ in general. For more details we refer to \cite{bunkeolbrich000}, pp. 108-109.
Equation (\ref{paul'}) now follows from (\ref{paul}) by meromorphy. This
finishes the proof of the theorem.
\hB

\begin{kor}\label{dense}
${}^\Gamma C^{-\infty}(\partial X,V(\sigma_\lambda,\varphi))$
is dense in ${}^\Gamma C^{-\omega}(\partial X,V(\sigma_\lambda,\varphi))$.
\end{kor}
\proof 
Choosing a holomorphic trivialization of the family of bundles
$$ \bigcup_z V_B(\sigma_{\lambda+z\alpha},\vp)$$
we can identify 
$C^{-\omega}(B,V(\sigma_\lambda,\varphi))$ with the space of constant families
$\mu\mapsto f_\mu \in  C^{-\omega}(B,V(\sigma_\mu,\varphi))$.
We set
$$ W_{reg}:=\{ f\in C^{-\omega}(B,V(\sigma_\lambda,\varphi))\:|\: ext_\mu f_\mu \mbox{ is regular at } \mu=\lambda\}\ .$$
Then $ext_\lambda: W_{reg}\rightarrow {}^\Gamma C^{-\omega}(\partial X,V(\sigma_\lambda,\varphi))$ is a well-defined continuous map. Since the singularity of $ext_\mu$ at $\mu=\lambda$ is at most finite-dimensional $W_{reg}\subset C^{-\omega}(B,V(\sigma_\lambda,\varphi))$
has finite codimension. Using that $C^{-\infty}(B,V(\sigma_\lambda,\varphi))\subset C^{-\omega}(B,V(\sigma_\lambda,\varphi))$ is dense we can find a finite-dimensional subspace $W_\infty\subset C^{-\infty}(B,V(\sigma_\lambda,\varphi))$ such that 
\begin{equation}\label{dec}
C^{-\omega}(B,V(\sigma_\lambda,\varphi))=W_{reg}\oplus W_\infty\ .
\end{equation}
Let now $f\in {}^\Gamma C^{-\omega}(\partial X,V(\sigma_\lambda,\varphi))$.
Choose a sequence $g_i\in C^{-\infty}(B,V(\sigma_\lambda,\varphi))$ converging
to $g=res(f)$. We decompose $g$ and $g_i$ according to (\ref{dec}), $g=g^{reg}+g^\infty$, $g_i=g_i^{reg}+g_i^\infty$. Then $g^\infty, g_i^\infty, g_i^{reg}\in C^{-\infty}(B,V(\sigma_\lambda,\varphi))$. We set
$$ f_i= f-ext_\lambda g^{reg}+ext_\lambda g^{reg}_i\ .$$
Because of the continuity of $ext_\lambda$ and of the splitting (\ref{dec}) 
the sequence $f_i$ converges to $f$. Moreover,
$$ res(f_i)=g^\infty+g_i^{reg}\in C^{-\infty}(B,V(\sigma_\lambda,\varphi))\ .$$
Hence $f_i\in {}^\Gamma C^{-\infty}(\partial X,V(\sigma_\lambda,\varphi))$ by
Theorem \ref{bart}. This proves the corollary.
\hB

\newpage
\section{Invariant distributions on the limit set}\label{round}

This section can be considered as a variation on the theme of \cite{bunkeolbrich000}, Section 6.
We retain the notation and assumptions of the previous section. 
Throughout this section we assume $\Omega\ne\emptyset$.
For a given representation $\sigma$ of $M$ we now consider
only such $\lambda\in\aca$, where $ext_\lambda$ is meromorphic with an at most
finite-dimensional singularity. By Theorem \ref{extmer} this is the case for all $\lambda\in\aca$ whenever
$X\ne \OO H^2$.

First we need the following result 
(compare \cite{bunkeolbrich990}, Proposition 3.4).

\begin{prop}\label{green}
If
$f\in {}^\Gamma C^{-\omega}(\Lambda,V(\sigma_\lambda,\varphi))$ and 
$\phi\in {}^\Gamma C^{-\omega}(\partial X,V(\tilde{\sigma}_\lambda,\tilde\varphi))$, 
then 
$$\langle res\circ \hat J_{\sigma,\lambda} (f),res(\phi)\rangle =0\ .$$
\end{prop}

Note that $res\circ \hat J_{\sigma,\lambda}(f)$ is regular since it
only sees the off-diagonal part of the intertwining operator. Moreover, $res\circ \hat J_{\sigma,\lambda} (f)\in C^{\omega}(B,V_B(\sigma_{-\lambda},\varphi))$ by Lemma \ref{offD}. 
Thus the above pairing is well-defined.

\noindent
\proof According to \cite{bunkeolbrich000}, Proposition 6.5, the assertion
is true for $f\in {}^\Gamma C^{-\infty}(\Lambda,V(\sigma_\lambda,\varphi))$ and 
$\phi\in {}^\Gamma C^{-\infty}(\partial X,V(\tilde{\sigma}_\lambda,\tilde\varphi))$. Now $f\in {}^\Gamma C^{-\infty}(\Lambda,V(\sigma_\lambda,\varphi))$ by Theorem \ref{bart}.
Because of Lemma \ref{bor} and Corollary \ref{dense} the assertion extends to $\phi\in {}^\Gamma C^{-\omega}(\partial X,V(\tilde{\sigma}_\lambda,\tilde\varphi))$ by continuity.
\hB

Next we introduce certain extensions of the bundles $V(\sigma_\lambda,\vp)$
with itselves. Let $\Pi$ be the space of polynomials on $\aaaa$. The group $A$ acts on it by translations.
This action extends to a representation
$ 1^+:MAN\rightarrow GL(\Pi)$ given by
$$ 1^+(man)f(H):=f(H-\log a)\ .$$
For any $k\in\nat$ the finite dimensional subspace $\Pi^k$ of polynomials of degree at most $k-1$
is invariant with respect to this action. We denote the restriction of
$1^+$ to $\Pi^k$ by $1^k$. 
Then we set $V^k(\sigma_\lambda,\vp):=V(\sigma_\lambda\otimes 1^k,\vp)$, $V^+(\sigma_\lambda,\vp):=V(\sigma_\lambda\otimes 1^+,\vp)$.  
There is a chain of inclusions
\begin{eqnarray*}
\{0\}=V^0(\sigma_\lambda,\vp)\subset V(\sigma_\lambda,\vp)&=&V^1(\sigma_\lambda,\vp)
\subset V^2(\sigma_\lambda,\vp)\subset\dots\\
\dots &\subset& V^k(\sigma_\lambda,\vp)\subset V^{k+1}(\sigma_\lambda,\vp)\subset\dots\subset V^+(\sigma_\lambda,\vp)\ .
\end{eqnarray*}
Again there is a restriction map between the corresponding spaces of sections
$$ res:{}^\Gamma C^{-\sharp}(\partial X, V^k(\sigma_\lambda,\vp))\rightarrow 
C^{-\sharp}(B, V^k_B(\sigma_\lambda,\vp))\ ,\qquad \sharp\in\{\infty,\omega\}\ .$$
We are especially interested in the space of invariant 
sections supported on the
limit set (recall Theorem \ref{bart}) 
$${}^\Gamma C^{-\infty}(\Lambda, V^k(\sigma_\lambda,\vp)):=\ker res\ .$$
We set
\begin{eqnarray*}
C^{\pm\sharp}(\partial X, V^+(\sigma_\lambda,\vp))&:=& \bigcup_{k\in\nat} C^{\pm\sharp}(\partial X, V^k(\sigma_\lambda,\vp))\ , \\ 
C^{\pm\sharp}(B, V_B^+(\sigma_\lambda,\vp))&:=& \bigcup_{k\in\nat} C^{\pm\sharp}(B, V_B^k(\sigma_\lambda,\vp))\ .
\end{eqnarray*}
The space $C^{\infty}(\partial X, V^+(\sigma_\lambda,\vp))\subset C^{-\sharp}(\partial X, V^+(\sigma_\lambda,\vp))$ can be described as follows
\begin{equation}\label{diep} 
\left\{f:G\times\aaaa\rightarrow
V_\sigma\otimes V_\vp\:|\: \begin{array}{l}f \mbox{ smooth and polynomial w.r.t. }\aaaa\:,\\  
f(gman,H)=a^{\lambda-\rho}(\sigma(m)^{-1}\otimes\id)f(g,H+\log a)
\end{array}\right\}\ .
\end{equation} 
 
We view the short root $\alpha$ as a coordinate on $\aaaa$. Differentiation
defines a $P$-equivariant operator 
$$ \frac{d}{d\alpha}: \Pi\rightarrow \Pi\ . $$
It induces a $\Gamma$-equivariant bundle homomorphism 
$$ \p: V^{+}(\sigma_\lambda,\vp)\rightarrow V^+(\sigma_\lambda,\vp) $$
which maps $V^{k+1}(\sigma_\lambda,\vp)$ onto $V^k(\sigma_\lambda,\vp)$.
$\p$ induces a kind of shift operators
\begin{eqnarray*}\p: C^{\pm\sharp}(\partial X, V^+(\sigma_\lambda,\vp))&\rightarrow& C^{\pm\sharp}(\partial X, V^+(\sigma_\lambda,\vp))\ ,\\ 
\p_\Gamma: {}^\Gamma C^{-\sharp}(\partial X, V^+(\sigma_\lambda,\vp))&\rightarrow& {}^\Gamma C^{-\sharp}(\partial X, V^+(\sigma_\lambda,\vp))\ ,\\
\p_B: C^{\pm\sharp}(B, V_B^+(\sigma_\lambda,\vp))&\rightarrow& C^{\pm\sharp}(B, V_B^+(\sigma_\lambda,\vp))\ .
\end{eqnarray*}
Note that $\p$ and $\p_B$ are surjective while the possible failure of surjectivity of $\p_\Gamma$ will turn out to have cohomological meaning
(see Corollary \ref{preth}). 
We also have the restriction map
$$ res:{}^\Gamma C^{-\sharp}(\partial X, V^+(\sigma_\lambda,\vp))\rightarrow 
C^{-\sharp}(B, V^+_B(\sigma_\lambda,\vp)$$
and its kernel 
${}^\Gamma C^{-\infty}(\Lambda, V^+(\sigma_\lambda,\vp))= \bigcup_{k\in\nat}{}^\Gamma C^{-\infty}(\Lambda, V^k(\sigma_\lambda,\vp))$.
We obtain
\begin{equation}\label{seni}
res\circ \p_\Gamma=\p_B\circ res\ .
\end{equation}

We denote the spaces of germs at $\lambda$ of holomorphic and meromorphic families $\mu\mapsto f_\mu\in C^{-\sharp}(\partial X,V(\sigma_\mu,\varphi))$ and $\mu\mapsto f_\mu\in C^{-\sharp}(B,V_B(\sigma_\mu,\varphi))$ by $\cO_\lambda C^{-\sharp}(\partial X,V(\sigma_.,\varphi))$, $\cM_\lambda C^{-\sharp}(\partial X,V(\sigma_.,\varphi))$ and
$\cO_\lambda C^{-\sharp}(B,V_B(\sigma_.,\varphi))$, $\cM_\lambda C^{-\sharp}(B,V_B(\sigma_.,\varphi))$, respectively. Let $\cM_\lambda^kC^{-\sharp}(\dots)$
be the space of germs of meromorphic families which have a pole of order at
most $k$ at $\lambda$.
Then 
$$\cO_\lambda C^{-\sharp}(\partial X,V(\sigma_.,\varphi))\subset \cM^k_\lambda C^{-\sharp}(\partial X,V(\sigma_.,\varphi))\subset\cM_\lambda C^{-\sharp}(\partial X,V(\sigma_.,\varphi))$$ 
become $\Gamma$-representations in a natural way. We consider the $\Gamma$-equivariant operator
$L_\lambda: \cM_\lambda C^{-\sharp}(\partial X,V(\sigma_.,\varphi))\rightarrow
\cM_\lambda C^{-\sharp}(\partial X,V(\sigma_.,\varphi))$
induced by $f_{\lambda+z\alpha}\mapsto zf_{\lambda+z\alpha}$.

\begin{lem}\label{fear}
There is a $\Gamma$-equivariant map
$$ ev: \cM_\lambda C^{-\sharp}(\partial X,V(\sigma_.,\varphi))\rightarrow
C^{-\sharp}(\partial X, V^+(\sigma_\lambda,\vp))$$
given by
$$ ev(f_\mu)(H):=\res_{z=0} (e^{\langle z\alpha,H\rangle}f_{\lambda+z\alpha})\ ,
\qquad H\in\aaaa\ .$$
Here we consider $f_\mu$ as a family of generalized functions on $G$ with values in $V_\sigma\otimes V_\vp$ satisfying 
$$f_\mu(.p)=(\sigma_\mu(p)^{-1}\otimes\id) f_\mu(.)\ , \qquad\forall p\in P \ . $$
$ev$ induces isomorphisms
\begin{eqnarray}\label{firth}
\cM^k_\lambda C^{-\sharp}(\partial X,V(\sigma_.,\varphi))/ \cO_\lambda C^{-\sharp}(\partial X,V(\sigma_.,\varphi))
\cong C^{-\sharp}(\partial X, V^k(\sigma_\lambda,\vp))\ ,\\
\cM_\lambda C^{-\sharp}(\partial X,V(\sigma_.,\varphi))/ \cO_\lambda C^{-\sharp}(\partial X,V(\sigma_.,\varphi))
\cong C^{-\sharp}(\partial X, V^+(\sigma_\lambda,\vp))\ . \label{firth1}
\end{eqnarray}
We have
\begin{equation}\label{eye} 
ev\circ L_\lambda =\p\circ ev\ .
\end{equation}
The analogous statements for $C^{-\sharp}(B,V^+_B(\sigma_\mu,\varphi))$ are also
true.
\end{lem}
\proof First one checks that $ev(f_\mu)$ satisfies the correct transformation law (\ref{diep}) with respect to the right action of $P$. Then one observes
that (\ref{firth}) is true for $k=1$. Indeed, the preimage $ev^{-1}(f)$ of $f\in C^{-\sharp}(\partial X, V(\sigma_\lambda,\vp))\cong C^{-\sharp}(K\times_M V_\sigma)\otimes V_\vp$ in $\cM^1_\lambda C^{-\sharp}(\partial X,V(\sigma_.,\varphi))/ \cO_\lambda C^{-\sharp}(\partial X,V(\sigma_.,\varphi))$
is given by $\lambda+z\alpha\mapsto \frac{1}{z} f\in  C^{-\sharp}(K\times_M V_\sigma)\otimes V_\vp\cong C^{-\sharp}(\partial X, V(\sigma_{\lambda+z\alpha},\vp))$. In order to complete the proof of the lemma
it suffices to show (\ref{eye}). Let $H_0\in\aaaa$ be the element determined by $\langle \alpha,H_0\rangle =1$. We compute
\begin{eqnarray*}
ev(L_\lambda f_\mu)(rH_0)&=& \res_{z=0} (e^{zr}zf_{\lambda+z\alpha})
=\res_{z=0} (\frac{d}{dr}e^{rz}f_{\lambda+z\alpha})\\
&=&\frac{d}{dr} ev(f_\mu)(rH_0)
= \p (ev(f_\mu))(rH_0)\ .
\end{eqnarray*}
\hB

In view of (\ref{eye}) the map $ev$ can be considered as a
version of taking the principal part of the Laurent expansion
of a meromorphic family $f_\mu\in \cM_\lambda C^{-\sharp}(\partial X,V(\sigma_.,\varphi))$ in a $\Gamma$-equivariant way.

For $l$ sufficiently large we are now able to define a pointwise shifted extension map
$$ ext[l] : C^{-\sharp}(B, V^+_B(\sigma_\lambda,\vp)) \rightarrow {}^\Gamma C^{-\sharp}(\partial X, V^+(\sigma_\lambda,\vp))\ ,\qquad \sharp\in\{\infty,\omega\}\ .$$

\begin{ddd}\label{crit} 
By $k_-=k_-(\sigma_\lambda,\vp)\in\nat_0$ we denote the order 
of the pole of $ext_\mu$ at $\mu=\lambda$.
Let $f\in C^{-\sharp}(B, V^+_B(\sigma_\lambda,\vp))$. We write $f=ev(f_\mu)$ for some $f_\mu\in \cM_\lambda C^{-\sharp}(B,V_B(\sigma_.,\varphi))$. For $l\ge k_-$ we define
$$
ext[l](f):=ev(ext_\mu (L_\lambda^l f_\mu))\in C^{-\sharp}(\partial X, V^+(\sigma_\lambda,\vp))\ .
$$
\end{ddd}
If $f_\mu\in \cO_\lambda C^{-\sharp}(B,V_B(\sigma_.,\varphi))$, then $ext_\mu(L_\lambda^l f_\mu)\in \cO_\lambda C^{-\sharp}(\partial X,V(\sigma_.,\varphi))$. Therefore the analog of (\ref{firth1}) for $B$ shows that $ext[l]$ is well-defined on $C^{-\sharp}(B,V^+_B(\sigma_\mu,\varphi))$. Moreover, $ext[l]$ maps
$C^{-\sharp}(B,V^k_B(\sigma_\lambda,\varphi))$ to ${}^\Gamma C^{-\sharp}(\partial X, V^k(\sigma_\lambda,\vp))$.

Equation (\ref{paul'}) now implies that
\begin{equation}\label{htak}
res\circ ext[l]=\p_B^l\ .
\end{equation}
Moreover, for $l\ge k_-$ and $r\in\nat$ we have by (\ref{eye})
\begin{equation}\label{sanoj}
ext[l+r]=\p_\Gamma^r\circ ext[l]=ext[l]\circ \p_B^r\ .
\end{equation}
Since $\p_B^r$ is surjective (\ref{sanoj}) implies that $\im\: ext[l]$ does not depend on the choice of $l\ge k_-$. 

\begin{ddd} We form
\begin{eqnarray*}
E^+(\sigma_\lambda,\vp)&:=&ext[l]\left(C^{-\omega}(B, V^+_B(\sigma_\lambda,\vp))\right) \subset {}^\Gamma  C^{-\omega}(\partial X, V^+(\sigma_\lambda,\vp))\ ,\\
E_\Lambda^+(\sigma_\lambda,\vp)&:=& E^+(\sigma_\lambda,\vp)\cap {}^\Gamma  C^{-\infty}(\Lambda, V^+(\sigma_\lambda,\vp))
\end{eqnarray*}
for some $l\ge k_-(\sigma_\lambda,\vp)$.
\end{ddd}

\begin{lem}\label{gm}
The space $E_\Lambda^+(\sigma_\lambda,\vp)$ is finite-dimensional, and we have
for any $l\ge k_-$
$$ E_\Lambda^+(\sigma_\lambda,\vp)=ext[l]\left(C^{-\omega}(B, V^l_B(\sigma_\lambda,\vp))\right) \subset {}^\Gamma C^{-\infty}(\Lambda, V^{k_-}(\sigma_\lambda,\vp))\ .$$
\end{lem}
\proof 
By (\ref{htak}) we have for $f\in C^{-\omega}(B, V^+_B(\sigma_\lambda,\vp))$
$$ res\circ ext[l](f)=\p_B^l (f)\ .$$
Thus $ext[l](f)\in \ker res$
if and only if $f\in C^{-\omega}(B, V^l_B(\sigma_\lambda,\vp))$. We conclude that 
$$E_\Lambda^+(\sigma_\lambda,\vp)=ext[l]\left(C^{-\omega}(B, V^l_B(\sigma_\lambda,\vp))\right)\ .$$ 
Hence the space $E_\Lambda^+(\sigma_\lambda,\vp)$
can be considered as the space of principal parts of Laurent expansions at
$\mu=\lambda$ of families of the form
$$ \mu\mapsto  ext_\mu f_\mu\ ,\quad f_\mu\in \cO_\lambda C^{-\omega}(B,V_B(\sigma_.,\varphi))\ .$$
Since the singularity of $ext_\mu$ at $\mu=\lambda$ is finite-dimensional we
therefore see that $E_\Lambda^+(\sigma_\lambda,\vp)$ is finite-dimensional. Using in addition that
$C^{-\infty}(B, V^l_B(\sigma_\lambda,\vp))\subset C^{-\omega}(B, V^l_B(\sigma_\lambda,\vp))$ is dense we conclude that
$$E_\Lambda^+(\sigma_\lambda,\vp)=ext[k_-]\left(C^{-\infty}(B, V^{k_-}_B(\sigma_\lambda,\vp))\right)\subset {}^\Gamma C^{-\infty}(\Lambda, V^{k_-}(\sigma_\lambda,\vp))\ .$$   
\hB
 
We choose a Cartan subalgebra $\taaa$ of $\maaa$. Then $\taaa\oplus\aaaa=:\haaa$ is a Cartan algebra of $\gaaa$. 
We consider the center $\cZ(\gaaa)$ of the universal enveloping algebra of
$\gaaa$. Via the Harish-Chandra isomorphism characters of $\cZ(\gaaa)$ are
parametrized by elements of $\haaa_\C^*/W(\gaaa_\C,\haaa_\C)$, where $W(\gaaa_\C,\haaa_\C)$ is the Weyl group of $(\gaaa_\C,\haaa_\C)$. We choose
a positive root system $\Delta^+$ of $(\gaaa_\C,\haaa_\C)$ such that all roots
with positive restriction to $\aaaa$ are positive. A character $\chi_\nu$, $\nu\in \haaa_\C^*$, is called integral,
if 
\begin{equation}\label{jemi}
2\frac{\langle \nu,\varepsilon\rangle}{\langle\varepsilon,\varepsilon\rangle}\in\Z
\end{equation}
for all $\varepsilon\in\Delta^+$. Here $\langle.,.\rangle$ is a $W(\gaaa_\C,\haaa_\C)$-invariant bilinear scalar product on $\haaa_\C^*$.  

By $\rho_\maaa\in \ii\taaa^*\subset\haaa_\C^*$ we denote half of the sum of the positive $\maaa_\C$-roots. Let $\hat M$ be the set of equivalence classes
of irreducible representations of $M$.
For $\sigma\in \hat{M}$ let $\mu_\sigma\in\ii\taaa^*$ be its highest weight. The infinitesimal character of the principal series representation $\pi^{\sigma,\lambda}$ of $G$, $\sigma\in \hat{M}$, $\lambda\in\aca$, is now given
by $\chi_{\sigma,\lambda}:=\chi_{\mu_\sigma+\rho_\maaa-\lambda}$.
If $X\not= \R H^2$, then we define for $\sigma\in\hat M$
$$I_\sigma:=\{\lambda\in\aaaa^*\:|\:\chi_{\sigma,\lambda}\mbox{ is integral}\}\ .$$
If $X=\R H^2$ and $G=SL(2,\R)$, then $M\cong \Z_2$. Let $\pm 1$ denote the trivial (+), resp. non-trivial (-) irreducible representation of $M$. We define 
$$I_{1}:=(\frac{1}{2}+\Z)\alpha,\quad I_{-1}:=\Z\alpha\ .$$
If $G=PSL(2,\R)$, then $M=\{1\}$, and we define $I_{1}:=(\frac{1}{2}+\Z)\alpha$.
Note that $I_{\sigma^w}=I_\sigma$, so the definition is compatible with 
our previous convention concerning the Weyl-invariance of $\sigma$. 
Let $I_\aaaa\subset\aaaa^*$ be the lattice generated by the short root $\alpha$,
if $2\alpha$ is a root, or by $\alpha/2$, if not. Note that $I_\sigma\subset I_\aaaa$. More precisely, we have either $I_\sigma=2I_\aaaa$ or 
$I_\sigma=\alpha^\prime+2I_\aaaa$, where $\alpha^\prime$ denotes the generator
of $I_\aaaa$.

A character $\chi_\nu$ is called regular, if none of the expressions (\ref{jemi})
vanishes. We call it weakly regular, if it is regular or if $X\ne\R H^n$ and
(\ref{jemi}) vanishes for at most one positive root $\ve$. We set
\begin{eqnarray*}
 I_\sigma^{r}&:=&\{\lambda\in I_\sigma\:|\: \chi_{\sigma,\lambda}\mbox{ is regular}\}\ ,\\
 I_\sigma^{wr}&:=&\{\lambda\in I_\sigma\:|\: \chi_{\sigma,\lambda}\mbox{ is weakly regular}\}\ .
\end{eqnarray*}
Note that $I_\sigma\setminus I_\sigma^r$ is finite.
 
We need the following irreducibility criterion for the principal series
representations \linebreak[4]
$C^{-\infty}(\partial X,V(\sigma_\lambda))$. Compare
related results in \cite{spehvogan80}, \cite{baldonisilvakraljevic80}, 
\cite{collingwood85}.

\begin{lem}\label{sauna}
Let $\sigma$ be irreducible, $\lambda\ne 0$. If $\lambda\not\in I_\sigma^{wr}$,
then $C^{-\infty}(\partial X,V(\sigma_\lambda))$ is irreducible. Assume in addition that $X\not=\R H^{2n+1}$ or that $p_\sigma(0)=0$.
Then for all $\lambda\in I_\sigma^r$ the representation $C^{-\infty}(\partial X,V(\sigma_\lambda))$ is reducible.
\end{lem}
\proof
For $\Ree(\lambda)=0$, $\lambda\ne 0$, the irreducibility of $C^{-\infty}(\partial X,V(\sigma_\lambda))$ is well-known (see e.g. \cite{knapp86}, Prop. 7.2). Thus we can assume that $\Ree(\lambda)\ne 0$.
For $\Ree(\lambda)>0$ ($\Ree(\lambda)<0$) there is a unique nonzero irreducible
$G$-submodule (quotient) of $C^{-\infty}(\partial X,V(\sigma_\lambda))$ which
is given by the image of $\hat J_{\sigma^{w},-\lambda}$ (of
$\hat J_{\sigma,\lambda}$) (see e.g. \cite{knapp86}, Prop. 7.24). (By
a $G$-submodule we always mean a representation on a closed subspace.) In particular, $C^{-\infty}(\partial X,V(\sigma_\lambda))$ is irreducible if and only if $C^{-\infty}(\partial X,V(\sigma^w_{-\lambda}))$ is, and this holds
precisely if $\hat J_{\sigma^{w},-\lambda}$ is surjective. Since $\hat J_{\sigma^{w},-\lambda}$ is always regular and the poles of $\hat J_{\sigma,\lambda}$ are of at most first order the functional
equation (\ref{rob}) now leads to following irreducibility criterion:

Let $\Ree(\lambda)>0$. Then $C^{-\infty}(\partial X,V(\sigma_\lambda))$
(and $C^{-\infty}(\partial X,V(\sigma^w_{-\lambda}))$) is irreducible
if and only if both $\hat J_{\sigma,\lambda}$ and $p_\sigma(\lambda)$ are regular or $p_\sigma(\lambda)=0$.

The poles of $\hat J_{\sigma,\lambda}$ are always contained in $I_\aaaa$ (\cite{knappstein71}, Thm. 3 and Prop. 43).
In addition, for the real hyperbolic spaces one knows that they are contained
$I_\aaaa\setminus I_\sigma$ in case of even dimension and in $I_\sigma$ in
case of odd dimension. For
a $SO(2n-1)$-representation $\sigma$ this follows from \cite{knappstein71}, Proposition 44. In any case, simple invariance considerations imply that the residue of $\hat J_{\sigma,\mu}$, which is a differential operator, vanishes
at the remaining points of $I_\aaaa$.

The Plancherel densities $p_\sigma$
are explicitly known (see e.g. \cite{knapp86}, Prop. 14.26):
\begin{equation}\label{kauz}
p_{\sigma}(\lambda)=c(\sigma)\left\{
\begin{array}{ccl}
\prod_{\ve\in \Delta^+}  \frac{\langle\mu_{\sigma}+\rho_\maaa-\lambda,\ve\rangle}{\langle \ve,\ve\rangle}\ ,&X=\R H^{2n+1}&\\
\tan \left(d\pi \frac{\langle\lambda,\alpha\rangle}{\langle\alpha,\alpha\rangle}\right)
\prod_{\ve\in \Delta^+}  \frac{\langle\mu_{\sigma}+\rho_\maaa-\lambda,\ve\rangle}{\langle \ve,\ve\rangle}, 
&X\ne\R H^{2n+1},&
0\not\in I_\sigma\\
\cot \left(d\pi \frac{\langle\lambda,\alpha\rangle}{\langle\alpha,\alpha\rangle}\right)
\prod_{\ve\in \Delta^+} \frac{\langle\mu_{\sigma}+\rho_\maaa-\lambda,\ve\rangle}{\langle \ve,\ve\rangle}, &X\ne\R H^{2n+1},&
0\in I_\sigma
\end{array} 
\right. \ , 
\end{equation}
where $d=1$ for real hyperbolic spaces, $d=\frac{1}{2}$ else, and $c(\sigma)$ is a positive constant depending on the normalization of the
Haar measure $dn$. 

An easy discussion of the poles and zeroes of $p_\sigma$ now shows that
the above irreducibility criterion implies the first assertion of the lemma.

If $X\ne\R H^{2n+1}$, then $p_\sigma$ has
poles at all $\lambda\in I_\sigma^r$. Thus in this case $C^{-\infty}(\partial X,V(\sigma_\lambda))$ is reducible. For $X=\R H^{2n+1}$ and $p_\sigma(0)=0$
the intertwining operator $\hat J_{\sigma,\lambda}$ has poles at all nonnegative $\lambda\in I_\sigma$
(\cite{knappstein71}, Proposition 44). If, in addition, $\lambda\in I_\sigma^r$,
then $p_\sigma(\lambda)\ne 0$. Thus by the above irreducibility criterion
$C^{-\infty}(\partial X,V(\sigma_\lambda))$ is reducible. This finishes
the proof of the lemma. 
\hB

One can show in addition, that for $X=\C H^n$ the intertwining operator $\hat J_{\sigma,\lambda}$ has a pole at all positive $\lambda\in I^{wr}_\sigma\setminus I_\sigma^{r}$. By definition, we have $I^{wr}_\sigma= I_\sigma^{r}$ for $X=\R H^n$. Thus for real and complex hyperbolic spaces and $p_\sigma(0)=0$ the
representations $C^{-\infty}(\partial X,V(\sigma_\lambda))$ are reducible
for all non-zero $\lambda\in I^{wr}_\sigma$. However, for $X=\HH H^n$ or $\OO H^2$ there are some exceptional pairs $(\sigma,\lambda)$ with $p_\sigma(0)=0$,
$\lambda\in I^{wr}_\sigma$,
and $C^{-\infty}(\partial X,V(\sigma_\lambda))$ irreducible (see \cite{baldonisilvakraljevic80}, \cite{baldonisilvabarbasch83}).

We now introduce a certain subset $I_\sigma^{wr,-}\subset I_\sigma^{wr}$ by 
$$I_\sigma^{wr,-} :=\{\lambda\in I_\sigma^{wr}\:|\: \mbox{there exist } \sigma^\prime\in\hat M, \lambda^\prime>\lambda\mbox{ s.th. } \chi_{\sigma^\prime,\lambda^\prime}=\chi_{\sigma,\lambda}\}\ .$$

\begin{ddd}\label{ere}
We call a pair $(\sigma,\lambda)$ special, if it satisfies one of
the following conditions:
\begin{enumerate}
\item $\Ree(\lambda)\ge 0$ and $\lambda\in I_\sigma^{wr,-}$. 
\item $\Ree(\lambda)<0$ and $\lambda\in I_\aaaa$. 
\end{enumerate}
A pair $(\sigma,\lambda)$ is called very special if $\lambda\in I_\sigma^{wr,-}$.
\end{ddd}

Recall that $I_\sigma^{wr,-}\subset I_\sigma \subset I_\aaaa$. Therefore a very special pair
is special. The following lemma provides some basic information on 
the set
$I_\sigma^{wr,-}$. We set $\aaaa^*_\pm:=\{\lambda\in\aaaa^*\:|\:\pm\lambda\ge 0\}$.

\begin{lem}\label{dnach}
\begin{enumerate}
\item Let $\lambda\in\aaaa^*$. If $\nu=\mu_\sigma+\rho_\maaa+\lambda$ belongs to the closed positive Weyl chamber, i.e., all the scalar products
in (\ref{jemi}) are non-negative, then $\lambda\not\in I_\sigma^{wr,-}$. Under
the additional assumption $\lambda\in I_\sigma^r$ also the opposite implication is true.
\item The set $I_\sigma^{wr,-}\cap\aaaa^*_+$ is finite.
\item $I_\sigma^{wr}\cap\aaaa^*_-\subset I_\sigma^{wr,-}$.
\item If $\sigma=1$ is the trivial representation, then $I_\sigma^{wr,-}\cap\aaaa^*_+=\emptyset$. In other words:
If $\sigma=1$ and $\Ree(\lambda)\ge 0$, then $(\sigma,\lambda)$ is not special.
The same is true for $\sigma=-1$ in case $G=SL(2,\R)$.
\end{enumerate}
\end{lem}
\proof
Assume that $\nu=\mu_\sigma+\rho_\maaa+\lambda$ belongs to the closed positive Weyl chamber. Let $(\sigma^\prime,\lambda^\prime)$ be such that
$\chi_{\sigma^\prime,\lambda^\prime}=\chi_{\sigma,\lambda}$. Then
$\mu_{\sigma^\prime}+\rho_\maaa-\lambda^\prime=w(\mu_\sigma+\rho_\maaa-\lambda)$
for some $w\in W(\gaaa_\C,\haaa_\C)$, hence
$$ w^\theta(\nu)= \mu_{\sigma^\prime}+\rho_\maaa+\lambda^\prime\ ,$$
where $\theta$ is the Cartan involution and $w^\theta:=\theta\circ w
\circ \theta$. Since $\nu$ belongs to the positive Weyl chamber the difference
$\nu-w\nu$ is a non-negative linear combination of positive roots.
Positive roots have non-negative restrictions to $\aaaa$. It follows that
$\lambda-\lambda^\prime=(\nu-w\nu)_{|\aaaa}\ge 0$, hence $\lambda\not\in I_\sigma^{wr,-}$. Vice versa, assume that $\lambda\in I_\sigma^{r}\setminus I_\sigma^{wr,-}$. Define $\nu$ as above. Let $\nu_0$ be the Weyl conjugate
of $\nu$ which belongs to the positive chamber. The condition $\lambda\in I_\sigma^{r}$ implies that 
$\nu_0=\mu_{\sigma^\prime}+\rho_\maaa+\lambda^\prime$ for
some pair $(\sigma^\prime,\lambda^\prime)$. Since $\lambda\not\in I_\sigma^{wr,-}$ we have $0\ge\lambda^\prime-\lambda=(\nu_0-\nu)_{|\aaaa}$. 
Since $\nu_0-\nu$ is a non-negative linear combination of positive roots,
we conclude that $\lambda^\prime=\lambda$ and $\tau:=\nu_0-\nu$ is a non-negative
linear combination of positive $\maaa$-roots. It follows that
$$ \|\mu_{\sigma}\|^2=\|\mu_{\sigma^\prime}\|^2
=\|\mu_\sigma+\tau\|^2
=\|\mu_{\sigma}\|^2+\|\tau\|^2+2\langle\mu_\sigma,\tau\rangle
\ge \|\mu_{\sigma}\|^2+\|\tau\|^2\ .$$
We conclude that $\tau=0$, hence $\nu=\nu_0$, i.e., $\nu$ belongs to the
positive chamber. This finishes the proof of the Assertion 1.

Fix $\sigma$. Then there exists $\lambda_0\in\aaaa^*$ such that
for $\lambda\ge \lambda_0$ the element $\mu_\sigma+\rho_\maaa+\lambda$
belongs to the positive Weyl chamber. Assertion 1 now implies that 
$I_\sigma^{wr,-}\cap\aaaa^*_+$ is contained in the finite set $[0,\lambda_0)\cap I_\sigma$. This proves 2.

Let $\lambda\in I_\sigma^{wr}\cap\aaaa^*_-$. If $\lambda<0$, then $\lambda\in I_\sigma^{wr,-}$, since $\chi_{\sigma^w,-\lambda}=\chi_{\sigma,\lambda}$. Here
$w$ denotes the non-trivial element of the small Weyl group $W(\gaaa,\aaaa)$.
In order to prove Assertion 3 it remains to discuss the case $\lambda=0$.
If $0\in I_\sigma^{r}$, then there exists $\ve\in\Delta^+$ such that 
$\langle\ve,\mu_\sigma+\rho\rangle <0$. Hence $0\in I_\sigma^{wr,-}$ by Assertion 1. 
If $0\in I_\sigma^{wr}\setminus I_\sigma^{r}$, then $X\ne \R H^n$ and $\langle \mu_\sigma+\rho_\maaa,\ve\rangle =0$ if and only if $\ve$ is the real root, i.e.,
$\ve_{|\taaa}=0$. A case by case check shows that there exists an element 
$w\in W(\gaaa_\C,\haaa_\C)$ such that $w\ve\ne \pm \ve$ and $w\ve_{|\aaaa}\ne 0$, i.e., $w\ve$ is not
a root of $\maaa$. This implies that $w(\mu_\sigma+\rho_\maaa)_{|\aaaa}\ne 0$ and
that $w(\mu_\sigma+\rho_\maaa)_{|\taaa}$ is $\maaa$-regular. Therefore we can find
an element $w_1\in W(\gaaa_\C,\haaa_\C)$ which commutes with $\theta$ such that
$w_1w(\mu_\sigma+\rho_\maaa)=\mu_{\sigma^\prime}+\rho_\maaa-\lambda^\prime$ for
some $\sigma^\prime\in\hat M$, $\lambda^\prime >0$. Thus $0\in I_\sigma^{wr,-}$.

We now prove 4. Let $\lambda\in I_1^{wr,-}$. Then $\|\lambda\|^2+\|\rho_\maaa\|^2   =\|\mu_{\sigma^\prime}+\rho_\maaa\|^2+\|\lambda^\prime\|^2$ for some pair $(\sigma^\prime,\lambda^\prime)$ with $\lambda^\prime>\lambda$. Now $\|\mu_{\sigma^\prime}+\rho_\maaa\|^2\ge \|\rho_\maaa\|^2$ implies
$\|\lambda^\prime\|^2\le \|\lambda\|^2$. It follows that $\lambda<0$. The same 
argument applies to the case $G=SL(2,\R)$, $\sigma=-1$.
\hB

The relevance of Definition \ref{ere} becomes manifest by the following lemma.

\begin{lem}\label{less}
Assume that $(\sigma,\lambda)$ is not special. Then 
$$res\circ \hat{J}_{\sigma,\lambda}:{}^\Gamma C^{-\infty}(\Lambda,V(\sigma_\lambda,\vp))\rightarrow C^{\omega}(B,V_B(\sigma_{-\lambda},\vp))$$
is injective.
\end{lem}
\proof
If we replace the condition $\lambda\in I_\sigma^{wr,-}$ in Definition \ref{ere} by $\lambda\in I_\sigma$, then the assertion is covered by Corollary 6.8 in \cite{bunkeolbrich000}. The proof of this statement in \cite{bunkeolbrich000} rests on Lemma 6.7 of
that paper which asserts that under certain conditions on $(\sigma,\lambda)$
for $0\ne f\in C^{-\infty}(\partial X, V(\sigma_\lambda))$ the sections $f$ and $\hat J_{\sigma,\lambda}f$ cannot vanish simultaneously on the same open subset
of $\partial X$. In particular, it is shown there for arbitrary $(\sigma,\lambda)$ that if there exists $f\in C^{-\infty}(\partial X, V(\sigma_\lambda))$ having this vanishing property,
then there are an integer $m>0$, a representation $\sigma^\prime\in\hat M$, and a nontrivial $G$-intertwining 
operator
\begin{equation}\label{city} 
D: C^{-\infty}(\partial X, V(\sigma^\prime_{\lambda+m\alpha}))\rightarrow C^{-\infty}(\partial X, V(\sigma_\lambda))\ .
\end{equation} 
Thus it remains to show that this cannot happen for $\lambda\not\in I_\sigma^{wr,-}$, $\Ree(\lambda)\ge 0$.

Indeed, if such an operator exists, then $C^{-\infty}(\partial X, V(\sigma^\prime_{\lambda+m\alpha}))$ has the same infinitesimal character
as $C^{-\infty}(\partial X, V(\sigma_\lambda))$. By the very definition of
$I_\sigma^{wr,-}$ this is impossible for
$\lambda\in I_\sigma^{wr}\setminus I_\sigma^{wr,-}$. Assume now that 
$\lambda\not\in I_{\sigma}^{wr}$.
Then $\lambda+m\alpha\not
\in I_{\sigma^\prime}^{wr}$. Since $\Ree(\lambda+m\alpha)>0$ it follows from Lemma \ref{sauna} that $C^{-\infty}(\partial X, V(\sigma^\prime_{\lambda+m\alpha}))$ is irreducible.
Thus the image of $D$ is isomorphic to the irreducible non-tempered representation $I^{\sigma^\prime,\lambda+m\alpha}_{-\infty}$. This is impossible
by Langlands classification (see e.g. \cite{wallach88}, Thm. 5.4.1):
Since $I^{\sigma^\prime,\lambda+m\alpha}_{-\infty}$ is non-tempered it follows
that $\Ree(\lambda)>0$, and for $\Ree(\lambda)>0$, $\Ree(\lambda^\prime)>0$,
the equivalence $I^{\sigma,\lambda}_{-\infty}\cong I^{\sigma^\prime,\lambda^\prime}_{-\infty}$ implies $\sigma=\sigma^\prime$,
$\lambda=\lambda^\prime$. This finishes the proof of the lemma.
\hB

\begin{lem}\label{moreless}
Let $\sigma$ be irreducible, $\lambda\ne 0$, and let $I\subset C^{-\infty}(\partial X, V(\sigma_\lambda))$ be a proper
$G$-submodule. If $f\in {}^\Gamma C^{-\infty}(\Lambda,V(\sigma_\lambda,\vp))\cap {}^\Gamma (I\otimes V_\vp)$,
then $res\circ \hat J_{\sigma,\lambda}(f)=0$.
\end{lem}
\proof The same argument which ensures the existence of the operator (\ref{city}) in the proof of Lemma 6.7 in \cite{bunkeolbrich000} shows that any $0\ne f\in I\subset C^{-\infty}(\partial X, V(\sigma_\lambda))$ vanishing on some open $U\subset\partial X$ gives
rise to an intertwining operator 
$D: C^{-\infty}(\partial X, V(\sigma_\lambda))
\rightarrow I$ which is nonzero whenever $(\hat J_{\sigma,\lambda}f)_{|U}\ne 0$.   If $\Ree(\lambda)>0$ ($\Ree(\lambda)<0$), then the $G$-module 
$C^{-\infty}(\partial X, V(\sigma_\lambda))$ has a unique irreducible
submodule (quotient), which, in addition, is not equivalent to any other subquotient of $C^{-\infty}(\partial X, V(\sigma_\lambda))$, whereas it is irreducible for $\Ree(\lambda)=0,
\lambda\ne 0$ (see e.g.
\cite{knapp86}, Theorem 7.2, Theorem 7.24, and Theorem 8.61). 
Thus there cannot exist a non-zero intertwining operator 
$D: C^{-\infty}(\partial X, V(\sigma_\lambda))
\rightarrow I$. The lemma follows.
\hB

In many situations, e.g. for $\Ree(\lambda)<0$, Lemma \ref{moreless} has a
much simpler proof. Indeed, in this case $\hat J_{\sigma,\lambda}$ is regular
and vanishes identically on each proper submodule. However, such a simple 
argument does not work if $\hat J_{\sigma,\mu}$ has a pole at $\mu=\lambda$. 

The following proposition can be considered as a refinement of \cite{bunkeolbrich000}, Proposition 6.11.

\begin{prop}\label{schnauz}
If $(\sigma,\lambda)$ is not special, then 
$${}^\Gamma  C^{-\omega}(\partial X, V^+(\sigma_\lambda,\vp))=E^+(\sigma_\lambda,\vp)\ .$$
In particular, in this case we have
$${}^\Gamma  C^{-\infty}(\Lambda, V^+(\sigma_\lambda,\vp))=E^+_\Lambda(\sigma_\lambda,\vp)=
{}^\Gamma C^{-\infty}(\Lambda, V^{k_-}(\sigma_\lambda,\vp))\ .$$
\end{prop}
\proof   
We set 
\begin{eqnarray*}
E^k(\sigma_\lambda,\vp)&:=&E^+(\sigma_\lambda,\vp)\cap {}^\Gamma  C^{-\omega}(\partial X, V^k(\sigma_\lambda,\vp))\ ,\\
E^k_\Lambda(\sigma_\lambda,\vp)&:=&E^k(\sigma_\lambda,\vp)\cap E^+_\Lambda(\sigma_\lambda,\vp)\ .
\end{eqnarray*}
We consider the restriction map
$$ res_1: {}^\Gamma  C^{-\omega}(\partial X, V^1(\sigma_\lambda,\vp))
\rightarrow C^{-\omega}(B, V^1(\sigma_\lambda,\vp)) \ .$$
The first step in the proof is to show the inequality $\dim E^1_\Lambda(\sigma_\lambda,\varphi)\ge \dim\coker\:res_1$, which is more or less obvious if $ext_\mu$ has a pole of at most first order at $\mu=\lambda$. The general case is more involved. For non-special $(\sigma,\lambda)$ we will 
conclude equality and
${}^\Gamma  C^{-\infty}(\partial X, V^1(\sigma_\lambda,\vp))=E^1(\sigma_\lambda,\vp)$. 
We then proceed by induction on $k$.

For some $l\ge k_-$ we form the space
$$B_\Lambda:= C^{-\omega}(B, V^l(\sigma_\lambda,\vp))/\ker ext[l] \ .$$
Note that because of (\ref{htak}) indeed $\ker ext[l]\subset C^{-\omega}(B, V^l(\sigma_\lambda,\vp))$. Then by Lemma \ref{gm} the map $ext[l]$ induces an isomorphism
$$ e_\Lambda: B_\Lambda\rightarrow E^+_\Lambda(\sigma_\lambda,\vp)$$
satisfying 
$$ \p_\Gamma\circ e_\Lambda = e_\Lambda \circ [\p_B]\ ,$$
where $[\p_B]: B_\Lambda\rightarrow B_\Lambda$ is induced by $\p_B$.
Then $E^1_\Lambda(\sigma_\lambda,\vp)\cong \ker\: [\p_B]$. 
We claim that 
\begin{equation}\label{wd}
\p_B^{l-1} (\ker ext[l])
= res_1\left(E^1(\sigma_\lambda,\vp)\right)\subset \im\: res_1\ .
\end{equation}
Indeed, let $\phi\in \ker ext[l]$
and choose $\psi \in C^{-\omega}(B, V^{l+1}(\sigma_\lambda,\vp))$ such that
$\p_B(\psi)=\phi$. Then $\p_\Gamma\circ ext[l](\psi)=ext[l]\circ \p_B(\psi)=0$, hence  
$ext[l](\psi)\in E^1(\sigma_\lambda,\vp)$,
and $res_1\circ ext[l](\psi)=\p_B^{l-1}(\phi)$. Vice versa, if $ext[l](\psi)\in E^1(\sigma_\lambda,\vp)$, then $\p_B(\psi)\in\ker ext[l]$, and $res_1\circ ext[l](\psi)=\p_B^{l-1}(\p_B(\psi))$. This proves (\ref{wd}).

Now (\ref{wd}) implies that $[\p_B]^{l-1}$ induces a surjective map from
$\coker\: [\p_B]$ onto $\coker\: res_1$. 
Summarizing the above discussion we obtain
\begin{eqnarray}
d:=\dim E^1_\Lambda(\sigma_\lambda,\varphi)&=&\dim \ker\:[\p_B]\nonumber\\
&=&\dim \coker\:[\p_B]\ge\dim\coker\:res_1\ .\label{surprise}
\end{eqnarray}
Set $\tilde d:= \dim E^1_\Lambda(\tilde\sigma_\lambda,\tilde\varphi)$. 
If $(\sigma,\lambda)$ is non-special, then so is 
$(\tilde\sigma,\lambda)$. Thus Proposition \ref{green} combined with Corollary \ref{less} implies 
$$\dim \coker\:res_1
\ge\dim {}^{\Gamma} C^{-\infty}(\Lambda,V(\tilde\sigma_\lambda,\tilde\varphi))\ .$$ 
It follows that
$$d\ge\dim\:\coker(res_1)
\ge \dim {}^{\Gamma} C^{-\infty}(\Lambda,V(\tilde\sigma_\lambda,\tilde\varphi))\ge\tilde d\ .$$
Changing the roles of $(\sigma,\varphi)$ and $(\tilde\sigma,\tilde\varphi)$
we obtain 
$d=\tilde d=\dim {}^{\Gamma} C^{-\infty}(\Lambda,V(\sigma_\lambda,\varphi))$,
i.e.,
\begin{equation}\label{1a}
{}^\Gamma  C^{-\infty}(\Lambda, V^1(\sigma_\lambda,\vp))=E^1_\Lambda(\sigma_\lambda,\vp)\ .
\end{equation}
Moreover, we conclude equality in (\ref{surprise}), hence in (\ref{wd}).
This together with (\ref{1a}) implies
$$
{}^\Gamma  C^{-\omega}(\partial X, V^1(\sigma_\lambda,\vp))=E^1(\sigma_\lambda,\vp)\ .
$$
Let us assume in addition that
$$
{}^\Gamma  C^{-\omega}(\partial X, V^{k-1}(\sigma_\lambda,\vp))=E^{k-1}(\sigma_\lambda,\vp)\ ,
$$
and let $f\in {}^\Gamma  C^{-\omega}(\partial X, V^k(\sigma_\lambda,\vp))$.
Then we can choose $\psi\in C^{-\omega}(B, V^+(\sigma_\lambda,\vp))$ such that
$\p_\Gamma(f)=ext[l](\p_B(\psi))$. This implies $\p_\Gamma(f-ext[l](\psi))=0$.
Thus 
$$f-ext[l](\psi)\in {}^\Gamma  C^{-\omega}(\partial X, V^1(\sigma_\lambda,\vp))= E^1(\sigma_\lambda,\vp)\ .$$ 
Hence $f\in E^k(\sigma_\lambda,\vp)$. 
The proposition now follows by induction.
\hB  

Now it is easy to derive the following

\begin{prop}\label{bartII}
For any pair $(\sigma,\lambda)$ the sequence of inclusions 
$$ {}^\Gamma C^{-\infty}(\Lambda, V^{1}(\sigma_\lambda,\vp))\subset\dots\subset {}^\Gamma C^{-\infty}(\Lambda, V^{k}(\sigma_\lambda,\vp))\subset{}^\Gamma C^{-\infty}(\Lambda, V^{k+1}(\sigma_\lambda,\vp))\subset \dots                             $$
stabilizes at some $k=:k_+(\sigma_\lambda,\vp)$, and 
${}^\Gamma C^{-\infty}(\Lambda, V^+(\sigma_\lambda,\vp))={}^\Gamma C^{-\infty}(\Lambda, V^{k_+}(\sigma_\lambda,\vp))$ is finite-dimensional.
\end{prop}
\proof
We can assume that $\sigma$ is irreducible. Then there exist finite-dimensional $G$-represen\-tations $\pi_{\sigma,\mu}$ satisfying (\ref{str}) and having arbitrary large highest weight $\mu\ge 0$.  We choose one such that $\Ree(\lambda)+\mu\ge 0$. By Lemma \ref{dnach} the pair $(1,\lambda+\mu)$
is not special. The embedding (\ref{str})
can be continued to the extended bundles. We obtain by Proposition \ref{schnauz} $$ i_{\sigma,\mu}\left({}^{\Gamma} C^{-\infty}(\Lambda,V^+(\sigma_\lambda,\varphi))\right)\subset {}^{\Gamma} C^{-\infty}(\Lambda,V^+(1_{\lambda+\mu},\varphi\otimes\pi_{\sigma,\mu}))=E_\Lambda^+(1_{\lambda+\mu},\varphi\otimes\pi_{\sigma,\mu})\ .$$
The proposition now follows since the space on the right hand side is finite-dimensional by Lemma \ref{gm}.
\hB

The assertion of Proposition \ref{schnauz} also holds for all pairs $(\sigma,\lambda)$ which are not very special. 
In order to prove this we first make
a couple of observations which are interesting in their own right.

\begin{kor}\label{iso1}
For fixed $\sigma$ and $\varphi$ the set $\{\lambda\in \aca\:|\:{}^{\Gamma} C^{-\infty}(\Lambda,V^+(\sigma_\lambda,\varphi))\not=\{0\}\}$ is a discrete subset
of the half-plane $\{\lambda\in \aca\:|\:\Ree(\lambda)\le \delta_\Gamma+\delta_\vp\}$.
\end{kor}
\proof By Proposition \ref{schnauz} we have
$$ \{\lambda\in \aca\:|\:{}^{\Gamma} C^{-\infty}(\Lambda,V^+(\sigma_\lambda,\varphi))\not=0\}\subset 
\{\lambda\in \aca\:|\:ext_\mu\mbox{ has a pole at }\mu=\lambda\mbox{ or $(\sigma,\lambda)$ is special} \}\ .$$
By the meromorphy of $ext_\mu$ the set on the right hand side is discrete.

If $f\in {}^{\Gamma} C^{-\infty}(\Lambda,V^k(\sigma_\lambda,\varphi))\setminus {}^{\Gamma} C^{-\infty}(\Lambda,V^{k-1}(\sigma_\lambda,\varphi))$, then
$\p_\Gamma^{k-1}f\in {}^{\Gamma} C^{-\infty}(\Lambda,V(\sigma_\lambda,\varphi))\setminus \{0\}$. Hence ${}^{\Gamma} C^{-\infty}(\Lambda,V^+(\sigma_\lambda,\varphi))\ne \{0\}$ if and only if ${}^{\Gamma} C^{-\infty}(\Lambda,V(\sigma_\lambda,\varphi))\ne \{0\}$. Thus by Corollary \ref{bartI} we have $\{\lambda\in \aca\:|\:{}^{\Gamma} C^{-\infty}(\Lambda,V^+(\sigma_\lambda,\varphi))\not=0\}\subset \{\lambda\in \aca\:|\:\Ree(\lambda)\le \delta_\Gamma+\delta_\vp\}$. This finishes the proof
of the corollary. 
\hB

\begin{kor}\label{iso2}
For $f_\mu \in  {}^{\Gamma} \cM_\lambda C^{-\omega}(\partial X,V(\sigma_.,\varphi))$ we have $ext_\mu\circ res(f_\mu)=f_\mu$.
\end{kor}
\proof We consider $h_\mu:=ext_\mu\circ res(f_\mu)-f_\mu$.
By (\ref{paul'}) we have $res(h_\mu)= 0$. Hence $h_\mu$ is
a $\Gamma$-invariant family which is supported on the limit set.
Corollary \ref{iso1} now implies that $h_\mu= 0$.
\hB

We can now give an alternative description of the space 
$E^+(\sigma_\lambda,\vp)$.

\begin{kor}\label{iso3}  
$
E^+(\sigma_\lambda,\vp)=ev\left({}^{\Gamma} \cM_\lambda C^{-\omega}(\partial X,V(\sigma_.,\varphi)) \right)\ .$
\end{kor}
\proof Let $f\in E^+(\sigma_\lambda,\vp)$. Then by definition of $ext[l]$, $l\ge k_-(\sigma_\lambda,\vp)$, and $E^+(\sigma_\lambda,\vp)$
there exists a family $h_\mu \in \cM_\lambda C^{-\omega}(B,V_B(\sigma_.,\varphi))$ such that $f=ev(ext_\mu(L_\lambda^l h_\mu))$.  
This shows
$E^+(\sigma_\lambda,\vp)\subset ev\left({}^{\Gamma} \cM_\lambda C^{-\omega}(\partial X,V(\sigma_.,\varphi)) \right)$.
The opposite inclusion follows from Corollary \ref{iso2}. Indeed, let 
$f_\mu\in {}^{\Gamma} \cM_\lambda C^{-\omega}(\partial X,V(\sigma_.,\varphi))$
and $h_\mu\in \cM_\lambda C^{-\omega}(B,V_B(\sigma_.,\varphi))$ such that
$L_\lambda^l h_\mu=res(f_\mu)$. Then 
$$
ev(f_\mu)=ev(ext_\mu\circ res(f_\mu))=ev(ext_\mu(L_\lambda^l h_\mu))
=ext[l] (ev(h_\mu))\in E^+(\sigma_\lambda,\vp)\ .
$$
\hB

\begin{prop}\label{zwirbel}
If $(\sigma,\lambda)$ is not very special, then the assertions of Proposition \ref{schnauz} hold true.
\end{prop}
\proof By Proposition \ref{schnauz} we can assume that $(\sigma,\lambda)$
is special, but not very special. This in particular implies that 
$\Ree(\lambda)<0$, hence $\hat J_{\sigma,\lambda}$ is regular, and $\lambda\not\in I_\sigma^{wr}$. By Lemma \ref{sauna} the principal series representations $C^{-\infty}(\partial X, V(\sigma^{\prime}_\mu))$ are irreducible
for each of the (at most two) irreducible components $\sigma^\prime$ of $\sigma$ and $\mu$ in a small neighbourhood of $\lambda$. Thus the intertwining operators $\hat J_{\sigma,\mu}$ are regular
and bijective for $\mu$ near $\lambda$, hence induce isomorphisms
\begin{eqnarray*} 
\hat J: \cM_\lambda C^{-\omega}(\partial X,V(\sigma_.,\varphi))
&\rightarrow&\cM_{-\lambda} C^{-\omega}(\partial X,V(\sigma_.,\varphi))\ ,\\
\hat J: \cO_\lambda C^{-\omega}(\partial X,V(\sigma_.,\varphi))
&\rightarrow&\cO_{-\lambda} C^{-\omega}(\partial X,V(\sigma_.,\varphi))\ ,\\
\hat J_+: {}^\Gamma  C^{-\omega}(\partial X, V^+(\sigma_\lambda,\vp))
&\rightarrow& {}^\Gamma  C^{-\omega}(\partial X, V^+(\sigma_{-\lambda},\vp))
\end{eqnarray*}
satisfying $\hat J_+\circ ev= -ev \circ \hat J$.
Now $-\lambda\not\in I_\sigma$, thus $(\sigma,-\lambda)$ is not special.
We conclude from Proposition \ref{schnauz} and Corollary \ref{iso3} that
\begin{eqnarray*}
\hat J_+\left( {}^\Gamma  C^{-\omega}(\partial X, V^+(\sigma_\lambda,\vp))
\right)&=&{}^\Gamma  C^{-\omega}(\partial X, V^+(\sigma_{-\lambda},\vp))
=E^+(\sigma_{-\lambda},\vp)\\
&=&ev\left( \hat J( \cM_\lambda C^{-\omega}(\partial X,V(\sigma_.,\varphi)))\right)\\
&=&\hat J_+\left(ev (\cM_\lambda C^{-\omega}(\partial X,V(\sigma_.,\varphi)))\right)\\
&=&\hat J_+\left(E^+(\sigma_\lambda,\vp)\right)\ .
\end{eqnarray*}
We obtain ${}^\Gamma  C^{-\omega}(\partial X, V^+(\sigma_\lambda,\vp))=   
E^+(\sigma_\lambda,\vp)$. This proves the proposition.
\hB

We now use Proposition \ref{schnauz} and Corollaries \ref{iso2} and \ref{iso3} in order
to generalize Proposition \ref{green} to the extended bundles. We introduce
a pairing between $C^\omega(B,V_B^k(\sigma_\lambda,\vp))$ and \linebreak[4] $C^{-\omega}(B,V_B^k(\tilde\sigma_{-\lambda},\tilde\vp))$ by
$$
\langle ev(f_\mu), ev(g_\nu)\rangle_k:= \res_{z=0}\langle  f_{\lambda+z\alpha},(L_{-\lambda}^k g)_{-\lambda-z\alpha}\rangle\ ,
$$
$f_\mu\in \cM^k_\lambda C^{\omega}(B,V_B(\sigma_.,\varphi))$, $g_\nu\in \cM^k_{-\lambda} C^{-\omega}(B,V_B(\tilde\sigma_.,\varphi))$.
It is easily checked that this pairing is well-defined and non-degenerate.
Since  
$$ res_\Omega\circ \hat J_{\sigma,\mu}: \cM_\lambda  C^{-\omega}(\Lambda, V(\sigma_.,\vp))
\rightarrow \cM_{-\lambda} C^{\omega}(\Omega, V_B(\sigma_.,\vp))$$ 
maps
germs of holomorphic families to germs of holomorphic families (see the remark following Proposition \ref{green}) it induces as in Proposition \ref{zwirbel}
an operator
$$(res_\Omega\circ\hat J)_+: C^{-\omega}(\Lambda, V^+(\sigma_\lambda,\vp))
\rightarrow C^{\omega}(\Omega, V^+(\sigma_{-\lambda},\vp))$$
and thus an operator 
$$(res\circ\hat J)_+: {}^\Gamma  C^{-\omega}(\Lambda, V^+(\sigma_\lambda,\vp))
\rightarrow C^{\omega}(B, V_B^+(\sigma_{-\lambda},\vp))\ .$$ 
Note that we have to use families $f_\mu\in \cM_\lambda  C^{-\omega}(\Lambda, V(\sigma_.,\vp))$ which cannot be $\Gamma$-invariant (see Corollary \ref{iso1}) in order to define
$(res\circ\hat J)_+$. However, if $\hat J_\mu$ is regular at $\mu=\lambda$,
then $(res\circ\hat J)_+ =res\circ\hat J_+$.
 
\begin{prop}\label{highgreen}
If $f\in {}^\Gamma C^{-\infty}(\Lambda,V^+(\sigma_\lambda,\vp))$ and $k$ is such
that 
$$(res\circ\hat J)_+(f)\in C^{\omega}(B, V_B^k(\sigma_{-\lambda},\vp))\ ,$$
then 
$$ \langle (res\circ \hat J)_+(f), res(\phi)\rangle_k = 0 $$
for all $\phi\in {}^\Gamma C^{-\omega}(\partial X,V^k(\tilde\sigma_{\lambda},\tilde\vp))$.
\end{prop}
\proof First we derive a formula for $(res\circ \hat J)_+(f)$ assuming that
$f\in E_\Lambda^+(\sigma_\lambda,\vp)$. Let $D_\lambda$ be the residue
of $\hat J_\mu$ at $\mu=\lambda$. It is a differential operator and thus
induces a differential operator
$$ D^B_\lambda: C^{-\omega}(B, V_B(\sigma_{\lambda},\vp))\rightarrow
C^{-\omega}(B, V_B(\sigma_{-\lambda},\vp))\subset C^{-\omega}(B, V_B^+(\sigma_{-\lambda},\vp))\ .$$
Choose $f_\mu \in \cO_\lambda C^{-\omega}(B, V_B(\sigma_.,\vp))$
such that $f=ev(ext_\mu f_\mu)$. We claim that
\begin{equation}\label{steif}
(res\circ \hat J)_+(f)= D^B_\lambda f_\lambda- ev(\hat S_\mu f_\mu)\ .
\end{equation}
Let $\tilde f_\mu\in \cM_\lambda  C^{-\omega}(\Lambda, V(\sigma_.,\vp))$
such that $ev(\tilde f_\mu)=f=ev(ext_\mu f_\mu)$, in other words
$$\tilde f_\mu-ext_\mu f_\mu \in \cO_\lambda  C^{-\omega}(\Lambda, V(\sigma_.,\vp))\ .$$
Using that $\hat J_\mu$ has a pole of at most first order at $\mu=\lambda$ we compute
$$
ev(res\circ\hat J_\mu(\tilde f_\mu-ext_\mu f_\mu))=
res\circ D_\lambda (\tilde f_\mu-ext_\mu f_\mu)_{|\mu=\lambda}
=-D^B_\lambda f_\lambda\ .
$$
Now we derive (\ref{steif}) as follows
\begin{eqnarray*}
(res\circ \hat J)_+(f)&=& -ev(res\circ\hat J_\mu(\tilde f_\mu)\\
&=& D^B_\lambda f_\lambda-ev(res\circ\hat J_\mu(ext_\mu f_\mu))\\
&=& D^B_\lambda f_\lambda- ev(\hat S_\mu f_\mu)\ .
\end{eqnarray*}

As in the proof of Proposition \ref{bartII} we consider an embedding
$$ i_{\sigma,\mu}: V(\sigma_\lambda,\vp)\hookrightarrow V(1_{\lambda+\mu},\pi_{\sigma,\mu}\otimes\vp)\ ,$$
$\mu$ such that $\Ree(\lambda)+\mu\ge 0$ (if $\sigma$ is the sum of two
irreducible Weyl-conjugate representations, then also $\pi_{\sigma,\mu}$
is the sum of two irreducible $G$-representations having the same highest
$\aaaa$-weight $\mu$). 
There is a corresponding embedding on the dual side
$$ i_{\tilde\sigma,\mu}: V(\tilde\sigma_\lambda,\vp)\hookrightarrow V(1_{\lambda+\mu},\pi_{\tilde\sigma,\mu}\otimes\tilde\vp)\ .$$
Now standard identities for intertwining operators 
imply (see the proof of Proposition 6.5 in \cite{bunkeolbrich000}, p. 116) 
$$ \langle (res\circ \hat J)_+(f), res(\phi)\rangle_k =
\langle (res\circ \hat J)_+\circ i_{\sigma,\mu}(f), res\circ i_{\tilde\sigma,\mu}(\phi)\rangle_k\ .$$
Therefore it is enough to prove the proposition for non-special pairs
like $(1,\lambda+\mu)$. We can thus assume $f\in E_\Lambda^+(\sigma_\lambda,\vp)$, $\phi\in E^k(\tilde\sigma_\lambda,\vp)$. 

Since $E_\Lambda^+(\sigma_\lambda,\vp)$ is finite-dimensional we can
choose $\tilde f_\mu\in \cO_\lambda  C^{\omega}(B, V_B(\sigma_.,\vp))$.
Then $D^B_\lambda f_\lambda$ and $ev(\hat S_\mu f_\mu)$ are analytic
sections. Write $\phi=ev(ext_\mu\psi_\mu)$ for some $\psi_\mu \in \cM_{\lambda} C^{-\omega}(B,V_B(\tilde\sigma_.,\tilde\varphi))$.
Using the adjointness relation (\ref{muh})
for the scattering matrix we obtain
\begin{eqnarray*}
\langle ev(\hat S_\mu f_\mu), res(\phi)\rangle_k 
&=&\res_{z=0}\langle \hat S_{\lambda-z\alpha}(f_{\lambda-z\alpha}), (L_\lambda^k\psi)_{\lambda-z\alpha}\rangle \\
&=&\res_{z=0}\langle  f_{\lambda-z\alpha}, \hat S_{\lambda-z\alpha}(L_\lambda^k\psi)_{\lambda-z\alpha}\rangle \\
&=&\res_{z=0}\langle f_{\lambda-z\alpha}, res\circ \hat J_{\lambda-z\alpha}
(L_\lambda^k\circ ext(\psi))_{\lambda-z\alpha}\rangle
\end{eqnarray*}
Since $\phi\in E^k(\tilde\sigma_\lambda,\vp)$ all families appearing in the
last line except $\hat J_{\lambda-z\alpha}$ are regular at $z=0$.
We obtain
\begin{eqnarray*}
\langle ev(\hat S_\mu f_\mu), res(\phi)\rangle_k 
&=&\langle f_\lambda, res\circ D_\lambda \circ
\p_\Gamma^{k-1}(\phi)\rangle\\
&=&\langle f_\lambda, D^B_\lambda
\circ\p_B^{k-1}(res(\phi))\rangle\\
&=&\langle D^B_\lambda f_\lambda,
\p_B^{k-1}(res(\phi))\rangle\ .
\end{eqnarray*}
Moreover, for any $f\in C^{\omega}(B, V_B^1(\sigma_{-\lambda},\vp))$,
$g\in C^{-\omega}(B, V_B^k(\tilde\sigma_{\lambda},\tilde\vp))$
we have 
$$ \langle f,g\rangle_k = \langle f, \p_B^{k-1} (g)\rangle\ .$$
We conclude that
$$ \langle ev(\hat S_\mu f_\mu), res(\phi)\rangle_k=
\langle D^B_\lambda f_\lambda, res(\phi)\rangle_k\ .$$
This together with (\ref{steif}) implies the proposition.
\hB

There is a characterization of $E^+(\sigma_\lambda,\vp)$ which is slightly weaker than Corollary \ref{iso3} (but strong enough to imply Propositions 
\ref{zwirbel} and \ref{highgreen}, too). Recall the definition of
$k_+(\sigma_\lambda,\vp)$ from Proposition \ref{bartII}. 

\begin{lem}\label{vorzug}
For any $k\ge k_+(\sigma_\lambda,\vp)$ one has
$$ E^+(\sigma_\lambda,\vp)=\im\p_\Gamma^{k}\ .$$
\end{lem}
\proof
It follows from (\ref{sanoj}) and the surjectivity of $\p^k_B$ that
$E^+(\sigma_\lambda,\vp)\subset \im\p^k_\Gamma$ for any $k\in \nat_0$.
We have to show that for $k\ge k_+$
\begin{equation}\label{kno}
\im\p_\Gamma^k\subset E^+(\sigma_\lambda,\vp)\ .
\end{equation}
Let $f=\p_\Gamma^k(\tilde f)$ for some $\tilde f\in {}^\Gamma C^{-\omega}(\partial X,V^+(\sigma_\lambda,\vp))$. Let $l\ge k_-$, and choose $\psi\in  C^{-\omega}(B,V^+_B(\sigma_\lambda,\vp))$ such that
$\p_B^l(\psi)=res(\tilde f)$. 
Using (\ref{htak}) we compute
$$ res(\tilde f-ext[l](\psi))=\p_B^{l}(\psi)-\p_B^{l}(\psi)=0\ .$$
Hence
$\tilde f-ext[l](\psi)$ is supported on the limit set. It is therefore contained
in the kernel of $\p_\Gamma^{k}$.
Thus by (\ref{sanoj})
$$ f=\p_\Gamma^k(\tilde f)=\p_\Gamma^k(ext[l](\psi))=ext[l](\p_B^k(\psi))\in E^+(\sigma_\lambda,\vp)\ .$$
This shows (\ref{kno}), and hence the lemma.
\hB

For $X=\R H^n$ we will see in Section \ref{sjp} (see also Corollary \ref{preth}) that  $\im\p_\Gamma^{k}=\im\p_\Gamma$ for all $k\in\nat$, hence $E^+(\sigma_\lambda,\vp)=\im\p_\Gamma$. It seems to be likely that this holds
in general.

For the rest of this section we assume that $\Ree(\lambda)\ge 0$ and that $\vp$ is unitary, which in particular implies $\delta_\vp=0$. In this
case much more information on the spaces ${}^{\Gamma} C^{-\infty}(\Lambda,V^+(\sigma_\lambda,\varphi))$ is available. We recall
the following definition from \cite{bunkeolbrich000}, Chapter 7.

\begin{ddd}\label{lenin}
For $\Ree(\lambda)\ge 0$ we define the space
of "stable" invariant distributions supported on the limit set as
$$
U_\Lambda(\sigma_\lambda,\varphi):=\{f\in {}^\Gamma C^{-\infty}(\Lambda,V(\sigma_\lambda,\vp))\:|
\:res\circ\hat J_{\sigma,\lambda}(f)=0\}\ .$$
\end{ddd}

By Lemma \ref{less} the space $U_\Lambda(\sigma_\lambda,\varphi)$ can be nontrivial
only if $\lambda\in I_\sigma^{wr,-}$. For irreducible $\sigma$ and $\Ree(\lambda)>0$ we denote by $I^{\sigma,\lambda}$
the underlying $(\gaaa,K)$-module of the unique irreducible submodule $I^{\sigma,\lambda}_{-\infty}$ of the principal
series
representation $C^{-\infty}(\partial X,V(\sigma_\lambda))$. We now restate
some of the results of \cite{bunkeolbrich000}.

\begin{prop}\label{oldstuff}
We have 
\begin{eqnarray}
{}^{\Gamma} C^{-\infty}(\Lambda,V(\sigma_\lambda,\varphi))&=&
E^1_\Lambda(\sigma_\lambda,\varphi)\oplus U_\Lambda(\sigma_\lambda,\varphi)\ ,\label{gaga}\\ E^+_\Lambda(\sigma_\lambda,\varphi)&=&E^1_\Lambda(\sigma_\lambda,\varphi)\ .
\label{gugu}
\end{eqnarray} 
If ${}^{\Gamma} C^{-\infty}(\Lambda,V^+(\sigma_\lambda,\varphi))$ is nonzero,
then $\lambda\in [0,\delta_\Gamma]\subset\aaaa^*$ and, for $\Ree(\lambda)>0$ and $\sigma$ irreducible,
the $(\gaaa,K)$-module $I^{\sigma,\lambda}$ is unitarizable. The latter condition implies
$p_\sigma(0)=0$.
\end{prop}
\proof Since ${}^{\Gamma} C^{-\infty}(\Lambda,V^+(\sigma_\lambda,\varphi))\ne\{0\}$ if and only if ${}^{\Gamma} C^{-\infty}(\Lambda,V(\sigma_\lambda,\varphi))\ne\{0\}$ (see the proof of Corollary \ref{iso1}) the proposition is a combination of a part of the statements of
Lemma 7.2, Propositions 7.3, 7.4, 7.6, and Proposition 9.2 in  \cite{bunkeolbrich000}.
\hB

For $\Ree(\lambda)>0$ we can show in addition

\begin{prop}\label{newer}
Let $\sigma$ be irreducible and $\Ree(\lambda)>0$. Then
\begin{enumerate}
\item ${}^{\Gamma} C^{-\infty}(\Lambda,V^+(\sigma_\lambda,\varphi))=
{}^{\Gamma} C^{-\infty}(\Lambda,V(\sigma_\lambda,\varphi))\subset {}^\Gamma (I^{\sigma,\lambda}_{-\infty}\otimes V_\vp)$.
\item If $C^{-\omega}(\partial X,V(\sigma_\lambda))$ is reducible, then 
\begin{equation}\label{phone}
{}^{\Gamma} C^{-\infty}(\Lambda,V(\sigma_\lambda,\varphi))=U_\Lambda(\sigma_\lambda,\varphi)\ .
\end{equation} 
In particular, $ext$ is regular at all $\lambda>0$, $\lambda\in I_\sigma^r$.
\end{enumerate}
\end{prop}
\proof Let $I^{\sigma,\lambda}_{-\omega}$ be the closure of $I^{\sigma,\lambda}_{-\infty}$ in $C^{-\omega}(\partial X,V(\sigma_\lambda))$. We claim that $I^{\sigma,\lambda}_{-\omega}$ is the only nonzero irreducible $G$-submodule of $C^{-\omega}(\partial X,V^+(\sigma_\lambda))$. For this we
consider the action of the shifted Casimir operator $\Omega_\lambda:=\Omega-\chi_{\sigma,\lambda}(\Omega)\in \cZ(\gaaa)$. One has
$\chi_{\mu_\sigma+\rho_\maaa-\mu}(\Omega)=\langle \mu,\mu\rangle - |\rho|^2
+|\mu_\sigma+\rho_\maaa|^2-|\rho_\maaa|^2$. Thus $\Omega_\lambda$ acts on
$\cM_\lambda C^{-\omega}(\partial X,V(\sigma_.))$ as multiplication by 
$\langle \mu,\mu\rangle-\langle \lambda,\lambda\rangle$, i.e. as 
$2\langle \lambda,\alpha\rangle L_\lambda+|\alpha|^2 L_\lambda^2$. We conclude from (\ref{eye}) that $\Omega_\lambda$ acts on 
$C^{-\omega}(\partial X,V^+(\sigma_\lambda))$ by
\begin{equation}\label{where}
2\langle \lambda,\alpha\rangle \p+|\alpha|^2 \p^2\ .
\end{equation}
Let 
$f\in C^{-\omega}(\partial X,V^k(\sigma_\lambda))\setminus  C^{-\omega}(\partial X,V^{k-1}(\sigma_\lambda))$ 
for some $k \ge 1$. Then by (\ref{where}) the element $\Omega_\lambda^{k-1} f$ belongs to $C^{-\omega}(\partial X,V^1(\sigma_\lambda)) \setminus\{0\}$. It follows that any nonzero submodule $I\subset C^{-\omega}(\partial X,V^+(\sigma_\lambda))$ has nonzero intersection with $C^{-\omega}(\partial X,V^1(\sigma_\lambda))$. This implies the claim.

We can now argue as in the proof of \cite{bunkeolbrich000}, Proposition 9.4. 
Let $W$ be the dual representation of $C^{-\omega}(\partial X,V^k(\sigma_\lambda))$ (it is isomorphic to $C^{\omega}(\partial X,V^k(\tilde\sigma_{-\lambda})))$. Let $L\subset G$ be a compact set such that $gP\in\Omega$ for all $g\in L$, $f\in {}^{\Gamma} C^{-\infty}(\Lambda,V^k(\sigma_\lambda,\varphi))$ and $\phi\in W$. Then
Lemma \ref{birth} provides a constant $C$ such that for all $a\in A_+$
$$|c_{f,\phi}(gak)|\le C (\alpha(\log a))^{k-1} a^{-(\lambda+\rho)}\ .$$
This implies that $c_{f,.}$ induces an intertwining operator from $W$ to the
unitary $G$-representation
$$ L^2(\Gamma\backslash G,\vp):=\{f:G\rightarrow V_\vp\:|\: f(gx)=\vp(g)f(x)\ \forall g\in\Gamma,\: x\in G,\:\ 
\int_{\Gamma\backslash G} |f(x)|^2\:dx<\infty\}\ .$$
By unitarity the image  of $c_{f,.}$ decomposes into a direct
sum of irreducible representations. But by the above $W$ has the unique irreducible quotient 
$W/(I^{\sigma,\lambda}_{-\omega})^\perp$. 
Thus $c_{f,.}$ factorizes over this quotient, and hence $f\in
{}^\Gamma(I^{\sigma,\lambda}_{-\omega}\otimes V_\vp)$. This finishes the
proof of assertion 1.

The equality 
(\ref{phone}) for reducible principal series representations follows from
Assertion 1 and Lemma \ref{moreless}. Principal series representations with regular integral
infinitesimal character and $p_\sigma(0)=0$ are  reducible by Lemma \ref{sauna}.
If $p_\sigma(0)\ne 0$, then (\ref{phone})
remains true by the vanishing result ${}^{\Gamma} C^{-\infty}(\Lambda,V(\sigma_\lambda,\varphi))=\{0\}$ of Proposition \ref{oldstuff}. The regularity of $ext_\lambda$ at these points now follows
from (\ref{gaga}). This finishes the proof of the proposition.
\hB

One knows that $\delta_\Gamma\le \rho$, where equality holds iff $\Omega=\emptyset$. By the theory of the Patterson-Sullivan measure one always has
\begin{equation}\label{ily}
\dim {}^{\Gamma} C^{-\infty}(\Lambda,V(1_{\delta_\Gamma}))=1\ .
\end{equation}
If $\delta_\Gamma >0$, this together with Proposition \ref{oldstuff} implies that $I^{1,\delta_\Gamma}$ is unitarizable. For $X=\HH H^n$, $n\ge 2$, and 
$X=\OO H^2$ this gives nontrivial restrictions on the possible values of 
$\delta_\Gamma$, an observation first made by Corlette \cite{corlette90}: 
If $\Omega\ne \emptyset$, then 
$\delta_\Gamma\le \rho-\beta$, where $\beta:=2\alpha$ for $X=\HH H^n$ and    
$\beta:=6\alpha$ for $X=\OO H^2$. Using now 
Proposition \ref{newer} we are able to sharpen this result.

\begin{kor}\label{corlette}
Let $\Gamma\subset G$ be a convex cocompact
non-cocompact discrete subgroup, where $G=Sp(1,n)$, $Sp(1,n)/\Z_2$ or 
$F_4^{-20}$, then $\delta_\Gamma<\rho-\beta$.
\end{kor}
\proof In any case, $\rho-\beta$ is the end of the spherical
complementary series of $G$, hence $C^{-\omega}(\partial X, V(1_{\rho-\beta}))$
is reducible. This result goes back to Kostant \cite{kostant69} (see also
\cite{helgason94}, Ch. VI, Thm. 3.6). Proposition \ref{newer}, 2 now says that
${}^{\Gamma} C^{-\infty}(\Lambda,V(1_{\rho-\beta}))=U_\Lambda(1_{\rho-\beta},1)$,
whereas $U_\Lambda(1_{\rho-\beta},1)=\{0\}$ by Lemma \ref{dnach}, 4, and Lemma \ref{less}.
It follows from (\ref{ily}) that $\delta_\Gamma\ne \rho-\beta$. Note that
we have not used the meromorphy of $ext_\lambda$ in the above argument,
thus it also works for the case $X=\OO H^2$.
\hB   

It remains to discuss the point $\lambda=0$.

\begin{prop}\label{new}
We have
\begin{eqnarray}\label{ktion} 
{}^{\Gamma} C^{-\infty}(\Lambda,V^+(\sigma_0,\varphi))&=&
{}^{\Gamma} C^{-\infty}(\Lambda,V(\sigma_0,\varphi))\ ,\\
{}^{\Gamma} C^{-\infty}(\Lambda,V(\sigma_0,\varphi))&=&\left\{\begin{array}{ccc}
E_\Lambda^1(\sigma_0,\vp)& &p_\sigma(0)=0\\
U_\Lambda(\sigma_0,\vp)& &p_\sigma(0)\ne 0
\end{array}\right. \ .\label{lation}
\end{eqnarray} 
In particular, if $p_\sigma(0)\ne 0$, then $E^+_\Lambda(\sigma_0,\varphi)=\{0\}$.
For $X=\R H^n$ we have in addition 
\begin{equation}\label{eja}
{}^{\Gamma} C^{-\infty}(\Lambda,V^+(\sigma_0,\varphi))=U_\Lambda(\sigma_0,\vp)=\{0\}
\end{equation}
whenever $\sigma$ is a faithful
representation of $Spin(n-1)$.
\end{prop}
\proof It was shown in \cite{bunkeolbrich000}, Prop. 7.4, that $ext_\lambda$ is regular at $\lambda=0$, whenever $p_\sigma(0)\ne 0$. If $p_\sigma(0)=0$, then
the explicit expression (\ref{kauz}) for $p_\sigma$ implies that $0\not\in I_\sigma^{wr}$. It follows
that in this case $U_\Lambda(\sigma_0,\vp)=\{0\}$. Equation (\ref{lation}) is now a consequence of (\ref{gaga}). For $p_\sigma(0)=0$ 
Equation (\ref{ktion}) follows from (\ref{gugu}).
The vanishing statement (\ref{eja}) follows from
the observation that in the case in question $p_\sigma(0)\ne 0$ and
$0\not\in I_\sigma^{wr}$.   

It remains to prove (\ref{ktion}) for $p_\sigma(0)\ne 0$.
Fix $k>1$, and let $W$ be the dual representation of $C^{-\omega}(\partial X,V^k(\sigma_0))$ (it is isomorphic to $C^{\omega}(\partial X,V^k(\tilde\sigma_0)))$.
As in the proof of Proposition \ref{newer} it is sufficient to show the following two assertions
\begin{enumerate}
\item[(i)] Let $f\in {}^{\Gamma} C^{-\infty}(\Lambda,V^k(\sigma_0,\varphi))$,
$\phi\in W$.
Then the matrix coefficient $c_{f,\phi}$ is square integrable.
\item[(ii)] All irreducible $G$-submodules of $C^{-\omega}(\partial X,V^k(\sigma_0))$ are contained in $C^{-\omega}(\partial X,V(\sigma_0))$.
\end{enumerate}  

Let $L\subset G$ be a compact set such that $gP\in\Omega$ for all $g\in L$, $f\in {}^{\Gamma} C^{-\infty}(\Lambda,V^k(\sigma_0,\varphi))$ and $\phi\in W$.
Let $\hat J_k$ be the restriction of $\hat J_+$ (see the proof of Proposition \ref{zwirbel}) to $C^{-\omega}(\partial X,V^k(\tilde{\sigma}_0,\tilde\varphi))$. As in \cite{bunkeolbrich000}, Lemma 6.2, (compare also (\ref{still})) one shows that for $a\to\infty$
\begin{equation}\label{blue} 
c_{f,\phi}(gak)= \langle (\hat J_k f)(ga),\tilde\sigma(w)\phi(k^{-1}w)\rangle_{V_\sigma} + O(a^{-\rho-\ve})\ .
\end{equation}
Since $p_{\tilde\sigma}(0)=p_\sigma(0)\ne 0$, the family $ext_\lambda$ is regular at $\lambda=0$, and the restriction map
$$ res: {}^\Gamma C^{-\omega}(\partial X,V^k(\tilde{\sigma}_0,\tilde\varphi))
\rightarrow C^{-\omega}(B,V^k_B(\tilde{\sigma}_0,\tilde\varphi))$$
is surjective. Thus Proposition \ref{highgreen} yields $res\circ\hat J_k (f)=0$. Now (\ref{blue}) implies Assertion (i).

The $G$-representation $Y:=C^{-\omega}(\partial X,V(\sigma_0))$ decomposes
into a direct sum of (two) irreducible representations. Any irreducible
composition factor of $Y^k:=C^{-\omega}(\partial X,V^k(\sigma_0))$ is isomorphic
to one of these. Thus any irreducible submodule of $Y^k$ is the image of an intertwining operator $I\in\Hom_G(Y,Y^k)\cong \Hom_{\gaaa,K}(Y_K,Y^k_K)
\cong \Hom_{MA}(Y_K/\naaa Y_K, V_{\sigma_\rho}\otimes\Pi^k)$. Here the subscript $K$ stands for the underlying $(\gaaa,K)$-module. The latter isomorphism is Casselman's Frobenius reciprocity (\cite{wallach88}, 4.2.2).
For $f\in Y_K$, $\phi\in W_K$ the matrix coefficients $c_{f,\phi}$ satisfy
for $a\to\infty$
$$ c_{f,\phi}(a)= p_{f,\phi} a^{-\rho} + O(a^{-\rho-\ve})\ ,$$
where $p_{f,\phi}$ is a constant not depending on $a$ (see e.g. \cite{bunkeolbrich000}, Lemma 9.1 and Lemma 6.2).
The main theorem in \cite{hechtschmid831} identifies the generalized $\rho$-eigenspace of $\aaaa$ in $Y_K/\naaa Y_K$ with the collection of functionals
$p_{f,.}$ on $W_K$. It follows that $\aaaa$ acts on it by scalars. This
implies that the image of $I^F$ is contained in $V_{\sigma_\rho}\subset V_{\sigma_\rho}\otimes\Pi^k$, where $I^F\in \Hom_{MA}(Y_K/\naaa Y_K, V_{\sigma_\rho}\otimes\Pi^k)$ corresponds to $I$ by Frobenius reciprocity.
Therefore the image of $I$ is contained in $Y$. Assertion (ii) follows.
\hB

We conclude this section by making the results of the preceding three propositions more explicit for the case $X=\R H^n$. If $\sigma\in\hat M$ is
a faithful representation of $Spin(n-1)$, then none of the $(\gaaa,K)$-modules
$I^{\sigma,\lambda}$, $\Ree(\lambda)>0$, is unitarizable (\cite{knappstein71}, 
Proposition 53). If $\sigma\in \hat M$ factorizes over $SO(n-1)$, then in standard coordinates of $i\taaa^*$ it has highest weight 
$$\mu_\sigma=(m_1,m_2,\dots,m_{[\frac{n-1}{2}]})\ ,\quad m_i\in\Z\ ,$$
where
\begin{eqnarray*}
m_1\ge m_2\ge\dots\ge m_\frac{n-2}{2}\ge 0\ ,&\ & n \mbox{ even}\ ,\\
m_1\ge m_2\ge\dots\ge m_\frac{n-3}{2}\ge |m_\frac{n-1}{2}|\ ,&\ & 
n \mbox{ odd}\ .
\end{eqnarray*}
Set
$$\lambda_\sigma:=\rho-(\max\{i\:|\: m_i\ne 0\})\alpha\in [0,\rho]\cap I_\sigma^r\ .$$
This maximum is understood to be zero for $\mu_\sigma=0$.
Then $C^\infty(\partial X,V(\sigma_\lambda))$, $\Ree(\lambda)>0$, is irreducible and unitarizable if and only if $0<\lambda<\lambda_\sigma$. The representations $C^\infty(\partial X,V(\sigma_\lambda))$ with $\lambda$ belonging to this interval are called representations of the complementary series. The only additional unitarizable Langlands submodule is $I^{\sigma,\lambda_\sigma}$, if $\lambda_\sigma>0$ (\cite{knappstein71}, 
Propositions 45, 49, 50).
Thus we obtain

\begin{prop}\label{knut}
Let $X=\R H^n$, $\Ree(\lambda)\ge 0$, $\sigma$ irreducible, and $\vp$ be unitary. If $\sigma$ is a faithful representation of $Spin(n-1)$, then
$$ {}^{\Gamma} C^{-\infty}(\Lambda,V^+(\sigma_\lambda,\varphi))=\{0\}\ .$$
If $\sigma$ factorizes over $SO(n-1)$, then
$$ {}^{\Gamma} C^{-\infty}(\Lambda,V^+(\sigma_\lambda,\varphi))=\left\{
\begin{array}{ccc}
E^1_\Lambda(\sigma_\lambda,\vp)&& 0\le\lambda<\lambda_\sigma,\  \lambda\le\delta_\Gamma\\
U_\Lambda(\sigma_\lambda,\vp)&& 0\le\lambda=\lambda_\sigma\le \delta_\Gamma\\
\{0\}&&\mbox{else}
\end{array}\ .\right. $$
\end{prop} 

\newpage
\section{Cohomology of real hyperbolic manifolds}\label{heart}

We retain the notation and assumptions of the previous two sections.
In addition, we assume that $X=\R H^n$, in particular $Y=\Gamma\backslash X$ is
a hyperbolic manifold of infinite volume and without cusps.

Any finite-dimensional representation $\vp$ of $\Gamma$ defines a flat vector
bundle $E(\vp):=\Gamma\backslash (X\times V_\vp)$ over $Y$. Here $\Gamma$ acts
diagonally on $X\times V_\vp$. Using the results of Section \ref{sur}
we want to relate the cohomology groups $H^p(Y,E(\vp))$
with spaces of invariant distributions on the limit set studied in the previous 
section. Since $X$ is contractible $H^p(Y,E(\vp))$ is isomorphic to the group cohomology $H^p(\Gamma,V_\vp)$ of $\Gamma$ with coefficients in $V_\vp$. The latter
can be defined as the $p$-th right derived functor of the left exact functor $V_\vp\mapsto {}^\Gamma V_\vp$ from the category of $\Gamma$-modules 
to the category of abelian groups. 

We have $\Gamma\subset G$, where $G$ is $SO(1,n)_0$ or its
double cover $Spin(1,n)$. Therefore we have maps $K\rightarrow SO(n)$, $M\rightarrow SO(n-1)$, which
are either isomorphisms or double covers. Let $\gaaa=\kaaa\oplus\paaa$
be the Cartan decomposition. There is a natural action of $K$ on $\paaa^*_\C$
(it is equivalent to the standard representation of $SO(n)$). Note
that $\aaaa\subset\paaa$. Let $\aaaa^\perp\subset \paaa^*_\C$ be its annihilator.
We have representations $\gamma^p$ of $K$ on $V_{\gamma^p}:=\Lambda^p\paaa^*_\C$ and $\sigma^p$ of $M$ on $V_{\sigma^p}:=\Lambda^p\aaaa^\perp$.
The embedding $T:\aaaa^\perp\hookrightarrow \paaa^*_\C$ induces $M$-equivariant
embeddings $T_p:=\Lambda^pT\in \Hom_M(V_{\sigma^p},V_{\gamma^p})$, $p=0,1,\dots,n-1$.
As $G$-homogeneous vector bundles we have
$$ \Lambda^p T_\C^*X=G\times_K V_{\gamma^p}\ ,\qquad  
\Lambda^p T_\C^*(\partial X)=V(\sigma^p_{\rho-p\alpha})\ .$$
Thus as $G$-representations $\Omega^p(X)=C^\infty(G,V_{\gamma^p})^K$,
$\Omega^p(\partial X)=C^\infty(\partial X,V(\sigma^p_{\rho-p\alpha}))$,
$\Omega^p_{-\omega}(\partial X)=C^{-\omega}(\partial X,V(\sigma^p_{\rho-p\alpha}))$. Here $\Omega^p_{-\omega}(\partial X)$ denotes
the space of $p$-hyperforms on $\partial X$, i.e., differential forms with
hyperfunction coefficients. In addition, we are interested in the 
$\Gamma$-module
$$\Omega^p_{-\infty}(\Lambda)=C^{-\infty}(\Lambda,V(\sigma^p_{\rho-p\alpha}))$$
of currents supported on the limit set.

Let $\tau$ be a finite-dimensional representation of $P=MAN$, $T\in\Hom_M(V_\tau,V_{\gamma^p})$. Then one defines
a $G$-equivariant map, called Poisson transform,
$$ P_\tau^T: C^{-\omega}(\partial X,V(\tau))
\rightarrow C^\infty(G,V_{\gamma^p})^K=\Omega^p(X)$$
by
\begin{equation}\label{akt} 
P_\tau^Tf(g):=\int_K \gamma^p(k)T f(gk) \:dk\:\in V_{\gamma^p}\ .
\end{equation}
The integral is a formal notation
meaning that the hyperfunction $\gamma^p(.)Tf(g.)$ on $K$ has to be applied to the constant
function $1$.

For $p\in\{0,1,\dots,n-1\}$ and $\lambda\in\aca$ we first consider the 
following Poisson transforms 
$$ P_{p,\lambda}:= P^{T_p}_{\sigma^p_\lambda}: C^{-\omega}(\partial X,V(\sigma_\lambda))\rightarrow \Omega^p(X)\ .$$
$P_{p,\lambda}$ coincides
with the Poisson transform $P^{T_p}_\lambda$ considered in \cite{bunkeolbrich000}. Moreover, for $z\in\C$ we have $P_{p,z\alpha}=\Phi^{\frac{n-1}{2}-z}_p$, where $\Phi^z_p$ is the transform
introduced and investigated in detail by Gaillard in \cite{Gai1}.

We choose a $G$-invariant Riemannian metric on $X$. It is unique up to a constant factor. We will remember this factor by the induced length of $\alpha\in\aaaa^*\subset\paaa^*\cong T_{eK}^*X$. The metric of constant sectional curvature $-1$ corresponds to $|\alpha|=1$. The Riemannian metric induces $p$-form Laplacians $\Delta$. 

We say that a differential form $h\in C^\infty(G,V_{\gamma^p})^K=\Omega^p(X)$
has moderate growth if there are constants $C$, $d\in \aaaa^*$ such that for all $g\in G$
$$ |h(g)|\le C a_g^{d}\ .$$   
Let $\Omega^p_{mg}(X)\subset\Omega^p(X)$ be the subspace of forms having moderate growth. We can now summarize the main results of \cite{Gai1} by
the following theorem. Compare
also the discussion in \cite{juhl88}.

\begin{theorem}[Gaillard \cite{Gai1}]\label{krah}
The Poisson transforms $P_{p,\lambda}$ define equivariant mappings
between the following $G$-representations
\begin{eqnarray}
P_{p,\lambda}: C^{-\omega}(\partial X,V(\sigma^p_\lambda))
\rightarrow \{\omega\in \Omega^p(X)\:|\: \Delta \omega =
(\langle \rho-p\alpha,\rho-p\alpha\rangle- \langle \lambda,\lambda\rangle)\omega, \delta\omega=0\}\ ,&& \label{sehn}\\
P_{p,\lambda}: C^{-\infty}(\partial X,V(\sigma^p_\lambda))
\rightarrow \{\omega\in \Omega^p_{mg}(X)\:|\: \Delta \omega =
(\langle \rho-p\alpha,\rho-p\alpha\rangle- \langle \lambda,\lambda\rangle)\omega, \delta\omega=0\}\ .&& \label{sucht}
\end{eqnarray}
(\ref{sehn}) and (\ref{sucht}) are isomorphisms if and only if 
$\lambda\not\in \{-\rho+p\alpha\}\cup -\rho-\nat\alpha$. In particular, if $p\ne\frac{n-1}{2}$, then
$$P_{p,\rho-p\alpha}:\Omega_{-\omega}^p(\partial X)\rightarrow  \Omega^p(X)_{\Delta,\delta}$$
is an isomorphism (see Section \ref{sur} for notation). Moreover,
\begin{eqnarray}
(n-1-p) dP_{p,\rho-p\alpha}= 
(n-1-2p)  P_{p+1,\rho-(p+1)\alpha}d\ ,&& p\le n-2\ ,\label{riki}\\
dP_{n-1,-\rho}f=(\int_{\partial X} f)(1-n)\frac{\vol_X}{|\alpha|^n}\ ,&&
\label{ruku}\\
P_{p+1,\rho-(p+1)\alpha}d=\frac{p}{|\alpha|}*P_{p,\rho-p\alpha}*\ ,&& p=\frac{n-1}{2}\ ,
\label{roko}
\end{eqnarray}
where $\vol_X$ denotes the volume form of $X=\R H^n$ with respect to the chosen normalization of the $G$-invariant metric and $*$ denotes the Hodge star operator on $X$ and $\partial X$, respectively.
\end{theorem} 

Note that the star operator on $X$ depends on the invariant metric on $X$, while the star operator on $\partial X$ in the middle dimension is conformally invariant.

Tensoring with $V_\vp$ and restricting to $\Gamma$-invariants $P_{p,\rho-p\alpha}\otimes\id_{V_\vp}$ induces a map
$$P_{p,\rho-p\alpha}: {}^\Gamma(\Omega_{-\omega}^p(\partial X)\otimes V_\vp)
\rightarrow \Omega^p(Y,E(\vp))_{\hat\Delta,\hat\delta}\ ,$$
which is an isomorphism for $p\ne\frac{n-1}{2}$.
Taking (\ref{riki}) and (\ref{ruku}) into account  
(note that $\vol_X$ spans the one-dimensional space $\Omega^n(X)_{\Delta,\delta}$) we see that in the case of even dimension suitably normalized Poisson transforms
provide an isomorphism between the complexes
\begin{equation}\label{kike}
0\rightarrow
{}^\Gamma(\Omega_{-\omega}^0(\partial X)\otimes V_\vp)
\stackrel{d}{\longrightarrow}
{}^\Gamma(\Omega_{-\omega}^1(\partial X)\otimes V_\vp)
\stackrel{d}{\longrightarrow}
\dots \stackrel{d}{\longrightarrow} 
{}^\Gamma(\Omega_{-\omega}^{n-1}(\partial X)\otimes V_\vp)
\stackrel{\int_{\partial X}}{\longrightarrow} {}^\Gamma V_\vp
\rightarrow 0
\end{equation}
and $(\Omega^*(Y,E(\vp))_{\hat\Delta,\hat\delta},d)$. It now follows from Theorem \ref{hop}
that (\ref{kike}) computes the cohomology groups $H^p(Y,E(\vp))$.
Note, however, that for non-cocompact $\Gamma$ all spaces appearing in (\ref{kike}) except ${}^\Gamma V_\vp$ are infinite-dimensional (\cite{bunkeolbrich000}, Corollary 9.7). In order to get something
finite-dimensional we want to reduce things to currents supported
on the limit set using the operators $ext_\lambda$. Such a reduction process
is not compatible with (\ref{kike}) (see, however, Subsection \ref{supp} for
the relation the subcomplex of (\ref{kike}) consisting of currents supported
on the limit set to cohomology with compact support).   
We will therefore use the description of cohomology given in Proposition \ref{infip} instead of the one given in Theorem
\ref{hop}. To do this have to describe the spaces $Z^p(Y,E(\vp))_{(\hat\Delta)}$
as images of suitably defined Poisson transforms.

For odd $n$ we first need a bijective transform into harmonic forms
also for $p=\frac{n-1}{2}$.  There is a decomposition $\sigma^\frac{n-1}{2}=\sigma^+\oplus\sigma^-$ into the $\pm i(-i)^\frac{n-3}{2}$-eigenspaces of
the star operator . Then
$$C^{-\omega}(\partial X,V(\sigma^\pm_0))=\Omega^\pm_{-\omega}(\partial X):=\{\omega\in\Omega^\frac{n-1}{2}_{-\omega}(\partial X)\:|\:
*\omega=\pm i (-i)^\frac{n-3}{2}\omega\}\ .$$
Let
$i^\pm:\sigma^\pm\rightarrow \sigma^\frac{n-1}{2}$ be the corresponding
embedding and $T^\pm=T^\frac{n-1}{2}i^\pm\in \Hom_M(V_\sigma^\pm,V_{\gamma^\frac{n-1}{2}})$. Set 
$$P_{\pm,\lambda}:=
P_{\sigma^\pm_\lambda}^{T^\pm}=
{P_{\frac{n-1}{2},\lambda}}_{|C^{-\omega}(\partial X,V(\sigma^\pm_\lambda))}\ .$$

\begin{kor}\label{raka}  
We have
\begin{eqnarray}
P_{\pm,\lambda}: C^{-\omega}(\partial X,V(\sigma^\pm_\lambda))
\rightarrow \{\omega\in \Omega^\frac{n-1}{2}(X)\:|\: *d \omega =
\pm i (-i)^\frac{n-3}{2}\langle \lambda,\frac{\alpha}{|\alpha|}\rangle\omega, \delta\omega=0\}\ ,&& \label{sehne}\\
P_{\pm,\lambda}: C^{-\infty}(\partial X,V(\sigma^\pm_\lambda))
\rightarrow \{\omega\in \Omega^\frac{n-1}{2}_{mg}(X)\:|\: *d \omega =
\pm  i (-i)^\frac{n-3}{2}\langle \lambda,\frac{\alpha}{|\alpha|}\rangle\omega, \delta\omega=0\}\ .&& \label{suechte}
\end{eqnarray}
(\ref{sehne}) and (\ref{suechte}) are isomorphisms for 
$\lambda\not\in -\rho-\nat\alpha$. In particular,
$$P_{\pm,0}:\Omega^\pm_{-\omega}(\partial X)\rightarrow  \ker d\cap\ker\delta\subset\Omega^\frac{n-1}{2}(X)$$
is an isomorphism.
\end{kor}
\proof One can identify $V_{\sigma^p}=\Lambda^p\naaa^*_\C$, $X=NA$, and,
using $NA$-invariant forms,
$\Omega^p(X)=C^\infty(NA)\otimes \Lambda^p(\aca\oplus\naaa^*_\C)$.
Let $\omega_p\in \Lambda^p\naaa^*_\C$, and $\delta_{\omega_p}\in C^{-\omega}(\partial X,V(\sigma^p_\lambda))$ be the corresponding
delta distribution at $eM\in K/M=\partial X$. Then $P_{p,\lambda}\delta_{\omega_p}(na)=a^{\lambda+\rho}\omega_p$, and one easily computes
$$*d(a^{\lambda+\rho}\omega_p)=\langle \lambda+\rho-p\alpha,\frac{\alpha}{|\alpha|}\rangle a^{\lambda+\rho}*_{\sigma^p}\omega_p\ .$$
Since the delta distributions at $eM$ with values in $V_{\sigma^\pm}$ generate
the $G$-representation \linebreak[4]
$C^{-\omega}(\partial X,V(\sigma^\pm_\lambda))$ it follows
that the image of $P_{\pm,\lambda}$ is contained in the eigenspace appearing
in (\ref{sehne}). In addition, since $(*d)^2=(-1)^\frac{n(n-1)}{2}\delta d$ we see that for $p=\frac{n-1}{2}$, $\lambda\ne 0$ the right hand side of (\ref{sehn}) is the
direct sum of the two eigenspaces appearing in (\ref{sehne}).
Thus for $\lambda\ne 0$ the corollary is a direct consequence of 
Theorem \ref{krah}.

For the case $\lambda=0$ we observe that the Hodge decomposition
$\Omega^\frac{n-1}{2}(\partial X)=\im\: d \oplus \im *d$ implies that
$d: C^{-\omega}(\partial X,V(\sigma^\pm_0))\rightarrow \ker D_\frac{n+1}{2}\subset
\Omega^\frac{n+1}{2}(\partial X)$ is bijective, where $D_\frac{n+1}{2}=d$ for $n>3$ and
$D_\frac{n+1}{2}=\int_{\partial X}$ for $n=3$. By Theorem \ref{krah} including
(\ref{riki}) and (\ref{ruku}) the Poisson transform $P_{_\frac{n+1}{2},-\alpha}$ provides an isomorphism between $\ker D_\frac{n+1}{2}$ and $\ker d\cap \ker \delta \subset \Omega^\frac{n+1}{2}(X)$. Now (\ref{roko}) shows that $P_{\pm,0}$
is an isomorphism.
\hB

Now we want to extend the Poisson transforms to the larger bundles
$V^k(\sigma^p_\lambda)$ and $V^+(\sigma^p_\lambda)$. Since we are only interested in forms annihilated
by a power of the Laplacian we restrict ourselves to the case $\lambda=\rho-p\alpha$. Set $S_p:=*_{\gamma^p} T_{n-p}\in \Hom_M(V_{\sigma^{n-p}},V_{\gamma^p})$. We define $q: \Pi\rightarrow \C$ by $\Pi\ni h\mapsto h(0)$. Note that $q$ is not $A$-equivariant with respect to the representation $1^+$. However,
it gives rise to an elements $S_p^k:=S_p\otimes q\in \Hom_M(V_{\sigma^{n-p}}\otimes \Pi^k,V_{\gamma^p})$, and hence to a Poisson
transforms
\begin{eqnarray*} 
\dot P^+_p: C^{-\omega}(\partial X,V^+(\sigma^{n-p}_{(p-1)\alpha-\rho}))
&\rightarrow &\Omega^p(X)\ ,\qquad p=1,\dots,n\ .\\
\dot P^+_{p|C^{-\omega}(\partial X,V^k(\sigma^{n-p}_{(p-1)\alpha-\rho}))}=\dot P^k_p&:=& P^{S_p^k}_{\sigma^{n-p}_{(p-1)\alpha-\rho}\otimes 1^k}\ .
\end{eqnarray*}
Observe that $\dot P^1_p=*P_{n-p,(p-1)\alpha-\rho}$. More generally,
using $\D q\circ ev(f_\mu)=\res_{z=0} f_{(p-1+z)\alpha-\rho}$, $f_\mu\in \cM_{(p-1)\alpha-\rho} C^{-\omega}(\partial X, V(\sigma^{n-p}_.))$, we see that
\begin{equation}\label{gugel}
\dot P^+_p\circ ev(f_\mu)
=\res_{z=0}
*(P_{n-p,(p-1+z)\alpha-\rho}f_{(p-1+z)\alpha-\rho})\ .
\end{equation} 
An analogous construction provides Poisson transforms $\dot P^+_\pm$, $\dot P^k_\pm$,
$$ \dot P^+_\pm:
C^{-\omega}(\partial X,V^+(\sigma^\pm_0)
\rightarrow \Omega^\frac{n+1}{2}(X)\ ,$$
characterized by
\begin{equation}\label{hupf}
\dot P^+_\pm\circ ev(f_\mu)
=\res_{z=0}
*(P_{\pm,z\alpha}f_{z\alpha})\ ,\qquad f_\mu\in \cM_0 C^{-\omega}(\partial X, V(\sigma^{\pm}_.))
\end{equation}
We now have

\begin{prop}\label{bow}
The Poisson transforms $\dot P^+_p$, $\dot P^k_p$, $p\ne\frac{n+1}{2}$, $\dot P^+_\pm$, $\dot P^k_\pm$, provide $G$-equivariant isomorphisms
\begin{eqnarray*}
\dot P^k_p&:& C^{-\omega}(\partial X, V^k(\sigma^{n-p}_{(p-1)\alpha-\rho}))\stackrel{\cong }\longrightarrow\{\omega\in\Omega^p(X)\:|\:\Delta^k\omega=0, d\omega=0\}\ ,\\
\dot P^k_\pm&:& C^{-\omega}(\partial X, V^k(\sigma^{\pm}_0))\stackrel{\cong }\longrightarrow \{\omega\in\Omega^\frac{n+1}{2}(X)\:|\: (d*)^k\omega=0, d\omega=0\}\ ,
\end{eqnarray*}
and
\begin{eqnarray*}
\dot P^+_p&:& C^{-\omega}(\partial X, V^+(\sigma^{n-p}_{(p-1)\alpha-\rho}))\stackrel{\cong }\longrightarrow Z^p(X)_{(\Delta)}\ ,\\
\dot P^+_\pm&:& C^{-\omega}(\partial X, V^+(\sigma^\pm_0))\stackrel{\cong }\longrightarrow Z^\frac{n+1}{2}(X)_{(\Delta)}
\end{eqnarray*}
satisfying 
\begin{eqnarray}\label{mupf}
\Delta\circ \dot P^+_p &=&|\alpha|^2 \dot P^+_p\circ \p\circ ((n+1-2p)\:\id-\p)\ ,\\   
\Delta\circ \dot P^+_\pm &=&-|\alpha|^2 \dot P^+_\pm\circ\p^2\ ,\label{mapf}\\
d*\circ \dot P^+_\pm &=& \pm i (-i)^\frac{n-3}{2}|\alpha|\dot P^+_\pm\circ\p \label{mipf} \ .
\end{eqnarray}
The corresponding statements for distribution sections on the boundary
and forms of moderate growth on $X$ are also true.
\end{prop}
\proof 
Using (\ref{gugel}), (\ref{sehn}), and (\ref{eye}) we compute
\begin{eqnarray*}
d \dot P^+_p\circ ev(f_\mu)&=& \res_{z=0}\:
d*(P_{n-p,(p-1+z)\alpha-\rho}f_{(p-1+z)\alpha-\rho})\\
&=& \res_{z=0}
*\delta (P_{n-p,(p-1+z)\alpha-\rho}f_{(p-1+z)\alpha-\rho})\\
&=&0\ ,\\  
\Delta \dot P^+_p\circ ev(f_\mu)&=& \res_{z=0}
\:\Delta*(P_{n-p,(p-1+z)\alpha-\rho}f_{(p-1+z)\alpha-\rho})\\
&=& \res_{z=0}
*\Delta(P_{n-p,(p-1+z)\alpha-\rho}f_{(p-1+z)\alpha-\rho})\\
&=&|\alpha|^2 \res_{z=0}
\: z(n+1-2p-z)*P_{n-p,(p-1+z)\alpha-\rho}
f_{(p-1+z)\alpha-\rho}\\
&=&|\alpha|^2\dot P^+_p\circ ev\circ L_{(p-1)\alpha-\rho}((n+1-2p)\id-L_{(p-1)\alpha-\rho})(f_\mu)\\
&=&\dot P^+_p\circ|\alpha|^2\p((n+1-2p)\id-\p)\circ ev(f_\mu)\ .
\end{eqnarray*}
This shows (\ref{mupf}) and that $\im \dot P^+_p$ consists of closed
forms. The desired mapping properties of $\dot P^k_p$ and $\dot P^+_p$
now follow from Theorem \ref{krah} by induction on $k$. (\ref{mipf}) and
the further mapping properties of $\dot P^k_\pm$ and $\dot P^+_\pm$
are proved in an analogous way using (\ref{hupf}) and Corollary \ref{raka}.
\hB

Tensoring with $V_\vp$ and restricting to $\Gamma$-invariants we obtain

\begin{kor}\label{preth}
The Poisson transforms $\dot P^+_p$, $\dot P^+_\pm$,
induce isomorphisms
\begin{eqnarray*}
\dot P^+_p: {}^\Gamma C^{-\omega}(\partial X, V^+(\sigma^{n-p}_{(p-1)\alpha-\rho},\vp))&\stackrel{\cong }\longrightarrow&
Z^p(Y,E(\vp))_{(\hat\Delta)}\ ,\qquad p\ne\frac{n+1}{2}\ ,\\
\dot P^+_\pm: {}^\Gamma C^{-\omega}(\partial X, V^+(\sigma^\pm_0,\vp))&\stackrel{\cong }\longrightarrow& 
Z^\frac{n+1}{2}(Y,E(\vp))_{(\hat\Delta)}\ ,
\end{eqnarray*}
and
\begin{eqnarray}
{}^\Gamma C^{-\omega}(\partial X, V^+(\sigma^{n-p}_{(p-1)\alpha-\rho},\vp))/\im\p_\Gamma
&\stackrel{\cong }\longrightarrow&
H^p(Y,E(\vp))\ ,\qquad p\ne\frac{n+1}{2}\ ,\label{vo}\\
{}^\Gamma C^{-\omega}(\partial X, V^+(\sigma^\pm_0,\vp))/\im\p_\Gamma
&\stackrel{\cong }\longrightarrow& 
H^\frac{n+1}{2}(Y,E(\vp))\ .\label{moe}
\end{eqnarray}
In all cases under consideration we have for all $k\in\nat$  
\begin{equation}\label{mu}
\im\p_\Gamma=\im\p_\Gamma^k=E^+(\sigma^._.,\vp)\ .
\end{equation}
Moreover, the maps
$$
{}^\Gamma C^{-\omega}(\partial X, V(\sigma^\pm_0,\vp))\rightarrow
{}^\Gamma C^{-\omega}(\partial X, V^+(\sigma^\pm_0,\vp))/E^+(\sigma^\pm_0,\vp)
\cong H^\frac{n+1}{2}(Y,E(\vp))$$
and
$$ Z^{n-p}(\vp)\rightarrow {}^\Gamma C^{-\omega}(\partial X, V^+(\sigma^{n-p}_{(p-1)\alpha-\rho},\vp))/E^+(\sigma^{n-p}_{(p-1)\alpha-\rho},\vp)\cong H^p(Y,E(\vp))\ ,$$
where $Z^{n-p}(\vp)\subset {}^\Gamma C^{-\omega}(\partial X, V(\sigma^{n-p}_{(p-1)\alpha-\rho},\vp))$ denotes the space of $(n-p)$-cocycles in (\ref{kike}), are surjective.
\end{kor}
\proof 
The first two isomorphisms are a direct consequence of Proposition \ref{bow}. 
(\ref{vo}), (\ref{moe}) and the first equation in (\ref{mu}) follow from Proposition \ref{infip} and Equations (\ref{mupf}) and (\ref{mipf}). Indeed, if $p\ne\frac{n+1}{2}$, then the map $(n+1-2p)\:\id-\p_\Gamma$
is invertible. Hence in this case
$\hat\Delta^k\left(Z^p(Y,E(\vp))_{(\hat\Delta)}\right)\cong \im \p_\Gamma^k$. 
For $p=\frac{n+1}{2}$ we have $\im\p_\Gamma^k\cong (d*)^k\left(Z^p(Y,E(\vp))_{(\hat\Delta)}\right)$. Since $d*\left(Z^p(Y,E(\vp))_{(\hat\Delta)}\right)$ consists of exact forms Proposition \ref{infip} implies 
$$d*\left(Z^p(Y,E(\vp))_{(\hat\Delta)}\right)\subset \hat\Delta\left(Z^p(Y,E(\vp))_{(\hat\Delta)}\right)= (d*)^2\left(Z^p(Y,E(\vp))_{(\hat\Delta)}\right)\ .$$
Now equality follows. 

The second equation in (\ref{mu}) is a consequence of Lemma \ref{vorzug}.

The surjectivity assertions rely on the fact that any cohomology class
has a closed and coclosed representative (see (\ref{ikp})). Indeed, by
the above discussion of (\ref{kike}) we have for $p\ne\frac{n+1}{2}$ 
$$ \dot P^1_p: Z^{n-p}(\vp)\stackrel{\cong}\longrightarrow \ker d\cap \ker\hat \delta\subset \Omega^p(Y,E(\vp))\ ,$$
while
$ \dot P^1_\pm$ maps ${}^\Gamma C^{-\omega}(\partial X, V(\sigma^\pm_0,\vp))$
bijectively to 
$\ker d\cap \ker\hat \delta\subset \Omega^\frac{n+1}{2}(Y,E(\vp))$.
The proof of the corollary is now complete.
\hB

Recall the definitions of the integers $0\le k_-(\sigma_\lambda,\vp)\le k_+(\sigma_\lambda,\vp)<\infty$ and of the space
$E^+_\Lambda(\sigma_\lambda,\vp)$ from Section \ref{round}. It was shown
there that
$$ E^+_\Lambda(\sigma_\lambda,\vp)\subset 
{}^\Gamma C^{-\infty}(\Lambda, V^{k_-}(\sigma_\lambda,\vp))
\subset {}^\Gamma C^{-\infty}(\Lambda, V^{k_+}(\sigma_\lambda,\vp))
={}^\Gamma C^{-\infty}(\Lambda, V^+(\sigma_\lambda,\vp)) $$
and that these spaces are finite-dimensional. We will also need the space
$$ Z^p_\Lambda(\vp):=Z^p(\vp)\cap (\Omega^{p}_{-\infty}(\Lambda)\otimes V_\vp)$$
of such cocycles of (\ref{kike}) which are supported on the limit set. 
Now we can state
the first main theorem of the paper.

\begin{theorem}\label{main1}
The Poisson transforms $\dot P^+_p$, $\dot P^+_\pm$ induce isomorphisms
\begin{eqnarray*} 
H^p(Y,E(\vp))&\cong& {}^\Gamma C^{-\infty}(\Lambda, V^+(\sigma^{n-p}_{(p-1)\alpha-\rho},\vp))/E^+_\Lambda(\sigma^{n-p}_{(p-1)\alpha-\rho},\vp)\ ,\quad 
p=1,\dots, n\ ,\ p\not=\frac{n+1}{2}\ ,\\
H^\frac{n+1}{2}(Y,E(\vp))&\cong& {}^\Gamma C^{-\infty}(\Lambda, V^+(\sigma^{\pm}_0,\vp))/E^+_\Lambda(\sigma^\pm_0,\vp)\ .
\end{eqnarray*}
Moreover, we have 
\begin{eqnarray*}
k_+(\sigma^{n-p}_{(p-1)\alpha-\rho},\vp))&\le&
k_-(\sigma^{n-p}_{(p-1)\alpha-\rho},\vp))+1\ ,\\
k_+(\sigma^{\pm}_0,\vp))&\le& k_-(\sigma^{\pm}_0,\vp))+1\ .
\end{eqnarray*}
If $\vp$ is unitary and $p\ge\frac{n+1}{2}$, then $k_-(\sigma^{n-p}_{(p-1)\alpha-\rho},\vp))=0$ and
\begin{eqnarray*}
H^p(Y,E(\vp))&\cong & \left\{\begin{array}{clc}
{}^\Gamma (\Omega^{n-p}_{-\infty}(\Lambda)\otimes V_\vp)=
U_\Lambda(\sigma^{n-p}_{(p-1)\alpha-\rho},\vp)\ ,&p>\frac{n+1}{2}\\    
{}^\Gamma (\Omega^\pm_{-\infty}(\Lambda)\otimes V_\vp)=
U_\Lambda(\sigma^\pm_0,\vp)\ ,&p=\frac{n+1}{2} \end{array}\right.\\
&\cong& Z^{n-p}_\Lambda(\vp)\ .
\end{eqnarray*}
\end{theorem}
\proof The theorem will be a consequence of Corollary \ref{preth}. The assertion concerning unitary $\vp$ will then follow from Proposition \ref{knut}.
In the
following one should replace $n-p$ by $\pm$, if $p=\frac{n+1}{2}$. First of all it is clear that
the natural map
$$ {}^\Gamma C^{-\infty}(\Lambda, V^+(\sigma^{n-p}_{(p-1)\alpha-\rho},\vp))/E^+_\Lambda(\sigma^{n-p}_{(p-1)\alpha-\rho},\vp)\rightarrow {}^\Gamma C^{-\omega}(\partial X, V^+(\sigma^{n-p}_{(p-1)\alpha-\rho},\vp))/E^+(\sigma^{n-p}_{(p-1)\alpha-\rho},\vp)$$
is injective. In view of (\ref{vo}), (\ref{moe}) and (\ref{mu}) it is therefore sufficient
to show that the natural map
$$ {}^\Gamma C^{-\infty}(\Lambda, V^{k_-+1}(\sigma^{n-p}_{(p-1)\alpha-\rho},\vp))/E^+_\Lambda(\sigma^{n-p}_{(p-1)\alpha-\rho},\vp)\rightarrow {}^\Gamma C^{-\omega}(\partial X, V^+(\sigma^{n-p}_{(p-1)\alpha-\rho},\vp))/E^+(\sigma^{n-p}_{(p-1)\alpha-\rho},\vp)$$
is surjective.

The last assertion of Corollary \ref{preth} tells us that any element in the
quotient
$${}^\Gamma C^{-\omega}(\partial X, V^+(\sigma^{n-p}_{(p-1)\alpha-\rho},\vp))/
E^+(\sigma^{n-p}_{(p-1)\alpha-\rho},\vp)$$ 
can be represented by an element $f\in {}^\Gamma C^{-\omega}(\partial X, V(\sigma^{n-p}_{(p-1)\alpha-\rho},\vp))$. Choose 
$$\psi \in
C^{-\omega}(B, V^{k_-+1}_B(\sigma^{n-p}_{(p-1)\alpha-\rho},\vp))$$
such that
$$\p_B^{k_-}(\psi)=res(f)\ .$$
Then $f^\prime:=f-ext[k_-](\psi)\in {}^\Gamma C^{-\infty}(\Lambda, V^{k_-+1}(\sigma^{n-p}_{(p-1)\alpha-\rho},\vp))$. Moreover,
$$ f^\prime \equiv f \quad \mod\: E^+(\sigma^{n-p}_{(p-1)\alpha-\rho},\vp)\ .$$
This proves the desired surjectivity, and hence the theorem.
\hB

For $p\ge \frac{n+1}{2}$ and trivial $\vp$ the theorem is a fairly direct consequence
of the results of Mazzeo/Phillips \cite{mazzeophillips90} concerning the
$L^2$-cohomology of $Y$. This was already noted by Patterson \cite{patterson93} and was also observed
by Lott \cite{lott00}. Indeed, this observation, among other things, led Patterson \cite{patterson93} to his general conjecture (\ref{02}) relating the cohomology
groups $H^p(\Gamma,F)$, where $F$ is a finite-dimensional representation of $G$,
to invariant distributions on the limit set. This conjecture will be treated
on the base of Theorem \ref{main1} in Section \ref{twist}. There we will also see that the isomorphism $H^p(Y,E(\vp))\cong Z^{n-p}_\Lambda(\vp)$ holds
for $p=\frac{n}{2}$, too. $L^2$-cohomology will be discussed in Section \ref{mazzeo}.

It would be more natural not to exclude the case of cocompact $\Gamma$ from
the considerations. Recall that in this case $\Lambda=\partial X$. In fact, for unitary $\vp$ the analog of Theorem \ref{main1} is a simple consequence of Theorem \ref{krah} and classical Hodge theory as has been already observed
in \cite{Gai1}. For nonunitary $\vp$ the assertion will be less
precise since Proposition \ref{infip} does not hold in this case. However,
combining Proposition \ref{comp} with Proposition \ref{bow} and using classical
Hodge theory for unitary $\vp$ we obtain

\begin{prop}\label{koks}
There are integers $k_+(\sigma^{n-p}_{(p-1)\alpha-\rho},\vp)$, $k_+(\sigma^\pm_0,\vp)$ such that
\begin{eqnarray*}
{}^\Gamma C^{-\infty}(\Lambda, V^+(\sigma^{n-p}_{(p-1)\alpha-\rho},\vp))&=&{}^\Gamma C^{-\infty}(\Lambda, V^{k_+}(\sigma^{n-p}_{(p-1)\alpha-\rho},\vp))\ ,\\
{}^\Gamma C^{-\infty}(\Lambda, V^+(\sigma^\pm_0,\vp))&=&{}^\Gamma C^{-\infty}(\Lambda, V^{k_+}(\sigma^\pm_0,\vp))
\end{eqnarray*}  
as in the noncocompact case. These spaces are finite-dimensional, and the
Poisson transforms $\dot P^+_p$, $\dot P^+_\pm$ induce surjections
$$ 
 {}^\Gamma C^{-\infty}(\Lambda, V^+(\sigma^{n-p}_{(p-1)\alpha-\rho},\vp))\rightarrow
H^p(Y,E(\vp))\ ,\quad 
{}^\Gamma C^{-\infty}(\Lambda, V^+(\sigma^{\pm}_0,\vp))\rightarrow H^\frac{n+1}{2}(Y,E(\vp)) \ .
$$
If $\vp$ is unitary, then Theorem \ref{main1} holds with the convention
$E^+_\Lambda(\sigma^._.,\vp)=\{0\}$, $k_-(\sigma^._.,\vp)=0$. Moreover,
in this case we have
for all $p$
$$ H^p(Y,E(\vp))\cong Z^{n-p}_\Lambda(\vp)\cong Z^p_\Lambda(\vp)\ .$$
\end{prop}

\newpage
\section{Finite-dimensional $G$-representations as coefficients}\label{twist}

In this section we continue the investigation of the
cohomology groups $H^p(\Gamma,V_{\psi})\cong H^p(Y, E(\psi))$ for
discrete subgroups $\Gamma$ acting convex cocompactly on $X=\R H^n$.
We want to understand more precisely the situation
that the $\Gamma$-representation 
$\psi$ has the form $\psi=\pi_{|\Gamma}\otimes\vp$, where
$(\pi,F)$ is a finite-dimensional representation of $G$. This understanding 
will
be achieved by applying to Theorem \ref{main1} a variant of the translation functor, i.e. tensoring
with $F$ followed by the projection to the generalized infinitesimal character
of $F$. It will bring into play
all the spaces ${}^{\Gamma} C^{-\infty}(\Lambda,V^+(\sigma_\lambda,\varphi))$,
$\sigma\in\hat M$, $\lambda\in I_\sigma^r$.
Since any finite-dimensional $G$-representation is
semisimple we can and will assume that $F$ is irreducible. Of course one is mainly interested in the case where the additional twist
$\vp$ is trivial or at least unitary.

For a while we can drop the assumption that $X=\R H^n$. In fact, Proposition \ref{trans} below, suitably interpreted, holds for any connected semisimple Lie group $G$ with
minimal parabolic subgroup $P$. We use the notation and conventions fixed in Section \ref{round} around Equation (\ref{jemi}). There are corresponding closed positive
Weyl chambers $\taaa_+^*\subset \ii\taaa^*$, $\haaa_+^*\subset \haaa_\C^+$.
They satisfy $\aaaa^*_+\subset \haaa_+^*\subset \aaaa^*_+\oplus \taaa_+^*$.
We set $\rho_\gaaa:=\rho+\rho_\maaa\in \haaa_+^*$. We consider the following subset of $W(\gaaa_\C,\haaa_\C)$:
$$ W^1:=\{w\in W(\gaaa_\C,\haaa_\C)\:|\: w(\haaa_+^*)\subset\aaaa^*\oplus \taaa_+^*\}\ . $$
To a finite-dimensional irreducible representation $(\pi_\nu,F)$ with highest weight 
$\nu\in\haaa_+^*$
and an element $w\in W^1$ we associate weights
\begin{eqnarray*}
\mu_w&:=&w(\nu+\rho_\gaaa)_{|\taaa}-\rho_\maaa\in\taaa^*_+\ ,\\
\lambda_w&:=&-w(\nu+\rho_\gaaa)_{|\aaaa}\in\aaaa^*\ .
\end{eqnarray*} 
If $M$ is connected, then there is exactly one irreducible $M$-representation
$\sigma^w$ with highest weight $\mu_w$. In general $\sigma^w$ is uniquely determined by the additional requirement that 
the center $Z(G)\subset M$ of $G$ acts
by the same charater as on $F$. Then $\lambda_w\in I_{\sigma^w}^r$. The only
linear rank one group, where $M$ is not connected, is $SL(2,\R)$. 
In this case $M=\{\pm\id\}$, and $\sigma^w(-\id)=(-1)^{\dim F-1}$ for both elements $w\in W^1$. 
We denote the resulting $P$-representation
$\sigma^w_{\lambda_w}$ by $\sigma^w_{F,\lambda_w}$. The set $$\{C^{-\omega}(\partial X, V(\sigma^w_{F,\lambda_w}))\:|\: w\in W^1\}$$ 
consists of all principal series representations having the infinitesimal character $\chi_F=\chi_{\nu+\rho_\gaaa}$ and the central character of $F$.

For a $\cZ(\gaaa)$-module $V$ and a character $\chi:\cZ(\gaaa)\rightarrow \C$
we consider the generalized eigenspace
$$ V^\chi:=\{v\in V\:|\:\exists k\in\nat\mbox{ such that } (z-\chi(z))^k v=0
\mbox{ for all } z\in \cZ(\gaaa)\}\ .$$
If $V$ is locally $\cZ(\gaaa)$-finite, i.e., $\dim \cZ(\gaaa)v<\infty$ for all $v\in V$, then we have 
$$ 
V=\bigoplus_\chi V^\chi\ .
$$
Since the underlying $(\gaaa,K)$-module of the $G$-representation
$$V:=C^{-\omega}(\partial X, V^k(\sigma_\lambda))\otimes F
=C^{-\omega}(\partial X, V^k(\sigma_\lambda,\pi_\nu)) $$
is admissible and finitely generated $V$ considered as a $Z(\gaaa)$-module
is locally finite, and the sum 
$$ 
C^{-\omega}(\partial X, V^k(\sigma_\lambda,\pi_\nu))=\bigoplus_\chi C^{-\omega}(\partial X, V^k(\sigma_\lambda,\pi_\nu))^\chi\ .
$$ 
is finite. Here the $\cZ(\gaaa)$-action comes from the diagonal $G$-action.
Restricting the attention to the subalgebra $\C[\Omega]\subset \cZ(\gaaa)$
we get a weaker finite decomposition into generalized eigenspaces of the 
Casimir operator
\begin{equation}\label{klu}
C^{-\omega}(\partial X, V^k(\sigma_\lambda,\pi_\nu))=\bigoplus_{\kappa\in\C}
C^{-\omega}(\partial X, V^k(\sigma_\lambda,\pi_\nu))^\kappa\ .
\end{equation}

The announced application of the translation functor will
rest on the the following

\begin{prop}\label{trans}
Let $(\pi_\nu,F)$ be a finite-dimensional irreducible representation of $G$
with highest weight $\nu\in\haaa^*_+$, and let $w\in W^1$. 
Let $\kappa_F:=\|\nu+\rho_\gaaa\|^2-\|\rho_\gaaa\|^2$ be the Casimir eigenvalue on $F$. Then for $\sharp\in\{\infty,\omega\}$
$$C^{-\sharp}(\partial X, V^+(\sigma^w_{F,\lambda_w}))\cong 
C^{-\sharp}(\partial X, V^+(\sigma^w_{\C,\lambda_w},\pi_\nu))^{\chi_F}
=C^{-\sharp}(\partial X, V^+(\sigma^w_{\C,\lambda_w},\pi_\nu))^{\kappa_F}\ .$$
More precisely, for any $k\in\nat$ there are $G$-equivariant differential
operators
\begin{eqnarray*} 
D_{k,w}: C^{-\omega}(\partial X, V^k(\sigma^w_{F,\lambda_w}))&\rightarrow&
C^{-\omega}(\partial X, V^k(\sigma^w_{\C,\lambda_w},\pi_\nu))\ ,\\  
D^{k,w}:
C^{-\omega}(\partial X, V^k(\sigma^w_{\C,\lambda_w},\pi_\nu))&\rightarrow& C^{-\omega}(\partial X, V^k(\sigma^w_{F,\lambda_w}))
\end{eqnarray*}
satisfying
\begin{enumerate}
\item $D_{k,w}: C^{-\omega}(\partial X, V^k(\sigma^w_{F,\lambda_w}))
\rightarrow
C^{-\omega}(\partial X, V^k(\sigma^w_{\C,\lambda_w},\pi_\nu))^{\kappa_F}$
is an isomorphism.
\item $D^{k,w}\circ D_{k,w}=\id$.
\item $D_{k+1,w\ |C^{-\omega}(\partial X, V^k(\sigma^w_{F,\lambda_w}))}
=D_{k,w}$, $\p\circ D_{k+1,w}=D_{k,w}\circ \p$. The analogous assertions hold for
$D^{k+1,w}$.
\end{enumerate}
\end{prop}
\proof The proposition, at least for $k=1$, is well-known among representation theorists (see e.g. \cite{spehvogan80}, Section 5). 
However, it is more convenient for us to give a direct proof here rather than to cite
all relevant results from several places in the literature.

Similar to the proof of Theorem \ref{extmer} we use the isomorphism of
$G$-bundles
$$ V^k(\sigma^w_{\C,\lambda_w},\pi_\nu))\cong V^k(\sigma^w_{\C,\lambda_w}\otimes\pi_{\nu |P})\ .$$
The $P$-representation $W=\sigma^w_{\C,\lambda_w}\otimes\pi_{\nu |P}$
has a Jordan-H\"older series
$$0=W_0\subset W_1\subset\dots\subset W_l=W $$ 
such that on each irreducible composition factor $Q_r:=W_r/W_{r-1}$ 
the group $MA$ acts irreducibly with highest weight of the form
\begin{equation}\label{glo} 
w(\rho_\gaaa)+\rho-\rho_\maaa+\mu_r\ ,
\end{equation}
where $\mu_r$ is a $MA$-weight of $F$. It follows that $C^{-\omega}(\partial X, V^k(\sigma^w_{\C,\lambda_w},\pi_\nu))$ has a composition series with
composition factors $R_r$ isomorphic to the induced representations
$$\mbox{Ind}^G_P(Q_r\otimes\Pi^k)=C^{-\omega}(\partial X, V^k(\sigma^r_{\lambda_r}))\ ,$$
where $\sigma^r$ is a $M$-representation having highest weight $(\mu_r+w(\rho_\gaaa))_{|\taaa\:}-\rho_\maaa$ and $\lambda_r=-(\mu_r+w(\rho_\gaaa))_{|\aaaa}$. Thus $R_r$
has generalized infinitesimal character $\chi_{\mu_r+w(\rho_\gaaa)}$.

We now assume that $R_r$ has generalized Casimir eigenvalue $\kappa_F$. This means
$$ \|\nu+\rho_\gaaa\|^2=\|\mu_r+w(\rho_\gaaa)\|^2\ .$$
We obtain
\begin{equation}\label{klux}
0\le\|\nu\|^2-\|\mu_r\|^2=-2\langle \nu-w^{-1}(\mu_r),\rho_\gaaa\rangle\  .
\end{equation}
Since $w^{-1}(\mu_r)$ is a weight of $F$, too, the difference $\nu-w^{-1}(\mu_r)$ is a nonnegative linear combination of positive roots. Hence the right hand side of
(\ref{klux}) cannot be positive. It follows that $\mu_r=w(\nu)$.
   
We claim that the $MA$-representation $\tau_{w,\nu}$ with highest weight
$$w(\rho_\gaaa)+\rho-\rho_\maaa+w(\nu) $$
occurs in $W$ exactly once. Indeed, $w(\nu)$ is an extremal weight of $F$,
thus it occurs in $F$ with multiplicity one. Hence $\tau_{w,\nu}$ can occur
in $W$ at most with multiplicity one. It remains to show that $w(\rho_\gaaa)+\rho-\rho_\maaa+w(\nu)$ is in fact a highest $MA$-weight.
If not, then $w(\rho_\gaaa)+\rho-\rho_\maaa+w(\nu)+\beta$, 
where $\beta$ is a sum
of positive $\maaa$-roots, would be a highest weight. Since $w\in W^1$
we have for any positive $\maaa$-root $\ve$
$$ \langle w^{-1}(\ve),\rho_\gaaa\rangle= \langle \ve,w(\rho_\gaaa)\rangle>0\ .$$
Hence $w^{-1}(\ve)\in \Delta^+$, and $\nu+w^{-1}(\beta)$ is not a weight of $F$.
Thus, also $w(\nu)+\beta$ does not occur in $F$. Therefore (see (\ref{glo}))
$w(\rho_\gaaa)+\rho-\rho_\maaa+w(\nu)+\beta$ cannot be a weight of $W$. 
The claim follows.

Summarizing the above discussion we see that exactly one of the composition
factors $R_r$ has generalized Casimir eigenvalue $\kappa_F$, and that this
composition factor is isomorphic to \linebreak[4]
$C^{-\omega}(\partial X, V^k(\sigma^w_{F,\lambda_w}))$ which has generalized infinitesimal character $\chi_F$. It then follows from (\ref{klu}) that
this composition factor occurs as a direct summand of $C^{-\omega}(\partial X, V^k(\sigma^w_{\C,\lambda_w},\pi_\nu))$. This gives the desired isomorphism.
In order to finish the proof we have to convince ourself that this
isomorphism is implemented by differential operators $D_{k,w}$ and $D^{k,w}$.
Let 
$$ P_k: C^{-\omega}(\partial X, V^k(\sigma^w_{\C,\lambda_w},\pi_\nu))
\rightarrow
C^{-\omega}(\partial X, V^k(\sigma^w_{\C,\lambda_w},\pi_\nu))$$
be the projection onto $C^{-\omega}(\partial X, V^k(\sigma^w_{\C,\lambda_w},\pi_\nu))^{\kappa_F}$ with respect to (\ref{klu}). Since the sum (\ref{klu}) is finite it is given by
the natural action of a certain polynomial of the Casimir operator on $C^{-\omega}(\partial X, V^k(\sigma^w_{\C,\lambda_w},\pi_\nu))$, and
is therefore a differential operator.
Let $r$ be such that $R_r\cong C^{-\omega}(\partial X, V^k(\sigma^w_{F,\lambda_w}))$ as above. By $\tau_r$ and $\tau_{r-1}$ we denote the $P$-representations on $W_r$ and $W_{r-1}$, respectively. Then $V^k(\tau_{r-1})\subset V^k(\tau_r)\subset V^k(\sigma^w_{\C,\lambda_w},\pi_\nu)$, and there is a projection
$$p_{k,r}: V^k(\tau_r)\rightarrow V^k(\sigma^w_{F,\lambda_w})\ .$$
Moreover,
\begin{equation}\label{rr}
\im P_k\subset C^{-\omega}(\partial X, V^k(\tau_r))\ ,
\quad C^{-\omega}(\partial X,V^k(\tau_{r-1}))\subset \ker P_k\ .
\end{equation}
Thus,
$$ D^{k,w}:=p_{k,r}\circ P_k $$
is a well-defined $G$-equivariant differential operator.
Choose an embedding
$$ i_r: V(\sigma^w_{F,\lambda_w})\rightarrow V(\tau_r)$$
such that
$$ p_{1,r}\circ i_{r}=\id\ .$$
We set
$$ i_{k,r}=i_r\otimes \id: V^k(\sigma^w_{F,\lambda_w})\rightarrow V^k(\tau_r)\subset V^k(\sigma^w_{\C,\lambda_w},\pi_\nu))\ . $$
While $i_{k,r}$ need not to be $G$-equivariant the composition
$$ D_{k,w}:=P_k\circ i_{k,r} $$
is $G$-equivariant because of (\ref{rr}).
It is now easily checked that the differential operators $D_{k,w}$ and
$D^{k,w}$ enjoy the Properties 1, 2 and 3. This completes the proof of the proposition.
\hB

Now let $\Gamma\subset G$ be convex cocompact, and let $(\vp, V_\vp)$
be a finite-dimensional representation of $\Gamma$. Let $\cZ(\gaaa)$
act on the first factor of 
$$C^{-\omega}(\partial X, V^k(\sigma_\lambda,\pi_\nu))\otimes V_\vp=C^{-\omega}(\partial X, V^k(\sigma_\lambda,\pi_\nu\otimes\vp))\ .$$
Then 
$${}^\Gamma C^{-\infty}(\Lambda, V^k(\sigma_\lambda,\pi_\nu\otimes\vp))
\subset {}^\Gamma C^{-\omega}(\partial X, V^k(\sigma_\lambda,\pi_\nu\otimes\vp))\subset C^{-\omega}(\partial X, V^+(\sigma_\lambda,\pi_\nu\otimes\vp))$$ 
are
$\cZ(\gaaa)$-submodules. Moreover, using e.g. Lemma \ref{vorzug}, we see that
$E^+(\sigma_\lambda,\pi_\nu\otimes\vp)$ and $E^+_\Lambda(\sigma_\lambda,\pi_\nu\otimes\vp)$ are $\cZ(\gaaa)$-submodules, too.
Here, if $\Gamma$ is cocompact, we set 
$E^+(\sigma_\lambda,\pi_\nu\otimes\vp):=\{0\}$.

\begin{kor}\label{beau}
The differential operators $D_{k,w}\otimes\id$ provide isomorphisms
\begin{eqnarray*}
{}^\Gamma C^{-\omega}(\partial X, V^k(\sigma^w_{F,\lambda_w},\vp))
&\cong& {}^\Gamma C^{-\omega}(\partial X, V^k(\sigma^w_{\C,\lambda_w},\pi_\nu\otimes\vp))^{\chi_F}\ ,\\
{}^\Gamma C^{-\infty}(\Lambda, V^+(\sigma^w_{F,\lambda_w},\vp))
&\cong& {}^\Gamma C^{-\infty}(\Lambda, V^+(\sigma^w_{\C,\lambda_w},\pi_\nu\otimes\vp))^{\chi_F}\ ,\\
E^+(\sigma^w_{F,\lambda_w},\vp)&\cong& E^+(\sigma^w_{\C,\lambda_w},\pi_\nu\otimes\vp)^{\chi_F}\ ,\\   
E^+_\Lambda(\sigma^w_{F,\lambda_w},\vp)&\cong& E^+_\Lambda(\sigma^w_{\C,\lambda_w},\pi_\nu\otimes\vp)^{\chi_F}\ .
\end{eqnarray*}
\end{kor}
\proof
The first isomorphism is obvious. The second one follows from the existence
of the left inverse $D^{k,w}$ which is again a local operator.
The third isomorphism is a consequence of Lemma \ref{vorzug} and Property 3
of the operators $D_{k,w}$ and $D^{k,w}$. Now the last isomorphism follows, too.
\hB

We now return to our assumption $X=\R H^n$, i.e., $G=Spin(1,n)$ or $G=SO(1,n)_0$. The generating set of simple reflections $\{s_\ve\:|\:\ve\in\Pi\}$, where $\Pi\subset \Delta^+$ is the set of simple roots, determines a length function
on $W(\gaaa_\C,\haaa_\C)$.
If $n$ is even, then 
$$W^1=\{w_0,w_1,\dots,w_{n-1}\}\ ,$$ 
where $w_i$ is the unique element of length $i$ in $W^1$. For odd $n$ we have $$W^1=\{w_0,\dots, w_{\frac{n-3}{2}},w_+,w_-,w_\frac{n+1}{2},\dots,w_{n-1}\}\ ,
$$ 
where $w_i$ and $w_\pm$ have length $i$ and $\frac{n-1}{2}$, respectively.
We abbreviate
\begin{eqnarray}
\mu_p:=\mu_w,&&  \lambda_p:=\lambda_w,\ \  \sigma^p_{F,\lambda_p}:=\sigma^w_{F,\lambda_w},\ \ \mbox{where } w=w_{n-1-p}\ ,\label{el}\\
\mu_\pm:=\mu_w,&&  \lambda_\pm:=\lambda_w,\ \   \sigma^\pm_{F,\lambda_\pm}:=\sigma^w_{F,\lambda_w},\ \ \mbox{where } w=w_\mp\ .
\label{end}
\end{eqnarray}
Occasionally, we will
also need a representation associated to $p=\frac{n-1}{2}$. We set  
$$\sigma^\frac{n-1}{2}_{F,\lambda_\frac{n-1}{2}}:=
\sigma^+_{F,\lambda_+}\oplus\sigma^-_{F,\lambda_-}\ .$$

We would like to write down things more explicitly. First let $n$ be even.
In standard coordinates $\aaaa^*\oplus \ii\taaa^*\cong \R\oplus \R^\frac{n-2}{2}$ with standard basis $e_0,e_1,\dots,e_\frac{n-2}{2}$ we have
\begin{eqnarray*}
\Pi&=&\{e_{i-1}-e_i\:|\:i=1,\dots,\frac{n-2}{2}\}\cup\{e_\frac{n-2}{2}\}\ ,\\
\haaa^*_+&=&\{(m_0,\dots,m_\frac{n-2}{2})\:|\: m_0\ge m_1\ge\dots\ge m_\frac{n-2}{2}\ge 0\}\ ,\\
\taaa^*_+&=&\{(m_1,\dots,m_\frac{n-2}{2})\:|\: m_1\ge m_2\ge\dots\ge m_\frac{n-2}{2}\ge 0\}\ ,\\
\rho_\gaaa&=&(\frac{n-1}{2},\frac{n-3}{2},\dots,\frac{1}{2})\ .
\end{eqnarray*}
The Weyl group $W(\gaaa_\C,\haaa_\C)$ consists of permutations of the
coordinates composed with an arbitrary number of sign chances. For $0\le p\le \frac{n-2}{2}$
we have
\begin{eqnarray*} 
w_p(m_0,\dots,m_\frac{n-2}{2})&=&(m_p, m_0,\dots, m_{p-1},m_{p+1},\dots,m_\frac{n-2}{2})\ ,\\
w_{n-1-p}(m_0,\dots,m_\frac{n-2}{2})&=&(-m_p, m_0,\dots, m_{p-1},m_{p+1},\dots,m_\frac{n-2}{2})\ .
\end{eqnarray*}
Associated to a highest weight $\nu=(m_0,\dots,m_\frac{n-2}{2})$ we get
for $0\le p\le \frac{n-2}{2}$
\begin{eqnarray*} 
\mu_p&=&(m_0+1,m_1+1,\dots, m_{p-1}+1,m_{p+1},\dots,m_\frac{n-2}{2})=\mu_{n-1-p}\ ,\\
\lambda_p&=&\rho+(m_p-p)\alpha\ ,\quad \lambda_{n-1-p}=-\lambda_p\ .
\end{eqnarray*}

Let now $n$ be odd.
In standard coordinates $\aaaa^*\oplus \ii\taaa^*\cong \R\oplus \R^\frac{n-1}{2}$ with standard basis $e_0,e_1,\dots,e_\frac{n-1}{2}$ we have
\begin{eqnarray*}
\Pi&=&\{e_{i-1}-e_i\:|\:i=1,\dots,\frac{n-1}{2}\}
\cup\{e_\frac{n-3}{2}+e_\frac{n-1}{2}\}\ ,\\
\haaa^*_+&=&\{(m_0,\dots,m_\frac{n-1}{2})\:|\: m_0\ge\dots\ge m_\frac{n-3}{2}\ge |m_\frac{n-1}{2}|\}\ ,\\
\taaa^*_+&=&\{(m_1,\dots,m_\frac{n-1}{2})\:|\: m_1\ge\dots\ge m_\frac{n-3}{2}\ge |m_\frac{n-1}{2}|\}\ ,\\
\rho_\gaaa&=&(\frac{n-1}{2},\frac{n-3}{2},\dots,0)\ .
\end{eqnarray*}
The Weyl group $W(\gaaa_\C,\haaa_\C)$ consists of permutations of the
coordinates composed with an even number sign chances. For $0\le p\le \frac{n-3}{2}$
we have
\begin{eqnarray*} 
w_p(m_0,\dots,m_\frac{n-1}{2})&=&(m_p, m_0,\dots, m_{p-1},m_{p+1},\dots,m_\frac{n-1}{2})\ ,\\
w_+(m_0,\dots,m_\frac{n-1}{2})&=&(m_\frac{n-1}{2}, m_0,\dots,m_\frac{n-3}{2})
\ ,\\
w_-(m_0,\dots,m_\frac{n-1}{2})&=&(-m_\frac{n-1}{2}, m_0,\dots,m_\frac{n-5}{2},-m_\frac{n-3}{2})
\ ,\\
w_{n-1-p}(m_0,\dots,m_\frac{n-1}{2})&=&(-m_p, m_0,\dots, m_{p-1},m_{p+1},\dots,m_\frac{n-3}{2},-m_\frac{n-1}{2})\ .
\end{eqnarray*}
Associated to a highest weight $\nu=(m_0,\dots,m_\frac{n-1}{2})$ we get
for $0\le p\le \frac{n-3}{2}$
\begin{eqnarray*} 
\mu_p&=&(m_0+1,m_1+1,\dots, m_{p-1}+1,m_{p+1},\dots,m_\frac{n-3}{2},-m_\frac{n-1}{2})\ ,\quad \lambda_p=\rho+(m_p-p)\alpha\ ,\\
\mu_\pm&=&(m_0+1,m_1+1,\dots,m_\frac{n-5}{2}+1,\mp (m_\frac{n-3}{2}+1))\ ,\quad 
\lambda_\pm=\pm m_\frac{n-1}{2}\alpha\ ,\\
\mu_{n-1-p}&=&(m_0+1,m_1+1,\dots, m_{p-1}+1,m_{p+1},\dots,m_\frac{n-1}{2})\ ,
\quad\lambda_{n-1-p}=-\lambda_p\ .
\end{eqnarray*}
In particular, the $P$-representations considered
in Section \ref{heart} are associated to the trivial $G$-representation:
\begin{equation}\label{not} 
\sigma^p_{\C,\lambda_p}=\sigma^p_{\rho-p\alpha}\ ,\quad \sigma^\pm_{\C,\lambda_\pm}=\sigma^\mp_0\ .
\end{equation}

As a next step we want to apply the translation functor to the complex
(\ref{kike}). Indeed, replacing there $\vp$ by $\pi_\nu\otimes \vp$ we
obtain a complex of $\cZ(\gaaa)$-modules which decomposes into
generalized $\cZ(\gaaa)$-eigenspaces. Taking the component with generalized
infinitesimal character $\chi_F$ and applying Proposition \ref{trans}
we obtain the complex
\begin{eqnarray}\label{zhelo}  
0\rightarrow {}^\Gamma C^{-\omega}(\partial X, V(\sigma^{0}_{F,\lambda_{0}},\vp))&\stackrel{d_F^0}\longrightarrow & {}^\Gamma C^{-\omega}(\partial X, V(\sigma^{1}_{F,\lambda_{1}},\vp))\stackrel{d_F^1}\longrightarrow \\
\dots&\stackrel{d_F^{n-2}}\longrightarrow &{}^\Gamma C^{-\omega}(\partial X, V(\sigma^{n-1}_{F,\lambda_{n-1}},\vp))\stackrel{d_F^{n-1}}\longrightarrow
{}^\Gamma (F\otimes V_\vp)\rightarrow 0\ ,\nonumber
\end{eqnarray}
where for $p<n-1$ the operator $d_F^p=D^{1,w_{n-1-p}}\circ d \circ D_{1,w_{n-1-p}}$ is a differential operator while 
$d_F^{n-1}=D^{1,w_0}\circ \int_{\partial X} \circ D_{1,w_0}$
can be identified with a multiple of the intertwining operator
$\hat J_{\sigma^{n-1}_{F},\lambda_{n-1}}$ if one considers $F\otimes V_\vp$
as sitting in $C^{-\omega}(\partial X, V(\sigma^{0}_{F,\lambda_{0}},\vp))$.
For trivial $\Gamma$ (and $\vp$) the complex (\ref{zhelo}) is exact and appears
in the literature under various names like BGG-resolution or \v Zelobenko complex (see e.g. \cite{zelo76}, \cite{bastoneastwood89}, \cite{juhl01}).
By $Z_F^p(\vp)$ and $Z^p_{F,\Lambda}(\vp)$ we denote the spaces of $p$-cocycles
of (\ref{zhelo}) and of $p$-cocycles supported on the limit set, respectively.

Now we can state
the second main theorem of the paper. Recall that by convention 
$E^+(\sigma_\lambda,\pi\otimes\vp)=\{0\}$ for cocompact $\Gamma$.

\begin{theorem}\label{main2}
Let $(\pi,F)$ be an irreducible finite-dimensional representation of $G$,
and let $(\vp,V_\vp)$ an irreducible representation of $\Gamma$. There
are surjections ${}^\Gamma C^{-\infty}(\Lambda, V^+(\sigma^{n-p}_{F,\lambda_{n-p}},\vp))\rightarrow H^p(\Gamma,F\otimes V_\vp)$, ${}^\Gamma C^{-\infty}(\Lambda, V^+(\sigma^{\pm}_{F,\lambda_\pm},\vp))\rightarrow H^\frac{n+1}{2}(\Gamma,F\otimes V_\vp)$. If $\Gamma$
is not cocompact or $\vp$ is unitary, then we have the following isomorphisms 
\begin{eqnarray*} 
H^p(\Gamma,F\otimes V_\vp)&\cong& {}^\Gamma C^{-\infty}(\Lambda, V^+(\sigma^{n-p}_{F,\lambda_{n-p}},\vp))/E^+_\Lambda(\sigma^{n-p}_{F,\lambda_{n-p}},\vp)\ ,\quad 
p=1,\dots, n\ ,\ p\not=\frac{n+1}{2}\ ,\\
H^\frac{n+1}{2}(\Gamma,F\otimes V_\vp)&\cong& {}^\Gamma C^{-\infty}(\Lambda, V^+(\sigma^{\pm}_{F,\lambda_\pm},\vp))/E^+_\Lambda(\sigma^\pm_{F,\lambda_\pm},\vp)\ .
\end{eqnarray*}
Moreover, we have 
\begin{eqnarray*}
k_+(\sigma^{n-p}_{F,\lambda_{n-p}},\vp)&\le&
k_-(\sigma^{n-p}_{F,\lambda_{n-p}},\vp)+1\ ,\\
k_+(\sigma^{\pm}_{F,\lambda_\pm},\vp))&\le& k_-(\sigma^{\pm}_{F,\lambda_\pm},\vp))+1\ .
\end{eqnarray*}
If $\vp$ is unitary and $p\ge\frac{n+1}{2}$, then $k_-(\sigma^{n-p}_{F,\lambda_{n-p}},\vp)=0$ and
\begin{eqnarray}
H^p(\Gamma,F\otimes V_\vp)&\cong & \left\{\begin{array}{clc}
{}^\Gamma C^{-\infty}(\Lambda, V^+(\sigma^{n-p}_{F,\lambda_{n-p}},\vp))=
U_\Lambda(\sigma^{n-p}_{F,\lambda_{n-p}},\vp)\ ,&p>\frac{n+1}{2}\\    
{}^\Gamma C^{-\infty}(\Lambda, V^+(\sigma^{\pm}_{F,\lambda_\pm},\vp))=
U_\Lambda(\sigma^\pm_{F,\lambda_\pm},\vp)\ ,&p=\frac{n+1}{2} \end{array}\right.
\label{ei}\\
&\cong& Z^{n-p}_{F,\Lambda}(\vp)\ .
\label{eiei}
\end{eqnarray}
If $\vp$ is unitary and $p=\frac{n}{2}$, then $k_+(\sigma^{n-p}_{F,\lambda_{n-p}},\vp)\le\max\{k_-(\sigma^{n-p}_{F,\lambda_{n-p}},\vp),1\}$ and
\begin{equation}\label{eieiei}
H^p(\Gamma,F\otimes V_\vp)\cong Z^{n-p}_{F,\Lambda}(\vp)\ .
\end{equation}
If, in addition, $\Gamma$ is cocompact, then  we have
for all $p$
$$H^p(\Gamma,F\otimes V_\vp)\cong Z^{n-p}_{F,\Lambda}(\vp)\cong Z^p_{F,\Lambda}(\vp)\ .$$
\end{theorem}
\proof
If $(\psi,V_\psi)$ is a finite-dimensional representation of $\Gamma$, then
we denote by $\tilde E(\psi)$ the flat vector bundle $X\times V_\psi$ over $X$.
It carries a $G$ action $(x,v)\mapsto (gx,v)$ and a $\Gamma$-action $(x,v)
\mapsto (\gamma x,\psi(\gamma)v)$. We equip the de Rham complex $\Omega^*(Y,E(\pi\otimes\vp))\cong {}^\Gamma \Omega^*(X,\tilde E(\pi\otimes\vp))$ with the $\cZ(\gaaa)$-action induced by the tensor product $G$-action on
\begin{equation}\label{Akt}
\Omega^*(X,\tilde E(\vp))\otimes F\cong\Omega^*(X,\tilde E(\pi\otimes\vp))\ .
\end{equation}
This induces a $\cZ(\gaaa)$-module structure on $H^p(\Gamma,F\otimes V_\vp)\cong H^p(Y,E(\pi\otimes\vp))$. Then the isomorphisms in Theorem \ref{main1} (with
$\vp$ replaced by $\pi\otimes\vp$) become isomorphisms of $\cZ(\gaaa)$-modules.
Now, by general principles, $\cZ(\gaaa)$ acts on $H^p(\Gamma,F\otimes V_\vp)$
by the infinitesimal character $\chi_F$. In order to see this we look at the category of
$(\C[\Gamma],\cZ(\gaaa))$-bimodules. Taking $\Gamma$-invariants
defines a left exact functor from this category to the category of $\cZ(\gaaa)$-modules. Its right derived functors coincide with $H^p(\Gamma,.)$,
but now provide also $\cZ(\gaaa)$-module structures on the cohomology
groups. We want to study the $\cZ(\gaaa)$-modules $H^p(\Gamma,F\otimes V_\vp)$, where $\cZ(\gaaa)$ acts on the first factor
of $F\otimes V_\vp$. For this
it is sufficient to look at resolutions of $F\otimes V_\vp$
by $\Gamma$-acyclic $(\C[\Gamma],\cZ(\gaaa))$-bimodules.
The $\Gamma$-modules $\Omega^p(X,\tilde E(\pi\otimes\vp))$ are acyclic by \cite{bunkeolbrich950}, Lemma 2.4. Taking $\Gamma$-invariants identifies
the cohomology of $\Omega^*(Y,E(\pi\otimes\vp))$ with $H^*(\Gamma,F\otimes V_\vp)$ as a $\cZ(\gaaa)$-module.
There is a second $\cZ(\gaaa)$-module
structure on this de Rham complex coming from the $G$-action on the second
factor in (\ref{Akt}) which gives rise to a resolution of $F\otimes V_\vp$,
too. It follows that $\cZ(\gaaa)$ acts on $H^*(\Gamma,F\otimes V_\vp)$ by $\chi_F$.

At least for noncocompact $\Gamma$ the theorem, except for (\ref{eieiei}), now follows from Corollary \ref{beau} by taking the components
with generalized infinitesimal character $\chi_F$ of all the spaces appearing
in Theorem \ref{main1}. For cocompact $\Gamma$ we use Proposition \ref{koks}
instead. That in case of unitary $\vp$ the resulting surjection is an isomorphism follows from classical Hodge theory with respect to the
Laplacian associated to a so called admissible scalar product on $F$. This
Laplacian coincides with the action of $-\Omega+\kappa_F$ (see 
\cite{borelwallach80}, Chapter II, \S 2, for all that), thus acts by zero
on the relevant spaces. These facts will be discussed further in Subsection \ref{nancy} and Section \ref{mazzeo}.

It remains to prove (\ref{eieiei}) for noncocompact $\Gamma$. For this we will need a couple of results which are of interest in their own right.

\begin{lem}\label{irr1}
Let $k\in\nat$, and let
$res_k: {}^\Gamma C^{-\omega}(\partial X,V^k(\sigma^{n-p}_{F,\lambda_{n-p}},\vp))\rightarrow
C^{-\omega}(B,V^k_B(\sigma^{n-p}_{F,\lambda_{n-p}},\vp))$ be the restriction map. 
Then $ext$ induces an isomorphism
$$ e_k: \coker\: res_k\stackrel{\cong}{\longrightarrow}
E^+_\Lambda(\sigma^{n-p}_{F,\lambda_{n-p}},\vp)/\p_\Gamma^k\left({}^\Gamma C^{-\infty}(\Lambda,V^+(\sigma^{n-p}_{F,\lambda_{n-p}},\vp))\right)\ .$$
In particular, if $k\ge k_+(\sigma^{n-p}_{F,\lambda_{n-p}},\vp)$, then
$$ \coker\: res_k\cong E^+_\Lambda(\sigma^{n-p}_{F,\lambda_{n-p}},\vp)\ .$$
\end{lem}
\proof 
Choose $l\ge\max\{k_-(\sigma^{n-p}_{F,\lambda_{n-p}},\vp),k\}$. Let $[f]\in\coker\: res_k$ be
represented by $\rho_B^{l-k} f_l$ for some $f_l \in C^{-\omega}(B,V^l_B(\sigma^{n-p}_{F,\lambda_{n-p}},\vp))$. Then we set
$$e_k([f]):=ext[l] f_l\ \mod \p_\Gamma^k\left({}^\Gamma C^{-\infty}(\Lambda,V^+(\sigma^{n-p}_{F,\lambda_{n-p}},\vp))\right)\ .$$
We have to check that $e_k$ is well-defined. Assume that 
$$\rho_B^{l-k} f_l=res_k\psi\quad \mbox{for some }
\psi\in {}^\Gamma C^{-\omega}(\partial X,V^k(\sigma^{n-p}_{F,\lambda_{n-p}},\vp))\ .$$ 
Choose $f_{k+l}\in C^{-\omega}(B,V^{k+l}_B(\sigma^{n-p}_{F,\lambda_{n-p}},\vp))$ such that
$\rho_B^k f_{k+l}=f_l$.
We obtain
$$ ext[l]f_l=ext[l](\p_B^k f_{k+l})=\p_\Gamma^k(ext[l] f_{k+l})
=\p_\Gamma^k(ext[l] f_{k+l}-\psi)\ .$$
Now $res(ext[l] f_{k+l}-\psi)=0$, hence $ext[l]f_l\in \p_\Gamma^k\left({}^\Gamma C^{-\infty}(\Lambda,V^+(\sigma^{n-p}_{F,\lambda_{n-p}},\vp))\right)$.
Thus $e_k$ is well-defined.

Vice versa, assume that 
$$ext[l] f_l=\p_\Gamma^k(ext[l] f_{k+l})=\p_\Gamma^k(\psi_0)\quad 
\mbox{for some }\psi_0\in {}^\Gamma C^{-\infty}(\Lambda,V^+(\sigma^{n-p}_{F,\lambda_{n-p}},\vp))\ .$$
Then 
$\psi:=ext[l] f_{k+l}-\psi_0\in C^{-\infty}(\partial X,V^k(\sigma^{n-p}_{F,\lambda_{n-p}},\vp))$ and $res_k\psi=f_l$.
Thus $e_k$ is injective.

Finally, it follows from Lemma \ref{gm} that $e_k$ is surjective.
\hB

\begin{prop}\label{irr2}
Let $\vp$ be unitary. Then
$$ \dim Z^{n-p}_{F,\Lambda}(\vp)\ge \dim H^p(\Gamma,F\otimes V_\vp)\ .$$
\end{prop}
\proof
For $p\ge\frac{n+1}{2}$ we have already shown equality. Thus we
can assume $p\le \frac{n}{2}$, in particular $\lambda_{F,n-p}<0$.
By Lemma \ref{irr1} we have that for $k\ge k_+(\sigma^{n-p}_{F,\lambda_{n-p}},\vp)$ 
\begin{equation}\label{libru}
\dim {}^\Gamma C^{-\infty}(\Lambda,V^+(\sigma^{n-p}_{F,\lambda_{n-p}},\vp))
= \dim H^p(\Gamma,F\otimes V_\vp)+\dim \:\coker\:res_k\ .
\end{equation}
As in the proof of Proposition \ref{zwirbel}, let $\hat J_+: C^{-\omega}(\partial X,V^+\sigma^{n-p}_{F,\lambda_{n-p}},\vp))\rightarrow C^{-\omega}(\partial X,V^+(\sigma^{p-1}_{F,\lambda_{p-1}},\vp))$ be the intertwining operator induced by $\hat J_{\sigma^{n-p}_F,\mu}$, which is regular at $\mu=\lambda_{F,n-p}<0$. The knowledge of the composition factors of
$C^{-\omega}(\partial X, V(\sigma^{n-p}_{F,\lambda_{n-p}}))$ (see e.g. 
\cite{collingwood85}, Ch. 5) gives for $p\ne \frac{n+1}{2}$
$$ Z^{n-p}_{F}(\vp)={}^\Gamma(\ker\hat J_+)\ .$$
Thus in view of (\ref{libru}) it is enough to show that
\begin{equation}\label{luck}
\dim \hat J_+\left({}^\Gamma C^{-\infty}(\Lambda,V^+(\sigma^{n-p}_{F,\lambda_{n-p}},\vp))\right) 
\le  \dim \:\coker\:res_k\ .
\end{equation}
As in the proof of Theorem \ref{main1} we see that any element
$h\in {}^\Gamma C^{-\infty}(\Lambda,V^k(\sigma^{n-p}_{F,\lambda_{n-p}},\vp))$
can be written as $h=\psi-ext[l]f_{l+1}$, where
$\psi\in Z^{n-p}_{F}(\vp)$ and $f_{l+1}\in C^{-\omega}(B,V_B^{l+1}(\sigma^{n-p}_{F,\lambda_{n-p}},\vp))$ such that $\p_B^{l} f_{l+1}=res(\psi)$.
This implies that $\hat J_+(h)= 
-\hat J_+(ext[l]f_{l+1})\in E^+(\sigma^{p-1}_{F,\lambda_{p-1}},\vp)$. Now $E^+_\Lambda(\sigma^{p-1}_{F,\lambda_{p-1}},\vp)=\{0\}$ by Proposition 
\ref{newer}, 2. It follows that
$$ \hat J_+\left({}^\Gamma C^{-\infty}(\Lambda,V^+(\sigma^{n-p}_{F,\lambda_{n-p}},\vp))\right)
\cap {}^\Gamma C^{-\infty}(\Lambda,V^+(\sigma^{p-1}_{F,\lambda_{p-1}},\vp))=\{0\}$$
and thus
$$\dim \hat J_+\left({}^\Gamma C^{-\infty}(\Lambda,V^+(\sigma^{n-p}_{F,\lambda_{n-p}},\vp))\right)
=\dim\: res\circ \hat J_+\left({}^\Gamma C^{-\infty}(\Lambda,V^+(\sigma^{n-p}_{F,\lambda_{n-p}},\vp))\right)\ .$$
Since $V(\tilde\sigma^{n-p}_{F,\lambda_{n-p}},\tilde\vp))$ is the complex
conjugate of $V(\sigma^{n-p}_{F,\lambda_{n-p}},\vp))$ and 
$res\circ \hat J_+=(res\circ \hat J)_+$ Proposition
\ref{highgreen} now implies that the dimension on the right hand side does not exceed
$\dim \:\coker\:res_k$. This proves (\ref{luck}), and hence the proposition.
\hB

\begin{lem}\label{dser}
Let $\vp$ be unitary and $n$ be even. Then
$$ Z^\frac{n}{2}_{F,\Lambda}(\vp)\cap E^1_\Lambda(\sigma^\frac{n}{2}_{F,\lambda_\frac{n}{2}},\vp)=\{0\}\ .$$
\end{lem}
\proof
The lemma relies on the fact that the $G$-representation $\ker d^\frac{n}{2}_F\subset C^{-\infty}(\partial X, V(\sigma^\frac{n}{2}_{F,\lambda_\frac{n}{2}}))$ is the direct sum of the distribution
globalization of the two irreducible discrete series representations having
infinitesimal character $\chi_F$ (see e.g. \cite{collingwood85}, Ch. 5).
Thus we can exploit unitarity in a similar way as for the corresponding
assertions for $p\ge\frac{n-1}{2}$. 

If $F$ has highest weight $(m_0,m_1,\dots,m_\frac{n-2}{2})$, then $\sigma^\frac{n}{2}_{F}$ has highest weight $(m_0+1,m_1+1,\dots,m_\frac{n-4}{2}+1)$. We consider the $K$-representation
$\gamma=\gamma^+\oplus\gamma^-$, where $\gamma^\pm$ has highest weight 
$(m_0+1,m_1+1,\dots,m_\frac{n-4}{2}+1,\pm(m_\frac{n-2}{2}+1))$. Let $\bar G$ be the
group $SO(1,n)$ or its double cover, respectively, having two connected components. Consider also the corresponding groups $\bar K$ and $\bar M$.
Then $\bar M\cong M\times \Z_2$. Letting the
$\Z_2$-factor act trivially we obtain an $\bar M$-action on $V_{\sigma^\frac{n}{2}_{F}}$. Thus $C^{-\infty}(\partial X, V^k(\sigma^\frac{n}{2}_{F,\lambda_\frac{n}{2}}))$, $\ker d^\frac{n}{2}_F\subset C^{-\infty}(\partial X, V(\sigma^\frac{n}{2}_{F,\lambda_\frac{n}{2}},\vp))$ become $\bar G$-representations. Then $\ker d^\frac{n}{2}_F$ is the unique irreducible $\bar G$-subrepresentation of $C^{-\infty}(\partial X, V(\sigma^\frac{n}{2}_{F,\lambda_\frac{n}{2}}))$.  The representation of $K$ on $V_\gamma$
can be uniquely extended to a $\bar K$ representation such that
$\dim \Hom_{\bar M}(V_{\sigma^\frac{n}{2}_{F}},V_\gamma)=1$.

Let $0\ne T\in \Hom_{\bar M}(V_{\sigma^\frac{n}{2}_{F}},V_\gamma)$. There
is a corresponding Poisson transform
$$ P^T:=P_{\sigma^\frac{n}{2}_{F,\lambda_\frac{n}{2}}}^T: C^{-\omega}(\partial X,V(\sigma^\frac{n}{2}_{F,\lambda_\frac{n}{2}}))
\rightarrow C^\infty(G,V_\gamma)^K\ ,$$
which becomes $\bar G$-equivariant. Since the $\bar K$-type $\gamma$ occurs in $\ker d^\frac{n}{2}_F$, it is injective. 

Let $\Phi\in Z^p_F(\vp)$, $f\in {}^\Gamma C^{-\infty}(\partial X,V^2(\sigma^\frac{n}{2}_{F,\lambda_\frac{n}{2}},\vp))$ such that $\p_\Gamma f \in Z^p_F(\vp)$ and 
$res(\Phi)$ and $\res(f)$ are smooth. We claim that $P^T\Phi, P^T\p_\Gamma f
\in [L^2(\Gamma\backslash G,\vp)\otimes V_\gamma]^K$ and that there is
a non-zero constant $c$ such that 
\begin{equation}\label{hellgruen}
(P^T\Phi, P^T\p_\Gamma f)=c(res(\Phi), res(\hat J_2 f))_B\ .
\end{equation}
Here $\hat J_2: C^{-\infty}(\partial X,V^2(\sigma^\frac{n}{2}_{F,\lambda_\frac{n}{2}},\vp))\rightarrow C^{-\infty}(\partial X,V^2(\sigma^\frac{n}{2}_{F,-\lambda_\frac{n}{2}},\vp))$
is the intertwining operator induced by the family $\hat J_{\sigma^\frac{n}{2}_{F},\mu}$ as in the proof of Proposition \ref{zwirbel},
and $(.,.)_B$ is the natural sesqui-linear pairing between $C^{\infty}(B,V_B(\sigma^\frac{n}{2}_{F,\lambda_\frac{n}{2}},\vp))$ and $C^{-\infty}(B,V_B(\sigma^\frac{n}{2}_{F,-\lambda_\frac{n}{2}},\vp))$.
Since $\rho_\Gamma f \in Z^p_F(\vp)\subset \ker\hat J_{\sigma^\frac{n}{2}_{F},\lambda_\frac{n}{2}}$ we have $\hat J_2 f \in C^{-\infty}(\partial X,V^1(\sigma^\frac{n}{2}_{F,-\lambda_\frac{n}{2}},\vp))$.

If $\Phi\in Z^p_{F,\Lambda}(\vp)\cap E^1_\Lambda(\sigma^\frac{n}{2}_{F,\lambda_\frac{n}{2}},\vp)$, 
then $\Phi=ext[k_-+1]\phi$ for some $\phi\in C^\infty(B,V^{k_-+1}(\sigma^\frac{n}{2}_{F,\lambda_\frac{n}{2}},\vp))$, thus
$\Phi=p_\Gamma f$ for $f=ext[k_-]\phi$. Assuming (\ref{hellgruen}) we obtain
$\|P^T\Phi\|^2=c(0,res(\hat J_2f))_B=0$, hence $\Phi=0$. This implies
the lemma.

It remains to prove (\ref{hellgruen}). Since $\ker d^\frac{n}{2}_F$ consists
of discrete series representations and has infinitesimal character $\chi_F$  we obtain for smooth $\psi\in \ker d^\frac{n}{2}_F$
$$P^T\psi(ka)=a^{\lambda_\frac{n}{2}-\rho}\tilde\psi(k)+O(a^{\lambda_\frac{n}{2}-\rho-\ve})\ ,$$
where $\tilde\psi$ is some smooth $V_\gamma$-valued function on $K$. Note that $\lambda_\frac{n}{2}<0$. The map $\psi\mapsto T^*\tilde\psi$ is a $\bar G$-intertwining operator from the smooth vectors
in $\ker d^\frac{n}{2}_F$ to $C^\infty(\partial X,V(\sigma^\frac{n}{2}_{F,\lambda_\frac{n}{2}}))$. Since $\ker d^\frac{n}{2}_F$ is the unique irreducible $\bar G$-subrepresentation of $C^{\infty}(\partial X, V(\sigma^\frac{n}{2}_{F,\lambda_\frac{n}{2}}))$ we conclude that
$T^*\tilde\psi$ is a multiple of $\psi$. We obtain as in \cite{bunkeolbrich000},
Lemma 6.2, 3, that for $kM\in\Omega$
$$P^T\Phi(ka)=a^{\lambda_\frac{n}{2}-\rho}\tilde\Phi(k)+O(a^{\lambda_\frac{n}{2}-\rho-\ve})$$
with $T^*\tilde\Phi$ a constant multiple of $res(\Phi)$.
As in Section \ref{heart}, $P^T$ can be extended to \linebreak[4]$C^{-\omega}(\partial X, V^2(\sigma^\frac{n}{2}_{F,\lambda_\frac{n}{2}}))$, and for $kM\in\Omega$
$$ P^Tf(ka)=a^{-\lambda_\frac{n}{2}-\rho}T\hat J_2 f(k)+O(a^{-\lambda_\frac{n}{2}-\rho-\ve})\ .$$
Since $(\Omega-\kappa_F)P^Tf$ is a constant multiple of $P^T\p_\Gamma f$ we 
can use Green's formula as in the proof of \cite{bunkeolbrich000}, Prop. 10.4, 
in order to conclude (\ref{hellgruen}).
\hB

Lemma \ref{dser} implies that the map 
$Z^\frac{n}{2}_{F,\Lambda}(\vp)\rightarrow H^\frac{n}{2}(\Gamma, F\otimes V_\vp)$ is injective. On the other hand
$\dim Z^\frac{n}{2}_{F,\Lambda}(\vp)\ge  H^\frac{n}{2}(\Gamma, F\otimes V_\vp)$ by Proposition \ref{irr2}.
Thus the map is also surjective. This finishes the proof of (\ref{eieiei}).
\hB 

That the cohomology groups $H^p(\Gamma,F)$ should be representable by 
currents on
the limit set was conjectured by Patterson \cite{patterson93}. As explained
in the introduction,
he suggested that, at least for trivial $\vp$, the isomorphism (\ref{eiei}), which we have now proved for $p\ge\frac{n}{2}$, should be true for all 
$p$.  
What we can say in general is the following

\begin{kor}\label{reg}
If $\vp$  is unitary and $ext_\mu$ is regular at $\mu=\lambda_{F,n-p}$, i.e. $k_-(\sigma^{n-p}_{F,\lambda_{n-p}},\vp)=0$, then
$${}^\Gamma C^{-\infty}(\Lambda,V^+(\sigma^{n-p}_{F,\lambda_{n-p}},\vp))
=Z^{n-p}_{F,\Lambda}(\vp)\cong H^p(\Gamma,F\otimes V_\vp)\ .$$
\end{kor}
\proof
By assumption $E^+_\Lambda(\sigma^{n-p}_{F,\lambda_{n-p}},\vp)=\{0\}$.
This implies injectivity of the map 
$${}^\Gamma C^{-\infty}(\Lambda,V^+(\sigma^{n-p}_{F,\lambda_{n-p}},\vp))\rightarrow
H^p(\Gamma,F\otimes V_\vp)\ .$$ 
We now apply Proposition \ref{irr2}.
\hB 

However, the map $Z^{n-p}_{F,\Lambda}(\vp)\rightarrow
H^p(\Gamma,F\otimes V_\vp)$ is not injective in general. 
Let $n$ be odd, and let $\Gamma$ be such that $B$ is disconnected (this
happens for example for cocompact subgroups of $SO(1,n-1)_0$ considered
as a subgroup of $SO(1,n)_0$). Let $f\in C^\infty(B,V(1_\rho))$ be the 
characteristic function of one connected component. The residue
$D_\rho$ of $\hat J_\lambda$ at $\lambda=\rho$ is a differential
operator. Its kernel is the one-dimensional space of constant functions while its image is equal to $\ker \hat J_{-\rho}$. Thus $0\ne h:=D_\rho (ext_\rho f)\in Z^{n-1}_\Lambda$.
Choosing a family $f_\mu \in C^\infty(B,V(1_\mu))$ such that $f_\rho=f$ we
obtain a family $h_\nu:=- \hat J_{-\nu}\circ ext_{-\nu} f_{-\nu}\in {}^\Gamma \cM^1_{-\rho}C^{-\infty}(\partial X, V(1_{.}))$ such that $ev(h_\nu)=h$. Thus 
$h\in E^1_\Lambda(1_{-\rho})$ by Corollary \ref{iso3} and hence
maps to zero in $H^1(\Gamma,\C)$. Using embedding we also get examples
for all even $n>2$.     
In spite of Lemma \ref{irr2} we do not know whether 
the map $Z^{n-p}_{F,\Lambda}(\vp)\rightarrow
H^p(\Gamma,F\otimes V_\vp)$ is surjective in general.

We have the following dimension formulas.

\begin{kor}\label{est}
Assume that $\Gamma$ is not cocompact. Then
\begin{eqnarray} 
\dim {}^\Gamma C^{-\infty}(\Lambda,V^k(\sigma^{n-p}_{F,\lambda_{n-p}},\vp))&=&
\dim H^p(\Gamma,F\otimes V_\vp)\nonumber\\
&&+\dim \left[E^+_\Lambda(\sigma^{n-p}_{F,\lambda_{n-p}},\vp)/\p_\Gamma^k\left({}^\Gamma C^{-\infty}(\Lambda,V(\sigma^{n-p}_{F,\lambda_{n-p}},\vp))\right)\right]\nonumber\\
&=&\dim H^p(\Gamma,F\otimes V_\vp)
+\dim\:\coker\: res_k\ .\label{rebru}
\end{eqnarray}
Here $res_k: {}^\Gamma C^{-\omega}(\partial X,V^k(\sigma^{n-p}_{F,\lambda_{n-p}},\vp))\rightarrow
C^{-\omega}(B,V^k_B(\sigma^{n-p}_{F,\lambda_{n-p}},\vp))$ is the restriction map.
\end{kor}
\proof
Let $\p_\Lambda$ be the map induced by $\p_\Gamma$ on
${}^\Gamma C^{-\infty}(\Lambda,V^+(\sigma^{n-p}_{F,\lambda_{n-p}},\vp))$. 
Since $\im \p_\Lambda^k\subset E^+_\Lambda(\sigma^{n-p}_{F,\lambda_{n-p}},\vp)$
(see (\ref{mu})) we obtain by Theorem \ref{main2}
\begin{eqnarray*}
\dim {}^\Gamma C^{-\infty}(\Lambda,V^k(\sigma^{n-p}_{F,\lambda_{n-p}},\vp))
&=&\dim \ker \p_\Lambda^k\\
&=&\dim\coker \p_\Lambda^k \\
&=&\dim {}^\Gamma C^{-\infty}(\Lambda,V^+(\sigma^{n-p}_{F,\lambda_{n-p}},\vp))/E^+_\Lambda(\sigma^{n-p}_{F,\lambda_{n-p}},\vp)\\
&&+ \dim E^+_\Lambda(\sigma^{n-p}_{F,\lambda_{n-p}},\vp)/\im\: \p_\Lambda^k\\
&=&\dim H^p(\Gamma,F\otimes V_\vp)+ \dim E^+_\Lambda(\sigma^{n-p}_{F,\lambda_{n-p}},\vp)/\im\:\p_\Lambda^k\ .
\end{eqnarray*}
Equation (\ref{rebru}) is now a consequence of Lemma \ref{irr1}.
\hB

In \cite{bunkeolbrich97} we obtained for $n=2$ by different methods that    
\begin{equation}\label{spargel}
\dim {}^\Gamma C^{-\infty}(\Lambda,V(1_{-\rho-l\alpha}))
=\left\{
\begin{array}{lc}
\dim H^{1}(\Gamma,F_{2l+1})+1&\Gamma\mbox{ elementary }\\
\dim H^{1}(\Gamma,F_{2l+1})&\mbox{ else }
\end{array}\right. \ .
\end{equation}
Here $l\in\nat_0$, and $F_{k}$ is the irreducible $SL(2,\R)$-representation of dimension $k$.
In that paper no natural map ${}^\Gamma C^{-\infty}(\Lambda,V(1_{-\rho-l\alpha}))\rightarrow H^{1}(\Gamma,F_{2l+1})$
was studied. This map is now provided by Theorem \ref{main2}. By (\ref{eieiei}) it is bijective
when restricted to $Z^1_{F_{2l+1},\Lambda}(1)$. Now (\ref{spargel}) implies that
for non-elementary $\Gamma$ we have 
\begin{equation}\label{mohr}
{}^\Gamma C^{-\infty}(\Lambda,V(1_{-\rho-n\alpha}))=Z^1_{F_{2l+1},\Lambda}(1)\ .
\end{equation}
In fact, (\ref{mohr}) has been proved directly by topological considerations
in \cite{bunkeolbrich97}, Lemma 5.2.

\begin{kor}
If $n=2$ and $\Gamma$ is non-elementary, then for all $k\in\frac{1}{2}\nat$
the extension map
$$ ext_\lambda: C^{-\omega}(B,V_B(1_\lambda))\rightarrow {}^\Gamma C^{-\omega}(\partial X,V(1_\lambda))$$
is regular at $\lambda=-k\alpha$.
\end{kor}
\proof
If $k$ is an integer, then ${}^\Gamma C^{-\infty}(\Lambda,V(1_{-k\alpha}))=\{0\}$ by \cite{bunkeolbrich97}, Proposition 4.3. If $k\in\frac{1}{2}+\nat$, 
then $E^1_\Lambda(1_{-k\alpha})=\{0\}$ by (\ref{mohr}) and Lemma \ref{dser}. Thus in any case $E^+_\Lambda(1_{-k\alpha})=\{0\}$.
\hB

\newpage
\section{Applications and related results}\label{apple}

\subsection{Vanishing results}\label{plum}

Let $\Gamma$ be a nontrivial torsion-free convex cocompact subgroup acting on $\R H^n$.
We define the real number $d_\Gamma\in [0,n-1]$ by 
$\delta_\Gamma+\rho=d_\Gamma\alpha$.
Now Theorem \ref{main1} combined with the vanishing result Corollary \ref{iso1}
has the immediate

\begin{kor}\label{van1}
Let $(\vp,V_\vp)$ be a finite-dimensional unitary representation of $\Gamma$. Then
$$ H^p(\Gamma,V_\vp)=\{0\} \quad \mbox{for all }\ p>d_\Gamma+1\ .$$
\end{kor}
\proof Indeed, if $p>d_\Gamma+1$, then $(p-1)\alpha-\rho>\delta_\Gamma$.
\hB

This result is neither new nor the best possible. 
By the celebrated result of Patterson \cite{patterson762} and Sullivan \cite{sullivan79} the number $d_\Gamma$ is equal to the Hausdorff
dimension of the limit set $\Lambda$. Using this and methods from algebraic
topology Izeki \cite{izeki95} proved

\begin{prop}[\cite{izeki95}, Prop. 4.13]\label{van2}
The cohomological dimension of $\Gamma$ is at most $d_\Gamma+1$.
In particular, for any (not necessarily finite-dimensional) representation
$(\vp,V_\vp)$ of $\Gamma$ we have
$$ H^p(\Gamma,V_\vp)=\{0\} \quad \mbox{for all } p>d_\Gamma+1\ .$$
\end{prop}

Nevertheless, it seems to be interesting that this topological statement
has an analytic proof at least in the interesting case of
the trivial $\Gamma$-representation which, moreover, does not use the relation between
the critical exponent and the Hausdorff dimension of the limit set.

We remark (compare \cite{izeki95}, Thm. 4.15) that in view of the exact 
sequence
$$  \dots\rightarrow H^p(Y)\rightarrow H^p(B)\rightarrow H^{p+1}(\bar Y,B)\rightarrow\dots $$
Proposition \ref{van2} as well as Corollary \ref{van1}
combined with Poincar\'e-Lefschetz duality implies Nayatani's vanishing
theorem \cite{nayatani97} 
$$ H^p(B)=0\quad \mbox{for all } p\in (d_\Gamma+1,n-2-d_\Gamma)\ .$$
This yields a cohomology free region if $\delta_\Gamma<-\alpha$ for odd $n$ 
and $\delta_\Gamma<-\frac{3}{2}\alpha$ for even $n$.

We can sharpen Proposition
\ref{van2} for $\Gamma$-representations of the form $F\otimes V_\vp$, where $\vp$ is unitary and $F$ is a finite-dimensional representation of $G$.

\begin{prop}\label{bowa}
Let $F$ be an irreducible finite-dimensional representation with highest
weight $(m_0,m_1,\dots,m_{[\frac{n-1}{2}]})$ as in Section \ref{twist}.
Set 
$$p_F:=\min\{i\:|\: m_i=0\}\in\left\{0,1,\dots\left[\frac{n+1}{2}\right]\right\}$$ 
(for this definition we set $m_i=0$ for
$i>[\frac{n-1}{2}]$). Let $(\vp,V_\vp)$ be a finite-dimensional unitary representation of $\Gamma$. Then
$$ H^p(\Gamma,F\otimes V_\vp)=\{0\} \quad \mbox{for all } p>\min\{d_\Gamma+1, n-p_F\}\ .$$
If $\Gamma$ is cocompact, then $H^p(\Gamma,F\otimes V_\vp)\ne\{0\}$ only if
$p\in [p_F,n-p_F]$.
\end{prop}
\proof 
Let $p=n-k$, where $k<p_F$. Then $m_k\ne 0$. We look at the associated 
$M$-representation $\sigma^k_F$ with highest weight $\mu_k=(\tilde m_1,\dots,\tilde m_{[\frac{n-1}{2}]})$ and the associated
element $\lambda_k\in\aaaa^*$ introduced in Section \ref{twist}. 
Before Proposition \ref{knut} we assigned to $\sigma^k_F$  an element $\lambda_{\sigma^k_F}\in\aaaa^*$ which specifies the end of the corresponding
complementary series. If $k=\frac{n-1}{2}$ we replace the index $k$ by $\sign(m_k)\in\{\pm\}$. If $k=0$, then $\lambda_k=\rho+m_0>\rho\ge \lambda_{\sigma^k_F}$. If $k\ge1$, then
$|\tilde m_k|=m_{k-1}+1\ne 0$. It follows that $\lambda_{\sigma^k_F}\le \rho-k\alpha<\rho+(|m_k|-k\alpha)=\lambda_k$. In any case $\lambda_k>\lambda_{\sigma^k_F}$. Proposition \ref{knut} now yields
${}^{\Gamma} C^{-\infty}(\Lambda,V^+(\sigma^k_{F,\lambda_k},\varphi))=\{0\}$.
Hence $H^p(\Gamma,F\otimes V_\vp)=\{0\}$ by Theorem \ref{main2}.

That $H^p(\Gamma,F\otimes V_\vp)$ vanishes for $p>d_\Gamma+1$ is clear from
Proposition \ref{van2}. Note, however, that for $\delta_\Gamma\ge -\frac{\alpha}{2}$ it is not necessary to refer to Proposition 
\ref{van2}. Indeed, in this case $p>d_\Gamma+1$ implies $p\ge\frac{n+1}{2}$
which yields
$$\lambda_{n-p}\ge \rho-(n-p)\alpha=(p-1)\alpha-\rho>\delta_\Gamma\ .$$
It again follows that ${}^{\Gamma} C^{-\infty}(\Lambda,V^+(\sigma^k_{F,\lambda_k},\varphi))=\{0\}= H^p(\Gamma,F\otimes V_\vp)$.
\hB  

Note that the above proof relies on two essential facts.
First, that the cohomology in degrees $p\ge \frac{n+1}{2}$ is related
to a specific irreducible representation $I^{\sigma^{n-p}_F,\lambda_{n-p}}$ 
of $G$ with infinitesimal character $\chi_F$, and second, that the cohomology
can only be non-zero if $I^{\sigma^{n-p}_F,\lambda_{n-p}}$ is unitary.
The point is that one precisely knows the parameters, where this happens.

For cocompact $\Gamma$ the result has been known before. This case
is usually investigated in the framework of $(\gaaa,K)$-cohomology of irreducible
unitary Harish-Chandra modules (see \cite{borelwallach80}, \cite{wallach88}, Ch. 9, in particular Theorem 9.5.8, \cite{voganzuckerman84};
compare Section \ref{mazzeo}). For cocompact $\Gamma$ and generic $F$,
i.e. $m_i\ne 0$ for all $i$, the proposition means that the cohomology is concentrated in the middle
degree or vanishes, if $n$ is even or odd, respectively.

\subsection{Cohomology with compact support}\label{supp}

As in the Section \ref{heart} we consider a hyperbolic manifold
$Y=\Gamma\backslash\R H^n$, where $\Gamma$ is convex cocompact,
and a finite-dimensional $\Gamma$-representation $(\vp,V_\vp)$.
Here we want to relate the cohomology with compact support $H^p_c(Y, E(\vp))$
with invariant currents supported on the limit set. Note that $E(\vp)$ extends
naturally to $\bar Y$: $E(\vp):= \Gamma\backslash ((X\cup\Omega)\times V_\vp)$. We use the
identification $H^p_c(Y, E(\vp))\cong H^p(\bar Y,B,E(\vp))$, where $H^p(\bar Y,B,E(\vp))$ denotes the relative cohomology of $\bar Y$ with respect to its
boundary $B$. For simplicity we only treat the case of even $n$. In
this case the complex (\ref{kike}) computes $H^p(Y,E(\vp))\cong H^p(\bar Y,E(\vp))$.
Thus $res$ induces a map $res_p: H^p(Y,E(\vp))\rightarrow H^p(B,E(\vp))$.
One expects that $res_p$ coincides with the natural map $i^*_p: H^p(\bar Y,E(\vp))\rightarrow H^p(B,E(\vp))$ induced by the embedding 
$B\hookrightarrow \bar Y$. This is indeed the case, see Proposition \ref{double}. Thus the cohomology of the complex of kernels of $res$, i.e.
$$ ({}^\Gamma[\Omega^*_{-\infty}(\Lambda)\otimes V_\vp], d) \mbox{ completed by } {}^\Gamma V_\vp\ ,$$
should be related to $H^p(Y,B,E(\vp))$. We will clarify this relation
in Proposition \ref{ancient}. That these two cohomology groups are not isomorphic in general
is caused by the lack of surjectivity of $res$. However, for unitary $\vp$
and $p\le \frac{n}{2}$ they are isomorphic.  

For any $\Gamma$-module $V$ we denote by $C^*_\Gamma V$ the standard 
group cohomology complex:
$$ C_\Gamma^p V:=\{f:\Gamma^{p+1}\rightarrow V\:|\: f(\gamma\gamma_0,\dots,\gamma\gamma_p)=\gamma f(\gamma_0,\dots,\gamma_p)\}$$
with differential
$$ \partial_pf(\gamma_0,\dots,\gamma_{p+1})=\sum_{i=0}^{p+1} (-1)^i f(\gamma_0,
\dots,\hat\gamma_i,\dots,\gamma_{p+1})\ .$$
Then $H^p(C^*_\Gamma V)=H^p(\Gamma,V)$.
If $V^*$ is a complex of $\Gamma$-modules, then $C_\Gamma^*V^*$ should
denote the total complex of the arising double complex.
There is a natural embedding of complexes
$$ c: {}^\Gamma V^*\hookrightarrow C_\Gamma^*V^* $$
sending ${}^\Gamma V^p$ to the constant functions in $C_\Gamma^0 V^p$.
If $V^*$ is a complex of acyclic $\Gamma$-modules, then $c$ is a quasi-isomorphism. If $\Phi: V^*\rightarrow W^*$ is a quasi-isomorphism of complexes of $\Gamma$-modules, then $C_\Gamma\Phi: C^*_\Gamma V^*\rightarrow C^*_\Gamma W^*$
is a quasi-isomorphism, too.

We will denote the completed de Rham complex (compare (\ref{kike}))
$$ 0\rightarrow
{}\Omega_{-\omega}^0(\partial X)
\stackrel{d}{\longrightarrow}
\Omega_{-\omega}^1(\partial X)
\stackrel{d}{\longrightarrow}
\dots \stackrel{d}{\longrightarrow} 
\Omega_{-\omega}^{n-1}(\partial X)
\stackrel{\int_{\partial X}}{\longrightarrow} \C
\rightarrow 0 $$
by $\dot\Omega_{-\omega}^*(\partial X)$.

\begin{lem}\label{fisch}
The embedding $c: {}^\Gamma(\dot\Omega^*_{-\omega}(\partial X)\otimes V_\vp)
\rightarrow C_\Gamma^*(\dot\Omega^*_{-\omega}(\partial X)\otimes V_\vp)$
is a quasi-isomorphism.
\end{lem}
\proof
Since the $\Gamma$-modules $\Omega^p(X)\otimes V_\vp$ are acyclic 
(\cite{bunkeolbrich950}, Lemma 2.4), the embedding
$$c_X: {}^\Gamma(\Omega^*(X)\otimes V_\vp)=\Omega^*(Y,E(\vp))
\hookrightarrow C_\Gamma^*(\Omega^*(X)\otimes V_\vp) $$
is a quasi-isomorphism. As already noted, by Theorem \ref{krah} suitably
normalized Poisson transforms define a quasi-isomorphisms
\begin{eqnarray*} 
P: \dot\Omega^*_{-\omega}(\partial X)\otimes V_\vp
&\rightarrow&  \Omega^*(X)\otimes V_\vp\ ,\\
{}^\Gamma P: {}^\Gamma(\dot\Omega^*_{-\omega}(\partial X)\otimes V_\vp)
&\rightarrow&  \Omega^*(Y,E(\vp))\ .
\end{eqnarray*}
In particular, $C_\Gamma P: C_\Gamma^*(\dot\Omega^*_{-\omega}(\partial X)\otimes V_\vp)
\rightarrow C_\Gamma^*(\Omega^*(X)\otimes V_\vp)$
is a quasi-isomorphism, too. We have
$C_\Gamma P\circ c=c_X\circ {}^\Gamma P$.
It follows that $c$ is a quasi-isomorphism. 
\hB

The precise normalization of the Poisson transform $P: \dot\Omega^*_{-\omega}(\partial X)
\otimes V_\vp
\rightarrow\Omega^*(X)\otimes V_\vp$ is given as follows. For
$p=0,\dots,n-1$ set $c_p:=\frac{(n-1)(n-3)\dots(n-2p+1)}{(n-1)(n-2)\dots
(n-p)}$. Then
$$ P_{|\Omega^p(X)\otimes V_\vp}:= c_p P_{p,\rho-p\alpha}\otimes\id\ ,\ p=0,1,\dots,n-1\ ;\quad P_{|\dot\Omega^n(X)\otimes V_\vp}:=c_{n-1}\frac{1-n}{|\alpha|^n}\vol_X\otimes\id\ .$$ 

\begin{prop}\label{double}
We consider the embeddings 
$i: B\rightarrow \bar Y$ and $j: Y\rightarrow \bar Y$ and the restriction map
$res: {}^\Gamma(\dot\Omega^*_{-\omega}(\partial X)\otimes V_\vp)\rightarrow
\Omega_{-\omega}^*(B,E(\vp))$. Here $res$ is defined to be zero on $\dot\Omega^n_{-\omega}(\partial X)\otimes V_\vp)$. Then the following
diagram commutes
$$\begin{array}{ccc}
H^p\left({}^\Gamma(\dot\Omega^*_{-\omega}(\partial X)\otimes V_\vp)\right)&\stackrel{({}^\Gamma P)_p}{\longrightarrow}&
H^p(Y,E(\vp))\\
&&\uparrow{\scriptstyle  j_p^*}\\
{\downarrow}{\scriptstyle res_p}&&H^p(\bar Y, E(\vp)) \\
&&\downarrow{\scriptstyle  i_p^*}\\
H^p\left(\Omega_{-\omega}^*(B,E(\vp))\right)&=&H^p(B,E(\vp))
\end{array}\ .
$$
$({}^\Gamma P)_p$ and $j^*_p$ are isomorphisms.
\end{prop}
\proof 
We also consider the embeddings $I: \Omega\rightarrow X\cup\Omega$ and $J: X\rightarrow X\cup\Omega$. $V_\vp$ can be embedded as the subcomplex
of constant $0$-forms into any $V_\vp$-valued de Rham complex. In particular,
the embeddings $V_\vp\hookrightarrow \Omega^*(X)\otimes V_\vp$,  $V_\vp\hookrightarrow \Omega^*(X\cup\Omega)\otimes V_\vp$, and $V_\vp\hookrightarrow \dot\Omega^*(\partial X)\otimes V_\vp$ are quasi-isomorphisms. We look at the following diagram
\begin{equation}\label{muenster}
\begin{array}{ccccc}
\Omega^*(Y, E(\vp))&\stackrel{{}^\Gamma P}{\longleftarrow}&{}^\Gamma(\dot\Omega^*_{-\omega}(\partial X)\otimes V_\vp)&\stackrel{res}{\longrightarrow}&\Omega^*_{-\omega}(B,E(\vp))\\
&&\downarrow {\scriptstyle  c}&& \downarrow {\scriptstyle  c_B}\\
\downarrow{\scriptstyle c_X}&&C^*_\Gamma(\dot\Omega^*_{-\omega}(\partial X)\otimes V_\vp)&\stackrel{C_\Gamma res_\Omega}{\longrightarrow}&C^*_\Gamma(\Omega^*_{-\omega}(\Omega)\otimes V_\vp)\\
&{\swarrow}{\scriptstyle C_\Gamma P}&\uparrow&\nearrow&\\
C^*_\Gamma(\Omega^*(X)\otimes V_\vp)&\longleftarrow&C^*_\Gamma V_\vp&&\uparrow\\
&\nwarrow{\scriptstyle C_\Gamma J^*}&\downarrow&\searrow&\\
\uparrow{\scriptstyle c_X}&&C^*_\Gamma(\Omega^*(X\cup\Omega)\otimes V_\vp)&\stackrel{C_\Gamma I^*}{\longrightarrow}&
C^*_\Gamma(\Omega^*(\Omega)\otimes V_\vp)\\
&&\uparrow{\scriptstyle  c}&&\uparrow{\scriptstyle  c_B}\\
\Omega^*(Y, E(\vp))&\stackrel{j^*}{\longleftarrow}&\Omega^*(\bar Y, E(\vp))&\stackrel{i^*}{\longrightarrow}&\Omega^*(B,E(\vp))
\end{array}\qquad .
\end{equation}
Using Lemma \ref{fisch} and the $\Gamma$-acyclicity of 
$\Omega^p(X)\otimes V_\vp$, $\Omega^p(X\cup\Omega)\otimes V_\vp$, $\Omega^p(\Omega)\otimes V_\vp$ (\cite{bunkeolbrich950}, Lemma 2.4), and
$\Omega^p_{-\omega}(\Omega)\otimes V_\vp$ (Lemma \ref{acB} below) we see that all vertical arrows and all arrows pointing to the left are quasi-isomorphisms.
The left and the right column induce the identity on $H^p(Y,E(\vp))$ and
$H^p(B,E(\vp))$, respectively. The proposition now follows from the observation
that all subdiagrams of (\ref{muenster}) commute.
\hB

The image $\im\:res$ of the restriction map $res: {}^\Gamma(\dot\Omega^*_{-\omega}(\partial X)\otimes V_\vp)\rightarrow
\Omega_{-\omega}^*(B,E(\vp))$ is a subcomplex of $\Omega_{-\omega}^*(B,E(\vp))$.
Let $\coker\:res$ denote the corresponding quotient complex.

\begin{prop}\label{ancient}
The cohomology of the finite-dimensional complex ${}^\Gamma(\dot\Omega^*_{-\infty}(\Lambda)\otimes V_\vp)$ of invariant currents
supported on the limit set fits into the following long exact sequence:
\begin{eqnarray*} 
&\dots&\rightarrow H^{p-2}(\coker\:res)\stackrel{\delta}{\longrightarrow}
H^p({}^\Gamma(\dot\Omega^*_{-\infty}(\Lambda)\otimes V_\vp))\rightarrow
H^p(\bar Y,B,E(\vp))\rightarrow H^{p-1}(\coker\:res)\\
&&\stackrel{\delta}   
{\longrightarrow}
H^{p+1}({}^\Gamma(\dot\Omega^*_{-\infty}(\Lambda)\otimes V_\vp))\rightarrow\dots\ .
\end{eqnarray*}
The map $\delta$ can be described as follows. Let $[\eta]\in H^{p-2}(\coker\:res)$ be represented by $\eta\in \Omega_{-\omega}^{p-2}(B,E(\vp))$ such that $d\eta=res(\omega)$ for some
$\omega\in {}^\Gamma(\dot\Omega^{p-1}_{-\omega}(\partial X)\otimes V_\vp)$.
Then $\delta [\eta]=[d\omega]\in H^p({}^\Gamma(\dot\Omega^*_{-\infty}(\Lambda)\otimes V_\vp))$. 
\end{prop}
\proof
We form the mapping cone $(C^*(res),d)$ which is given by
$$ C^p(res)={}^\Gamma(\dot\Omega^p_{-\omega}(\partial X)\otimes V_\vp)\oplus
\Omega_{-\omega}^{p-1}(B,E(\vp))\ ,\quad d(\omega,\eta)=(d\omega,res(\omega)-d\eta)\ .$$
Proposition \ref{double} implies that it is quasi-isomophic to the mapping
cone $C^*(i^*)$,\linebreak[4] $i^*: \Omega^*(\bar Y,E(\vp))\rightarrow \Omega^*(B,E(\vp))$.
The cohomology of $C^*(i^*)$ is equal to $H^*(\bar Y,B,E(\vp))$. Thus
\begin{equation}\label{konrad}
H^p(C^*(res))\cong H^p(\bar Y,B,E(\vp))\ .
\end{equation}
We consider the exact sequence of complexes
\begin{equation}\label{lorenz}
0\rightarrow {}^\Gamma(\dot\Omega^*_{-\infty}(\Lambda)\otimes V_\vp)) 
\rightarrow C^*(res)\rightarrow  D^*(res) \rightarrow 0\ ,
\end{equation}
where $D^p(res):=\im\: res^{p}\oplus \Omega^{p-1}(B,E(\vp))$, 
$d(\omega,\eta)=(d\omega,\omega-d\eta)$. It is now easy to check that
the map $D^p(res)\ni (\omega,\eta)\mapsto (-1)^p [\eta]\in\coker\:res^{p-1}$
defines a quasi-isomorphism between $D^*(res)$ and $(\coker\:res)[1]$.
This together 
with (\ref{konrad}) and the long exact sequence associated to (\ref{lorenz})
implies the proposition.
\hB

Let $Z^p_\Lambda(\vp)$ be the space of $p$-cocycles of ${}^\Gamma(\dot\Omega^*_{-\infty}(\Lambda)\otimes V_\vp))$ as in Section \ref{heart}. 

\begin{kor}\label{helga}
If $\vp$ is unitary and $p\le\frac{n}{2}$, then the map
$Z^p_\Lambda(\vp)\rightarrow H^p(\bar Y,B,E(\vp))$ given by Proposition
\ref{ancient} is an isomorphism.
\end{kor}
\proof
Proposition \ref{newer} implies for all $q<p$ that $\coker\:res^q=0$ and
that ${}^\Gamma(\dot\Omega^q_{-\infty}(\Lambda)\otimes V_\vp))=Z^q_\Lambda(\vp)$. Thus $H^p({}^\Gamma(\dot\Omega^*_{-\infty}(\Lambda)\otimes V_\vp))\rightarrow
H^p(\bar Y,B,E(\vp))$ is an isomorphism and $H^p({}^\Gamma(\dot\Omega^*_{-\infty}(\Lambda)\otimes V_\vp))=Z^p_\Lambda(\vp)$.
\hB

For $p\le\frac{n}{2}$ and general $\vp$ one can describe the map
$Z^p_\Lambda(\vp)\rightarrow H^*(\bar Y,B,E(\vp))$ more explicitly.
The standard relative cohomology complex $\Omega^*(\bar Y,B,E(\vp))$ for the pair  
$(\bar Y,B)$ is given by
$$ \Omega^p(\bar Y,B,E(\vp)):=\{\omega\in \Omega^p(\bar Y,E(\vp))\:|\:i^*\omega=0\}\ .$$
Indeed, the natural embedding $\Omega^*(\bar Y,B,E(\vp))\hookrightarrow
C^*(i^*)$ is a quasi-isomorphism. Using the asymptotic behaviour for $a\to\infty$ of $P_{p,\rho-p\alpha}\omega(ka)$, $\omega\in Z^p_\Lambda(\vp)$, $kM\in\Omega$ (\cite{bunkeolbrich950}, Lemma 6.2, Equation (37)) it is
not difficult to show that $c_pP_{p,\rho-p\alpha}\omega$ defines an element
in $\Omega^p(\bar Y,B,E(\vp))$ (this is not true for $p>\frac{n}{2}$), which
is closed and represents the correct cohomology class in $H^p(\bar Y,B,E(\vp))$.

Let $(\pi,F)$ be an irreducible finite-dimensional representation of $G$. 
Proposition \ref{ancient} in particular yields a map
$$ r_p: Z^p_{\Lambda}(\pi\otimes\vp)\rightarrow H^p(\bar Y,B,E(\pi\otimes\vp))\ .$$

\begin{kor}\label{comptwist}
Let $D_\pi: Z^p_{F,\Lambda}(\vp)\rightarrow Z^p_{\Lambda}(\pi\otimes\vp)$ be the embedding constructed in Proposition \ref{trans}. If $n$ is even, $\vp$ is unitary, and
$p\le\frac{n}{2}$, then 
$$r_p\circ D_\pi: Z^p_{F,\Lambda}(\vp)\rightarrow H^p(\bar Y,B,E(\pi\otimes\vp))$$ is an isomorphism.
\end{kor}
\proof
By the argument at the beginning of the proof of Theorem \ref{main2} and since the complex of sheaves
$\Omega_{-\omega}^*(.,E(\pi\otimes\vp))^{\chi_F}=C^{-\omega}(.,V_B(\sigma^*_{F,\lambda_*},\vp))$ on $B$ is a flabby resolution of the sheaf of
parallel sections of $E(\pi\otimes\vp)$ on $B$ Proposition \ref{double} applied $\pi\otimes\vp$ remains true if we
replace ${}^\Gamma(\dot\Omega^*_{-\omega}(\partial X)\otimes F\otimes V_{\vp})$
by (\ref{zhelo}) and $\Omega_{-\omega}^*(B,E(\pi\otimes\vp))$ by $C^{-\omega}(B,V_B(\sigma^*_{F,\lambda_*},\vp))$. Thus we can also replace in Proposition \ref{ancient} the complexes ${}^\Gamma(\dot\Omega^*_{-\infty}(\Lambda)\otimes F\otimes V_{\vp})$ and $\coker\:res$
by their subcomplexes having infinitesimal character $\chi_F$. Now one can argue as in the proof of Corollary
\ref{helga}.
\hB

\subsection{The Harder-Borel conjecture}\label{nancy}

Let $G$ be a connected real reductive group, $K\subset G$ be a maximal compact
subgroup, $\Gamma\subset G$ be a finitely generated discrete subgroup. Let 
$(\pi,F)$ and $(\vp,V_\vp)$ be finite-dimensional representations
of $G$ and $\Gamma$, respectively. We assume that $F$ is irreducible, hence has
infinitesimal character $\chi_F$. If $\Gamma$ is cocompact, we require $\vp$ to
be unitary. 
One would like to study the cohomology groups
$$ H^p(\Gamma, F\otimes V_\vp) $$
by means of complexes of {\em automorphic} differential forms. Of course, the de Rham complex
$$\Omega^*(Y,E(\pi\otimes\vp)):={}^\Gamma(\Omega^*(X)\otimes F\otimes V_\vp)\ $$
computes these cohomology groups.
If $\Gamma$ is torsion-free, then $Y:=\Gamma\backslash X$ is a manifold, and
this complex is the de Rham complex associated to the flat vector bundle
$E(\pi\otimes\vp)$ as considered in Section \ref{sur}. By Selberg's Lemma
(see e.g. \cite{ratcliffe94}) there is always a torsion-free normal subgroup 
$\Gamma_0\subset \Gamma$ of finite index. Then the finite group
$\Gamma/\Gamma_0$ acts on the manifold $Y_0:=\Gamma_0\backslash X$, and we have 
\begin{eqnarray}
\Omega^*(Y,E(\pi\otimes\vp))={}^{\Gamma/\Gamma_0}\Omega^*(Y_0,E(\pi\otimes\vp)) \ ,\label{kung}\\
H^p(\Gamma, F\otimes V_\vp)={}^{\Gamma/\Gamma_0} H^p(\Gamma_0, F\otimes V_\vp)\ .
\label{fu}
\end{eqnarray}
These equations show that up to an action of a finite group, which is harmless,
we are still in the situation of the preceding sections.
 
Fix a $G$-invariant Riemannian metric
on $X=G/K$ and an admissible positive definite scalar product on $F$: it has
to be $K$-invariant and the Iwasawa $\aaaa$ has to act by selfadjoint endomorphisms. By the usual twisting isomorphism we look at $X\times F\cong
G\times_K F$ as a $G$-homogeneous vector bundle on $X$ which now comes with
a $G$-invariant Hermitian metric. It induces a Hermitian metric on 
$\Lambda^pT^*Y\otimes E(\pi)$ and thus gives rise to a codifferential
$$ \delta=\delta_F: \Omega^{p+1}(Y,E(\pi))\rightarrow \Omega^p(Y,E(\pi))$$
and a corresponding Laplacian  
$$ \Delta_F=\delta_F d+d\delta_F: \Omega^{p}(Y,E(\pi))\rightarrow \Omega^p(Y,E(\pi))$$
as in Section \ref{sur}. The codifferential $\delta_F$ can be twisted with the flat
connection of $E(\vp))$ as in (\ref{Hut}), and we obtain operators
$$ \delta_F: \Omega^{p+1}(Y,E(\pi\otimes\vp))\rightarrow \Omega^p(Y,E(\pi\otimes\vp))\ ,\ \  \Delta_F=\delta_F d+d\delta_F: \Omega^{p}(Y,E(\pi\otimes\vp))\rightarrow \Omega^p(Y,E(\pi\otimes\vp))\ .$$
  
A form $\omega\in \Omega^*(Y,E(\pi\otimes\vp))$ is said to have
moderate growth, if for all $k\in\nat_0$ there exist
constants $C,r$ such that with respect to some chosen norm on $V_\vp$ 
$$ |\Delta^k_F\omega (x)|\le C e^{r\dist(x,eK)}\ .$$

The $G$-action on $\Omega^*(X)\otimes F$ induces an action of the center $\cZ(\gaaa)$ of the universal enveloping algebra of $\gaaa$ on $\Omega^*(Y,E(\pi\otimes\vp))$. In particular, it makes sense to say that
a form has generalized infinitesimal character. Note that the Casimir operator
acts by $-\Delta_F+\kappa_F$. Let
us introduce the following subcomplexes of $\Omega^*(Y,E(\pi\otimes\vp))$:
\begin{itemize}
\item $\Omega_{mg}^*(Y,E(\pi\otimes\vp))$ (forms of moderate growth)
\item $\cA^*_F(Y,E(\pi\otimes\vp)):=\Omega_{mg}^*(Y,E(\pi\otimes\vp))^{\chi_F}$ (automorphic forms)
\item $\cA_{F,ch}^*(Y,E(\pi\otimes\vp)):=\cA^*_F(Y,E(\pi\otimes\vp))\cap \Omega^*(Y,E(\pi\otimes\vp))_{\Delta_F,\delta_F}$ (coclosed harmonic automorphic forms)
\end{itemize}
Here, as in Section \ref{twist}, the superscript $\chi_F$ means forms having
generalized infinitesimal character $\chi_F$. Some comments concerning these
definitions are in order. First, the definition of moderate growth given
here is slightly stronger than the one given in Section \ref{heart}. This is necessary in order to make $\Omega_{mg}^*(Y,E(\pi\otimes\vp))$ into a complex.
$\Omega^*(Y,E(\pi\otimes\vp))^{\chi_F}$ is a subcomplex of $\Omega^*(Y,E(\pi\otimes\vp))_{(\Delta_F)}$. This shows that the weaker definition of moderate growth would lead to the same space of automorphic forms.  
The Harder-Borel conjecture,
following Gaillard \cite{gaillard01}, now asserts

\begin{con}\label{HBC}
The inclusions $\cA_{F,ch}^*(Y,E(\pi\otimes\vp))
\subset \cA^*_F(Y,E(\pi\otimes\vp))\subset \Omega^*(Y,E(\pi\otimes\vp))$
are quasi-isomorphisms.
\end{con}

For cocompact $\Gamma$ we observe $\Omega_{mg}^*(Y,E(\pi\otimes\vp))=\Omega^*(Y,E(\pi\otimes\vp))$ and
thus Conjecture \ref{HBC} is a direct consequence of classical Hodge theory.

The conjecture arose in the work of Harder, Borel and others (see e.g. \cite{borel83} and the survey article \cite{schwermer90}) on
the cohomology
of arithmetic subgroups of reductive algebraic groups which is in fact the most
interesting case. Here one usually considers the case of trivial $\vp$, only.
Indeed, in this situation (or more generally, if $\Gamma\backslash G$ has finite
volume) $\cA_{F,ch}^*(Y,E(\pi))$ in contrast to $\Omega^*(Y,E(\pi))_{\Delta_F,\delta_F}$ is finite-dimensional
and carries a lot of additional structure. However, it is at least of theoretical interest to understand the generality in which statements like
Conjecture \ref{HBC} hold. Gaillard \cite{gaillard01} formulates and discusses
a corresponding conjecture for trivial $\pi$ and $\vp$, but for a discrete subgroup of an arbitrary connected
Lie group.

The most far reaching result concerning Conjecture \ref{HBC} known up
to now is the following

\begin{theorem}[Franke \cite{franke98}]
Let $G$ be reductive algebraic defined over $\Q$, $\Gamma\subset G$ be a congruence subgroup. Let $(\pi,F)$ be a rational representation of $G$. Then the inclusion $\cA^*_F(Y,E(\pi))$ $\subset \Omega^*(Y,E(\pi))$ is a quasi-isomorphism.
\end{theorem}

For $G$ of $\Q$-rank one the result was previously obtained by Speh and
Casselman, see \cite{casselman84} for a discussion of the case $G=SL(2,\R)$.
Compare also Harder's early results \cite{harder75}.
For $G$ of general rank it is not known whether the inclusion  
$\cA_{F,ch}^*(Y,E(\pi))
\subset \Omega^*(Y,E(\pi))$ is a quasi-isomorphism. 
However, \cite{bunkeolbrich980}, Theorem 6.3 implies

\begin{prop}\label{smom}
If $G$ has real rank one, $\Gamma\subset G$ has finite covolume, and $\vp$ is trivial, then Conjecture \ref{HBC} holds true.
\end{prop}
 
We will shortly discuss the proof of the proposition in order to indicate
what kind of results one could try to prove if one wants to establish
Conjecture \ref{HBC}. First of all one needs that $\Omega_{mg}^*(Y,E(\pi))$
is quasi-isomorphic to $\Omega^*(Y,E(\pi))$. This was shown by Borel \cite{borel80} (compare also \cite{borel83}, \cite{bunkeolbrich980}, Theorem 5.6, \cite{franke98}, Section 2.3). Now one looks at the image of the 
Laplacian
$$\Delta_F=-\Omega+\kappa_F: \Omega_{mg}^{p}(Y,E(\pi))\rightarrow \Omega_{mg}^p(Y,E(\pi))\ .$$
\cite{bunkeolbrich980}, Theorem 6.3 tells us that it is the orthogonal complement of rapidly decreasing harmonic forms, 
which are automatically closed and coclosed. Arguing as in the proof of Theorem
\ref{ho} we find that the complexes $\Omega_{mg}^*(Y,E(\pi))_{(\Delta_F)}$ and $\Omega_{mg}^*(Y,E(\pi))_{\Delta_F,\delta_F}$ are quasi-isomorphic to $\Omega_{mg}^*(Y,E(\pi))$. The rank one assumption now ensures that the action of $\cZ(\gaaa)$ on these complexes is locally finite. Thus taking the components
with generalized infinitesimal character $\chi_F$ does not change the cohomology
(see the proof of Proposition \ref{main2}). This completes the proof of Proposition \ref{smom}.

One could try to extend this approach to, say, convex cocompact or even
geometrically finite discrete subgroups of rank one Lie groups. It seems to
be not too difficult to establish that $\Omega_{mg}^*(Y,E(\pi))$
is quasi-isomorphic to $\Omega^*(Y,E(\pi))$. However, the determination of
the image
of the Laplacian on forms of moderate growth seems to be out of reach in the
moment. The crucial point is to show that the Laplacian has closed range.

However, if $G=SO(1,n)_0$ or $G=Spin(1,n)$ and $\Gamma$ is convex cocompact, then
we can use Theorem
\ref{main2} and or rather the tools developed for its proof in order to establish Conjecture \ref{HBC}.

\begin{theorem}\label{blub}
If $G=SO(1,n)_0$ or $G=Spin(1,n)$ and $\Gamma$ is convex cocompact, then
Conjecture \ref{HBC} is true.
\end{theorem}

For the special case of elementary $\Gamma$ (and $\pi=\vp=1$) the theorem was previously obtained by Delacroix \cite{delacroix00}. The proof of the theorem will proceed in several steps. 
For cocompact $\Gamma$ there is nothing to show. We thus assume $\Gamma$
to be non-cocompact. In view of (\ref{kung}) and (\ref{fu}) we may also assume
$\Gamma$ to be torsion-free.

We introduce close relatives of the complexes $\cA_{F}^*(Y,E(\pi\otimes\vp))$ and $\cA_{F,ch}^*(Y,E(\pi\otimes\vp))$
which can be directly investigated by means of Theorem \ref{main2} and which coincide with the original ones in case of trivial $\pi$. Recall the definitions
of the codifferential $\hat\delta$ and the corresponding Laplacian $\hat\Delta$
from Section \ref{sur}. We set
\begin{eqnarray*}
\hat\cA_{F}^*(Y,E(\pi\otimes\vp))&:=&\cA_{F}^*(Y,E(\pi\otimes\vp))\cap
\Omega^*(Y,E(\pi\otimes\vp))_{(\hat\Delta)}\ ,\\
\hat\cA_{F,ch}^*(Y,E(\pi\otimes\vp))&:=&\cA_{F}^*(Y,E(\pi\otimes\vp))\cap
\Omega^*(Y,E(\pi\otimes\vp))_{\hat\Delta,\hat\delta}\ .
\end{eqnarray*}
These complexes of automorphic forms are as natural as the corresponding
unhatted versions. We will first prove the analog of Conjecture \ref{HBC}
for them (see Propositions \ref{HBC1} and \ref{HBC2} below).

The following observation holds for any discrete subgroup $\Gamma$ of a real reductive group $G$.

\begin{lem}\label{xxx}
The inclusion $\hat\cA_{F}^*(Y,E(\pi\otimes\vp))\subset \cA_{F}^*(Y,E(\pi\otimes\vp))$ is a quasi-isomorphism.
\end{lem}
\proof $\hat\Delta$ defines a $G$-invariant differential operator acting
on $\Omega^*(X,E(\pi))=\Omega^*(X)\otimes F$. The algebra of $G$-invariant
differential operators acting on a homogeneous vector bundle on $X$ is
a finitely generated module over $\cZ(\gaaa)$. This implies that $\C[\hat\Delta]$ acts locally finitely on $\Omega^*(X,E(\pi))^{\chi_F}$, thus also on
$\cA_{F}^*(Y,E(\pi\otimes\vp))$. As in (\ref{klu}) we obtain a 
direct sum decomposition 
\begin{equation}\label{reinhold}
\cA_{F}^*(Y,E(\pi\otimes\vp))=\hat\cA_{F}^*(Y,E(\pi\otimes\vp))\oplus
\cA_{F}^*(Y,E(\pi\otimes\vp))_R\ ,
\end{equation}
where $\cA_{F}^*(Y,E(\pi\otimes\vp))_R$ is the sum of all generalized eigenspaces of $\hat\Delta$ corresponding to non-zero eigenvalues. From this the lemma follows easily.
\hB

Note that there is no similar simple relation between $\hat\cA_{F,ch}^*(Y,E(\pi\otimes\vp))$ and $\cA_{F,ch}^*(Y,E(\pi\otimes\vp))$.
In particular, $\hat\Delta$ and $\delta_F$ do not commute, hence $\hat\Delta$
does not act on $\cA_{F,ch}^*(Y,E(\pi\otimes\vp))$.

We now return to the case of convex cocompact $\Gamma$ acting on $\R H^n$.

\begin{prop}\label{HBC1}
The inclusion $\hat\cA_{F}^*(Y,E(\pi\otimes\vp))\subset \Omega^*(Y,E(\pi\otimes\vp))$ is a quasi-isomorphism.
\end{prop}
\proof It is clear that this inclusion induces an isomorphism on
$H^0$. Let $p>0$. Then
by Proposition \ref{bow} and Proposition \ref{trans} the space 
$Z^p(\hat\cA_{F}^*(Y,E(\pi\otimes\vp)))$ of $p$-cocycles of $\hat\cA_{F}^*(Y,E(\pi\otimes\vp))$ is isomorphic to ${}^\Gamma C^{-\infty}(\partial X, V^+(\sigma^{n-p}_{F,\lambda_{n-p}},\vp))$. For $l\ge k_-(\sigma^{n-p}_{F,\lambda_{n-p}},\vp))$ we set
\begin{eqnarray}\label{endlich}
E^+_{-\infty}(\sigma^{n-p}_{F,\lambda_{n-p}},\vp)&:=&ext[l]\left( 
C^{-\infty}(B, V_B^+(\sigma^{n-p}_{F,\lambda_{n-p}},\vp))\right)\nonumber\\
&=&E^+(\sigma^{n-p}_{F,\lambda_{n-p}},\vp) \cap C^{-\infty}(\partial X, V^+(\sigma^{n-p}_{F,\lambda_{n-p}},\vp))\ .\nonumber
\end{eqnarray}
Arguing as in the proof of Theorem \ref{main1} we find that
$$ {}^\Gamma 
C^{-\infty}(\partial X, V^+(\sigma^{n-p}_{F,\lambda_{n-p}},\vp))/
E^+_{-\infty}(\sigma^{n-p}_{F,\lambda_{n-p}},\vp) \cong 
{}^\Gamma C^{-\infty}(\Lambda, V^+(\sigma^{n-p}_{F,\lambda_{n-p}},\vp))/
E^+_\Lambda(\sigma^{n-p}_{F,\lambda_{n-p}},\vp)\ .$$
We obtain the following commutative diagram
$$
\begin{array}{ccc}
{}^\Gamma 
C^{-\infty}(\partial X, V^+(\sigma^{n-p}_{F,\lambda_{n-p}},\vp))/
E^+_{-\infty}(\sigma^{n-p}_{F,\lambda_{n-p}},\vp)&\stackrel{a}{\longleftarrow}&
{}^\Gamma C^{-\infty}(\Lambda, V^+(\sigma^{n-p}_{F,\lambda_{n-p}},\vp))/
E^+_\Lambda(\sigma^{n-p}_{F,\lambda_{n-p}},\vp)\\
\downarrow{\scriptstyle b}&&\downarrow{\scriptstyle c}\\
H^p(\hat\cA_{F}^*(Y,E(\pi\otimes\vp)))&\stackrel{i}{\longrightarrow}&
H^p(\Gamma,F\otimes V_\vp)\ .
\end{array}
$$
By Theorem \ref{main2} the map $c$ is an isomorphism. Since $a$ is an isomorphism, $b$ is surjective, and $c=i\circ b \circ a$, we conclude
that $i$ (as well as $b$) is an isomorphism. The proposition follows.
\hB 

Proposition \ref{HBC1} together with Lemma \ref{xxx} shows that the 
inclusion $\cA^*_F(Y,E(\pi\otimes\vp))\subset \Omega^*(Y,E(\pi\otimes\vp))$
is a quasi-isomorphism which proves the first half of Theorem \ref{blub}.

\begin{lem}\label{otto}
$E^+_{-\infty}(\sigma^{n-p}_{F,\lambda_{n-p}},\vp)=\p_\Gamma\left({}^\Gamma 
C^{-\infty}(\partial X, V^+(\sigma^{n-p}_{F,\lambda_{n-p}},\vp))\right)\ .$
\end{lem}
\proof We set $E^\prime:=\p_\Gamma\left({}^\Gamma 
C^{-\infty}(\partial X, V^+(\sigma^{n-p}_{F,\lambda_{n-p}},\vp))\right)$.
As in Lemma \ref{vorzug} we obtain for $k\ge k_-(\sigma^{n-p}_{F,\lambda_{n-p}},\vp))$
$$ E^+_{-\infty}(\sigma^{n-p}_{F,\lambda_{n-p}},\vp)=\p_\Gamma^k\left({}^\Gamma 
C^{-\infty}(\partial X, V^+(\sigma^{n-p}_{F,\lambda_{n-p}},\vp))\right)
\subset E^\prime\ .$$
On the other hand we have by Corollary \ref{preth} and Corollary \ref{beau}
$$ E^\prime\subset \im\p_\Gamma\cap C^{-\infty}(\partial X, V^+(\sigma^{n-p}_{F,\lambda_{n-p}},\vp))= E^+_{-\infty}(\sigma^{n-p}_{F,\lambda_{n-p}},\vp)\ .$$
The lemma follows.
\hB

\begin{kor}\label{paulernst}
$H^p(\hat\cA_{F}^*(Y,E(\pi\otimes\vp)))\cong Z^p(\hat\cA_{F}^*(Y,E(\pi\otimes\vp)))/\hat\Delta\left( Z^p(\hat\cA_{F}^*(Y,E(\pi\otimes\vp)))\right)\ .$
\end{kor}
\proof For $p=0$ the assertion is obvious. If $p>0$, then the space on the 
right hand side is isomorphic to (compare Corollary \ref{preth})
$${}^\Gamma 
C^{-\infty}(\partial X, V^+(\sigma^{n-p}_{F,\lambda_{n-p}},\vp))/\p_\Gamma\left({}^\Gamma 
C^{-\infty}(\partial X, V^+(\sigma^{n-p}_{F,\lambda_{n-p}},\vp))\right)\ ,$$
which in turn is equal to
$$ {}^\Gamma 
C^{-\infty}(\partial X, V^+(\sigma^{n-p}_{F,\lambda_{n-p}},\vp))/
E^+_{-\infty}(\sigma^{n-p}_{F,\lambda_{n-p}},\vp)$$
by Lemma \ref{otto}. The latter quotient is isomorphic to $H^p(\hat\cA_{F}^*(Y,E(\pi\otimes\vp)))$ via the map $b$ appearing in the proof
of Proposition \ref{HBC1}.
\hB

\begin{prop}\label{HBC2}
The inclusion $\hat\cA_{F,ch}^*(Y,E(\pi\otimes\vp))
\subset \hat\cA^*_F(Y,E(\pi\otimes\vp))$ is a quasi-isomorphism.
\end{prop}
\proof
We first prove injectivity. Let $\omega\in Z^p(\hat\cA_{F,ch}^*(Y,E(\pi\otimes\vp)))$ such that $[\omega]=0$ in
$H^p(\hat\cA_{F}^*(Y,E(\pi\otimes\vp)))$. By Corollary \ref{paulernst}
there exists an element $\psi \in Z^p(\hat\cA_{F}^*(Y,E(\pi\otimes\vp)))$
such that
$$ \omega=\hat\Delta \psi =d(\hat\delta\psi)\ .$$
Now $\hat\delta\psi\in \hat\cA_{F,ch}^{p-1}(Y,E(\pi\otimes\vp))$, hence
$[\omega]=0$ in $H^p(\hat\cA_{F,ch}^*(Y,E(\pi\otimes\vp)))$.

$Y$ is orientable. Therefore there is a Hodge-$*$-operator on $\Omega^*(Y)$. It induces an operator
$$ *:=*\otimes\id_{E(\pi\otimes\vp)} : \hat\cA_{F}^p(Y,E(\pi\otimes\vp)))\rightarrow \hat\cA_{F}^{n-p}(Y,E(\pi\otimes\vp)))$$
satisfying
$$*^2=(-1)^{p(n-1)}\id\ ,\quad \hat\delta=(-1)^{p(n-1)+1}*d*\ ,
\quad *\hat\Delta=\hat\Delta*\ .$$

Now we can prove surjectivity. Let $\omega\in Z^p(\hat\cA_{F}^p(Y,E(\pi\otimes\vp)))$. We want to find an element $\eta\in
\hat\cA_{F}^{p-1}(Y,E(\pi\otimes\vp)))$ such that $\hat\delta(\omega-d\eta)=0$.
The element $\phi=*\hat\delta\omega$ is exact in $\hat\cA_{F}^*(Y,E(\pi\otimes\vp))$. Thus by Corollary \ref{paulernst} we have $\phi=\hat\Delta \psi$ for some $\psi\in Z^{n+1-p}(\hat\cA_{F}^*(Y,E(\pi\otimes\vp)))$. Now $\eta=(-1)^{(p-1)(n-1)}*\psi$
does the job. This finishes the proof of the proposition.
\hB

Note that by the results of Section \ref{twist} the space of cocyles
$Z^p(\hat\cA_{F,ch}^*(Y,E(\pi\otimes\vp)))$ is isomorphic to
$Z^{n-p}_F(\vp)\cap C^{-\infty}(\partial X, V^+(\sigma^{n-p}_{F,\lambda_{n-p}},\vp))$. The surjectivity assertion could
also have been proved by showing that 
$Z^{n-p}_F(\vp)\cap C^{-\infty}(\partial X, V^+(\sigma^{n-p}_{F,\lambda_{n-p}},\vp))$ is dense in $Z^{n-p}_F(\vp)$.
This is indeed possible by a refinement of the proof of Corollary \ref{dense}.
 
Proposition \ref{HBC2}, while interesting in its own
right, has no direct impact for the proof of Theorem \ref{blub} (except for the case $\pi=1$).
However, we will use the technique employed in its proof in order to
show that the inclusion $\cA_{F,ch}^*(Y,E(\pi\otimes\vp))
\subset \cA^*_F(Y,E(\pi\otimes\vp))$ is a quasi-isomorphism. For this
we need the following analogue of Corollary \ref{paulernst}

\begin{lem}\label{BS&T}
$H^p(\cA_{F}^*(Y,E(\pi\otimes\vp)))\cong Z^p(\cA_{F}^*(Y,E(\pi\otimes\vp)))/\Delta_F\left( Z^p(\cA_{F}^*(Y,E(\pi\otimes\vp)))\right)\ .$
\end{lem}
\proof 
We employ the decomposition (\ref{reinhold}) which is stable under the action
of $\C[\Delta_F]$. For $p>0$ the space $Z^p(\hat\cA_{F}^*(Y,E(\pi\otimes\vp)))$ is isomorphic to 
${}^\Gamma C^{-\infty}(\partial X, V^+(\sigma^{n-p}_{F,\lambda_{n-p}},\vp))$. Using that $\Delta_F=-\Omega+\kappa_F$ we find as in the proof of Corollary \ref{preth} that this isomorphism
sends $\Delta_F\left( Z^p(\hat\cA_{F}^*(Y,E(\pi\otimes\vp)))\right)$ to
$\p_\Gamma\left({}^\Gamma 
C^{-\infty}(\partial X, V^+(\sigma^{n-p}_{F,\lambda_{n-p}},\vp))\right)$
which is isomorphic to \linebreak[4]
$\hat\Delta\left( Z^p(\hat\cA_{F}^*(Y,E(\pi\otimes\vp)))\right)$. Now Lemma \ref{xxx} and
Corollary \ref{paulernst} imply
$$H^p(\cA_{F}^*(Y,E(\pi\otimes\vp)))\cong Z^p(\hat\cA_{F}^*(Y,E(\pi\otimes\vp)))/\Delta_F\left( Z^p(\hat\cA_{F}^*(Y,E(\pi\otimes\vp)))\right)\ .$$
Thus it remains to show that
$$ \Delta_F: Z^p(\cA_{F}^*(Y,E(\pi\otimes\vp))_R)\rightarrow Z^p(\cA_{F}^*(Y,E(\pi\otimes\vp))_R) $$
is surjective. 

In the same way as Corollary \ref{preth} has been derived from Theorem \ref{krah} one
shows that
\begin{equation}\label{eif}  
Z^p(\cA_{F}^*(Y,E(\pi\otimes\vp))_R)\cong \bigoplus_{\lambda\in \aaaa^*_+\setminus\{\pm ((p-1)\alpha-\rho)\}} 
{}^\Gamma C^{-\infty}(\partial X, V^+(\sigma^{n-p}_\lambda,\pi\otimes\vp))^{\chi_F}\ .
\end{equation}
The sum on the right hand side is finite. It is therefore sufficient to show that
\begin{equation}\label{fel}
A:=-\Omega_F+\kappa_F: {}^\Gamma C^{-\infty}(\partial X, V^+(\sigma^{n-p}_\lambda,\pi\otimes\vp))\rightarrow {}^\Gamma C^{-\infty}(\partial X, V^+(\sigma^{n-p}_\lambda,\pi\otimes\vp)) 
\end{equation}
is surjective for all $\lambda$ appearing in (\ref{eif}). 

As in in the proof of Proposition \ref{trans} we see that the 
$G$-representations $C^{-\infty}(\partial X, V^+(\sigma^{n-p}_\mu,\pi))$ have composition
series with composition factors of the form $C^{-\infty}(\partial X, V^+(\sigma^\prime_{\mu^\prime(\mu)}))$. Therefore there exist pairwise different non-zero polynomials $p_i$ on 
$\aca$ and natural numbers $k_i$, $i=1,\dots,r$, such that
$\prod_{i=1}^{r} (A-p_i(\mu))^{k_i}$ acts by zero on $C^{-\infty}(\partial X, V(\sigma^{n-p}_\mu,\pi\otimes\vp))$. 
 
Then we can find non-zero polynomials $q_i\in \C[\aca][x]$, $b_i\in \C[\aca]$ such that 
$$xq_i(x,\mu)\prod_{j\ne i} (x-p_j(\mu))^{k_j}\equiv b_i(\mu)\quad \mod  \ (x-p_i(\mu))^{k_i}\ .$$ 
We now introduce
a rational family $Q_\mu$ on $\aca$ of operators in $\C[A]$ by
$$Q_\mu=\sum_{i=1}^{r} \frac{q_i(A,\mu)}{b_i(\mu)} \prod_{j\ne i} (A-p_j(\mu))^{k_j}\ .$$
We now fix $\lambda$ appearing in (\ref{eif}).
For $\mu$ in a pointed neighbourhood of $\lambda$ the family $Q_\mu$ is regular and  we have
\begin{equation}\label{spitz} 
AQ_\mu f_\mu=f_\mu\ ,\quad f_\mu\in C^{-\infty}(\partial X, V(\sigma^{n-p}_\mu,\pi\otimes\vp))\ .
\end{equation}

Recall the definition of $I_\sigma^{wr,-}$ from Section \ref{round}.
For $\sigma=\sigma^{n-p}$ one checks that
$$ I_\sigma^{wr,-}\cap \aaaa^*_+=\left\{\begin{array}{ccl}
\emptyset&& p=1,n\\
\{(p-1)\alpha-\rho\}&& \frac{n+1}{2}\le p \le n-1\\
\{\rho-(p-1)\alpha\}&& 2\le p\le \frac{n}{2}\quad .
\end{array}\right.
$$
Thus all the pairs $(\sigma^{n-p},\lambda)$ occurring in (\ref{eif}) are not special. Proposition \ref{schnauz} implies that
$${}^\Gamma C^{-\infty}(\partial X, V^+(\sigma^{n-p}_\lambda,\pi\otimes\vp))=
E^+_{-\infty}(\sigma^{n-p}_\lambda,\pi\otimes\vp)\ .$$
Let now $f\in {}^\Gamma C^{-\infty}(\partial X, V^+(\sigma^{n-p}_\lambda,\pi\otimes\vp))$. Arguing as in the proof of Corollary \ref{iso3} we find 
a family $f_\mu \in  {}^{\Gamma} \cM_\lambda C^{-\infty}(\partial X,V(\sigma_.,\pi\otimes\varphi))$ such that $f=ev(f_\mu)$.
Then 
$$ g:=ev(Q_\mu f_\mu)\in {}^\Gamma C^{-\infty}(\partial X, V^+(\sigma^{n-p}_\lambda,\pi\otimes\vp))$$
satisfies $Ag=f$ in view of (\ref{spitz}). This proves
the desired surjectivity.
\hB

In order to finish the proof of Theorem \ref{HBC} it remains to show

\begin{prop}\label{HBC3}
The inclusion $\cA_{F,ch}^*(Y,E(\pi\otimes\vp))
\subset \cA^*_F(Y,E(\pi\otimes\vp))$ is a quasi-isomorphism.
\end{prop}
\proof
We will argue similarly as in the proof of Proposition \ref{HBC2}. 
The proof of injectivity is straightforward. Indeed, let $\omega\in Z^p(\cA_{F,ch}^*(Y,E(\pi\otimes\vp)))$ such that $[\omega]=0$ in
$H^p(\cA_{F}^*(Y,E(\pi\otimes\vp)))$. By Lemma \ref{BS&T}
there exists an element $\psi \in Z^p(\cA_{F}^*(Y,E(\pi\otimes\vp)))$
such that
$$ \omega=\Delta_F \psi =d(\delta_F\psi)\ .$$
Now $\delta_F\psi\in \cA_{F,ch}^{p-1}(Y,E(\pi\otimes\vp))$, hence
$[\omega]=0$ in $H^p(\cA_{F,ch}^*(Y,E(\pi\otimes\vp)))$.

The proof of surjectivity needs a little preparation. Let $\gaaa=\kaaa\oplus\paaa$ be the Cartan decomposition. Let $(r, C^\infty(\Gamma\backslash G,\vp))$ be the $G$-representation given by
$$C^\infty(\Gamma\backslash G,\vp):=\{f:G\rightarrow V_\vp\:|\: f(gx)=\vp(g)f(x)\ \forall g\in\Gamma,\: x\in G\}\ ,\quad (r(g)f)(x):=f(xg)\ .$$ Then one can identify
$\Omega^*(Y,E(\pi\otimes\vp)))$ with $[C^\infty(\Gamma\backslash G,\vp)\otimes
F\otimes \Lambda^*\paaa^*]^K$. Any element $X\in\paaa$ defines operators
$\ve(X):\Lambda^p\paaa^*\rightarrow \Lambda^{p+1}\paaa^*$ and $i(X):\Lambda^p\paaa^*\rightarrow \Lambda^{p-1}\paaa^*$ given by
$\ve(X)(\omega)=X^\sharp\wedge\omega$, where $X^\sharp$ is the $1$-form
corresponding to $X$ via the invariant bilinear form on $\paaa$, i.e. the
Riemannian metric on $Y$, and $i(X)(\omega)(Y_1,\dots,Y_{p-1})=\omega(X,Y_1,\dots,Y_{p-1})$. 
Let $X_i$, $i=1,\dots,n$ be an orthogonal
basis of $\paaa$. Using the above identification we have
\begin{eqnarray*}
d&=& \sum_{i=1}^{n} r(X_i)\otimes\id\otimes\ve(X_i)
+\id\otimes\pi(X_i)\otimes\ve(X_i)\ ,\\
\delta_F&=& \sum_{i=1}^{n} -r(X_i)\otimes\id\otimes i(X_i)
+\id\otimes\pi(X_i)\otimes i(X_i)\ ,\\
\hat\delta&=&\sum_{i=1}^{n} -r(X_i)\otimes\id\otimes i(X_i)
-\id\otimes\pi(X_i)\otimes i(X_i)\ .
\end{eqnarray*}
Let $\theta: G\rightarrow G$ be the Cartan involution. We define a new
representation $\pi^\theta$ on $F$ by $\pi^\theta(g):=\pi(\theta(g))$.
$F$ equipped with the representation $\pi^\theta$ will sometimes be denoted
by $F^\theta$.
In particular, we have ${\pi^\theta}_{|K}=\pi_{|K}$, $\pi^\theta(X_i)=-\pi(X_i)$.
Thus we can identify the bundle $E(\pi^\theta\otimes\vp)$ with $E(\pi\otimes\vp)$. This identification induces on $E(\pi\otimes\vp)$
a second flat connection, and hence on $\Omega^*(Y,E(\pi\otimes\vp)))$ operators
$d_\theta$, $\delta_{F^\theta}$, and $\hat\delta_\theta$. The above formulas
show that $\delta_F=\hat\delta_\theta$ and $\delta_{F^\theta}=\hat\delta$.
In particular, $\delta_F=(-1)^{p(n-1)+1}*d_\theta *$.
Note that for the corresponding Laplacians we have $\Delta_F=\Delta_{F^\theta}$, whereas in general $\hat\Delta\ne \hat\Delta_\theta$.

Let $B^p_\theta(\cA_{F}^*(Y,E(\pi\otimes\vp)))$, $Z^p_\theta(\cA_{F}^*(Y,E(\pi\otimes\vp)))$ denote the spaces of $p$-coboundaries and $p$-cocycles
of the complex $(\cA_{F}^*(Y,E(\pi\otimes\vp)), d_\theta)$, respectively. 
We will need the following

\begin{lem}\label{notleast}
$B^p_\theta(\cA_{F}^*(Y,E(\pi\otimes\vp)))
=\Delta_{F^\theta}(Z^p_\theta(\cA_{F}^*(Y,E(\pi\otimes\vp))))\ .$
\end{lem}
\proof 
If $\pi^\theta\cong\pi$, i.e. $\chi_{F^\theta}=\chi_F$, then there
exists a bundle automorphism of $E(\pi\otimes\vp)$ which intertwines
$d_\theta$ with $d$, $\delta_{F^\theta}$ with $\delta_F$. Thus in this case
the lemma follows immediately from Lemma \ref{BS&T}. In order to deal
with the opposite case (which
can occur for odd $n$, only) we observe that similarly to (\ref{eif})
\begin{eqnarray*}
Z^p_\theta(\cA_{F}^*(Y,E(\pi\otimes\vp)))&\cong& 
{}^\Gamma C^{-\infty}(\partial X, V^+(\sigma^{n-p}_{(p-1)\alpha-\rho},\pi^\theta\otimes\vp))^{\chi_F}\\
&&\oplus \bigoplus_{\lambda\in \aaaa^*_+\setminus\{\pm((p-1)\alpha-\rho)\}} 
{}^\Gamma C^{-\infty}(\partial X, V^+(\sigma^{n-p}_\lambda,\pi^\theta\otimes\vp))^{\chi_F}
\end{eqnarray*}
for $p\ne\frac{n+1}{2}$. If $p\ne\frac{n+1}{2}$ one has to replace
$C^{-\infty}(\partial X, V^+(\sigma^{n-p}_{(p-1)\alpha-\rho},\pi^\theta\otimes\vp))^{\chi_{F}}$
by \linebreak[4] $C^{-\infty}(\partial X, V^+(\sigma^{\pm}_{0},\pi^\theta\otimes\vp))^{\chi_{F}}$.

We claim that if $\chi_{F^\theta}\ne \chi_F$, then $C^{-\infty}(\partial X, V^+(\sigma^{n-p}_{(p-1)\alpha-\rho},\pi^\theta\otimes\vp))^{\chi_{F}}=\{0\}$.
Indeed, using Proposition \ref{trans} we obtain
\begin{eqnarray*}
C^{-\infty}(\partial X, V^+(\sigma^{n-p}_{(p-1)\alpha-\rho},\pi^\theta\otimes\vp))^{\chi_{F}}   
&\subset& C^{-\infty}(\partial X, V^+(\sigma^{n-p}_{(p-1)\alpha-\rho},\pi^\theta\otimes\vp))^{\kappa_{F}}\\
&=& C^{-\infty}(\partial X, V^+(\sigma^{n-p}_{(p-1)\alpha-\rho},\pi^\theta\otimes\vp))^{\kappa_{F^\theta}}\\
&=& C^{-\infty}(\partial X, V^+(\sigma^{n-p}_{(p-1)\alpha-\rho},\pi^\theta\otimes\vp))^{\chi_{F^\theta}}\ .
\end{eqnarray*}  
Therefore the surjectivity of the operator $A$ in (\ref{fel}) implies
the lemma in this case, too.
\hB  
  
Now we can prove surjectivity. Let $\omega\in Z^p(\cA_{F}^*(Y,E(\pi\otimes\vp)))$. We want to find an element $\eta\in
\cA_{F}^{p-1}(Y,E(\pi\otimes\vp)))$ such that $\delta_F(\omega-d\eta)=0$.
We have $*\delta_F\omega\in B^{n+1-p}_\theta(\cA_{F}^*(Y,E(\pi\otimes\vp)))$. Thus by Lemma \ref{notleast} there exists $\psi\in Z_\theta^{n+1-p}(\cA_{F}^*(Y,E(\pi\otimes\vp)))$ such that 
$*\delta_F\omega=\Delta_{F^\theta} \psi=d_\theta\hat\delta\psi$. Now $\eta=(-1)^{(p-1)(n-1)}*\psi$
does the job. This finishes the proof of the proposition.
\hB

\newpage
\section{Hyperfunctions on the limit set as coefficients}\label{sjp}

We return to the assumptions of Sections \ref{hyper} and \ref{round}. 
That is, $\Gamma$ is a discrete (torsion-free) convex cocompact subgroup of a linear rank one Lie group
$G$. In the present section we want to study the cohomology of $\Gamma$
with coefficients in the infinite-dimensional $\Gamma$-representations
$C^{-\omega}(\Lambda, V^k(\sigma_\lambda,\vp))$ and $C^{-\omega}(\Lambda, V^+(\sigma_\lambda,\vp))$. Their invariants, i.e. the cohomology groups
$H^0(\Gamma,.)$, have been already investigated in Sections \ref{round}-\ref{twist}. In this section we want to understand the higher cohomology groups, too. 
At all places, where we use the extension operator
$ext_\lambda$ and the spaces $E^+(\sigma_\lambda,\vp)$ we have to assume
that $X\ne \OO H^2$ (or that $\Ree(\lambda)>\delta_\Gamma+\delta_\vp$).

This study is motivated by the theory of the
Selberg zeta function associated to an irreducible $M$-representation $\sigma$ and a finite-dimensional $\Gamma$-representation $\vp$. Note that
any element $\gamma\in\Gamma$ is conjugated in $G$ to an element $m_\gamma a_\gamma\in MA^+$. Then for $\mu\in\aca$, $\Ree(\mu)>\delta_\Gamma+\delta_\vp$,
the Selberg zeta function can be defined by the infinite product
$$ Z_{S,\sigma,\vp}(\mu):=\prod_{[\gamma]}\prod_{k=0}^{\infty}
\det\left(\id_{V_\sigma\otimes S^k(\bar \naaa)\otimes V_\vp}-a_\gamma^{-2\rho}
(\sigma_\mu\otimes S^k(\Ad_{|\bar\naaa}))(m_\gamma a_\gamma)\otimes \vp(\gamma)\right)\ .$$
Here the first product runs over all non-trivial primitive conjugacy
classes of $\Gamma$, and $S^k$ denotes the $k$-fold symmetric power.
For a discussion of the geometric meaning of this formula we refer to
\cite{fried86}, \cite{bunkeolbrich955}, \cite{juhl01}. It is a folklore
theorem that $Z_{S,\sigma,\vp}$ admits a meromorphic continuation
to all of $\aca$. For cocompact $\Gamma$ this has been
proved by dynamical methods in \cite{fried86} and at
various places by trace formula methods (see 
\cite{juhl01}, \cite{bunkeolbrich955} and the literature cited therein). In the latter approach one has to assume $\vp$ to be unitary. The dynamical
methods also work for the general convex cocompact case as has been 
explained in \cite{pattersonperry01}, Sect. 2, 
while the trace formula methods developed up to now 
(see \cite{pattersonperry01}, \cite{bunkeolbrich00}) imply only that
the logarithmic derivative of $Z_{S,\sigma,\vp}$ is meromorphic.
In \cite{patterson93} Patterson conjectured
a precise relationship between the divisor of $Z_{S,\sigma,\vp}$ (again for trivial $\vp$) and the cohomology groups 
$H^*(\Gamma,C^{-\infty}(\Lambda, V(\sigma_\lambda,\vp)))$. Guided by the
experience gained in the work \cite{bunkeolbrich950},
\cite{bunkeolbrich990} on this conjecture we now state a slightly modified conjecture. Recall the definition of the natural number $k_+(\sigma_\lambda,\vp)$ from Proposition \ref{bartII}.

\begin{con}\label{gopatt}\mbox{ }\vspace{-0.8cm}\\ 
\begin{enumerate}
\item[(i)] The cohomology groups $H^*(\Gamma,C^{-\omega}(\Lambda, V^k(\sigma_\lambda,\vp)))$ and $H^*(\Gamma,C^{-\omega}(\Lambda, V^+(\sigma_\lambda,\vp)))$ are finite-dimensional.
\item[(ii)] The Euler characteristic satisfies 
$\chi(\Gamma,C^{-\omega}(\Lambda, V^k(\sigma_\lambda,\vp)))=0$ for all $k\in\nat$.
\item[(iii)] If $k\ge k_+(\sigma_\lambda,\vp)$, then $\dim H^*(\Gamma,C^{-\omega}(\Lambda, V^k(\sigma_\lambda,\vp)))=\dim H^*(\Gamma,C^{-\omega}(\Lambda, V^{k_+}(\sigma_\lambda,\vp)))$.
\item[(iv)] If $k\ge k_+(\sigma_\lambda,\vp)$, then 
$$ \chi(\Gamma, C^{-\omega}(\Lambda, V^+(\sigma_\lambda,\vp)))=
-\chi_1(\Gamma, C^{-\omega}(\Lambda, V^k(\sigma_\lambda,\vp)))\ ,$$
where for any $\Gamma$-module $V$ with $\dim H^*(\Gamma,V)<\infty$ the first derived Euler characteristic
$\chi_1(\Gamma,V)$ is defined by
$$\chi_1(\Gamma,V):=\sum_{p=1}^{n}(-1)^p p\dim H^p(\Gamma,V)\ .$$
\item[(v)] The order of the singularity of $Z_{S,\sigma,\vp}$ at $\mu=\lambda$
is given by
$$\ord_\lambda (Z_{S,\sigma,\vp})=\chi(\Gamma,C^{-\omega}(\Lambda, V^+(\sigma_\lambda,\vp)))\ .$$
\end{enumerate}
\end{con}

For a further discussion of this conjecture we refer to the introduction
of \cite{bunkeolbrich990}. There the role of the $\Gamma$-module $C^{-\omega}(\Lambda, V^+(\sigma_\lambda,\vp))$ is played by the 
$\Gamma$-module $\cO_\lambda C^{-\omega}(\Lambda, V(\sigma_.,\vp))$.
It will turn out that the module $C^{-\omega}(\Lambda, V^+(\sigma_\lambda,\vp))$
is much easier to deal with, but will have the same cohomology as $\cO_\lambda C^{-\omega}(\Lambda, V(\sigma_.,\vp))$ shifted by degree one (at least in the cases, where the latter has been computed). The conjecture (in slightly different
but equivalent formulation) has been proved in the following two cases:
\begin{itemize}
\item $\Gamma$ cocompact, $\vp=1$ (\cite{bunkeolbrich950}).
\item $X=\R H^n$, $\sigma=\vp=1$ (\cite{bunkeolbrich990}).
\end{itemize}
The interested reader will find a discussion of these topics and related
results and conjectures in the recent monograph \cite{juhl01}.

In this section we will establish Assertions (i)-(iv) of Conjecture \ref{gopatt}
for pairs $(\sigma,\lambda)$ which are not very special (see Definition \ref{ere}) and for all $(\sigma,\lambda)$ in case $X=\R H^n$. This will be done
by an explicit calculation of the cohomology groups in question in terms
of the finite-dimensional spaces ${}^\Gamma C^{-\infty}(\Lambda, V^+(\sigma_\lambda,\vp))$, $E^+_\Lambda(\sigma_\lambda,\vp)$ and 
$H^p(\Gamma,F\otimes V_\vp)$, where $F$ is an irreducible finite-dimensional
representation of $G$. These results give additional support to Assertion (v) of Conjecture \ref{gopatt} since the trace formula approach (if it succeeds)
usually provides a description of the divisor of $Z_{S,\sigma,\vp}$ in terms of scattering data like
$E^+_\Lambda(\sigma_\lambda,\vp)$ and topological data like $H^p(\Gamma,F\otimes V_\vp)$. Indeed, the results obtained here are a direct generalization
of the corresponding results for $X=\R H^n$, $\sigma=1$, obtained in \cite{bunkeolbrich990}, which were comparable to the trace formula
results of \cite{pattersonperry01}. However, the trace formula for the
general situation (see \cite{bunkeolbrich00} for trivial $\vp$) is up to now
not in a sufficiently explicit shape in order to perform the analogous comparison, except for $\lambda\not\in I_\aaaa\cup
[-\delta_\Gamma,\delta_\Gamma]$, where we have
$$
\ord_\lambda(Z_{S,\sigma,1})=\dim {}^\Gamma C^{-\infty}(\Lambda, V^+(\sigma_\lambda))
=\chi(\Gamma,C^{-\omega}(\Lambda, V^+(\sigma_\lambda)))
$$
(see Remark 3 at the end of \cite{bunkeolbrich00} and Theorem \ref{impo}
below).

We will first investigate the cohomology of the modules $C^{-\omega}(\partial X, V^k(\sigma_\lambda,\vp)))$ and \linebreak[4]
$C^{-\omega}(\partial X, V^+(\sigma_\lambda,\vp)))$. Therefore we have to look for suitable acyclic resolutions of these modules. This can be done for arbitrary
discrete subgroups $\Gamma\subset G$.
Recall that a $\Gamma$-module $V$ is called acyclic, if $H^p(\Gamma,V)=\{0\}$
for all $p\ge 1$.

Let $(\gamma,V_\gamma)$ be a finite-dimensional representation of $K$.
By $V(\gamma):= G\times_K V_\gamma$ we denote the associated homogeneous
vector bundle on $X=G/K$. We equip $V(\gamma,\vp):=V(\gamma)\otimes V_\vp$
with the tensor product $\Gamma$-action. The $G$-action on the first factor
of $V(\gamma,\vp)$ induces a $\Gamma$-equivariant action of $\cZ(\gaaa)$ on
$C^\infty(X,V(\gamma,\vp))$. In particular, for $\kappa\in\C$ and $k\in\nat$
the kernel $C^\infty(X,V(\gamma,\vp))_{A_\kappa^k}$ of the operator $A_\kappa^k:=(\Omega-\kappa)^k\in\cZ(\gaaa)$ becomes a $\Gamma$-module. Because of our rank one assumption $\cZ(\gaaa)$
acts locally finitely on $C^\infty(X,V(\gamma,\vp))_{A_\kappa^k}$, and we
obtain a finite direct sum decomposition 
\begin{equation}\label{dali} C^\infty(X,V(\gamma,\vp))_{A_\kappa^k}=\bigoplus_{\{\chi\:|\:\chi(\Omega)=\kappa\}}
C^\infty(X,V(\gamma,\vp))_{A_\kappa^k}^\chi
\end{equation}
into $\cZ(\gaaa)$-modules with generalized infinitesimal character.
This decomposition respects the $\Gamma$-action. We also consider the $\Gamma$-modules
\begin{eqnarray*}
C^\infty(X,V(\gamma,\vp))_{(A_\kappa)}&:=&\bigcup_{k\in\nat} C^\infty(X,V(\gamma,\vp))_{A_\kappa^k}\ ,\\ 
C^\infty(X,V(\gamma,\vp))^\chi&=&\bigcup_{k\in\nat} C^\infty(X,V(\gamma,\vp))_{A_\kappa^k}^\chi\ .
\end{eqnarray*}
We form the bundle $V_Y(\gamma,\vp):=\Gamma\backslash V(\gamma,\vp)$ over $Y=\Gamma\backslash X$. Then we have the corresponding
spaces of sections $C^\infty(Y,V_Y(\gamma,\vp))$, $C^\infty(Y,V_Y(\gamma,\vp))_{A_\kappa^k}$ etc.
 
\begin{lem}\label{ac}
Let $\Gamma\subset G$ be a torsion-free discrete subgroup.
Then
\begin{enumerate}
\item
If $\Gamma$ is not cocompact, then the $\Gamma$-modules $C^\infty(X,V(\gamma,\vp))_{A_\kappa^k}$ and $C^\infty(X,V(\gamma,\vp))_{A_\kappa^k}^\chi$ are acyclic.
\item If $\Gamma$ is cocompact, then
\begin{eqnarray*} 
\dim H^0\left(\Gamma,C^\infty(X,V(\gamma,\vp))_{A_\kappa^k}\right)&=& \dim H^1\left(\Gamma,C^\infty(X,V(\gamma,\vp))_{A_\kappa^k}\right)\\
&=&\dim C^\infty(Y,V_Y(\gamma,\vp))_{A_\kappa^k}<\infty\ ,\\
\dim H^0\left(\Gamma,C^\infty(X,V(\gamma,\vp))_{A_\kappa^k}^\chi\right)&=& \dim H^1\left(\Gamma,C^\infty(X,V(\gamma,\vp))_{A_\kappa^k}^\chi\right)\\
&=&\dim C^\infty(Y,V_Y(\gamma,\vp))_{A_\kappa^k}^\chi<\infty\ ,
\end{eqnarray*}
and
$$H^p\left(\Gamma,C^\infty(X,V(\gamma,\vp))_{A_\kappa^k}\right)
=H^p\left(\Gamma,C^\infty(X,V(\gamma,\vp))_{A_\kappa^k}^\chi\right)=\{0\}  \quad \mbox{for all }
p\ge 2\ .$$
\item If $\Gamma$ is cocompact, then there exists a minimal number $k_0\in\nat_0$ such that
$$ H^0\left(\Gamma,C^\infty(X,V(\gamma,\vp))_{(A_\kappa)}\right)=C^\infty(Y,V_Y(\gamma,\vp))_{(A_\kappa)}=C^\infty(Y,V_Y(\gamma,\vp))_{A_\kappa^{k_0}}\ .$$
If, in addition, $\vp$ is unitary, then $k_0\le 1$. 
\item In any case, the modules $C^\infty(X,V(\gamma,\vp))_{(A_\kappa)}$ and $C^\infty(X,V(\gamma,\vp))^\chi$ are acyclic.
\end{enumerate}
\end{lem}
\proof 
Let $\C[A_\kappa]$ be the ring of all polynomials in $A_\kappa$.
We consider functors $\Fin_k$, $k\in\nat$, and $\Fin$ on the category
of $\C[A_\kappa]$-modules defined by
$$
\Fin_k(V):=\ker_V A^k_\kappa\ ,\qquad \Fin(V):=\bigcup_{k\in\nat} \Fin_k(V)\ .
$$
They are left exact. Let us denote their higher derived functors by $\Fin_k^q$
and $\Fin^q$. Then (see \cite{bunkeolbrich980}, p. 52) $\Fin^q_k=\Fin^q=0$
for all $k$ and $q\ge 2$, and
\begin{equation}\label{finfin} 
\Fin_k^1(V):=\coker_V A^k_\kappa\ ,\qquad \Fin^1(V):=\lim_{\stackrel{\longrightarrow}{k}} \coker_V A^k_\kappa\ ,
\end{equation}
where the limit has to be taken with respect to the map 
$A_\kappa: \coker_V A^k_\kappa\rightarrow \coker_V A^{k+1}_\kappa$.

Assume we are given a $(\C[\Gamma], \C[A_\kappa])$-bimodule
which is acyclic as a $\Gamma$-module and satisfies $\Fin^1_1(V)=0$. Then $\Fin_k^1(V)=\Fin(V)=0$ for all $k$.
It follows that
\begin{equation}\label{holger} 
H^p(\Gamma, \Fin_k(V))=\Fin_k^p({}^\Gamma V)\ ,\quad  H^p(\Gamma, \Fin(V))=\Fin^p({}^\Gamma V)\ .
\end{equation}

We want to apply (\ref{holger}) to the module $V=C^\infty(X,V(\gamma,\vp))$.
In \cite{bunkeolbrich950}, Lemma 2.4., it was shown that $C^\infty(X,V(\gamma))$ is $\Gamma$-acyclic. The argument given
their also works with the additional twist $V_\vp$. Hence $C^\infty(X,V(\gamma,\vp))$ is acyclic. Moreover, the operator
$$A_\kappa: C^\infty(X,V(\gamma,\vp))\rightarrow C^\infty(X,V(\gamma,\vp))$$
is surjective by Theorem \ref{masur}. In other words, $\Fin^1_1(C^\infty(X,V(\gamma,\vp)))=0$.
Now (\ref{holger}) and (\ref{finfin}) yield
\begin{eqnarray*}
H^p\left(\Gamma,C^\infty(X,V(\gamma,\vp))_{A_\kappa^k}\right)
&=&H^p\left(\Gamma,C^\infty(X,V(\gamma,\vp))_{(A_\kappa^k)}\right)=\{0\}  \quad \mbox{for all }
p\ge 2\ ,\\
H^1\left(\Gamma,C^\infty(X,V(\gamma,\vp))_{A_\kappa^k}\right)
&=&\coker\left(A_\kappa^k: C^\infty(Y,V_Y(\gamma,\vp))\rightarrow
C^\infty(Y,V_Y(\gamma,\vp))\right)\ ,\\ 
H^1\left(\Gamma,C^\infty(X,V(\gamma,\vp))_{(A_\kappa)}\right)
&=&\lim_{\stackrel{\longrightarrow}{k}} \coker\left(A_\kappa^k: C^\infty(Y,V_Y(\gamma,\vp))\rightarrow
C^\infty(Y,V_Y(\gamma,\vp))\right)\ .
\end{eqnarray*} 
In particular, if $\Gamma$ is not cocompact, then $$H^1\left(\Gamma,C^\infty(X,V(\gamma,\vp))_{A_\kappa^k}\right)=
H^1\left(\Gamma,C^\infty(X,V(\gamma,\vp))_{(A_\kappa)}\right)=\{0\}$$ by Theorem \ref{masur}. 

We claim that for cocompact $\Gamma$ there exists $k_0\in\nat_0$
such that
$$C^\infty(Y,V_Y(\gamma,\vp))_{(A_\kappa)}=C^\infty(Y,V_Y(\gamma,\vp))_{A_\kappa^{k_0}}$$
is finite-dimensional and 
$$ C^\infty(X,V(\gamma,\vp))=C^\infty(Y,V_Y(\gamma,\vp))_{A_\kappa^{k_0}}
\oplus A_\kappa^{k_0}\left(C^\infty(Y,V_Y(\gamma,\vp))\right)\ .
$$
Indeed, if $\vp$ is unitary, then this assertion with $k_0=1$ is a consequence
of the spectral theory of elliptic selfadjoint operators on a
compact manifold.
If $\vp$ fails to be unitary, we choose a Hermitian metric
on $V_Y(\vp)$ which gives rise to an $L^2$ scalar product on $C^\infty(X,V(\gamma,\vp))$ and
a corresponding Bochner-Laplace operator $\nabla^*\nabla$ on $L^2(X,V(\gamma,\vp))$. It has the same principal symbol as $-A_\kappa$.
Now we can argue as in the proof of Proposition \ref{comp}.

The claim implies in particular that 
$\D H^1\left(\Gamma,C^\infty(X,V(\gamma,\vp))_{(A_\kappa)}\right)
=\lim_{\stackrel{\longrightarrow}{k}} \coker A_\kappa^k=0$. 
Thus all assertions of the lemma concerning
$H^p\left(\Gamma,C^\infty(X,V(\gamma,\vp))_{A_\kappa^k}\right)$
and $H^p\left(\Gamma,C^\infty(X,V(\gamma,\vp))_{(A_\kappa)}\right)$
are now proved.

The assertions for the components with generalized infinitesimal
character $\chi$ are now immediate consequences of the direct sum decomposition 
(\ref{dali}).
\hB

Let $\gamma$ be a finite-dimensional representation of $K$, and let $T\in\Hom_M(V_\sigma,V_\gamma)$. Then Formula (\ref{akt}) 
(with $\gamma^p$ replaced by $\gamma$) defines a Poisson transform
$$ P^T_{\sigma_\lambda}: C^{-\omega}(\partial X, V(\sigma_\lambda))\rightarrow
C^\infty(X,V(\gamma))\ .$$
By $G$-equivariance it has values in $C^\infty(X,V(\gamma))^{\chi_{\sigma,\lambda}}_{A_{\sigma,\lambda}}$, where
$$ A_{\sigma,\lambda}=\Omega-\chi_{\sigma,\lambda}(\Omega)\ .$$
Recall that $\chi_{\sigma,\lambda}(\Omega)=\langle \lambda,\lambda\rangle -\|\rho_\gaaa\|^2+\|\mu+\rho_\maaa\|^2$.

We call an element $T\in\Hom_M(V_\sigma,V_\gamma)$ generating if the elements
$v_T\in C^{\omega}(\partial X, V(\tilde\sigma_{-\lambda}))$, $v\in V_{\tilde\gamma}$, defined by
$v_T(k)=T^*(\tilde\gamma(k^{-1})v)$, generate the $G$-module $C^{-\omega}(\partial X, V(\tilde\sigma_{-\lambda}))$. It follows from the
interpretation of the Poisson transform as a linear combination of matrix coefficients
$$\langle P^T_{\sigma_\lambda}f(g), v\rangle= c_{f,v_T}(g)$$ 
that $P^T_{\sigma_\lambda}$ is injective if and only if $T$ is generating.
For any pair $(\sigma,\lambda)$ one can always find a (not necessarily irreducible) $\gamma$ such that $\Hom_M(V_\sigma,V_\gamma)$ contains
generating elements. 

\begin{prop}\label{chicago}
Assume that $\lambda\not\in I^{wr,-}_\sigma$. In case $p_\sigma(0)\ne 0$ we require
in addition that $\lambda\ne 0$. Then for any generating $T\in\Hom_M(V_\sigma,V_\gamma)$ the Poisson transform $P^T_{\sigma_\lambda}$
identifies the $G$-module $C^{-\omega}(\partial X, V(\sigma_\lambda))$
with a direct summand of $C^\infty(X,V(\gamma))^{\chi_{\sigma,\lambda}}_{A_{\sigma,\lambda}}$.
Moreover, if $\lambda\not\in I^{wr}_\sigma$, then any non-zero $T\in\Hom_M(V_\sigma,V_\gamma)$ is generating.
\end{prop}
\proof The proposition is a consequence of the theory of asymptotic expansions
of matrix coefficients of Harish-Chandra modules (see e.g. \cite{wallach88}, Chapter 4; \cite{knapp86}, Chapter VIII; compare also \cite{olbrichdiss}). We denote by $Y^k_{\sigma,\lambda}$ and $Z$
the underlying $(\gaaa,K)$-modules of $C^{-\omega}(\partial X, V^k(\sigma_\lambda))$ and $C^\infty(X,V(\gamma))^{\chi_{\sigma,\lambda}}_{A_{\sigma,\lambda}}$, respectively. They
are admissible and finitely generated (such modules are usually called Harish-Chandra modules). Then $C^{-\omega}(\partial X, V^k(\sigma_\lambda))$ and $C^\infty(X,V(\gamma))^{\chi_{\sigma,\lambda}}_{A_{\sigma,\lambda}}$ constitute
the maximal globalizations in the sense of Schmid \cite{schmid85}, \cite{kashiwaraschmid94} of $Y^k_{\sigma,\lambda}$ and $Z$, respectively. By the main theorem of
\cite{schmid85} or \cite{kashiwaraschmid94} it is therefore enough to prove
that $P^T_{\sigma_\lambda}$
identifies the $(\gaaa,K)$-module $Y_{\sigma,\lambda}$ with a direct summand of $Z$.

There is a countable set $E\subset \aca$ such that any $f\in Z$ has
an asymptotic expansion for $a\to \infty$
\begin{equation}\label{ae} 
f(a)\sim \sum_{\mu\in E} p_{\mu,f}(\log a)\: a^\mu\ ,
\end{equation}
where $p_{\mu,f}$ are polynomials on $\aaaa$ with values in $V_\gamma$. We call $\mu\in E$ a leading
exponent of $Z$, if $p_{\mu+m\alpha,f}=0$ for all $m\in\nat$ and $f\in Z$.
Let $E^l\subset E$ be the set of leading exponents of $Z$. If $\sigma^\prime \in \hat M$, $S\in \Hom_M(V_\gamma,V_{\sigma^\prime})$ and $\mu\in E^l$, then there is a well-defined $MAN$-equivariant map
$$ \hat\beta^S_{\sigma^\prime,\mu+\rho}: Z/\naaa Z\rightarrow V_{\sigma^\prime_{\mu+\rho}}\otimes \Pi$$
given by $\hat\beta^S_{\sigma^\prime,\mu+\rho}([f]):=S\circ p_{\mu,f}$.
Using that $Z$ consists of eigensections of $\Omega$
it follows as in the proof of \cite{bunkeolbrich950}, Prop. 4.1, that
for any leading exponent $\mu\ne -\rho$ of $Z$ the polynomial $p_{\mu,f}$
is in fact a constant. In addition, if $-\rho$ is a leading exponent, then
$\deg(p_{-\rho,f})\le 1$. By Casselman's Frobenius reciprocity (\cite{wallach88}, 4.2.2) we obtain $(\gaaa,K)$-equivariant maps
$$ \beta^S_{\sigma^\prime,\nu}: Z \rightarrow Y_{\sigma^\prime,\nu}\ ,\qquad
\nu-\rho\in E^l\setminus\{-\rho\}, \ S\in \Hom_M(V_\gamma,V_{\sigma^\prime})\ .$$
If $-\rho\in E^l$, then we obtain maps
$\tilde\beta^S_{\sigma^\prime,0}: Z \rightarrow Y^2_{\sigma^\prime,0}$.
We then set
$$ \beta^S_{\sigma^\prime,0}:=\p\circ \tilde\beta^S_{\sigma^\prime,0}\ .$$
We call the maps $\beta^S_{\sigma^\prime,\nu}$ boundary value maps.  

For $\Ree(\lambda)>0$, $f\in Y_{\sigma,\lambda}$, we have the following limit formula (see \cite{vanderven94} or \cite{olbrichdiss}, also
\cite{wallach88}, Thm. 5.3.4)
\begin{equation}\label{warz} 
\lim_{a\to\infty} a^{\rho-\lambda}P^T_{\sigma_\lambda}f(a)=c_\gamma(\lambda)Tf(1)\ .
\end{equation}
For a discussion of the meromorphic function $c_\gamma:\aca\rightarrow \End_M(V_\gamma)$ and its relation to the Knapp-Stein intertwining operators 
we refer to \cite{bunkeolbrich000}, Section 5, in particular Lemma 5.5.
It is regular for $\Ree(\lambda)\ge 0$, $\lambda\ne 0$, and in this region it has the property
that $c_\gamma(\lambda)T\ne 0$ for generating $T$.

If $\nu-\rho$ is a leading exponent of $Z$ such that the map $\beta^S_{\sigma^\prime,\nu}$ is non-zero, then $Y_{\sigma^\prime,\nu}$ has infinitesimal
character $\chi_{\sigma,\lambda}$.
Thus, if $\lambda\in I_\sigma^{wr}\setminus I_\sigma^{wr,-}$,
then $\lambda-\rho$ is a leading
exponent of $Z$. Hence $\beta^S_{\sigma,
\lambda}$ is defined on $Z$. Moreover, by Lemma \ref{dnach}, 3, we have
$\lambda>0$.
Choose $S_0\in
\Hom_M(V_\gamma, V_\sigma)$ such that $S_0\circ c_\gamma(\lambda) T=\id$. Then
$(\ref{warz})$ implies that $\beta^{S_0}_{\sigma,
\lambda}\circ P^T_{\sigma_\lambda}=\id$. It follows that
\begin{equation}\label{schwein} 
Z=\im P^T_{\sigma_\lambda}\oplus \ker \beta^{S_0}_{\sigma,
\lambda}\ ,
\end{equation}
which proves the proposition for $\lambda\in I_\sigma^{wr}\setminus I_\sigma^{wr,-}$.

We now assume $\lambda\not\in I_\sigma^{wr}$, $\Ree(\lambda)\ge 0$. Choose
$m\in\nat_0$ minimal such that $\lambda-\rho+m\alpha\in E^l$. We first
assume that $m>0$. By the above we obtain a bunch of boundary value maps
$$ \beta^S_{\sigma^\prime,\lambda+m\alpha}: Z\rightarrow Y_{\sigma^\prime,\lambda+m\alpha}\ .$$
Let us denote the irreducible $M$-representations $\sigma^\prime$ which appear in such non-zero boundary value maps by $\sigma^m_1,\dots,\sigma^m_k$. 
The maps $\beta^S_{\sigma_l^m,\lambda+m\alpha}$, $S\in\Hom_M(V_\gamma, V_{\sigma^m_l})$, fit together
to a map
$$ \beta_{\lambda+m\alpha}=(\beta_{\sigma^m_1,\lambda+m\alpha},\dots,
\beta_{\sigma^m_k,\lambda+m\alpha}): Z\rightarrow 
\bigoplus_{l=1}^k \Hom_M(V_{\sigma^m_l}, V_\gamma)\otimes Y_{\sigma^m_l,\lambda+m\alpha}
$$
determined by the condition $(S\otimes\id)\circ \beta_{\sigma^m_l,\lambda+m\alpha}=\beta^S_{\sigma^m_l,\lambda+m\alpha}$, $S\in \Hom_M( V_\gamma,V_{\sigma^m_l})$. Here 
we have identified $\End_M(V_{\sigma^m_l})$ with
$\C$ by Schur's Lemma.
We also have a Poisson transform
$$ P_{\lambda+m\alpha}: \bigoplus_{l=1}^k \Hom_M(V_{\sigma^m_l}, V_\gamma)\otimes Y_{\sigma^m_l,\lambda+m\alpha}\rightarrow Z $$
given by $P_{\lambda+m\alpha}(T\otimes f):=P^T_{\sigma^m_{l,\lambda+m\alpha}}f$,
$T\in \Hom_M(V_{\sigma^m_l}, V_\gamma)$, $f\in Y_{\sigma^m_l,\lambda+m\alpha}$. Let $c_\gamma(\sigma,\lambda)$ denote the operator on $\Hom_M(V_{\sigma},V_\gamma)$ given by $c_\gamma(\sigma,\lambda)(T):=c_\gamma(\lambda)\circ T$. Now (\ref{warz}) implies that
$$ \beta_{\lambda+m\alpha}\circ P_{\lambda+m\alpha}
=(c_\gamma(\sigma_1,\lambda+m\alpha)\otimes\id,\dots, c_\gamma(\sigma_k,\lambda+m\alpha)\otimes\id)\ .$$
Since the $(\gaaa,K)$-modules $Y_{\sigma^m_l,\lambda+m\alpha}$ have infinitesimal character $\chi_{\sigma,\lambda}$ they are irreducible by Lemma \ref{sauna}. This implies that $c_\gamma(\sigma^m_l,\lambda+m\alpha)$
is a non-zero multiple of the identity. Thus $\beta_{\lambda+m\alpha}\circ P_{\lambda+m\alpha}$ is bijective, and we obtain a direct
sum decomposition
$$ 
Z=\im P_{\lambda+m\alpha}\oplus \ker \beta_{\lambda+m\alpha}\ .
$$

Set $Z_1:=\ker \beta_{\lambda+m\alpha}$. Then $\lambda-\rho+(m-1)\alpha$ is a leading exponent of $Z_1$. Therefore we
can define a boundary value $\beta_{\lambda+(m-1)\alpha}$ on $Z_1$. Moreover, by
(\ref{warz}) or a corresponding formula for $\Ree(\lambda)=0$ (see \cite{bunkeolbrich000}, Lemma 6.2) we have $\im P^T_{\sigma_\lambda}\subset Z_1$. Arguing inductively
we find a direct sum decomposition
\begin{equation}\label{wolf} 
Z= W\oplus Z_m\ ,
\end{equation}
such that $\lambda-\rho$ is a leading exponent of $Z_m$ and $\im P^T_{\sigma_\lambda}\subset Z_m$. Again using (\ref{warz}) or the analogous
formulas for $\Ree(\lambda)=0$ we find as in the proof of (\ref{schwein})
an element $S_0 \in \Hom_M(V_\gamma,V_\sigma)$ such that 
\begin{equation}\label{wein} 
Z_m=\im P^T_{\sigma_\lambda}\oplus \ker \beta^{S_0}_{\sigma,
\lambda}\ .
\end{equation}
Note that this last argument does not work in the case $p_\sigma(0)\ne 0$, $\lambda=0$,
since in this case $\deg p_{-\rho,P^T_{\sigma_0}f}=0$, and hence $\beta^S_{\sigma,0}\circ P^T_{\sigma_0}=0$ for any choice
of $S$.
Now for $\Ree(\lambda)\ge 0$, $\lambda\not\in I^{wr}_\sigma$ the proposition follows from (\ref{wolf}) and
(\ref{wein}). That any non-zero $T\in\Hom_M(V_\sigma,V_\gamma)$ is generating 
is a consequence
of the irreducibility of $Y_{\sigma,\lambda}$ (see Lemma \ref{sauna} and \cite{knapp86}, Cor. 14.30).

It remains to discuss the case $\Ree(\lambda)<0$, $\lambda\not\in I^{wr}_\sigma$.
The irreducibility of $Y_{\sigma,\lambda}$ and the functional equation of the Poisson transform (\cite{bunkeolbrich000}, Eq. (18)) now
implies that $\im P^T_{\sigma_\lambda}=\im P^{T^-}_{\sigma^w_{-\lambda}}$,
where $\D T^-:=\lim_{\mu\to\lambda}p_\sigma(\mu)\gamma(w)c_\gamma(\mu)T\sigma(w^{-1})\in \Hom_M(V_{\sigma^w},V_\gamma)$. Moreover, $T^-\ne 0$ whenever $T\ne 0$. Note that $-\lambda\not\in I^{wr}_{\sigma^w}$. We have therefore reduced the assertion to the
case $\Ree(\lambda)\ge 0$. This finishes the proof of the proposition.     
\hB

Let $q:\Pi\rightarrow \C$, $q(h)=h(0)$, be the map taking
the constant term. For $T\in \Hom_M(V_\sigma,V_\gamma)$ we set
$T^+:=T\otimes q \in \Hom_M(V_\sigma\otimes\Pi, V_\gamma)$. Let $T^k$
be the restriction of $T^+$ to $V_\sigma\otimes\Pi^k$. As in Section \ref{round}
we obtain Poisson transforms
\begin{eqnarray*} 
P^{T^k}_{\sigma_\lambda}: C^{-\omega}(\partial X, V^k(\sigma_\lambda))&\rightarrow& C^\infty(X,V(\gamma))\ ,\\
P^{T^+}_{\sigma_\lambda}: C^{-\omega}(\partial X, V^+(\sigma_\lambda))&\rightarrow& C^\infty(X,V(\gamma))\ .
\end{eqnarray*}
Using (\ref{where}) one derives that
\begin{equation}\label{schnurz}
A_{\sigma,\lambda}\circ P^{T^+}_{\sigma_\lambda}=P^{T^+}_{\sigma_\lambda}
\circ \p(2\langle \lambda,\alpha \rangle + |\alpha|^2\p)\ .
\end{equation}
For $\lambda\ne 0$ this implies that $P^{T^+}_{\sigma_\lambda}$ is injective
whenever $P^{T}_{\sigma_\lambda}$ is injective. Moreover, we see that
$P^{T^k}_{\sigma_\lambda}$ and $P^{T^+}_{\sigma_\lambda}$ have values in
$C^\infty(X,V(\gamma))^{\chi_{\sigma,\lambda}}_{A^k_{\sigma,\lambda}}$ and
$C^\infty(X,V(\gamma))^{\chi_{\sigma,\lambda}}$, respectively.

\begin{prop}\label{newyork}
Assume that $\lambda\not\in I^{wr,-}_\sigma\cup\{0\}$. Then for any generating $T\in\Hom_M(V_\sigma,V_\gamma)$ the Poisson transforms $P^{T^k}_{\sigma_\lambda}$ and $P^{T^+}_{\sigma_\lambda}$
identify the $G$-modules $C^{-\omega}(\partial X, V^k(\sigma_\lambda))$
and \linebreak[4]
$C^{-\omega}(\partial X, V^+(\sigma_\lambda))$ with a direct summand of $C^\infty(X,V(\gamma))^{\chi_{\sigma,\lambda}}_{A^k_{\sigma,\lambda}}$ and $C^\infty(X,V(\gamma))^{\chi_{\sigma,\lambda}}$, respectively.
\end{prop}
\proof We will follow the lines of the proof of Proposition \ref{chicago}.
Again, it is sufficient to check the assertion for $\Ree(\lambda)\ge 0$ 
on the level of
the underlying $(\gaaa,K)$-modules. 
Let $Z^k$ be the underlying $(\gaaa,K)$-module of $C^\infty(X,V(\gamma))^{\chi_{\sigma,\lambda}}_{A^k_{\sigma,\lambda}}$.
Then any $f\in Z^k$ has an asymptotic expansion of the form (\ref{ae}), 
where the
set of (leading) exponents of $Z^k$ coincides with that of $Z$. For
any leading exponent $\mu\ne -\rho$ we have $\deg p_{\mu,f}\le k-1$.
Thus we obtain boundary value maps 
$$ \beta^S_{\sigma^\prime,\nu}: Z^k \rightarrow Y^k_{\sigma^\prime,\nu}\ ,\qquad
\nu-\rho\in E^l\setminus\{-\rho\}, \ S\in \Hom_M(V_\gamma,V_{\sigma^\prime})\ ,$$
and, if $m$ is as in the proof of  Proposition \ref{chicago},
$$ \beta_{\lambda+m\alpha}: Z^k\rightarrow 
\bigoplus_{l=1}^{r}\Hom_M(V_{\sigma_l^m}, V_\gamma)\otimes Y^k_{\sigma_l^m,\lambda+m\alpha}\ .
$$
We define
$$ P^k_{\lambda+m\alpha}: \bigoplus_{l=1}^{r}\Hom_M(V_{\sigma_l^m}, V_\gamma)\otimes Y^k_{\sigma_l^m,\lambda+m\alpha}\rightarrow Z^k $$
by $P^k_{\sigma^m_{l,\lambda+m\alpha}}(T\otimes f):=P^{T^k}_{\sigma^m_{l,\lambda+m\alpha}}f$.
Since $\beta_{\lambda+m\alpha}\circ P^k_{\lambda+m\alpha}$ commutes with the action of the Casimir operator and with the projections on the direct summands it commutes with $\id\otimes\p$ (here one uses $\lambda+m\alpha\ne 0$, compare (\ref{schnurz})).
Thus the bijectivity of $\beta_{\lambda+m\alpha}\circ P_{\lambda+m\alpha}$ on $\bigoplus_{l=1}^{r}\Hom_M(V_{\sigma_l^m}, V_\gamma)\otimes Y_{\sigma_l^m,\lambda+m\alpha}$ implies the bijectivity of $\beta_{\lambda+m\alpha}\circ P^k_{\lambda+m\alpha}$ on \linebreak[4]
$\bigoplus_{l=1}^{r}\Hom_M(V_{\sigma_l^m}, V_\gamma)\otimes Y^k_{\sigma_l^m,\lambda+m\alpha}$. It
follows that
$$ Z^k=\im P^k_{\lambda+m\alpha}\oplus Z^k_1\ ,$$
where $Z^k_1:=\ker\beta_{\lambda+m\alpha}$. Inductively we find a
decomposition
$$ Z^k=W^k\oplus Z_m^k $$
such that $\lambda-\rho$ is a leading exponent of $Z_m^k$ and $\im P^{T^k}_{\sigma_\lambda}\subset Z^m_k$. Let $S_0$ be as in the proof of
Proposition \ref{chicago}. Because of $\lambda\ne 0$ we see that
$\beta_{\sigma,\lambda}^{S_0}\circ P^{T^k}_{\sigma_\lambda}$ is bijective, too.
It follows that $Z^k=W^k\oplus\im P^{T^k}_{\sigma_\lambda}\oplus \ker\beta_{\sigma,\lambda}^{S_0}$.
\hB

It remains to discuss the case $\lambda=0$, $0\not\in I^{wr}_\sigma$.
Recall that it is impossible that $0\in I^{wr}_\sigma\setminus I^{wr,-}_\sigma$.
 
\begin{prop}\label{harlem}
Assume that $0\not\in I^{wr}_\sigma$ and that $p_\sigma(0)\ne 0$. Then for any generating $T\in\Hom_M(V_\sigma,V_\gamma)$ the Poisson transforms $P^{T^k}_{\sigma_0}$ and $P^{T^+}_{\sigma_0}$ are injective.
Moreover, $P^{T^{2k}}_{\sigma_0}$ and $P^{T^+}_{\sigma_0}$ identify the $G$-modules $C^{-\omega}(\partial X, V^{2k}(\sigma_0))$
and $C^{-\omega}(\partial X, V^+(\sigma_0))$ with a direct summand of $C^\infty(X,V(\gamma))^{\chi_{\sigma,0}}_{A^k_{\sigma,0}}$ and $C^\infty(X,V(\gamma))^{\chi_{\sigma,0}}$, respectively.
\end{prop}
\proof If $P^{T^k}_{\sigma_0}$ would have a nontrivial kernel then
it would intersect $C^{-\omega}(\partial X, V(\sigma_0))$ nontrivially
by Assertion (ii) appearing in the proof of Proposition \ref{new}. This
is impossible since $T$ is generating.
That $\im P^{T^{2k}}_{\sigma_0}\subset C^\infty(X,V(\gamma))^{\chi_{\sigma,\lambda}}_{A^k_{\sigma,\lambda}}$ is a
consequence of (\ref{schnurz}).
In order to prove that $\im P^{T^{2k}}_{\sigma_0}$ is a direct summand
of $C^\infty(X,V(\gamma))^{\chi_{\sigma,\lambda}}_{A^k_{\sigma,\lambda}}$
one proceeds as in the proof of Proposition \ref{newyork}. Arguing exactly
as there we arrive at a decomposition
$$ Z^k=W^k\oplus Z^k_m$$
such that $-\rho$ is a leading exponent of $Z_m^k$ and $\im P^{T^{2k}}_{\sigma_0}\subset Z^k_m$. For $f\in Z^k_m$ the polynomials 
$p_{-\rho,f}$ have degree at most $2k-1$. For $S\in\Hom_M(V_\gamma,V_\sigma)$
we obtain boundary value maps
$$ \beta^S_{\sigma,0}: Z^k_m \rightarrow Y^{2k}_{\sigma,0}\ .$$
We now consider the limit formula (\cite{bunkeolbrich000}, Eq. (36))
\begin{equation}\label{bue} 
\lim_{a\to\infty} a^{\rho}P^T_{\sigma_0}f(a)=c_\gamma(0)Tf(1)+
T^w(\hat J_{\sigma,0}f)(1)\ ,
\end{equation}
where $T^w:=\gamma(w)T\sigma(w^{-1})$.

If $\sigma\not\cong\sigma^w$, then $c_\gamma(0)T$ and $T^w$ have values
in the different $M$-isotypic components $V_\gamma(\sigma)$ and $V_\gamma(\sigma^w)$, respectively. If $\sigma\cong\sigma^w$, then
$Y$ splits into two eigenspaces of $\hat J_{\sigma,0}$. 
Therefore the relation between $c_\gamma$ and $\hat J_{\sigma,0}$ (see e.g. \cite{bunkeolbrich000}, Lemma 5.5) yields a decomposition $T=T_++T_-$ such
that $T^w=T_+^w+T_-^w$ and $c_\gamma(0)T=\frac{1}{\sqrt{p_\sigma(0)}}(T^w_+-T^w_-)$. Moreover, $T$ is generating if and only if both components $T_+$ and $T_-$ are non-zero.
It follows that $\im\, c_\gamma(0)T\cap \im T^w=\{0\}$.
Thus in any case we can find
$S_0\in\Hom_M(V_\gamma,V_\sigma)$ such that $S_0\circ c_\gamma(0)T=\id$ and
$S_0\circ T^w=0$. Then
$\beta^{S_0}_{\sigma,0}\circ P^T_{\sigma_0}=\id$. As at the beginning of the proof it follows
that $\beta^{S_0}_{\sigma,0}\circ P^{T^{2k}}_{\sigma_0}$ is injective, and hence bijective. Therefore,
$$ Z_m^k=\im P^{T^{2k}}_{\sigma_0}\oplus\ker \beta^{S_0}_{\sigma,0}\ ,$$
which implies the proposition.
\hB

We now assume $p_\sigma(0)=0$. In particular, $0\not\in I_\sigma^{wr}$ and
$\sigma\cong\sigma^w$. Choose an irreducible representation $\gamma$ of $K$ such
that $\dim\Hom_M(V_\sigma,V_\gamma)=1$. Such a representation always exists (see
e.g. \cite{knapp86}, Ch. XV). Then one can define a meromorphic function
$c:\aca\rightarrow \C$ by
$$ c(\lambda)T:=(c_\gamma(\lambda)T)^w\ ,\qquad T\in\Hom_M(V_\sigma,V_\gamma)\ .$$
The superscript $w$ has the same meaning as in (\ref{bue}). As in \cite{bunkeolbrich000},
Sect. 5, we introduce a normalized intertwining operator
$$ J_{\sigma,\lambda}:=\frac{1}{c(-\lambda)}\hat J_{\sigma,\lambda}\ .$$
It satisfies the functional equation
\begin{equation}\label{erna}  
J_{\sigma,-\lambda}\circ J_{\sigma,\lambda}=\id\ .
\end{equation}
Moreover, we have the functional equation of the Poisson transform
(compare \cite{bunkeolbrich000}, (18))
\begin{equation}\label{fut}
P^T_{\sigma_\lambda}\circ J_{\sigma,-\lambda}=P^T_{\sigma_{-\lambda}}\ .
\end{equation}
By Proposition \ref{chicago} the Poisson transform $P^T_{\sigma_\lambda}$ is
injective for $T\ne 0$. It follows that $J_{\sigma_0}=\id$.
The normalized intertwining operators induce an operator
$$ J_+: C^{-\omega}(\partial X,V^+(\sigma_0)) \rightarrow C^{-\omega}(\partial X ,V^+(\sigma_0))$$
given by 
$$J_+ev(f_\mu):=-ev({J}_{\sigma,-\mu}f_{-\mu})\ .$$
Here $f_\mu\in \cM_0 C^{-\omega}(\partial X,V(\sigma_.))$. Because of (\ref{erna}) we have $J_+^2=\id$.
We set 
$$ C_{\pm}^{-\omega}(\partial X ,V^+(\sigma_0)):=\{f\in C^{-\omega}(\partial X ,V^+(\sigma_0))\:|\: J_+ f=\pm f\}$$ 
and $C_{\pm}^{-\omega}(\partial X ,V^k(\sigma_0)):=C_{\pm}^{-\omega}(\partial X ,V^+(\sigma_0))\cap C^{-\omega}(\partial X ,V^k(\sigma_0))$.
Observe that
\linebreak[4] 
$\p\left(C_{\pm}^{-\omega}(\partial X ,V^+(\sigma_0))\right)=
C_{\mp}^{-\omega}(\partial X ,V^+(\sigma_0))$. We obtain the following decompositions of $G$- and 
$\C[\p^2]$-modules
\begin{eqnarray}
C^{-\omega}(\partial X ,V^k(\sigma_0))&=&C_+^{-\omega}(\partial X ,V^k(\sigma_0))\oplus C_-^{-\omega}(\partial X ,V^k(\sigma_0))\ ,\label{futII}\\
C^{-\omega}(\partial X ,V^+(\sigma_0))&=&C_+^{-\omega}(\partial X ,V^+(\sigma_0))\oplus C_-^{-\omega}(\partial X ,V^+(\sigma_0))\ .\nonumber
\end{eqnarray}
We have
$C_+^{-\omega}(\partial X ,V^{2k}(\sigma_0))=C_+^{-\omega}(\partial X ,V^{2k-1}(\sigma_0))$ and $C_+^{-\omega}(\partial X ,V^1(\sigma_0))=C^{-\omega}(\partial X ,V^1(\sigma_0))$.

\begin{prop}\label{brooklyn}
Assume that $p_\sigma(0)=0$ and that $\dim\Hom_M(V_\sigma,V_\gamma)=1$.
Then $\p$ induces an isomorphism
$$C_-^{-\omega}(\partial X ,V^+(\sigma_0))\cong C_{+}^{-\omega}(\partial X ,V^+(\sigma_0))\ .$$ 
If $0\ne T\in\Hom_M(V_\sigma,V_\gamma)$, then $\ker P^{T^+}_{\sigma_0}=C_-^{-\omega}(\partial X ,V^+(\sigma_0))$. 
Moreover, $P^{T^{2k-1}}_{\sigma_0}$ and $P^{T^+}_{\sigma_0}$ identify the $G$-modules $C_{+}^{-\omega}(\partial X, V^{2k-1}(\sigma_0))$
and $C_{+}^{-\omega}(\partial X, V^+(\sigma_0))$ with a direct summand of $C^\infty(X,V(\gamma))^{\chi_{\sigma,0}}_{A^k_{\sigma,0}}$ and $C^\infty(X,V(\gamma))^{\chi_{\sigma,0}}$, respectively.
\end{prop}
\proof
The first assertion is obvious. As in (\ref{gugel}) we obtain for 
$f\in \cM_0 C^{-\omega}(\partial X,V(\sigma_.))$ 
$$ P^{T^+}_{\sigma_0}\circ ev(f_\mu)
=\res_{z=0}
P^T_{\sigma,z\alpha}f_{z\alpha}\ .$$
Assume now that $ev(f_\mu)\in C_-^{-\omega}(\partial X ,V^+(\sigma_0))$.
Then $J_{\sigma,-z\alpha}f_{-z\alpha}=f_{z\alpha}\ \mod\: \cO_0 C^{-\omega}(\partial X,V(\sigma_.))$.
Equation (\ref{fut}) implies that
$$ P^T_{\sigma,z\alpha}f_{z\alpha}=P^T_{\sigma,-z\alpha}f_{-z\alpha}\ \mod\: 
\cO_0 C^{-\omega}(\partial X,V(\sigma_.))\ .$$
It follows that $P^{T^+}_{\sigma_0}\circ ev(f_\mu)=0$, hence
$\ker P^{T^+}_{\sigma_0}\supset C_-^{-\omega}(\partial X ,V^+(\sigma_0))$.

By (\ref{schnurz}) any nontrivial $G$-submodule of $C_{+}^{-\omega}(\partial X, V^+(\sigma_0))$ has nontrivial intersection with $C^{-\omega}(\partial X ,V^1(\sigma_0))$.
We know from Proposition \ref{chicago} that $\ker P^{T^+}_{\sigma_0}\cap
C^{-\omega}(\partial X ,V^1(\sigma_0))=\{0\}$. It follows
that $\ker P^{T^+}_{\sigma_0}=C_-^{-\omega}(\partial X ,V^+(\sigma_0))$.

In order to show that $P^{T^{2k-1}}_{\sigma_0}\left(C_{+}^{-\omega}(\partial X, V^{2k-1}(\sigma_0))\right)$ is a direct summand of $C^\infty(X,V(\gamma))^{\chi_{\sigma,0}}_{A^k_{\sigma,0}}$ one proceeds as in the proof of
Proposition \ref{harlem}. Instead of the boundary value map
$$ \beta^S_{\sigma,0}: Z^k_m \rightarrow Y^{2k}_{\sigma,0}$$
one uses the map
$$ p_+\circ\p\circ \beta^S_{\sigma,0}: Z^k_m \rightarrow Y^{2k-1}_{\sigma,+}\ ,$$
where $Y^{2k-1}_{\sigma,+}$ is the underlying $(\gaaa,K)$-module of 
$C_{+}^{-\omega}(\partial X, V^{2k-1}(\sigma_0))$ and $p_+: Y^{2k-1}_{\sigma,0}\rightarrow Y^{2k-1}_{\sigma,+}$ is the projection corresponding to (\ref{futII}).
\hB

We now combine the previous four propositions with Lemma \ref{ac}.

\begin{prop}\label{acb}
Let $\Gamma\subset G$ be a torsion-free discrete subgroup, and let
$(\vp, V_\vp)$ be a finite-dimensional representation of $\Gamma$.
Assume that $(\sigma,\lambda)$ is not very special, i.e., $\lambda\not\in I_\sigma^{wr,-}$. Then
\begin{enumerate}
\item The $\Gamma$-module $C^{-\omega}(\partial X, V^+(\sigma_\lambda,\vp))$
is acyclic.
\item
If $\Gamma$ is not cocompact, then for any $k\in\nat$ the $\Gamma$-module $C^{-\omega}(\partial X, V^k(\sigma_\lambda,\vp))$ is acyclic.
\item Let $\Gamma$ be cocompact. Then there exists a minimal number $k_+=k_+(\sigma,\lambda)\in\nat_0$ such that
$$ H^0\left(\Gamma,C^{-\omega}(\partial X, V^+(\sigma_\lambda,\vp))\right)=H^0\left(\Gamma,C^{-\omega}(\partial X, V^{k_+}(\sigma_\lambda,\vp))\right)\ .$$
The spaces $H^0\left(\Gamma,C^{-\omega}(\partial X, V^+(\sigma_\lambda,\vp))\right)$ and $H^0\left(\Gamma,C^{-\omega}(\partial X, V^{k}(\sigma_\lambda,\vp))\right)$ are finite-dimen\-sional.
If, in addition, $\vp$ is unitary, then $k_+(\sigma,\lambda)\le 1$ except for $p_\sigma(0)=0$, $\lambda=0$, in which case $k_+(\sigma,0)\in\{0,2\}$.
\item If $\Gamma$ is cocompact, then
\begin{eqnarray} 
\dim H^0\left(\Gamma,C^{-\omega}(\partial X, V^k(\sigma_\lambda,\vp))\right)&=& \dim H^1\left(\Gamma,C^{-\omega}(\partial X, V^k(\sigma_\lambda,\vp))\right)<\infty\ ,\nonumber\\
H^p\left(\Gamma,C^{-\omega}(\partial X, V^k(\sigma_\lambda,\vp))\right)
&=&\{0\}  \quad \mbox{for all }
p\ge 2\ .\label{past}
\end{eqnarray}
In particular, $\chi\left(\Gamma,C^{-\omega}(\partial X, V^k(\sigma_\lambda,\vp))\right)=0$.
\item If $\Gamma$ is cocompact, then for any $k\ge k_+(\sigma,\lambda)$  
$$ \chi\left(\Gamma,C^{-\omega}(\partial X, V^+(\sigma_\lambda,\vp))\right)
=-\chi_1\left(\Gamma,C^{-\omega}(\partial X, V^k(\sigma_\lambda,\vp))\right)
=\dim{}^\Gamma C^{-\infty}(\partial X, V^+(\sigma_\lambda,\vp))\ .$$
\end{enumerate}
Assertion 3 also holds without the assumption $\lambda\not\in I_\sigma^{wr,-}$.
\end{prop}
\proof For any pair $(\sigma,\lambda)$ and a generating $T\in\Hom_M(V_\sigma,V_\gamma)$ the Poisson transform $P^{T^k}_{\sigma_\lambda}:=P^{T^k}_{\sigma_\lambda}\otimes\id_{V_\vp}$ injects $C^{-\omega}(\partial X, V^k(\sigma_\lambda,\vp))$ 
($C_+^{-\omega}(\partial X, V^k(\sigma_0,\vp))$, if $\lambda=0$ and $p_\sigma(0)=0$) into $C^\infty(X,V(\gamma,\vp))_{A_{\sigma,\lambda}^k}$. Thus the first statement of Assertion 3 is a consequence of Lemma \ref{ac}, 3. Moreover, for $\lambda\ne 0$
we have $\im P^{T^{k+r}}_{\sigma_\lambda}\cap C^\infty(X,V(\gamma,\vp))_{A_{\sigma,\lambda}^k}=\im P^{T^k}_{\sigma_\lambda}$.
It follows, that $k_+(\sigma,\lambda)\le 1$ for unitary $\vp$ and $\lambda\ne 0$. That this
also holds for $p_\sigma(0)\ne 0$, $\lambda=0$, can be proved as in Proposition \ref{new}. The case $p_\sigma(0)=0$, $\lambda=0$ is covered by Proposition \ref{brooklyn}.

Assertions 1,2, and Equation (\ref{past}) are immediate consequences
of Lemma \ref{ac}, and Propositions \ref{chicago}, \ref{newyork}, \ref{harlem},
and \ref{brooklyn}, respectively.

The long exact cohomology sequence associated to
$$ 0\rightarrow C^{-\omega}(\partial X, V^k(\sigma_\lambda,\vp))\rightarrow
C^{-\omega}(\partial X, V^+(\sigma_\lambda,\vp))\stackrel{\p^k}{\longrightarrow}
C^{-\omega}(\partial X, V^+(\sigma_\lambda,\vp))\rightarrow 0 $$
degenerates to
\begin{eqnarray*}
0\rightarrow {}^\Gamma C^{-\omega}(\partial X, V^k(\sigma_\lambda,\vp))&\rightarrow&
{}^\Gamma C^{-\omega}(\partial X, V^+(\sigma_\lambda,\vp))\stackrel{\p_\Gamma^k}{\longrightarrow}
{}^\Gamma C^{-\omega}(\partial X, V^+(\sigma_\lambda,\vp))\\
&\rightarrow&
H^1\left(\Gamma,C^{-\omega}(\partial X, V^k(\sigma_\lambda,\vp))\right)\rightarrow 0\ .
\end{eqnarray*}
Since for cocompact $\Gamma$ the space ${}^\Gamma C^{-\omega}(\partial X, V^+(\sigma_\lambda,\vp))$ is finite-dimensional the equality
$$\dim H^0\left(\Gamma,C^{-\omega}(\partial X, V^k(\sigma_\lambda,\vp))\right)=
\dim H^1\left(\Gamma,C^{-\omega}(\partial X, V^k(\sigma_\lambda,\vp))\right)$$
as well as the remaining assertions of the proposition follow. For Assertion
5 we have also used Theorem \ref{bart}.
\hB

We now consider the decomposition $\partial X=\Lambda\cup\Omega$. Since the
sheaf of hyperfunction sections of a vector bundle is flabby we obtain
the following exact sequences of $\Gamma$-modules
\begin{eqnarray}
0\rightarrow C^{-\omega}(\Lambda, V^k(\sigma_\lambda,\vp))&\rightarrow&C^{-\omega}(\partial X, V^k(\sigma_\lambda,\vp))\stackrel{res_\Omega}{\longrightarrow}
C^{-\omega}(\Omega, V^k(\sigma_\lambda,\vp))\rightarrow 0\ ,\label{bolek}\\
0\rightarrow C^{-\omega}(\Lambda, V^+(\sigma_\lambda,\vp))&\rightarrow&C^{-\omega}(\partial X, V^+(\sigma_\lambda,\vp))\stackrel{res_\Omega}{\longrightarrow}
C^{-\omega}(\Omega, V^+(\sigma_\lambda,\vp))\rightarrow 0\ .\label{lolek}
\end{eqnarray}

\begin{lem}\label{acB}
The $\Gamma$-modules $C^{-\omega}(\Omega, V^k(\sigma_\lambda,\vp))$ and
$C^{-\omega}(\Omega, V^+(\sigma_\lambda,\vp))$ are acyclic.
\end{lem}
\proof 
Since the action of $\Gamma$ on $\Omega$ is properly discontinuous
and $C^{-\omega}(\Omega, V^k(\sigma_\lambda,\vp))$, \linebreak[4] $C^{-\omega}(\Omega, V^+(\sigma_\lambda,\vp))$ are the spaces of global sections of a flabby $\Gamma$-equivariant sheaves on $\Omega$ the lemma can be shown in the same 
way as \cite{bunkeolbrich950}, Lemma 2.6.
\hB

\begin{theorem}\label{impo}
Let $\Gamma\subset G$ be a torsion-free convex cocompact discrete subgroup
of $\Gamma$, and let
$(\vp, V_\vp)$ be a finite-dimensional representation of $\Gamma$.
Assume that $(\sigma,\lambda)$ is not very special, i.e., $\lambda\not\in I_\sigma^{wr,-}$. Then Assertions (i)-(iv) of Conjecture \ref{gopatt}
hold true. In addition, we have
\begin{enumerate}
\item The $\Gamma$-module $C^{-\omega}(\Lambda, V^+(\sigma_\lambda,\vp))$
is acyclic.
\item For all $k\in\nat$
\begin{eqnarray*} 
\dim H^0\left(\Gamma,C^{-\omega}(\Lambda, V^k(\sigma_\lambda,\vp))\right)&=& \dim H^1\left(\Gamma,C^{-\omega}(\Lambda, V^k(\sigma_\lambda,\vp))\right)\\
&=&\dim E^k_\Lambda(\sigma_\lambda,\vp)\ ,\qquad\mbox{if $\Gamma$ is not cocompact}\ ,\\
H^p\left(\Gamma,C^{-\omega}(\Lambda, V^k(\sigma_\lambda,\vp))\right)
&=&\{0\}  \quad \mbox{for all }
p\ge 2\ .
\end{eqnarray*}
\item For any $k\ge k_+(\sigma,\lambda)$  
$$ \chi\left(\Gamma,C^{-\omega}(\Lambda, V^+(\sigma_\lambda,\vp))\right)
=-\chi_1\left(\Gamma,C^{-\omega}(\Lambda, V^k(\sigma_\lambda,\vp))\right)
=\dim{}^\Gamma C^{-\infty}(\Lambda, V^+(\sigma_\lambda,\vp))\ .$$
\end{enumerate}
\end{theorem}
\proof
The case of cocompact $\Gamma$ is already covered by Proposition \ref{acb}.
We thus assume $\Omega\ne\emptyset$. By Proposition \ref{acb} and Lemma
\ref{acB} the sequence (\ref{lolek}) provides an acyclic
resolution of 
$C^{-\omega}(\Lambda, V^+(\sigma_\lambda,\vp))$. We thus
have to determine the cohomology of the complex
$$ 0\rightarrow {}^\Gamma C^{-\omega}(\partial X, V^+(\sigma_\lambda,\vp))\stackrel{res}{\longrightarrow}
C^{-\omega}(B, V_B^+(\sigma_\lambda,\vp))\rightarrow 0\ .$$
The surjectivity of $\p_B$ combined with (\ref{htak}) implies the
surjectivity of
$$ res: {}^\Gamma C^{-\omega}(\partial X, V^+(\sigma_\lambda,\vp))\rightarrow
C^{-\omega}(B, V_B^+(\sigma_\lambda,\vp))\ .$$
It follows that all higher cohomology groups of $C^{-\omega}(\Lambda, V^+(\sigma_\lambda,\vp))$ vanish. By Proposition \ref{zwirbel}
we have the following equality of finite-dimensional spaces
$$ \dim H^0\left(\Gamma,C^{-\omega}(\Lambda, V^+(\sigma_\lambda,\vp))\right)
=E^+_\Lambda(\sigma_\lambda,\vp)=E^{k_+}_\Lambda(\sigma_\lambda,\vp)\ .$$
The long exact cohomology sequence associated to 
$$ 0\rightarrow C^{-\omega}(\Lambda, V^k(\sigma_\lambda,\vp))\rightarrow
C^{-\omega}(\Lambda, V^+(\sigma_\lambda,\vp))\stackrel{\p^k}{\longrightarrow}
C^{-\omega}(\Lambda, V^+(\sigma_\lambda,\vp))\rightarrow 0 $$
now reads as
$$
0\rightarrow E^{k}_\Lambda(\sigma_\lambda,\vp)\rightarrow
E^{+}_\Lambda(\sigma_\lambda,\vp)\stackrel{\p_\Gamma^k}{\longrightarrow}
E^{+}_\Lambda(\sigma_\lambda,\vp)
\rightarrow
H^1\left(\Gamma,C^{-\omega}(\Lambda, V^k(\sigma_\lambda,\vp))\right)\rightarrow 0\ .
$$
From this it is easy to derive the remaining assertions of the theorem.
\hB

Theorem \ref{impo} says that interesting phenomena, i.e., higher cohomology
groups of \linebreak[4]
$C^{-\omega}(\Lambda, V^+(\sigma_\lambda,\vp))$, occur for $\lambda\in
I_\sigma^{wr,-}$, only. In the rest of this section we will determine these
cohomology groups in case $X=\R H^n$. This is based on the results of Sections
\ref{heart} and \ref{twist}.  

From now on let $G$ be one of the groups $Spin(1,n)$ or $SO(1,n)_0$. Recall that in this case $I_\sigma^{wr}=I_\sigma^r$.
We have seen in Section \ref{twist} that for
$\lambda\in I_\sigma^r$ there exist an irreducible $G$-representation $(\pi,F)$
having infinitesimal character $\chi_{\sigma,\lambda}$ and an element
$p\in\{0,1,\dots,n-1\}\cup\{+,-\}$ such that $\sigma_\lambda=\sigma^p_{F,\lambda_p}$.
Then we set
$$ (\pi_{\sigma_\lambda},F_{\sigma_\lambda}):=(\pi,F)\ ,\quad  l_{\sigma_\lambda}=n-p\ \mbox{ for } p\ne\pm\ \ \mbox{ and }\  l_{\sigma_\lambda}:=\frac{n+1}{2} \ \mbox{ if } p=\pm\ .$$
Using the formulas in Section \ref{twist} this definition can be made more
explicit. Let $\lambda=r\alpha$, and let $\sigma$ have highest weight 
$\mu_\sigma=(m_1,m_2,\dots,m_{[\frac{n-1}{2}]})$. Let $i$ be the smallest
number such that $|r|>m_i+\frac{n-1}{2}-i$. Then $$l_{\sigma_\lambda}=\left\{\begin{array}{ccc}
n+1-i&&\lambda\ge 0\\
i&& \lambda\le 0
\end{array}\right.\ .$$
Moreover, if $n$ is even, then $F_{\sigma_\lambda}$ has highest weight
$$ (m_1-1,\dots, m_{i-1}-1, |r|+i-\frac{n+1}{2}, m_i,\dots, m_{\frac{n-2}{2}})\ .$$
For odd $n$ this highest weight is given by
\begin{eqnarray*}
&&(m_1-1,\dots, m_{i-1}-1, |r|+i-\frac{n+1}{2}, m_i,\dots, m_{\frac{n-3}{2}},
\sign(r)m_{\frac{n-1}{2}})\ ,\quad i\le\frac{n-1}{2}\ ,\\
&&(m_1-1,\dots, m_{\frac{n-3}{2}}-1, |m_{\frac{n-1}{2}}|-1,\sign(m_{\frac{n-1}{2}})r)\ ,\quad i=\frac{n+1}{2}\ .
\end{eqnarray*}
By $\tilde Z^p_{F_{\sigma_\lambda}}$ we denote the spaces of cocycles of
(\ref{zhelo}) for trivial $\Gamma$ and $\vp$. If $l_{\sigma_\lambda}\ne \frac{n+1}{2}$, then $\tilde Z^{n-l_{\sigma_\lambda}}_{F_{\sigma_\lambda}}$ is a $G$-submodule of $C^{-\omega}(\partial X, V(\sigma_\lambda))$, whereas
for $l_{\sigma_\lambda}=\frac{n+1}{2}$ we have $\tilde Z^\frac{n-1}{2}_{F_{\sigma_\lambda}}\cong C^{-\omega}(\partial X, V(\sigma_\lambda))$. We let $G$ act on both factors the flat vector bundle $\tilde E(\pi_{\sigma_\lambda}):=X\times F_{\sigma_\lambda}$ over $X$.
This $G$-action induces a $\cZ(\gaaa)$-action on $\Omega^*(X,\tilde E(\pi_{\sigma_\lambda}))\cong\Omega^*(X)
\otimes F_{\sigma_\lambda}$.

\begin{prop}\label{dRh}
If $\lambda\in I_\sigma^r$, then the following sequences are exact:
\begin{eqnarray*}
0&\rightarrow& C^{-\omega}(\partial X, V^+(\sigma_\lambda))\rightarrow
\Omega^{l_{\sigma_\lambda}}(X,\tilde E(\pi_{\sigma_\lambda}))_{(\hat\Delta)}
^{\chi_{\sigma,\lambda}}
\stackrel{d}{\longrightarrow} \Omega^{l_{\sigma_\lambda}+1}(X,\tilde E(\pi_{\sigma_\lambda}))_{(\hat\Delta)}
^{\chi_{\sigma,\lambda}}
\\
&&\hspace{6cm}\stackrel{d}{\longrightarrow}\dots 
\stackrel{d}{\longrightarrow}   
\Omega^n(X,\tilde E(\pi_{\sigma_\lambda}))_{(\hat\Delta)}
^{\chi_{\sigma,\lambda}}\rightarrow 0\ ,
\end{eqnarray*}
\begin{eqnarray*}
0&\rightarrow& \tilde Z^{n-l_{\sigma_\lambda}}_{F_{\sigma_\lambda}}\rightarrow \Omega^{l_{\sigma_\lambda}}(X,\tilde E(\pi_{\sigma_\lambda}))_{(\hat\Delta)}
^{\chi_{\sigma,\lambda}}
\stackrel{{\scriptsize \left(\begin{array}{c}d\\ d*\end{array}\right)}} {\longrightarrow}
\begin{array}{c}\Omega^{l_{\sigma_\lambda}+1}(X,\tilde E(\pi_{\sigma_\lambda}))
_{(\hat\Delta)}^{\chi_{\sigma,\lambda}}\\ \oplus\\ \Omega^{n+1-l_{\sigma_\lambda}}(X,\tilde E(\pi_{\sigma_\lambda}))_{(\hat\Delta)}
^{\chi_{\sigma,\lambda}}\end{array}
\stackrel{ \left({\scriptsize\begin{array}{cc}d&0\\0&d\end{array}} \right)}{\longrightarrow}\\
&\dots&\stackrel{ \left({\scriptsize\begin{array}{cc}
d&0\\ 0&d\end{array}}\right)}{\longrightarrow}
\begin{array}{c}\Omega^{n}(X,\tilde E(\pi_{\sigma_\lambda}))_{(\hat\Delta)}
^{\chi_{\sigma,\lambda}}\\ \oplus\\ \Omega^{2n-2l_{\sigma_\lambda}}(X,\tilde E(\pi_{\sigma_\lambda}))_{(\hat\Delta)}
^{\chi_{\sigma,\lambda}}\end{array}
\stackrel{ \left({\scriptsize\begin{array}{cc}
0&d\end{array}}\right)}{\longrightarrow}
\Omega^{2n-2l_{\sigma_\lambda}+1}(X,\tilde E(\pi_{\sigma_\lambda}))_{(\hat\Delta)}
^{\chi_{\sigma,\lambda}}\stackrel{d}{\longrightarrow}\\
&\dots&\stackrel{d}{\longrightarrow}
\Omega^{n}(X,\tilde E(\pi_{\sigma_\lambda}))^{\chi_{\sigma,\lambda}}_{(\hat\Delta)}\rightarrow 0\ .
\end{eqnarray*}
\end{prop}
\proof The proposition is a consequence of Proposition \ref{trans}, Proposition \ref{bow}, and the surjectivity of $\hat\Delta:\Omega^{p}(X,\tilde E(\pi_{\sigma_\lambda}))\rightarrow
\Omega^{p}(X,\tilde E(\pi_{\sigma_\lambda}))$. In fact, the exactness of the pieces of the de Rham complexes
$\Omega^*(X,\tilde E(\pi_{\sigma_\lambda}))
_{(\hat\Delta)}^{\chi_{\sigma,\lambda}}$ appearing in these complexes is a direct consequence of Theorem \ref{hop} and the fact that $\cZ(\gaaa)$ acts locally finitely on
$\Omega^*(X,\tilde E(\pi_{\sigma_\lambda}))
_{(\hat\Delta)}$. It remains to show that for a pair of closed forms $(\omega,\eta)\in
\Omega^{p+1}(X,\tilde E(\pi_{\sigma_\lambda}))
_{(\hat\Delta)}^{\chi_{\sigma,\lambda}}\oplus \Omega^{n-p+1}(X,\tilde E(\pi_{\sigma_\lambda}))_{(\hat\Delta)}
^{\chi_{\sigma,\lambda}}$ there exist $\psi\in\Omega^{p}(X,\tilde E(\pi_{\sigma_\lambda}))
_{(\hat\Delta)}^{\chi_{\sigma,\lambda}}$ such that $d\psi=\omega$, $d*\psi=\eta$. Since $H^{k}(X,\tilde E(\pi_{\sigma_\lambda}))=\{0\}$ for $k>0$
we find $\psi_0,\psi_1\in\Omega^{p}(X,\tilde E(\pi_{\sigma_\lambda}))_{(\hat\Delta)}
^{\chi_{\sigma,\lambda}}$ such that
$d\psi_0=\omega$, $d*\psi_1=\eta$. Now by Corollary \ref{ellen} the equation
$\hat\Delta \psi_2=*(\psi_1-\psi_0)$ is solvable in $\Omega^{n-p}(X,\tilde E(\pi_{\sigma_\lambda}))_{(\hat\Delta)}
^{\chi_{\sigma,\lambda}}$. We set $\psi:=\psi_0+(-1)^{(p-1)(n-1)+1}d*d\psi_2$.
Then $d\psi=\omega$ and
$d*\psi=d*\psi_0+d\hat\Delta \psi_2
=d*\psi_1=\eta$.
\hB

\begin{prop}\label{molo}
Let $\Gamma\subset G$ be a torsion-free discrete subgroup of $G$, and let
$(\vp,V_\vp)$ be a finite-dimensional representation of $\Gamma$. Let $\lambda\in I_\sigma^r$. Then
we have for all $p\ge 1$
\begin{equation}\label{ba}
H^p\left(\Gamma,C^{-\omega}(\partial X, V^+(\sigma_\lambda,\vp))\right)\cong
H^{l_{\sigma_\lambda}+p}(\Gamma,F_{\sigma_\lambda}\otimes V_\vp)\ .
\end{equation}
If $\Gamma$ is not cocompact or $\vp$ is unitary, 
then for all $k\in\nat$, $p\ge 1$ there is a natural isomorphism
\begin{equation}\label{nane}
H^p\left(\Gamma,C^{-\omega}(\partial X, V^k(\sigma_\lambda,\vp))\right)\cong
H^{l_{\sigma_\lambda}+p}(\Gamma,F_{\sigma_\lambda}\otimes V_\vp)\oplus 
H^{l_{\sigma_\lambda}+p-1}(\Gamma,F_{\sigma_\lambda}\otimes V_\vp)\ .
\end{equation}
For cocompact $\Gamma$ and general $\vp$ the isomorphism (\ref{nane})
holds for all $p\ge 2$, while for $p=1$ we have 
\begin{eqnarray}\label{utsch}
\dim H^1\left(\Gamma,C^{-\omega}(\partial X, V^k(\sigma_\lambda,\vp))\right)
&=&\dim H^{l_{\sigma_\lambda}+1}(\Gamma,F_{\sigma_\lambda}\otimes V_\vp)\\
&&+\dim H^0\left(\Gamma,C^{-\omega}(\partial X, V^k(\sigma_\lambda,\vp))\right) \  .\nonumber
\end{eqnarray}
\end{prop} 
\proof
By Lemma \ref{ac} the $\Gamma$-module $\Omega^p(X,\tilde E(\pi_{\sigma_\lambda}\otimes\vp))_{(\Delta)}$ is acyclic.
Since $\cZ(\gaaa)$ acts locally finitely on this module $\Omega^p(X,\tilde E(\pi_{\sigma_\lambda}\otimes\vp))_{(\Delta)}^{\chi_{\sigma,\lambda}}$ is a direct summand of it, hence is acyclic, too.
Now Proposition \ref{dRh} provides an acyclic resolution of $C^{-\omega}(\partial X, V^+(\sigma_\lambda,\vp))$, namely
\begin{eqnarray}
0\rightarrow C^{-\omega}(\partial X, V^+(\sigma_\lambda,\vp))&\rightarrow&
\Omega^{l_{\sigma_\lambda}}(X,\tilde E(\pi_{\sigma_\lambda}\otimes\vp))_{(\hat\Delta)}
^{\chi_{\sigma,\lambda}}
\stackrel{d}{\longrightarrow} \Omega^{l_{\sigma_\lambda}+1}(X,\tilde E(\pi_{\sigma_\lambda}\otimes\vp))_{(\hat\Delta)}
^{\chi_{\sigma,\lambda}}\nonumber
\\
&&\stackrel{d}{\longrightarrow}\dots \stackrel{d}{\longrightarrow}   
\Omega^n(X,\tilde E(\pi_{\sigma_\lambda}\otimes\vp))_{(\hat\Delta)}
^{\chi_{\sigma,\lambda}}\rightarrow 0\ .\label{reso}
\end{eqnarray}
It follows that $H^*\left(\Gamma,C^{-\omega}(\partial X, V^+(\sigma_\lambda,\vp))\right)$ is given as the cohomology of the complex
\begin{eqnarray}
0\rightarrow
\Omega^{l_{\sigma_\lambda}}(Y,E(\pi_{\sigma_\lambda}\otimes\vp))_{(\hat\Delta)}
^{\chi_{\sigma,\lambda}}
&\stackrel{d}{\longrightarrow}& \Omega^{l_{\sigma_\lambda}+1}(Y, E(\pi_{\sigma_\lambda}\otimes\vp))_{(\hat\Delta)}
^{\chi_{\sigma,\lambda}}
\stackrel{d}{\longrightarrow}\nonumber\\
\dots &\stackrel{d}{\longrightarrow}&   
\Omega^n(Y,E(\pi_{\sigma_\lambda}\otimes\vp))_{(\hat\Delta)}
^{\chi_{\sigma,\lambda}}\rightarrow 0\ .\label{invcom}
\end{eqnarray}
(\ref{ba}) now follows from Theorem \ref{hop}
combined with the observation that taking the component having the correct
infinitesimal character does not change the cohomology (compare the proof
of Theorem \ref{main2}).

We now look at the long exact sequence
associated to
$$ 0\rightarrow C^{-\omega}(\partial X, V^k(\sigma_\lambda,\vp))
\stackrel{i}{\longrightarrow} C^{-\omega}(\partial X, V^+(\sigma_\lambda,\vp))
\stackrel{\p^k}{\longrightarrow} 
C^{-\omega}(\partial X, V^+(\sigma_\lambda,\vp))\rightarrow 0\ . $$
We need the following

\begin{lem}\label{lazy}
Let $\p_p$ be the operator induced by $\p$ on $H^p\left(\Gamma,C^{-\omega}(\partial X, V^+(\sigma_\lambda,\vp))\right)$.
Then $\p_p=0$ for all $p\ge 1$.
\end{lem}
\proof We look at the operator $|\alpha|^2\p((n+1-2l_{\sigma_\lambda})\id-\p)$. By (\ref{mupf}) the
Laplacians $\hat\Delta$ on $\Omega^{k}(X,\tilde E(\pi_{\sigma_\lambda}\otimes\vp))_{(\hat\Delta)}
^{\chi_{\sigma,\lambda}}$, $k\ge l_{\sigma_\lambda}$ extend it to the
acyclic resolution (\ref{reso}) of $C^{-\omega}(\partial X, V^+(\sigma_\lambda,\vp))$. However, the Laplacians act by zero on the higher cohomology groups  of (\ref{invcom}). It follows that $|\alpha|^2\p_p((n+1-2l_{\sigma_\lambda})\id-\p)_p=0$ for $p\ge 1$. Since for $l_{\sigma_\lambda}\ne\frac{n+1}{2}$ the operator $(n+1-2l_{\sigma_\lambda})\id-\p$ is invertible $((n+1-2l_{\sigma_\lambda})\id-\p)_p$ is invertible, too. Hence $\p_p=0$.
For $l_{\sigma_\lambda}=\frac{n+1}{2}$ we extend $\p$ by a constant multiple 
of $d*$ on $\Omega^{l_{\sigma_\lambda}}(X,\tilde E(\pi_{\sigma_\lambda}\otimes\vp))_{(\hat\Delta)}
^{\chi_{\sigma,\lambda}}$ (see (\ref{mipf})) and by the zero operator
on $\Omega^{k}(X,\tilde E(\pi_{\sigma_\lambda}\otimes\vp))_{(\hat\Delta)}
^{\chi_{\sigma,\lambda}}$ for $k>l_{\sigma_\lambda}$. Again $\p_p=0$ follows.
\hB

We continue the proof of Proposition \ref{molo}. By Lemma \ref{lazy} the above
mentioned long exact sequence gives rise to short exact sequences
\begin{eqnarray}
0&\rightarrow& \coker \p_0^k \rightarrow H^1\left(\Gamma,C^{-\omega}(\partial X, V^k(\sigma_\lambda,\vp))\right)
\stackrel{i_1}{\longrightarrow} H^1\left(\Gamma,C^{-\omega}(\partial X, V^+(\sigma_\lambda,\vp))\right)\rightarrow 0\ ,\nonumber\\ 
0&\rightarrow& H^{p-1}\left(\Gamma,C^{-\omega}(\partial X, V^+(\sigma_\lambda,\vp))\right)\rightarrow  H^{p}\left(\Gamma,C^{-\omega}(\partial X, V^k(\sigma_\lambda,\vp))\right)
\label{gurke}\\
&&\hspace{5cm}\stackrel{i_p}{\longrightarrow} H^p\left(\Gamma,C^{-\omega}(\partial X, V^+(\sigma_\lambda,\vp))\right)\rightarrow 0\ ,\ p\ge 2\ .\nonumber
\end{eqnarray}
By Proposition \ref{acb}, 3, the space $H^0\left(\Gamma,C^{-\omega}(\partial X, V^+(\sigma_\lambda,\vp))\right)$ is finite-dimensional for cocompact $\Gamma$. It follows
that $\dim\coker\p_0^k=\dim\ker\p_0^k=\dim H^0\left(\Gamma,C^{-\omega}(\partial X, V^k(\sigma_\lambda,\vp))\right)$. This implies (\ref{utsch}). 

It remains to construct a split of $i_p$ in the sequences (\ref{gurke}) except for
$p=1$, $\Gamma$ cocompact, $\vp$ not unitary. The second complex in Proposition \ref{dRh} gives rise to an acyclic resolution of $\tilde Z^{n-l_{\sigma_\lambda}}_{F_{\sigma_\lambda}}\otimes V_\vp$.
Sending a cocycle $\omega\in
Z^{l_{\sigma_\lambda}+p}(Y,E(\pi_{\sigma_\lambda}\otimes\vp))_{(\hat\Delta)}
^{\chi_{\sigma,\lambda}}$ of (\ref{invcom}) to $(\omega,0)\in \Omega^{l_{\sigma_\lambda}+p}(Y,E(\pi_{\sigma_\lambda}))
_{(\hat\Delta)}^{\chi_{\sigma,\lambda}}\oplus \Omega^{n-l_{\sigma_\lambda}+p}(Y,E(\pi_{\sigma_\lambda}))_{(\hat\Delta)}
^{\chi_{\sigma,\lambda}}$ defines for $p\ge 2$ a map
\begin{equation}\label{kaese} 
q_p: H^p\left(\Gamma,C^{-\omega}(\partial X, V^+(\sigma_\lambda,\vp))\right)
\rightarrow H^p\left(\Gamma,\tilde Z^{n-l_{\sigma_\lambda}}_{F_{\sigma_\lambda}}\otimes V_\vp))\right)\ .
\end{equation}
Using Hodge theoretic arguments as for example in the proof of Proposition \ref {dRh} we see that $q_p$ is also well-defined for $p=1$ unless $\Gamma$ is cocompact and $\vp$ not unitary. Let $j:\tilde Z^{n-l_{\sigma_\lambda}}_{F_{\sigma_\lambda}}\otimes V_\vp\rightarrow
C^{-\omega}(\partial X, V^k(\sigma_\lambda,\vp))$ be the inclusion. Then by
construction $(i\circ j)_p\circ q_p=\id$. Thus $j_p\circ q_p$ provides the
desired splitting. 
\hB 

Using the second sequence in Proposition \ref{dRh} one can also
compute the higher cohomology groups of $\tilde Z^{n-l_{\sigma_\lambda}}_{F_{\sigma_\lambda}}\otimes V_\vp$. We find
for $p\ge 1$ ($p\ge 2$ if $\Gamma$ is cocompact and $\vp$ not unitary)
that
$$ H^p\left(\Gamma,\tilde Z^{n-l_{\sigma_\lambda}}_{F_{\sigma_\lambda}}\otimes V_\vp\right)\cong
H^{l_{\sigma_\lambda}+p}(\Gamma,F_{\sigma_\lambda}\otimes V_\vp)\oplus 
H^{n-l_{\sigma_\lambda}+p}(\Gamma,F_{\sigma_\lambda}\otimes V_\vp)\ .
$$

Now we assume $\Gamma$ to be convex cocompact. In order to get
a uniform description of the cohomology for all cases including
cocompact $\Gamma$ we let $\tilde E^+_\Lambda(\sigma_\lambda,\vp)$ be the kernel
of the surjection $C^{-\omega}(\partial X, V^+(\sigma_\lambda,\vp))\rightarrow
H^{l_{\sigma_\lambda}}(\Gamma,F_{\sigma_\lambda}\otimes V_\vp)$ appearing
in Theorem \ref{main2}. Then for noncocmpact $\Gamma$ we have
$\tilde E^+_\Lambda(\sigma_\lambda,\vp)=E^+_\Lambda(\sigma_\lambda,\vp)$, 
whereas $\tilde E^+_\Lambda(\sigma_\lambda,\vp)=\{0\}$, if $\Gamma$ is cocompact and $\vp$ is unitary. 

\begin{theorem}\label{louvre}
Let $\Gamma\subset G$ be a torsion-free convex cocompact discrete subgroup
of $\Gamma$, and let
$\lambda\in I_\sigma^r$. Then Assertions (i)-(iv) of Conjecture \ref{gopatt}
hold true. More precisely, we have
\begin{enumerate}
\item For all 
$p\ge 1$ 
$$H^p\left(\Gamma,C^{-\omega}(\Lambda, V^+(\sigma_\lambda,\vp))\right)\cong
H^{l_{\sigma_\lambda}+p}(\Gamma,F_{\sigma_\lambda}\otimes V_\vp)\ ,$$  
and there is
an exact sequence
$$0\rightarrow \tilde E^+_\Lambda(\sigma_\lambda,\vp) \rightarrow H^0\left(\Gamma,C^{-\omega}(\Lambda, V^+(\sigma_\lambda,\vp))\right)\rightarrow H^{l_{\sigma_\lambda}}(\Gamma,F_{\sigma_\lambda}\otimes V_\vp)\rightarrow 0\ .$$
For $k\ge k_+(\sigma_\lambda,\vp)$ we have 
$$ H^0\left(\Gamma,C^{-\omega}(\Lambda, V^k(\sigma_\lambda,\vp))\right)=H^0\left(\Gamma,C^{-\omega}(\Lambda, V^+(\sigma_\lambda,\vp))\right)\ .$$
\item For all $k\in\nat$
\begin{eqnarray*}
H^p\left(\Gamma,C^{-\omega}(\Lambda, V^k(\sigma_\lambda,\vp))\right)&\cong&
H^{l_{\sigma_\lambda}+p}(\Gamma,F_{\sigma_\lambda}\otimes V_\vp)\\
&&\oplus H^{l_{\sigma_\lambda}+p-1}(\Gamma,F_{\sigma_\lambda}\otimes V_\vp)\ ,\quad p\ge 2\ ,\\
\dim H^1\left(\Gamma,C^{-\omega}(\Lambda, V^k(\sigma_\lambda,\vp))\right)&=&
\dim H^0\left(\Gamma,C^{-\omega}(\Lambda, V^k(\sigma_\lambda,\vp))\right)\\
&&+\dim H^{l_{\sigma_\lambda}+1}(\Gamma,F_{\sigma_\lambda}\otimes V_\vp)\ . 
\end{eqnarray*}
\item For any $k\ge k_+(\sigma,\lambda)$  
\begin{eqnarray*}
\lefteqn{\chi\left(\Gamma,C^{-\omega}(\Lambda, V^+(\sigma_\lambda,\vp))\right)
=-\chi_1\left(\Gamma,C^{-\omega}(\Lambda, V^k(\sigma_\lambda,\vp))\right)}\\
&=&\dim{}^\Gamma C^{-\infty}(\Lambda, V^+(\sigma_\lambda,\vp))+
\sum_{p=l_{\sigma_\lambda}+1}^{n} (-1)^{p-l_{\sigma_\lambda}}\dim H^p(\Gamma,F_{\sigma_\lambda}\otimes V_\vp)\\
&=&\dim \tilde E^+_\Lambda(\sigma_\lambda,\vp)+
\sum_{p=l_{\sigma_\lambda}}^{n} (-1)^{p-l_{\sigma_\lambda}}\dim H^p(\Gamma,F_{\sigma_\lambda}\otimes V_\vp)  \ .
\end{eqnarray*}
For $l_{\sigma_\lambda}=1$ this sum can be rewritten as
\begin{eqnarray*}
\chi\left(\Gamma,C^{-\omega}(\Lambda, V^+(\sigma_\lambda,\vp))\right)&=& \dim \tilde E^+_\Lambda(\sigma_\lambda,\vp)+\dim {}^\Gamma(F_{\sigma_\lambda}\otimes V_\vp)
-\chi(\Gamma,F_{\sigma_\lambda}\otimes V_\vp)\\
&=& \dim \tilde E^+_\Lambda(\sigma_\lambda,\vp)+\dim {}^\Gamma(F_{\sigma_\lambda}\otimes V_\vp)-
\dim (F_{\sigma_\lambda})\dim(V_\vp)\chi(Y)\ .
\end{eqnarray*}
Here $\chi(Y)$ is the topological Euler characteristic of $Y$.
\end{enumerate}
\end{theorem}
\proof It is sufficient to prove Assertions 1 and 2 and that $\chi\left(\Gamma,C^{-\omega}(\Lambda, V^k(\sigma_\lambda,\vp))\right)=0$.
All facts claimed concerning $H^0\left(\Gamma,C^{-\omega}(\Lambda, V^k(\sigma_\lambda,\vp))\right)$ and $H^0\left(\Gamma,C^{-\omega}(\Lambda, V^+(\sigma_\lambda,\vp))\right)$ are already known from Theorem \ref{main2}.
For cocompact $\Gamma$ the theorem now follows immediately from Proposition
\ref{molo}. We thus assume $\Gamma$ to be noncocompact.

The long exact sequence associated to (\ref{lolek}) together with Lemma \ref{acB}
and the surjectivity of 
$$ res: {}^\Gamma C^{-\omega}(\partial X, V^+(\sigma_\lambda,\vp))\rightarrow
C^{-\omega}(B, V_B^+(\sigma_\lambda,\vp))$$
(compare the proof of Theorem \ref{impo}) provides isomorphisms of $\C[\p]$-modules
\begin{equation}\label{cheat}
H^p\left(\Gamma,C^{-\omega}(\Lambda, V^+(\sigma_\lambda,\vp))\right)
\cong H^p\left(\Gamma,C^{-\omega}(\partial X, V^+(\sigma_\lambda,\vp))\right)\ ,\qquad p\ge 1\ .
\end{equation}
Now Assertion 1 follows from Proposition \ref{molo}. In particular, all
cohomology groups of \linebreak[4]$C^{-\omega}(\Lambda, V^+(\sigma_\lambda,\vp))$ are
finite-dimensional. Let $\p_\Lambda$ be the restriction of $\p$ to
$C^{-\omega}(\Lambda, V^+(\sigma_\lambda,\vp))$. Now the exact sequence
\begin{equation}\label{tango}
0\rightarrow C^{-\omega}(\Lambda, V^k(\sigma_\lambda,\vp))
\stackrel{i}{\longrightarrow} C^{-\omega}(\Lambda, V^+(\sigma_\lambda,\vp))
\stackrel{\p_\Lambda^k}{\longrightarrow} 
C^{-\omega}(\Lambda, V^+(\sigma_\lambda,\vp))\rightarrow 0
\end{equation} 
implies that $H^p\left(\Gamma,C^{-\omega}(\Lambda, V^k(\sigma_\lambda,\vp))\right)$ is finite-dimensional and that
$$\chi\left(\Gamma,C^{-\omega}(\Lambda, V^k(\sigma_\lambda,\vp))\right)=
\chi\left(\Gamma,C^{-\omega}(\Lambda, V^+(\sigma_\lambda,\vp))\right)-\chi\left(\Gamma,C^{-\omega}(\Lambda, V^+(\sigma_\lambda,\vp))\right)=0\ .$$
$\p_{\Lambda}$ induces an operator $\p_{\Lambda,0}$ on $H^0\left(\Gamma,C^{-\omega}(\Lambda, V^+(\sigma_\lambda,\vp))\right)$.
We look at the long exact sequence
associated to (\ref{tango}).
Because of (\ref{cheat}) and Lemma \ref{lazy} we obtain the exact sequence   
$$0\rightarrow \coker (\p_{\Lambda,0})^k \rightarrow H^1\left(\Gamma,C^{-\omega}(\Lambda, V^k(\sigma_\lambda,\vp))\right)
\stackrel{i_1}{\longrightarrow} H^1\left(\Gamma,C^{-\omega}(\Lambda, V^+(\sigma_\lambda,\vp))\right)\rightarrow 0\ .$$
The space $H^0\left(\Gamma,C^{-\omega}(\Lambda, V^+(\sigma_\lambda,\vp))\right)$ is finite-dimensional. It follows
that $\dim\coker(\p_{\Lambda,0})^k=\dim\ker(\p_{\Lambda,0})^k=\dim H^0\left(\Gamma,C^{-\omega}(\Lambda, V^k(\sigma_\lambda,\vp))\right)$.
Moreover, 
$$H^1\left(\Gamma,C^{-\omega}(\Lambda, V^+(\sigma_\lambda,\vp))\right)
\cong H^{l_{\sigma_\lambda}+1}(\Gamma,F_{\sigma_\lambda}\otimes V_\vp)$$ 
by Assertion 1. It follows that
$$ \dim H^1\left(\Gamma,C^{-\omega}(\Lambda, V^k(\sigma_\lambda,\vp))\right)=
\dim H^0\left(\Gamma,C^{-\omega}(\Lambda, V^k(\sigma_\lambda,\vp))\right)+\dim H^{l_{\sigma_\lambda}+1}(\Gamma,F_{\sigma_\lambda}\otimes V_\vp)\ .$$

The direct sum decomposition of $H^p\left(\Gamma,C^{-\omega}(\Lambda, V^k(\sigma_\lambda,\vp))\right)$ claimed in Assertion 2 follows from the long exact sequence
associated to (\ref{bolek}), Lemma \ref{acB} and Proposition \ref{molo}.
The proof of the theorem is now complete.
\hB

We conclude this section by a couple of remarks on Theorem \ref{louvre}.

For any $\sigma$ there exists a constant $c(\sigma)\ge 0$ such that
for $\lambda\in I_\sigma^r$ with $|\lambda|> c_\sigma$ we have
$$ l_{\sigma_\lambda}=\left\{\begin{array}{ccc}
n&& \lambda>0\\1&& \lambda<0
\end{array}\right.\ .$$
Assume $\vp$ to be unitary. For $l_{\sigma_\lambda}=n$ the vanishing results Corollary \ref{iso1} and Proposition \ref{bowa} imply  that $H^p\left(\Gamma,C^{-\omega}(\Lambda, V^+(\sigma_\lambda,\vp))\right)=\{0\}$
for all $p$ with the obvious exception $\sigma=1$ and $\Gamma$ cocompact.
Thus the generic non-zero contribution to the expected "topological part"
of the divisor of the Selberg zeta function
$$ \{\left(\lambda, \chi\left(\Gamma,C^{-\omega}(\Lambda, V^+(\sigma_\lambda,\vp))\right)\right)\:|\: \lambda\in I^r_\sigma\}\ $$
is given by
$$ \{(\lambda, \dim E^+_\Lambda(\sigma_\lambda,\vp)+\dim {}^\Gamma(F_{\sigma_\lambda}\otimes V_\vp)-
\dim (F_{\sigma_\lambda})\dim(V_\vp)\chi(Y))\:|\:\lambda<-c_\sigma\} . $$
For cocompact $\Gamma$ and $\sigma_\lambda\ne 1_{-\rho}$ we have 
$E^+_\Lambda(\sigma_\lambda,\vp)=\{0\}$, 
${}^\Gamma(F_{\sigma_\lambda}\otimes V_\vp)=\{0\}$.
Thus the above multiplicity simplifies to $-\dim (F_{\sigma_\lambda})\dim(V_\vp)\chi(Y)$. This coincides with the formula for
the order of the singularity of $Z_{S,\sigma,\vp}$ at very negative integer points obtained in \cite{bunkeolbrich955}, Thm. 3.15.

If $\sigma$ is trivial, then we can choose $c(\sigma)=0$. If in addition, $\Gamma$ is not cocompact and $\vp$ is trivial, then $\dim  E^+_\Lambda(\sigma_\lambda,\vp)+\dim {}^\Gamma(F_{\sigma_\lambda}\otimes V_\vp)$
can be identified with the multiplicity $n_\lambda$ of the pole
of a normalized version of the scattering matrix $\hat S_\mu$ at $\mu=\lambda$ (see \cite{bunkeolbrich990}, Section 5). We obtain for negative $\lambda\in I^r_1$
$$\chi\left(\Gamma,C^{-\omega}(\Lambda,V^+(1_\lambda))\right)=n_\lambda-\chi(Y)\dim F_{1_\lambda}\ .$$
It was proved by Patterson and Perry \cite{pattersonperry01} that the right hand side coincides with $\ord_\lambda (Z_{S,1,1})$. This eventually
led to the proof of Conjecture \ref{gopatt}, (v), in this special case \cite{bunkeolbrich990}.
One of the difficulties one is confronted with, if one tries to extend this
approach to general $\sigma$, is that for $l_{\sigma_\lambda}\not\in\{1,n\}$ the normalized scattering matrix has
an infinite-dimensional pole at $\lambda$, and thus there is no easy definition
of the multiplicity $n_\lambda$. 

One can also get more precise information on the
spaces $H^p\left(\Gamma,C^{-\omega}(\Lambda, V^k(\sigma_\lambda,\vp))\right)$, $p=0,1$. Assume $\Gamma$ to be noncocompact. Let 
$$ res_k: {}^\Gamma C^{-\omega}(\partial X, V^k(\sigma_\lambda,\vp))\rightarrow
C^{-\omega}(B, V_B^k(\sigma_\lambda,\vp))$$
be the restriction map. The long exact sequence associated to (\ref{bolek})
combined with (\ref{nane}) yields an exact sequence
$$ 0\rightarrow \coker\: res_k\rightarrow H^1\left(\Gamma,C^{-\omega}(\Lambda, V^k(\sigma_\lambda,\vp))\right)\rightarrow H^{l_{\sigma_\lambda}+1}(\Gamma,F_{\sigma_\lambda}\otimes V_\vp)\oplus 
H^{l_{\sigma_\lambda}}(\Gamma,F_{\sigma_\lambda}\otimes V_\vp)
\rightarrow 0\ .$$
Combining this with Assertion 2 of Theorem \ref{louvre} we obtain an alternative
proof of
Corollary \ref{est}:
$$\dim H^0\left(\Gamma,C^{-\omega}(\Lambda, V^k(\sigma_\lambda,\vp))\right)
=\dim \coker\  res_k +\dim H^{l_{\sigma_\lambda}}(\Gamma,F_{\sigma_\lambda}\otimes V_\vp)\ .$$

For $n=2$ and noncocompact $\Gamma$ the dimensions of the cohomology groups $H^p\left(\Gamma,C^{-\omega}(\Lambda, V(1_\lambda))\right)$ have been computed
in \cite{bunkeolbrich97} by slightly different methods. Note that in this
case $\Gamma$ has cohomological dimension 1. Thus 
$$\dim H^0\left(\Gamma,C^{-\omega}(\Lambda, V(1_\lambda))\right)=
\dim H^1\left(\Gamma,C^{-\omega}(\Lambda, V(1_\lambda))\right)\ .$$
We refer to the end of Section \ref{twist} for a discussion of the dimension of $H^0$ at negative half integers $\lambda$.

\newpage
\section[$L^2$-cohomology of Kleinian manifolds]{The discrete spectrum and $L^2$-cohomology of Kleinian manifolds}\label{mazzeo}

In this section we consider a torsion-free convex cocompact non-cocompact
subgroup of a linear rank one Lie group $G$. Let $K\subset G$ be a maximal compact subgroup, and let $\gaaa=\kaaa\oplus\paaa$ be the corresponding Cartan decomposition of the Lie algebra. We freely use all the notation introduced in
previous sections. Whenever spaces of the form $E^1_\Lambda(\sigma_\lambda,\varphi)$ occur we tacitly assume that
$X=G/K\ne\OO H^2$ or that $\delta_\Gamma<0$. Let $(\vp,V_\vp)$
be a finite-dimensional unitary representation of $\Gamma$. 
$G$ acts unitarily on the Hilbert space
$$ L^2(\Gamma\backslash G,\vp):=\{f:G\rightarrow V_\vp\:|\: f(gx)=\vp(g)f(x)\ \forall g\in\Gamma,\: x\in G,\:\ 
\int_{\Gamma\backslash G} |f(x)|^2\:dx<\infty\}\ ,$$
where the action is given by right translation. One of the main
results of \cite{bunkeolbrich000} was the decomposition of this representation
into irreducible ones, the so called Plancherel Theorem. Based on the results of Section \ref{round} we will give 
a more precise description of the irreducible representations
which can occur discretely in this decomposition.

The second theme will be $L^2$-cohomology of Kleinian manifolds. Let $(\pi,F)$ be an irreducible finite-dimensional representation of $G$. The admissible
scalar product on $F$ (see Subsection \ref{nancy}) induces an $L^2$ scalar
product on the space $\Omega^*_c(Y,E(\pi\otimes\vp))$ of compactly supported
differential forms with values in the flat vector bundle $E(\pi\otimes\vp)$ 
such that the corresponding Hilbert space of square-integrable forms satisfies
\begin{equation}\label{hihi} 
\Omega^*_{(2)}(Y,E(\pi\otimes\vp))\cong [L^2(\Gamma\backslash G,\vp)\otimes F\otimes \Lambda^*\paaa^*]^K\ .
\end{equation}
Using results on $(\gaaa,K)$-cohomology we will compute the $L^2$-cohomology groups $H^p_{(2)}(Y,E(\pi\otimes\vp))$ and their reduced versions $\cH^p_{(2)}(Y,E(\pi\otimes\vp))$ in terms of multiplicities of the unitary
representations with infinitesimal character $\chi_{\tilde F}$ in the discrete part of $L^2(\Gamma\backslash G,\vp)$, i.e., in terms of invariant distributions
supported on the limit set. Here $\tilde F$
is the dual representation of $F$. A comparison of these results with Theorem \ref{main2} then shows that for $X=\R H^n$ and $p\ge \frac{n+1}{2}$ the natural
map
$$\cH^p_{(2)}(Y,E(\pi\otimes\vp))\rightarrow H^p(Y,E(\pi\otimes\vp))$$
is an isomorphism. This generalizes the main result of \cite{mazzeophillips90}
for convex cocompact $\Gamma$ to non-trivial $\pi$ and $\vp$.

In order to state our refined version of the Plancherel Theorem for $L^2(\Gamma\backslash G,\vp)$ we will introduce a couple of $(\gaaa,K)$-modules.
By $H^{\sigma,\lambda}$ we denote the underlying $(\gaaa,K)$-module of $C^\infty(\partial X,V(\sigma_\lambda))$. Recall that for $\Ree(\lambda)>0$
it has the unique nontrivial irreducible submodule $I^{\sigma,\lambda}$. 

Let $\sigma\in\hat M$ such that $p_\sigma(0)=0$. Let $0<\lambda_\sigma\in\aaaa^*$ be the corresponding end of the complementary series, i.e., the interval $(0,\lambda_\sigma)$ consists of all $\lambda$
with $\Ree(\lambda)>0$ such that $H^{\sigma,\lambda}$ is irreducible and
unitarizable. $\lambda_\sigma$ can be characterized as the smallest $\lambda>0$
such that $H^{\sigma,\lambda}$ becomes reducible. In particular, 
$\lambda_\sigma\in I^{wr}_\sigma$ (see Lemma \ref{sauna}). The value of $\lambda_\sigma$ is explicitly known in all cases (see \cite{baldonisilvabarbasch83} and the literature cited therein).
For $X=\R H^n$ it is given before Proposition \ref{knut}. Note that $\lambda_{\tilde\sigma}=\lambda_\sigma$.
We set   
$$ \cH_\res(\vp):=
\bigoplus_{\{\sigma\in\hat M\:|\: p_\sigma(0)=0\}}
\bigoplus_{\lambda\in (0,\lambda_\sigma)\cap (0,\delta_\Gamma]} E^1_\Lambda(\tilde\sigma_\lambda,\varphi)\otimes H^{\sigma,\lambda}\ .$$
For any $\sigma$ the sum over $\lambda$ is actually finite by Corollary \ref{iso1}. By (38) and Lemma \ref{less} we have 
$E^1_\Lambda(\tilde\sigma_\lambda,\varphi)={}^{\Gamma} C^{-\infty}(\Lambda,V(\tilde\sigma_\lambda,\varphi))$ unless $\lambda\in I^{wr,-}_\sigma$. By the remark following Lemma \ref{sauna} this cannot
happen for $\lambda\in (0,\lambda_\sigma)$ and $X=\R H^n$ or $X=\C H^n$.
We have seen in the proof of Proposition \ref{newer} that the matrix
coefficient map
$$ E^1_\Lambda(\tilde\sigma_\lambda,\varphi)\otimes H^{\sigma,\lambda}\ni
f\otimes\phi \mapsto c_{f,\phi} $$
injects $E^1_\Lambda(\tilde\sigma_\lambda,\varphi)\otimes H^{\sigma,\lambda}$
in a $(\gaaa,K)$-equivariant way into $L^2(\Gamma\backslash G,\vp)_K$, where
$L^2(\Gamma\backslash G,\vp)_K$ denotes the $(\gaaa,K)$-module of $K$-finite
smooth vectors in $L^2(\Gamma\backslash G,\vp)$. We obtain an embedding
of $\cH_\res(\vp)$ into $L^2(\Gamma\backslash G,\vp)_K$.
Let $L^2(\Gamma\backslash G,\vp)_\res$ be the closure in $L^2(\Gamma\backslash G,\vp)$ of the image of $\cH_\res(\vp)$ under the matrix coefficient map.

We further introduce
$$ \cH_{U,+}(\vp):=
\bigoplus_{\{\sigma\in\hat M\setminus\{1\}\:|\: p_\sigma(0)=0\}}
\bigoplus_{\{\lambda\in I^{wr,-}_\sigma\cap (0,\delta_\Gamma]\:|\: I^{\sigma,\lambda}{\mathrm unitarizable}\}} U_\Lambda(\tilde\sigma_\lambda,\varphi)\otimes I^{\sigma,\lambda}\ .$$
Again by the remark following Lemma \ref{sauna} we have in case $X=\R H^n$ or 
$X=\C H^n$ that $U_\Lambda(\tilde\sigma_\lambda,\varphi)={}^{\Gamma} C^{-\infty}(\Lambda,V(\tilde\sigma_\lambda,\varphi))$ for $\lambda\in I^{wr,-}_\sigma$. Moreover, it is known that in these two cases
\begin{equation}\label{tery}
\{\lambda\in \aaaa^*\:|\: \Ree(\lambda)>0, I^{\sigma,\lambda}{\mathrm unitarizable}\}=(0,\lambda_\sigma]\ .
\end{equation}
It follows that
$$\{\lambda\in I^{wr,-}_\sigma\cap (0,\delta_\Gamma]\:|\: I^{\sigma,\lambda}{\mathrm unitarizable}\}\subset\{\lambda_\sigma\}\ .$$

Last but not least we collect the contributions of the principal
series modules $H^{\sigma,0}$ for $p_\sigma(0)\ne 0$.
For $X=\R H^n$, $n$ odd, let $\hat M_0\subset \hat M$ a set of representatives
of $W(\gaaa,\aaaa)$-orbits of representations which factorize over $SO(n-1)$
and which are not Weyl-invariant. We set  
$$ \cH_{U,0}(\vp):=
\bigoplus_{\sigma\in\hat M_0} 
U_\Lambda(\tilde\sigma_0,\varphi)\otimes H^{\sigma,0}\ .$$
If $\dim X$ is even, then $\hat M_0:=\{\sigma\in\hat M\:|\: 0\in I^{wr}_\sigma\}$. Note that for $X=\R H^n$, $n$ even, $\hat M_0=\emptyset$.
The module $H^{\sigma,0}$ splits into two eigenspaces of $\hat J_{\sigma,0}$
$$ H^{\sigma,0}=H^{\sigma,+}\oplus H^{\sigma,-}\ .$$
Similarly, we obtain
$$ U_\Lambda(\tilde\sigma_0,\varphi)=U_\Lambda^+(\tilde\sigma_0,\varphi)
\oplus U^-_\Lambda(\tilde\sigma_0,\varphi)\ .$$
Then
$$ \cH_{U,0}(\vp):=
\bigoplus_{\sigma\in\hat M_0} 
U^+_\Lambda(\tilde\sigma_0,\varphi)\otimes H^{\sigma,+}
\oplus U^-_\Lambda(\tilde\sigma_0,\varphi)\otimes H^{\sigma,-}\ .$$

We set
$$ \cH_{U}(\vp):=\cH_{U,+}(\vp)\oplus \cH_{U,0}(\vp)$$
Let $L^2(\Gamma\backslash G,\vp)_U$ be the closure in $L^2(\Gamma\backslash G,\vp)$ of the image of $\cH_U(\vp)$ under the matrix coefficient map.

Recall the definition and classification of discrete series representations
of $G$ from e.g. \cite{knapp86}, Chapters IX and XII. We can now formulate the following refinement of \cite{bunkeolbrich000}, 
Theorem 11.1.

\begin{theorem}\label{plaff}
If $X\ne\OO H^2$ or $\delta_\Gamma<0$, then there is a $G$-equivariant decomposition
$$ L^2(\Gamma\backslash G,\vp) = L^2(\Gamma\backslash G,\vp)_{ac}\oplus L^2(\Gamma\backslash G,\vp)_d
\ ,$$
where the discrete subspace $L^2(\Gamma\backslash G,\vp)_d$, being the 
closure of the sum of all irreducible subrepresentations of $L^2(\Gamma\backslash G,\vp)$,  
has a further decomposition
$$L^2(\Gamma\backslash G,\vp)_d=L^2(\Gamma\backslash G,\vp)_{cusp}\oplus L^2(\Gamma\backslash G,\vp)_{U}\oplus L^2(\Gamma\backslash G,\vp)_{\res}\ .
$$
$L^2(\Gamma\backslash G,\vp)_{cusp}$ decomposes into discrete series representations of $G$, each
discrete series representation of $G$ occurs with infinite multiplicity. It is zero
iff $X=\R H^n$, $n$ odd. The remaining part of $L^2(\Gamma\backslash G,\vp)_d$
is zero, if $\delta_\Gamma<0$.

$L^2(\Gamma\backslash G,\vp)_{ac}$ decomposes into a sum of direct integrals
corresponding to the unitary principal series representations of $G$, each
occurring with infinite multiplicity. The corresponding Plancherel measures
are absolutely continuous with respect to the Lebesgue measure of $i\aaaa^*$.
\end{theorem}
\proof The only point where Theorem \ref{plaff} goes beyond \cite{bunkeolbrich000}, 
Theorem 11.1, is that certain contributions of the form $U_\Lambda(\tilde\sigma_\lambda,\varphi)\otimes I^{\sigma,\lambda}$ and
$E^1_\Lambda(\tilde\sigma_\lambda,\varphi)\otimes H^{\sigma,\lambda}$ which
appear in \cite{bunkeolbrich000} do not occur in our definition of
$L^2(\Gamma\backslash G,\vp)_{U}$ and
$L^2(\Gamma\backslash G,\vp)_{\res}$. This is justified by Proposition \ref{newer} and Proposition \ref{new}.
\hB

If one does not know that $ext_\lambda$ is meromorphic in the region 
$\{\Ree(\lambda)>-\ve\}$ for some $\ve>0$, then it is not possible to 
determine $L^2(\Gamma\backslash G,\vp)_{ac}$ and to show that there is
no singular continuous spectrum. However, if one traces back the arguments
in \cite{bunkeolbrich000} and the additional information obtained in Section \ref{round} which led to the determination of 
$L^2(\Gamma\backslash G,\vp)_d$, one obtains
a decomposition of $L^2(\Gamma\backslash G,\vp)_d$ without using $ext_\lambda$.
In particular, it was indicated in \cite{bunkeolbrich000}, pp. 121-122,
how one can prove the finite-dimensionality of ${}^\Gamma C^{-\infty}(\Lambda, V(\sigma_\lambda,\varphi))$ without referring to $ext_\lambda$.
Indeed, also for $X=\OO H^2$ and $\delta_\Gamma\ge 0$ the following is true.

\begin{prop}\label{exc}
Let $L^2(\Gamma\backslash G,\vp)_d$ be the 
closure of the sum of all irreducible subrepresentations of $L^2(\Gamma\backslash G,\vp)$. Let 
$L^2(\Gamma\backslash G,\vp)_{U_0}\subset L^2(\Gamma\backslash G,\vp)_{U}$
be the orthogonal complement of the contributions of those $(\sigma,\lambda)$
such that $\lambda\in (0,\lambda_\sigma)$. Then
\begin{eqnarray*}
L^2(\Gamma\backslash G,\vp)_d&=&L^2(\Gamma\backslash G,\vp)_{cusp}\oplus L^2(\Gamma\backslash G,\vp)_{U_0}\\ 
&&\oplus\  {\mathrm cl}\left(\bigoplus_{\{\sigma\in\hat M\:|\: p_\sigma(0)=0\}}
\bigoplus_{\lambda\in (0,\lambda_\sigma)\cap (0,\delta_\Gamma]} {}^\Gamma C^{-\infty}(\Lambda, V(\tilde\sigma_\lambda,\varphi))\otimes H^{\sigma,\lambda}
\right)\ .
\end{eqnarray*}
Here cl means the closure in $L^2(\Gamma\backslash G,\vp)$.
$L^2(\Gamma\backslash G,\vp)_{cusp}$ decomposes into discrete series representations of $G$, each
discrete series representation of $G$ occurs with infinite multiplicity.
While for fixed $\sigma$ the sum over $\lambda$
might be actually infinite, the individual multiplicity spaces ${}^\Gamma C^{-\infty}(\Lambda, V(\tilde\sigma_\lambda,\varphi))$ and 
$U_\Lambda(\tilde\sigma_0,\varphi)$ are finite-dimensional.
\end{prop}  

We now want to compute the $L^2$-cohomology groups $H^p_{(2)}(Y,E(\pi\otimes\vp))$. They are defined as follows
$$ H^p_{(2)}(Y,E(\pi\otimes\vp)):=\frac{\{\omega\in \Omega^p_{(2)}(Y,E(\pi\otimes\vp))\:|\: d\omega=0\}}{\{d\eta\:|\: \eta \in \Omega^{p-1}_{(2)}(Y,E(\pi\otimes\vp)),d\eta\in \Omega^p_{(2)}(Y,E(\pi\otimes\vp))\}}\ .$$
Here the differential $d$ is taken in the distributional sense. Since $\vp$ is
unitary the operators $\delta_F$ and $\Delta_F$ defined in Subsection \ref{nancy} coincide with the operators $\delta$ and $\Delta$ associated to
the Hermitian metric on $E(\pi\otimes\vp))$ which were introduced in Section \ref{sur}.
We also consider the space of square-integrable harmonic forms
\begin{eqnarray*}
\cH^p_{(2)}(Y,E(\pi\otimes\vp))&=&\{\omega\in \Omega^p_{(2)}(Y,E(\pi\otimes\vp))\:|\: \Delta_F\omega=0\}\\
&=&
\{\omega\in \Omega^p_{(2)}(Y,E(\pi\otimes\vp))\:|\: d\omega=0,\delta_F\omega=0\}
\ .
\end{eqnarray*}
The second equality uses the completeness of $Y$. It follows that the natural map $\cH^p_{(2)}(Y,E(\pi\otimes\vp))\rightarrow H^p_{(2)}(Y,E(\pi\otimes\vp))$ is injective. It is an isomorphism
if and only if the linear subspace $B^p_{(2)}:=\{d\eta\:|\: \eta \in \Omega^{p-1}_{(2)}(Y,E(\pi\otimes\vp)),d\eta\in \Omega^p_{(2)}(Y,E(\pi\otimes\vp))\}\subset \Omega^p_{(2)}(Y,E(\pi\otimes\vp))$
is closed. In general, one has an orthogonal decomposition
\begin{equation}\label{bloed} 
H^p_{(2)}(Y,E(\pi\otimes\vp)):= \cH^p_{(2)}(Y,E(\pi\otimes\vp))
\oplus \overline{B^p_{(2)}}/B^p_{(2)}\ ,
\end{equation}
where $\bar{}$ denotes the closure, and the second summand is infinite-dimensional if not zero. \linebreak[4] $\cH^*_{(2)}(Y,E(\pi\otimes\vp))$
considered as a summand
of $H^*_{(2)}(Y,E(\pi\otimes\vp))$ is often called the reduced $L^2$-cohomology.
For these and more generalities on $L^2$-cohomology we refer to \cite{brylinskizucker90}, Chapter 3. Moreover, we have the following
general fact

\begin{lem}\label{doof}
$\overline{B^p_{(2)}}=B^p_{(2)}$ if and only if $0\not\in \spec(\Delta_{F\:|\overline{B^p_{(2)}}})$.
\end{lem}
\proof Assume that $0\not\in \spec(\Delta_{F\:|\overline{B^p_{(2)}}})$.
This means that $\Delta_{F\:|\overline{B^p_{(2)}}}=d\delta_{F\:|\overline{B^p_{(2)}}}: {\mathrm dom}(\Delta_F)\cap \overline{B^p_{(2)}}\rightarrow \overline{B^p_{(2)}}$ is 
surjective.
It follows that $\im\: d=\overline{B^p_{(2)}}$.

Vice versa, if $d$ has closed range, then by general principles
$\delta_F=d^*$ has closed range, too. 
It follows that $\im \delta_F=\ker d^\perp
\subset \Omega^{p-1}_{(2)}(Y,E(\pi\otimes\vp))$. By assumption any element
in $\overline{B^p_{(2)}}$ can be written as $d\eta$ for some $\eta$. Here
we can assume that $\eta\in \ker d^\perp$. It follows that $d\delta_F$ is surjective on $\overline{B^p_{(2)}}$. By (\ref{bloed}) the Laplacian is always injective on ${\mathrm dom}(\Delta_F)\cap \overline{B^p_{(2)}}$. This shows
that $0\not\in \spec(\Delta_{F\:|\overline{B^p_{(2)}}})$. 
\hB

We want to employ the theory of $(\gaaa,K)$-cohomology of irreducible
unitarizable $(\gaaa,K)$-modules in a similar manner as it has been done
in the computation of the $L^2$-cohomology for symmetric spaces in
\cite{borel85} (compare also \cite{olbrich00}) and for locally 
symmetric spaces of
finite volume in \cite{borelcasselman83}.

The functor of $(\gaaa,K)$-cohomology $V\mapsto H^*(\gaaa,K,V)$ which
goes from the category of $(\gaaa,K)$-modules to the category of vector spaces
is the right derived functor of the left exact functor taking $(\gaaa,K)$-invariants. $H^*(\gaaa,K,V)$ can be computed using the standard relative
Lie algebra cohomology complex $([V\otimes \Lambda^*\paaa^*]^K, d)$, where
$$ d\omega(X_0,\dots,X_p)=\sum_{i=0}^{p}(-1)^i \pi(X_i)\omega(X_0,\dots,\hat X_i,\dots, X_p)\ ,\quad \omega\in [V\otimes \Lambda^p\paaa^*]^K,\  X_i\in\paaa\ .$$
Note that for $V=C^\infty(\Gamma\backslash G,\vp)_K\otimes F$ this complex is isomorphic to the de Rham
complex $\Omega^{*}(Y,E(\pi\otimes\vp))$.
The following basic result can be considered as an algebraic version of Hodge theory.

\begin{prop}[\cite{borelwallach80}, II.3.1. and I.5.3.]\label{char}
Let $(\pi,F)$ be an irreducible finite-dimensional $(\gaaa,K)$-module.
Then
\begin{enumerate}
\item 
For any irreducible unitarizable $(\gaaa,K)$-module $(\xi,V)$
$$ H^p(\gaaa,K,V\otimes F)=\left\{
\begin{array}{cc} 
[V\otimes F\otimes \Lambda^p\paaa^*]^K&\xi(\Omega)=\pi(\Omega)\\
\{0\}&\xi(\Omega)\ne \pi(\Omega)
\end{array}\right. \ .$$
\item
If $V$ is a locally $\cZ(\gaaa)$-finite $(\gaaa,K)$-module, then
$$ H^p(\gaaa,K,V\otimes F)=H^p(\gaaa,K,V^{\chi_{\tilde F}}\otimes F)\ ,$$
where $\tilde F$ is the dual representation of $F$.
In particular, if $V$ has (generalized) infinitesimal character $\chi_V\ne \chi_{\tilde F}$, then
$H^p(\gaaa,K,V\otimes F)=\{0\}$.
\end{enumerate}
\end{prop}

We now give the list of the cohomology groups for all irreducible $V$
such that $\chi_V=\chi_{\tilde F}$ for some finite-dimensional irreducible $F$.
Note that discrete series modules always have this property. Moreover, iff
$X\ne \R H^{2n+1}$, then then for any $F$ there are exactly
$k_X:=|W(\gaaa_\C,\haaa_\kaaa)/W(\kaaa_\C,\haaa_\kaaa)|$ inequivalent discrete
series modules having infinitesimal character $\chi_F$. Here $\haaa_\kaaa$
is a Cartan subalgebra of $\kaaa$ which is also a Cartan subalgebra of $\gaaa$.
 We have \vspace{0.2cm}\\  
\centerline{\begin{tabular}{|c||c|c|c|c|}
\hline
$X$&$\R H^{2n}$&$\C H^n$&$\HH H^n$&$\OO H^2$\\
\hline
$k_X$&
$2$&$n+1$&$n+1$&$3$\\
\hline
\end{tabular} . }\vspace{0.3cm}
In the odd-dimensional case $X\ne \R H^{2n+1}$ there are no discrete
series representations.

\begin{prop}[\cite{borelwallach80},II.5.3.]\label{ds}
Let $V$ be a discrete series module with infinitesimal character
$\chi_{\tilde F}$. Then
$$H^p(\gaaa,K,V\otimes F)=\left\{
\begin{array}{cc} 
\C&p=\frac{\dim X}{2}\\
\{0\}&\mbox{otherwise}
\end{array}\right. \ .$$
\end{prop}

By Langlands classification (see e.g. \cite{collingwood85}, Sect. 2.1)
an irreducible $(\gaaa,K)$-modules with
infinitesimal character $\chi_F$ which is not a discrete series module
is equivalent to exactly one of the modules $I_F^w:=I^{\sigma,\lambda}$, where
$$ w\in W_0:=\{w\in W(\gaaa_\C,\haaa_\C)\:|\: w(\haaa_+^*)\subset\aaaa_-^*\oplus \taaa_+^*\}\subset W^1 $$
and $\sigma_\lambda=\sigma^w_{F,\lambda_w}$. 
Here we have used the notation of Section \ref{twist}, and $\aaaa_-^*:=-\aaaa_+^*$. If $\lambda_w=0$,
which can only happen for $X=\R H^{2n+1}$, then $H^{\sigma^w_F,0}$ is
irreducible, and $I_F^w:=H^{\sigma^w_F,0}$.

If $X=\R H^n$, then we have in the notation of Section \ref{twist}
\begin{eqnarray*}
W^0=\{w_\frac{n}{2},w_{\frac{n}{2}+1},\dots,w_{n-1}\}&& n\mbox{ even },\\
W^0=\{w_-,w_\frac{n+1}{2},\dots,w_{n-1}\}&& n\mbox{ odd }.
\end{eqnarray*}
For $p\in\{0,1,\dots,[\frac{n}{2}]-1\}$ we set $I_F^p:=I_F^{w_{n-1-p}}$,
$I_F^+:=I_F^{w_-}$. Thus $I_F^p$, $I_F^+$ are the underlying $(\gaaa,K)$-modules of the $G$-representations $\tilde Z^p_F$ (introduced before Proposition \ref{dRh}) and $C^{-\omega}(\partial X, V(\sigma^{+}_{F,\lambda_{+}}))$,
respectively.

In the remaining cases the set $W^0$ can be labeled in a natural way as follows
(compare the so-called Hasse diagrams in \cite{collingwood85}, pp. 177-180) \begin{eqnarray*}
\{w_{ij}\:|\: i,j\in \nat_0, i+j\le n-1\}\ ,&& X=\C H^n\ ,\\
\{w_{ij}\:|\: i\le j\in \nat_0, i+j\le 2n-1\}\ ,&& X=\HH H^n\ ,\\
\{w_{00},\dots,w_{04},w_{13},\dots,w_{16},w_{23},w_{24},w_{25}\}\ ,&& X=\OO H^2\ .
\end{eqnarray*}
Note that $w_{ij}$ has length $\dim \partial X-(i+j)$. We will denote
the modules $I_F^{w_{ij}}$ by $I^{ij}_F$.

\begin{prop}\label{koller}
We have for
\begin{eqnarray*}
X=\R H^n:&&
H^p(\gaaa,K,I^q_{\tilde F}\otimes F)=\left\{
\begin{array}{cc} 
\C&p=q, n-q\\
\{0\}&\mbox{otherwise}
\end{array}\right. \ ,\\
&&H^p(\gaaa,K,I^+_{\tilde F}\otimes F)=\left\{
\begin{array}{cc} 
\C&p=\frac{n-1}{2}, \frac{n+1}{2}\\
\{0\}&\mbox{otherwise}
\end{array}\right. \ ,\\
X=\C H^n:
&&H^p(\gaaa,K,I^{ij}_{\tilde F}\otimes F)=\left\{
\begin{array}{cc} 
\C&p=i+j, i+j+2, \dots, 2n-(i+j)\\
\{0\}&\mbox{otherwise}
\end{array}\right. \ ,\\
X=\HH H^n:&&
H^p(\gaaa,K,I^{ii}_{\tilde F}\otimes F)=\left\{
\begin{array}{cc} 
\C&p=2i, 2i+4, \dots, 4n-2i\\
\{0\}&\mbox{otherwise}
\end{array}\right. \ ,\\
&&H^p(\gaaa,K,I^{ij}_{\tilde F}\otimes F)=\left\{
\begin{array}{cc} 
\C&p=i+j, i+j+2, \dots 4n-(i+j)\\
\{0\}&\mbox{otherwise}
\end{array}\right.\ ,\ i\ne j\ ,\\
X=\OO H^2:
&&H^p(\gaaa,K,I^{00}_{\tilde F}\otimes F)=\left\{
\begin{array}{cc} 
\C&p=0,8,16\\
\{0\}&\mbox{otherwise}
\end{array}\right.\ ,\\
&&H^p(\gaaa,K,I^{01}_{\tilde F}\otimes F)=\left\{
\begin{array}{cc} 
\C&p=1,7,9,15\\
\{0\}&\mbox{otherwise}
\end{array}\right.\ ,\\ 
&&H^p(\gaaa,K,I^{02}_{\tilde F}\otimes F)=\left\{
\begin{array}{cc} 
\C&p=2,6,8,10,14\\
\{0\}&\mbox{otherwise}
\end{array}\right.\ ,\\
&&H^p(\gaaa,K,I^{03}_{\tilde F}\otimes F)=\left\{
\begin{array}{cc} 
\C&p=3,5,7,9,11,13\\
\{0\}&\mbox{otherwise}
\end{array}\right.\ ,\\
&&H^p(\gaaa,K,I^{04}_{\tilde F}\otimes F)=\left\{
\begin{array}{cc} 
\C&p=4,8,12\\
\{0\}&\mbox{otherwise}
\end{array}\right.\ ,\\
&&H^p(\gaaa,K,I^{13}_{\tilde F}\otimes F)=\left\{
\begin{array}{cc} 
\C&p=4,6,10,12\\
\{0\}&\mbox{otherwise}
\end{array}\right.\ ,\\
&&H^p(\gaaa,K,I^{14}_{\tilde F}\otimes F)=\left\{
\begin{array}{cc} 
\C&p=5,7,9,11\\
\{0\}&\mbox{otherwise}
\end{array}\right.\ ,\\
&&H^p(\gaaa,K,I^{15}_{\tilde F}\otimes F)=\left\{
\begin{array}{cc} 
\C&p=6,8,10\\
\{0\}&\mbox{otherwise}
\end{array}\right.\ ,\\
&&H^p(\gaaa,K,I^{16}_{\tilde F}\otimes F)=\left\{
\begin{array}{cc} 
\C&p=7,9\\
\{0\}&\mbox{otherwise}
\end{array}\right.\ ,\\
&&H^p(\gaaa,K,I^{2j}_{\tilde F}\otimes F)=\left\{
\begin{array}{cc} 
\C&p=2+j, 14-j\\
\{0\}&\mbox{otherwise}
\end{array}\right.\ ,\ j=3,4,5\ .
\end{eqnarray*}
\end{prop}
\proof
For the trivial representation $F=\C$ the result is given in \cite{collingwood85}, Chapter 7. It is not difficult to derive from Proposition \ref{trans} that $I^w_{\tilde F}\cong (I^w_\C\otimes\tilde F)^{\chi_{\tilde F}}$. Using Proposition \ref{char}, 2., we obtain
$$ H^p(\gaaa,K,I^w_{\tilde F}\otimes F)
=H^p(\gaaa,K,I^w_\C\otimes\tilde F\otimes F)
=H^p(\gaaa,K,I^w_\C\otimes(\tilde F\otimes F)^{\chi_\C})
=H^p(\gaaa,K,I^w_\C)\ .$$
This finishes the proof of the proposition.
\hB

Here we are only interested in the cohomology groups for
unitarizable modules. $I^w_F$ is called isolated unitary, if the corresponding
representation of $G$ is isolated in the unitary dual, i.e., 
$\lambda_w>\lambda_{\sigma^w_F}$. Such modules exist for $X=\HH H^n$ and
$X=\OO H^2$, only. 

\begin{prop}\label{uni}
Among the modules $I^w_F$ precisely the following are unitarizable: 
\begin{itemize}
\item $X=\R H^n$: $I^q_F$ for $p_F\le q$ and $I^+_F$ for $p_F\le\frac{n-1}{2}$.
Here $p_F$ is as in Proposition \ref{bowa}. In these cases we have 
$\lambda_q=\rho-q\alpha$, $\lambda_+=0$. 
\item $X=\C H^n$: $I^{ij}_F$ for those $F$ having highest weight
$\nu=(m_0,m_1,\dots,m_n)\in\haaa^*_+$ with $m_i=m_{i+1}=\dots=m_{n-j}$.
In this case we have $\lambda_{ij}=(n-(i+j))\alpha$.
\item $X=\HH H^n$: $I^{ij}_F$ is isolated unitary iff $i=j$ and the highest weight $\nu=(m_0,m_1,\dots,m_n)\in\haaa^*_+$ of $F$ satisfies 
$m_i=m_{i+1}=\dots=m_n=0$. Then we have $\lambda_{ii}=(2n+1-2i)\alpha$. The remaining unitarizable modules are $I^{ij}_F$
for $j\ge n$ and $m_i=m_{i+1}=\dots=m_{2n-j}$. In these
cases $\lambda_{ij}=(2n-(i+j))\alpha$.
\item $X=\OO H^2$:
\begin{tabular}{|c|c|c|c|}
\hline
$ij$&$\nu$&$\lambda_{ij}$&isolated\\
\hline
$00$&$m_0=\dots=m_3=0$&$11\alpha$&$\times$\\
$04$&$m_0=m_1, m_2=m_3=0$&$5\alpha$&$\times$\\
$15$&$m_0=m_1+m_2, m_3=0$&$2\alpha$&\\
$16$&$m_0=m_1+m_2+m_3$&$\alpha$&\\
$23$&$m_1=m_2=m_3=0$&$5\alpha$&\\
$24$&$m_2=m_3=0$&$3\alpha$&\\
$25$&$m_3=0$&$\alpha$&\\
\hline
\end{tabular}\\ 
Here the highest weight of $F$ is given by $\nu=(m_0,m_1,m_2,m_3)$,
where $m_i\in\frac{1}{2}\Z, m_i-m_j\in\Z, m_0\ge m_1+m_2+m_3, m_1\ge m_2\ge m_3\ge 0$.  \end{itemize}
\end{prop}
\proof The assertion for $X=\R H^n$ was already obtained in the proof of
Proposition \ref{bowa}. For $X=\C H^n$ the module $H^{\sigma,\lambda}$ is irreducible if and only if $\lambda\not\in I^{wr}_\sigma$ (see Lemma \ref{sauna} and the remark following it).  
Now (\ref{tery}) implies that $I^{ij}_F$ is unitarizable if and only if
$\lambda_{ij}-2k\alpha\not\in  I^{wr}_{\sigma_F^{ij}}$
for all $k\in\nat$ such that $\lambda_{F,ij}-2k\alpha\ge 0$. Let $\nu=(m_0,m_1,\dots,m_n)\in\haaa^*_+$ be the highest weight of $F$. $2\alpha$ considered as an element of $\haaa^*_+$ is given by $(1,0,\dots,0,-1)$, $\rho_\gaaa=(\frac{n}{2},\frac{n-2}{2},\dots,-\frac{n-2}{2},-\frac{n}{2})$. 
Then
\begin{eqnarray*} 
-(\lambda_{F,ij}-2k\alpha)+\mu_{\sigma_F^{ij}}+\rho_\maaa
&=&w_{ij}(\nu+\rho_\gaaa)+2k\alpha\\
&=&(m_{n-j}+j+k-\frac{n}{2},m_0+\frac{n}{2},\dots,m_{n}-\frac{n}{2},
m_i+\frac{n}{2}-i-k)\ .
\end{eqnarray*}
Such an element is not weakly regular if and only if at least two coordinates coincide. Now the assertion can be easily checked.

The list for $\HH H^n$ is extracted from the classification of the
unitary dual of $Sp(1,n)$ given in \cite{baldonisilva81}, Thm. 7.1.
Here we have used that, if $\haaa^*_+=\{(m_0,\dots,m_n)\:|\:m_0\ge m_1\ge\dots\ge m_n\ge 0\}$, then $\rho_\gaaa=(n+1,n,\dots,1)$, $2\alpha=(1,1,0,\dots,0)$, and
$$ w_{ij}(\nu+\rho_\gaaa)=\left\{\begin{array}{ll}
(-(m_{j+1}+n-j),-(m_i+n+1-i),m_0+n+1,\dots,m_{n}+1)& j\le n-1\\
(m_{2n-j}+j+1-n,-(m_i+n+1-i),m_0+n+1,\dots,m_{n}+1)& j\ge n
\end{array}\right. \ .$$

For $X=\OO H^2$ and $F=\C$ the list is taken from \cite{collingwood85}, p. 155,
where also the real semisimple part ${}^0\laaa_0$ of the Levi subalgebra of the $\theta$-stable
parabolic $\qaaa\subset\gaaa_\C$ is given such that $I^{ij}_\C=A_\qaaa(0)$. 
Here 
$A_\qaaa(0)$ denotes the Zuckerman module associated to $\qaaa$ 
(for this notion see e.g. \cite{voganzuckerman84}, 
\cite{wallach88}, Chapters 6 and 9).
The reader has to be aware that in this list the algebra 
${}^0\laaa_0={\frak sp}(1,2)$
is missing. Therefore $I^{15}_\C$ is unitarizable 
though stated conversely. Now it follows from the non-vanishing of the
cohomology of all modules in question (see Proposition \ref{koller}) and 
the Vogan-Zuckerman classification of unitarizable modules with cohomology (\cite{voganzuckerman84}, \cite{wallach88}, Thm. 9.7.1)
that $I^{ij}_F$ is unitary iff $I^{ij}_\C$ is so and the corresponding
${}^0\laaa_0$ annihilates the highest weight vector of $F$. This leads to the
conditions on $\nu$. Compare also the list in terms of lowest 
$K$-types given in \cite{baldonisilvabarbasch83}.
\hB

We remark that we could have used the above mentioned Vogan-Zuckerman classification in terms of the modules $A_\qaaa(\nu)$ (together with the a 
priori knowledge that all unitarizable modules with integral regular infinitesimal character have cohomology \cite{salamancariba99}) 
in order to obtain
the list of all unitarizable modules $I^w_F$ together with their
cohomology at once. For this one has to relate the parametrization by $W^0$ to
the parametrization by $\theta$-stable
parabolic subalgebras $\qaaa\subset\gaaa_\C$ (modulo certain equivalence relations). A short indication how this rather involved correspondence
works can be found in \cite{baldonisilvabarbasch83}, pp. 30-31. 
However, Theorems \ref{koller} and \ref{uni} together
contain more information which might be of some value
for possible generalizations of Theorem \ref{main1} and Theorem \ref{main2} to arbitrary
rank one groups.

For each discrete series module $(\xi_{F,i},V_{\xi_{F,i}})$, 
$i=0,\dots,k_X-1$, with
infinitesimal character $\chi_F$ we introduce the Hilbert
space of multiplicities of the dual module $$D_{F,i}(\vp):=\Hom_{\gaaa,K}(V_{\tilde\xi_{F,i}},L^2(\Gamma\backslash G,\vp)_K)\ .$$ 
It can be identified with a subspace of $\Gamma$-invariant elements in the
tensor product of the distribution globalization of $V_{\xi_{F,i}}$ 
with $V_\vp$ (see \cite{bunkeolbrich000}, Ch. 8). 
By Theorem \ref{plaff} we have $\dim D_{F,i}(\vp)=\infty$.
As in Subsection \ref{plum} we consider the non-negative real number 
$d_\Gamma$ defined
by $d_\Gamma\alpha:=\delta_\Gamma+\rho$. Note that $d_\Gamma$ is equal to
the Hausdorff dimension of $\Lambda$ with respect to the natural class of sub-Riemannian metrics on $\partial X$ (\cite{corlette90},\cite{patterson762}, \cite{sullivan79}). We have (see Corollary \ref{corlette})  
\begin{equation}\label{Hd}
\mbox{\begin{tabular}{|c||c|c|c|c|}
\hline
$X$&$\R H^n$&$\C H^n$&$\HH H^n$&$\OO H^2$\\
\hline
$d_\Gamma<$&$n-1$&$2n$&$4n$&$16$\\
\hline
\end{tabular} . }
\end{equation}

Now we can give a complete description of the $L^2$- cohomology
groups $H^p_{(2)}(Y,E(\pi\otimes\vp))$ and $\cH^p_{(2)}(Y,E(\pi\otimes\vp))$
(recall the decomposition (\ref{bloed})).

\begin{theorem}\label{l2co}
Let $(\pi,F)$ be an irreducible finite-dimensional representation of $G$,
$(\vp, V_\vp)$ a finite-dimensional unitary representation of $\Gamma$.
We use the same notation for the highest weight of $F$ as in Proposition \ref{uni}.

If $X\ne\OO H^2$ or if $\delta_\Gamma<0$, then $H^p_{(2)}(Y,E(\pi\otimes\vp))\ne\cH^p_{(2)}(Y,E(\pi\otimes\vp))$
if and only if $X=\R H^n$, $n$ odd, $p=\frac{n+1}{2}$, and 
$p_F\le\frac{n-1}{2}$ (i.e. $m_\frac{n-1}{2}=0$).
In general, we have for
\begin{itemize}
\item $X=\R H^n$: 
\begin{eqnarray*}
\cH^p_{(2)}(Y,E(\pi\otimes\vp))&\cong& 
U_\Lambda(\sigma^p_{F,\lambda_p},\vp)\ ,
\quad p\le\frac{n-2}{2}\ ,\\ 
\cH^\frac{n-1}{2}_{(2)}(Y,E(\pi\otimes\vp))&\cong& U_\Lambda(\sigma^{\pm}_{F,\lambda_\pm},\vp)\ ,\\       
\cH^\frac{n}{2}_{(2)}(Y,E(\pi\otimes\vp))&\cong& D_{F,0}(\vp)\oplus D_{F,1}(\vp)\ .
\end{eqnarray*}
If $p<\min\{p_F,n-1-d_\Gamma\}$, then $\cH^p_{(2)}(Y,E(\pi\otimes\vp))=\{0\}$.
\item $X=\C H^n$:
\begin{eqnarray*}
\cH^p_{(2)}(Y,E(\pi\otimes\vp))&\cong&\bigoplus_{l=0}^{[\frac{p-1}{2}]}
\bigoplus_{i+j=p-2l} U_\Lambda(\sigma^{ij}_{F,\lambda_{ij}},\vp)\ ,
\quad p\le n-1\ ,\\ 
\cH^n_{(2)}(Y,E(\pi\otimes\vp))&\cong& \bigoplus_{i=0}^{n} D_{F,i}(\vp) \oplus \bigoplus_{l=1}^{[\frac{n-1}{2}]}\bigoplus_{i+j=n-2l} U_\Lambda(\sigma^{ij}_{F,\lambda_{ij}},\vp)\ .
\end{eqnarray*}
If $m_i>m_{n-j}$ or $i+j<2n-d_\Gamma$, then $U_\Lambda(\sigma^{ij}_{F,\lambda_{ij}},\vp)=\{0\}$.
\item $X=\HH H^n$:
\begin{eqnarray*}
\cH^p_{(2)}(Y,E(\pi\otimes\vp))&\cong&\bigoplus_{l=0}^{[\frac{p-n}{2}]}
\bigoplus_{i+j=p-2l,j\ge n} U_\Lambda(\sigma^{ij}_{F,\lambda_{ij}},\vp)\ ,
\quad p\le 2n-1,\ p \mbox{ odd}\ ,\\
\cH^p_{(2)}(Y,E(\pi\otimes\vp))&\cong&\bigoplus_{i=3,5,\dots,\frac{p}{2}} U_\Lambda(\sigma^{ii}_{F,\lambda_{ii}},\vp)\oplus\bigoplus_{l=0}^{[\frac{p-n}{2}]}
\bigoplus_{i+j=p-2l,j\ge n} U_\Lambda(\sigma^{ij}_{F,\lambda_{ij}},\vp)\ ,\\
&&\hspace{5.6cm} p\le 2n-1,\ p\equiv 2\:\mod\: 4\ ,\\ 
\cH^p_{(2)}(Y,E(\pi\otimes\vp))&\cong&\bigoplus_{i=2,4,\dots,\frac{p}{2}} U_\Lambda(\sigma^{ii}_{F,\lambda_{ii}},\vp)\oplus\bigoplus_{l=0}^{[\frac{p-n}{2}]}
\bigoplus_{i+j=p-2l,j\ge n} U_\Lambda(\sigma^{ij}_{F,\lambda_{ij}},\vp)\ ,\\
&&\hspace{5.6cm} p\le 2n-1,\ p\equiv 0\:\mod\: 4\ ,\\
\cH^{2n}_{(2)}(Y,E(\pi\otimes\vp))&\cong& \bigoplus_{i=0}^{n} D_{F,i}(\vp) \oplus \bigoplus_{i=3,5,\dots,n} U_\Lambda(\sigma^{ii}_{F,\lambda_{ii}},\vp)\\
&&\hspace{1.5cm}\oplus\bigoplus_{l=0}^{[\frac{p-n}{2}]}
\bigoplus_{i+j=p-2l,j\ge n} U_\Lambda(\sigma^{ij}_{F,\lambda_{ij}},\vp)\ ,
\quad n\mbox{ odd}\ ,\\
\cH^{2n}_{(2)}(Y,E(\pi\otimes\vp))&\cong& \bigoplus_{i=0}^{n} D_{F,i}(\vp) \oplus \bigoplus_{i=2,4,\dots,n} U_\Lambda(\sigma^{ii}_{F,\lambda_{ii}},\vp)\\
&&\hspace{1.5cm}\oplus\bigoplus_{l=0}^{[\frac{p-n}{2}]}
\bigoplus_{i+j=p-2l,j\ge n} U_\Lambda(\sigma^{ij}_{F,\lambda_{ij}},\vp)\ ,
\quad n\mbox{ even}\ .
\end{eqnarray*}
If $m_i\ne 0$ or $2i<4n+2-d_\Gamma$, then $U_\Lambda(\sigma^{ii}_{F,\lambda_{ii}},\vp)=\{0\}$. If $m_i>m_{2n-j}$ or
$i+j<4n+1-d_\Gamma$, then $U_\Lambda(\sigma^{ii}_{F,\lambda_{ij}},\vp)=\{0\}$.
\item $X=\OO H^2$: 
\begin{eqnarray*}
\cH^p_{(2)}(Y,E(\pi\otimes\vp))&=&\{0\}\ ,
\quad p\le 5\ ,\\
\cH^6_{(2)}(Y,E(\pi\otimes\vp))&\cong&
U_\Lambda(\sigma^{15}_{F,\lambda_{15}},\vp)
\oplus U_\Lambda(\sigma^{24}_{F,\lambda_{24}},\vp)\ ,\\
\cH^7_{(2)}(Y,E(\pi\otimes\vp))&\cong&
U_\Lambda(\sigma^{16}_{F,\lambda_{16}},\vp)
\oplus U_\Lambda(\sigma^{25}_{F,\lambda_{25}},\vp)\ ,\\        
\cH^8_{(2)}(Y,E(\pi\otimes\vp))&\cong& D_{F,0}(\vp)\oplus D_{F,1}(\vp)\oplus D_{F,2}(\vp)
\oplus U_\Lambda(\sigma^{15}_{F,\lambda_{15}},\vp)\ .
\end{eqnarray*}
We have $U_\Lambda(\sigma^{ij}_{F,\lambda_{ij}},\vp)=\{0\}$, if $\lambda_{ij}>\delta_\Gamma$ or if the highest weight $\nu=(m_0,m_1,m_2,m_3)$ of $F$ does not satisfy the condition
given by the table in Proposition \ref{uni}. 
\end{itemize}
By Poincar\'e duality $\cH^{\dim Y-p}_{(2)}(Y,E(\pi\otimes\vp))\cong \cH^p_{(2)}(Y,E(\pi\otimes\vp))$. For $p\ne\frac{\dim Y}{2}$
the spaces $\cH^p_{(2)}(Y,E(\pi\otimes\vp))$ are finite-dimensional.
In case of even dimension $\cH^\frac{\dim Y}{2}_{(2)}(Y,E(\pi\otimes\vp))$ is always infinite-dimensional.
\end{theorem}
\proof We employ the isomorphism (\ref{hihi}) which sends
$\Delta_F$ to $-\Omega\otimes \id\otimes\id+\id\otimes\pi(\Omega)\otimes\id$. 

Lemma \ref{doof}
and the decomposition (\ref{bloed}) tell us that $\overline{B^p_{(2)}}\ne B^p_{(2)}$ if and only if $0$ belongs to the continuous spectrum of $\Delta_F$
restricted to $p$-cocycles.
By Theorem \ref{plaff} this can only occur if
$$ 0\in \{\chi_{\sigma,\lambda}(\Omega)-\pi(\Omega)\:|\: \lambda\in i\aaaa^*, \sigma\in\hat M \mbox{ s.th. }
[H^{\sigma,\lambda}\otimes F\otimes \Lambda^p\paaa^*]^K\ne\{0\}\}\ .$$
If $(\sigma,\lambda)\in \hat M\times i\aaaa^*$ satifies $\chi_{\sigma,\lambda}(\Omega)-\pi(\Omega)=0$, then Proposition \ref{char}
implies that 
\begin{equation}\label{hundert}
[H^{\sigma,\lambda}\otimes F\otimes \Lambda^p\paaa^*]^K=H^p(\gaaa,K,H^{\sigma,\lambda}\otimes F)
\end{equation}
and that $\chi_{\sigma,\lambda}=\chi_{\tilde F}$. 
By Propositions \ref{koller} and \ref{uni} the space $H^p(\gaaa,K,H^{\sigma,\lambda}\otimes F)$ is non-zero if and only if
$X=\R H^n$, $n$ odd, $p=\frac{n\pm 1}{2}$, 
$p_{\tilde F}\le\frac{n-1}{2}$, $\lambda=0$, and $\sigma=\sigma^{\pm}_{\tilde F}$. Note that $p_{\tilde F}=p_F$. Fix such a $\sigma$. 
It remains to check for which $p\in\{\frac{n-1}{2},\frac{n+1}{2}\}$ the
corresponding direct integral over $i\aaaa^*$ contributes to cocycles
in $\Omega^p_{(2)}(Y,E(\pi\otimes\vp))$. This happens iff the space
of $p$-cocycles in the relative Lie algebra cohomology complex
$([H^{\sigma,\lambda}\otimes F\otimes \Lambda^p\paaa^*]^K, d)$ is non-zero
for a nonempty open subset of $i\aaaa^*$. The left hand side of (\ref{hundert})
does not depend on $\lambda$. It follows that $[H^{\sigma,\lambda}\otimes F\otimes \Lambda^p\paaa^*]^K=\{0\}$ for $p\not\in\{\frac{n-1}{2},\frac{n+1}{2}\}$ and all $\lambda$. Since for
$\lambda\ne 0$ the cohomology $H^p(\gaaa,K,H^{\sigma,\lambda}\otimes F)$
vanishes for all $p$ we conclude that $d$ is injective on $[H^{\sigma,\lambda}\otimes F\otimes \Lambda^\frac{n-1}{2}\paaa^*]^K$ and
that $[H^{\sigma,\lambda}\otimes F\otimes \Lambda^\frac{n+1}{2}\paaa^*]^K$
consists entirely of cocyles.
This proves the first assertion of the theorem.

Let $\hat G$ be the unitary dual of $G$. If we write
$$ L^2(\Gamma\backslash G,\vp)_d=\bigoplus_{\xi\in\hat G}^{\mathrm Hilbert}
\Hom_G(W_\xi,L^2(\Gamma\backslash G,\vp))\hat\otimes W_\xi\ ,$$
then
$$\cH^p_{(2)}(Y,E(\pi\otimes\vp))=\bigoplus_{\{\xi\in\hat G\:|\:\xi(\Omega)=\pi(\Omega)\}}^{\mathrm Hilbert}
\Hom_G(W_\xi,L^2(\Gamma\backslash G,\vp))\otimes [W_{\xi,K}\otimes F\otimes\Lambda^p\paaa^*]^K\ .$$
Here $W_{\xi,K}$ denotes the underlying $(\gaaa,K)$-module of $W_\xi$.
By Proposition \ref{char} the right hand side is equal to the finite
sum
$$ \bigoplus_{\{\xi\in\hat G\:|\:\chi_\xi=\chi_{\tilde F}\}}
\Hom_G(W_\xi,L^2(\Gamma\backslash G,\vp))\otimes 
H^p(\gaaa,K,W_{\xi,K}\otimes F)\ .$$
If $W_\xi$ is a discrete series representation with $\chi_\xi=\chi_{\tilde F}$, then $\Hom_G(W_\xi,L^2(\Gamma\backslash G,\vp))=D_{F,i}(\vp)$ for some
$i\in\{0,\dots,k_X-1\}$. The underlying $(\gaaa, K)$-modules of the remaining
unitary representations with $\chi_\xi=\chi_{\tilde F}$ are listed in Proposition \ref{uni}. Theorem \ref{plaff} and Proposition \ref{exc} then assert that $\dim D_{F,i}(\vp)=\infty$ and
that for the remaining representations $\Hom_G(W_\xi,L^2(\Gamma\backslash G,\vp))$ is one of the finite-dimensional spaces $U_\Lambda(\sigma^{w}_{F,\lambda_{w}},\vp)$, $w\in W^0$. Note that $\widetilde
{\sigma^{w}_{\tilde F}}=\sigma^{w}_{F}$ unless $X=\R H^n$, $n$ odd, $w=w_-$. The restrictions for the possible values of 
$\delta_\Gamma$ indicated in table (\ref{Hd}) force some of these multiplicity spaces to be zero.
The cohomology groups $H^p(\gaaa,K,W_{\xi,K}\otimes F)$ have been determined in Propositions \ref{ds} and \ref{koller}.
This finishes the proof of the theorem.
\hB

We would like to emphasize the following two vanishing results which are contained in Theorem \ref{l2co}.

\begin{kor}\label{ver1}
Define the critical strip $S_X$ by the following table:
$$ \mbox{\begin{tabular}{|c|c|}
\hline
$X$&$S_X$\\
\hline\hline
$\R H^n$&$[n-1-d_\Gamma, d_\Gamma+1]$\\
\hline
$\C H^n$&$[2n-d_\Gamma,d_\Gamma]$\\
\hline
$\HH H^n$&$[2+(4n-d_\Gamma), d_\Gamma-2]\cap 2\Z\ \cup\  [1+(4n-d_\Gamma), d_\Gamma-1]\cap [n,3n]$\\
\hline
$\OO H^2$
&$[6,10]\cap [19-d_\Gamma, d_\Gamma-3]$\\
\hline
\end{tabular} }\ .
$$
If $p\not\in S_X$ and $p\ne\frac{\dim Y}{2}$, then $\cH^p_{(2)}(Y,E(\pi\otimes\vp))=\{0\}$.
\end{kor}

\begin{kor}\label{ver2}
If $F$ is generic, i.e., its highest weight $\nu$ lies in the interior of
the positive Weyl chamber $\haaa_+^*$, then $\cH^p_{(2)}(Y,E(\pi\otimes\vp))=\{0\}$ for $p\ne\frac{\dim Y}{2}$.
\end{kor}

We also conclude

\begin{kor}\label{maphi}
If $X=\R H^n$ and $p\ge\frac{n+1}{2}$, then the natural map  
$$\cH^p_{(2)}(Y,E(\pi\otimes\vp))\rightarrow H^p(Y,E(\pi\otimes\vp))$$
is an isomorphism. For $p=\frac{n}{2}$ it is surjective.
\end{kor}
\proof Let $p\ge\frac{n+1}{2}$. By Theorem \ref{main2} and Theorem \ref{l2co} both sides
are isomorphic to $U_\Lambda(\sigma^{n-p}_{F,\lambda_{n-p}},\vp)$
(or to $U_\Lambda(\sigma^\pm_{F,0},\vp)$, respectively). We leave the simple verification that these
isomorphisms are compatible with the map $\cH^p_{(2)}(Y,E(\pi\otimes\vp))\rightarrow H^p(Y,E(\pi\otimes\vp))$ to the reader. Similarly, for $p=\frac{n}{2}$ we use the fact that $Z^\frac{n}{2}_{F,\Lambda}(\vp)\subset D_{F,0}(\vp)\oplus D_{F,1}(\vp)$
(compare the proof of Lemma \ref{dser}).
\hB

It was asked in the introduction of \cite{bunkeolbrich000} what
the significance of the space $L^2(\Gamma\backslash G,\vp)_U$, 
that is of
the multiplicity spaces $U_\Lambda(\sigma_\lambda,\vp)$, could be. For regular infinitesimal character $\chi_{\sigma,\lambda}$ Theorem
\ref{l2co} provides an answer. The spaces $U_\Lambda(\sigma_\lambda,\vp)$ appear as a kind of primitive parts of the $L^2$-cohomology spaces $\cH^p_{(2)}(Y,E(\pi\otimes\vp))$, $\chi_{\sigma,\lambda}=\chi_F$.
Indeed, for $X=\C H^n$ it is not difficult to verify that 
$$ \cH^{p,q}_{0,(2)}(Y,E(\pi\otimes\vp))=\left\{\begin{array}{ccc}
U_\Lambda(\sigma^{qp}_{F,\lambda_{qp}},\vp)&& p+q\le n-1\\
D_{F,q}&& p=n-q
\end{array}\right.
\ ,$$
where $\cH^{p,q}_{0,(2)}(Y,E(\pi\otimes\vp))$ denotes the space of
square integrable harmonic primitive $(p,q)$-forms. (It would be interesting
to find an analogous interpretation of the spaces $U_\Lambda(\sigma_\lambda,\vp)$ for weakly regular, singular
infinitesimal character $\chi_{\sigma,\lambda}$.)
We see that the spaces $U_\Lambda(\sigma_\lambda,\vp)$
play a similar role as the multiplicity spaces of the discrete series.
This phenomenon can also be observed from a pure harmonic analysis
point of view. Let $\cC(\Gamma\backslash G,\vp)\subset L^2(\Gamma\backslash G,\vp)$ be the Schwartz space (see \cite{bunkeolbrich000}, Ch. 8) which is
the analog for $\Gamma\backslash G$ of the Harish-Chandra Schwartz space of $G$.
Its subspace
$$ {}^0\cC(\Gamma\backslash G,\vp)=\{f\in \cC(\Gamma\backslash G,\vp)\:|\: 
\int_N f(gnh)\:dn =0 \mbox{ for all } g,h\in G \mbox{ s.th. } gP\in\Omega\}$$
may be considered as the analog of the Harish-Chandra space of cusp forms on $G$.
The closure of ${}^0\cC(\Gamma\backslash G,\vp)$ in $L^2(\Gamma\backslash G,\vp)$ can be shown to be equal to
$$ L^2(\Gamma\backslash G,\vp)_{cusp}\oplus L^2(\Gamma\backslash G,\vp)_U\ .$$

\newpage

\bibliographystyle{plain}

\end{document}